\documentclass[12pt,authoryear]{elsarticle}
\journal{}

\makeatletter
\def\ps@pprintTitle{%
 \let\@oddhead\@empty
 \let\@evenhead\@empty
 \def\@oddfoot{\today\hfill\thepage}%
 \let\@evenfoot\@oddfoot}
\makeatother

\usepackage{amsmath,amssymb,amsfonts,amsthm,graphicx,appendix}
\usepackage[bookmarksnumbered,colorlinks=true]{hyperref}
\usepackage[labelfont=bf]{caption}

\providecommand{\doi}[1]{\href{https://doi.org/#1}{DOI:#1}}
\usepackage{xurl} 
\renewcommand{\doi}[1]{%
 \href{https://doi.org/#1}{\nolinkurl{DOI:#1}}%
}

\usepackage[OT1]{fontenc}
\usepackage{natbib,enumitem,mathtools,mathrsfs,dsfont}
\usepackage{tikz,float}
\usetikzlibrary{arrows.meta}
\usepackage{geometry}
\newtheorem{theorem}{Theorem}[section]
\newtheorem{lemma}[theorem]{Lemma}
\newtheorem{corollary}[theorem]{Corollary}
\newtheorem{proposition}[theorem]{Proposition}
\newtheorem{conjecture}[theorem]{Conjecture}

\theoremstyle{definition}

\newtheorem{remark}[theorem]{Remark}

\numberwithin{equation}{section}

\newcommand{\N}{\mathbb{N}}
\newcommand{\Z}{\mathbb{Z}}
\newcommand{\R}{\mathbb{R}}
\newcommand{\C}{\mathbb{C}}
\newcommand{\PP}{\mathsf{P}}
\newcommand{\EE}{\mathsf{E}}
\newcommand{\Var}{\mathsf{Var}}
\newcommand{\Cov}{\mathsf{Cov}}
\newcommand{\bb}[1]{\boldsymbol{#1}}
\newcommand{\ind}{\mathds{1}}
\newcommand{\rd}{\mathrm{d}}
\newcommand{\e}{\varepsilon}
\newcommand{\ii}{\mathrm{i}}
\renewcommand{\Re}{\mathrm{Re}\, }
\renewcommand{\Im}{\mathrm{Im}\, }
\newcommand{\OO}{\mathcal{O}}
\newcommand{\oo}{\mathrm{o}}
\newcommand{\leqdef}{\vcentcolon=}

\newcommand{\sig}{\sigma}
\newcommand{\logp}{\log_{+}}

\allowdisplaybreaks

\begin{document}

\begin{frontmatter}

\title{Poisson--Dirichlet statistics for the log-modulus of the Riemann zeta function on expanding intervals near the critical line}

\author[a1,a2]{Fr\'{e}d\'{e}ric Ouimet}

\address[a1]{California Institute of Technology, USA*\footnote{*A rough draft of this paper, containing all the main results stated in Section~\ref{sec:main.results}, was written in 2019--2021 during the author's postdoctoral appointment at the California Institute of Technology.}}
\address[a2]{Universit\'{e} du Qu\'{e}bec \`{a} Trois-Rivi\`{e}res, Canada. \ead{frederic.ouimet2@uqtr.ca}\vspace{-5mm}}

\begin{abstract}
For $T > 0$ large, let $\tau$ be uniformly distributed on $[T,2T]$ and let $\sig = \sig(T)$ be a deterministic real part approaching $1/2$ from the right. For the Riemann zeta function $\zeta$, we study the Gibbs measure with density proportional to $|\zeta(\sig + \ii\tau + \ii h)|^{\beta}$ on $I_{\theta} \leqdef [ - (\log{T})^{\theta},(\log{T})^{\theta}]$, for every $\theta > 0$ and $\beta > 0$. As $\beta$ increases, the Gibbs measure increasingly favors positions where $|\zeta|$ is large, providing a way to study the spatial organization of extreme values. Under explicit conditions on the rate at which $\sig$ approaches the critical line, a phase transition is established at the critical inverse temperature $\beta_{c}(\theta) \leqdef 2\sqrt{1 + \theta}$. The limiting distribution of the correlation (or overlap) between the field values at two positions sampled from the Gibbs measure is concentrated at $0$ when $\beta \leq \beta_{c}(\theta)$, and assigns masses $\beta_{c}(\theta)/\beta$ and $1 - \beta_{c}(\theta)/\beta$ to $0$ and $1$, respectively, when $\beta > \beta_{c}(\theta)$. Joint convergence of the overlaps is also established, with the limiting arrays described by Poisson--Dirichlet weights in the low-temperature regime. The Ghirlanda--Guerra identities, which provide recursive relations for the joint overlap distributions, are established in the limit, uniformly for bounded deterministic observables of the sampled positions. All our results are established independently of the Riemann hypothesis.
\end{abstract}


\end{frontmatter}

\setcounter{tocdepth}{2}
\tableofcontents

\section{Introduction}\label{sec:intro}

\subsection{Background}\label{sec:background}

For a complex argument $s$, the Riemann zeta function is defined by
\[
\zeta(s) \leqdef \sum_{n = 1}^{\infty} \frac{1}{n^{s}} = \prod_{p\text{ prime}} \left(1 - \frac{1}{p^{s}}\right)^{-1}, \qquad \Re(s) > 1,
\]
where $\Re(s)$ is the real part of $s$ and $p$ ranges over the prime numbers. It extends meromorphically to the complex plane $\C$, with a simple pole at $s = 1$.

Let $T > 0$ be the large height parameter and let $\tau$ be uniformly distributed on $[T,2T]$ under the probability law $\PP$, with corresponding expectation denoted by $\EE$. Selberg's central limit theorem \citep{MR0020594,MR1220477} showed that $\log|\zeta(1/2 + \ii \tau)|$ has an asymptotically normal distribution with mean zero and variance $(1/2)\log \log T$. A short elegant proof was given by \citet{MR3832861}, and a multidimensional extension was proved by \citet{MR2678896}. More precisely, for every fixed $a < b$,
\[
\PP\left(\frac{\log|\zeta(1/2 + \ii \tau)|}{\sqrt{(1/2)\log \log T}} \in [a,b]\right) \to \frac{1}{\sqrt{2\pi}}\int_{a}^{b} e^{-u^{2}/2}\rd{u}, \qquad T\to \infty.
\]
Thus the typical fluctuations of $\log|\zeta(1/2 + \ii \tau)|$ have size $\sqrt{\log{\log{T}}}$.

Recent work studied the Gaussian decay suggested by Selberg's theorem for the upper-tail probability $\PP(\log|\zeta(1/2 + \ii \tau)| > V)$ when $V \sim \gamma\log \log T$. Unconditional upper bounds of Gaussian order were proved for $0 < \gamma < 2$ by \citet{ArguinBailey2023}. Complementary unconditional lower bounds of the same order were subsequently obtained for every $\gamma > 0$ on the line $\Re(s) = 1/2 + \delta/\log T$, with $\delta$ sufficiently large depending on $\gamma$ \citep{ArguinBailey2025}, and then on the critical line itself \citep{ArguinCreighton2026}. Assuming the Riemann hypothesis, \citet{ArguinBaileyRoberts2026} established the corresponding upper bound for every fixed $\gamma > 0$, together with a new proof of the sharp-order conditional moment upper bounds of \citet{Harper2013Sharp}. These estimates quantify the tail regime that feeds both high moments and short-interval extremes.

As shown by \citet{MR2678896}, the covariance between the values $\log|\zeta(1/2 + \ii\tau + \ii h)|$ and $\log|\zeta(1/2 + \ii\tau + \ii h')|$ grows logarithmically as $|h - h'|$ decreases toward the microscopic scale $(\log T)^{-1}$. This logarithmic correlation structure suggested an analogy with characteristic polynomials of random unitary matrices. Building on that analogy, and on the corresponding predictions from statistical mechanics for log-correlated fields, \citet{FyodorovHiaryKeating2012} and \citet{MR3151088} studied the extreme-value behavior of $\log|\zeta(1/2 + \ii\tau + \ii h)|$ and the moments of $|\zeta(1/2 + \ii\tau + \ii h)|$ on intervals of bounded length, such as $I_0 = [-1,1]$. In particular, they predicted tightness of the local maximum $\max_{h\in I_{0}} \log|\zeta(1/2 + \ii\tau + \ii h)|$ after recentering by $\log \log T - (3/4)\log \log \log T$, as well as a change in the leading growth rate of the $\beta$-moment
\begin{equation}\label{eq:zeta.moment}
\int_{I_{0}} |\zeta(1/2 + \ii\tau + \ii h)|^{\beta}\rd{h}, \qquad \beta > 0,
\end{equation}
when $\beta$ crosses the critical value $\beta_{c}(0) = 2$. The leading order of the local maximum was subsequently proved for $\theta = 0$ in \citet{MR3851835}, conditionally on the Riemann hypothesis for the lower bound, and unconditionally in \citet{ABBRS_2019}. The phase transition of the $\beta$-moments was proved by \citet{AOR2019}. The complementary bounds of \citet{ABR-2020} and \citet{ArguinBourgadeRadziwill2023} established tightness of the recentered local maximum and the sharp order of its positive right tail.

The length of the interval over which the shift parameter $h$ varies also affects the extremes and the growth of the moments. For a fixed interval-length exponent $\theta > - 1$, define
\begin{equation}\label{eq:I.theta}
I_{\theta} \leqdef [ - (\log{T})^{\theta},(\log{T})^{\theta}], \qquad q_{0} \leqdef |\hspace{-0.4mm}\min\{\theta,0\}| \equiv |\theta \wedge 0|.
\end{equation}
The role of $\theta$ can already be seen from the Euler product. A formal expansion of its logarithm, followed by truncation of the resulting Dirichlet series, suggests the approximation
\begin{equation}\label{eq:log.zeta.approximation}
\log\zeta(1/2 + \ii\tau + \ii h) \approx \sum_{k = 1}^{\lfloor \log \log T\rfloor}\sum_{\substack{p \leq T \\ e^{k - 1} < \log{p} \leq e^{k}}}\frac{p^{-\ii\tau}}{p^{1/2 + \ii h}},
\end{equation}
up to lower-order contributions from prime powers and the truncation. Here $k$ indexes logarithmic prime scales, since
\[
e^{k} = (\log T)^{k/\log \log T},
\]
so the normalized coordinate $k/\log \log T \in [0,1]$ records the fraction of the logarithmic prime scales that has been accumulated. The dependence on $h$ is encoded, scale by scale and prime by prime, by the phases $p^{-\ii h}$. The contribution of a given range of primes therefore depends on how much these phases vary over $I_{\theta}$.

When $ - 1 < \theta < 0$, for every fixed $0 < \e < q_{0} = |\theta|$,
\[
\sup_{\substack{h \in I_{\theta} \\ \log{p} \leq (\log T)^{|\theta| - \e}}}|h\log{p}| \leq (\log T)^{-\e} \to 0, \qquad T \to \infty.
\]
Thus $p^{-\ii h} = 1 + \oo(1)$ uniformly over these primes and over $h \in I_{\theta}$, so their contribution is approximately the same throughout the interval. In terms of the logarithmic prime-scale coordinate $k/\log \log T$, this means that the prime contributions are approximately common at every fixed scale below $|\theta|$, whereas scales above $|\theta|$ produce nontrivial variation in $h$.

\begin{figure}[!t]
\centering
\includegraphics[width = \textwidth]{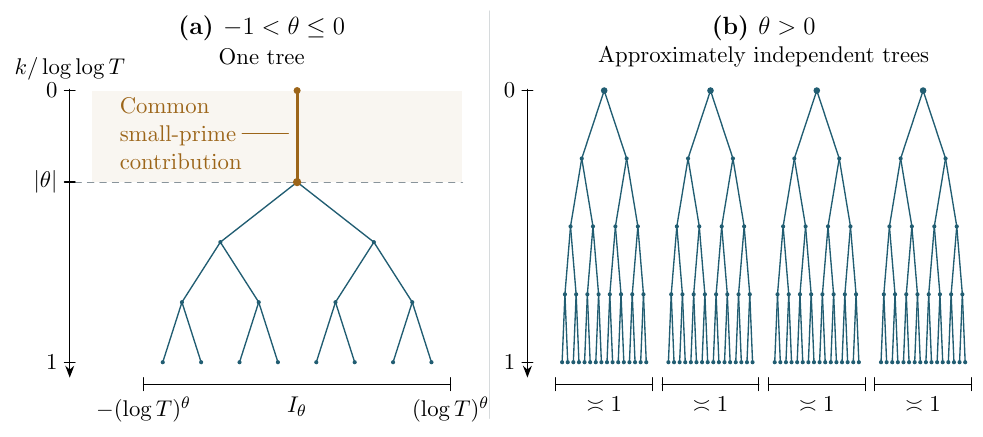}
\caption{Schematic logarithmic prime-scale description on $I_{\theta}$. The vertical coordinate is $k/\log \log T$; leaves represent spatial points separated on the scale $(\log T)^{-1}$. For $ - 1 < \theta \leq 0$, the common trunk represents the approximately common prime contributions at scales below $|\theta|$, followed by branching over the remaining depth $1 + \theta$. The common trunk disappears when $\theta = 0$. For $\theta > 0$, about $(\log T)^{\theta}$ bounded subintervals are represented by separate trees of depth $1$, each with about $\log T$ leaves; this depicts the vanishing correlation between field values in well-separated subintervals. Both panels represent about $(\log T)^{1 + \theta}$ spatial points.}\label{fig:prime.scale.regimes}
\end{figure}

When $\theta > 0$, the interval $I_{\theta}$ contains about $(\log T)^{\theta}$ subintervals of bounded length, and values in well-separated subintervals have correlation tending to zero, as made explicit for the Euler field later in \eqref{eq:overlaps.Re.Q.asymptotics}. The trees in Figure~\ref{fig:prime.scale.regimes} illustrate the leading covariance structure in these two regimes.

These facts explain both the smallest possible limiting correlation $|\theta \wedge 0|$ and the dependence on $\theta$ of the critical inverse temperature
\[
\beta_{c}(\theta) \leqdef
\begin{cases}
2, & \theta \leq 0, \\
2 \sqrt{1 + \theta}, & \theta > 0.
\end{cases}
\]
The regimes $\beta < \beta_{c}(\theta)$, $\beta = \beta_{c}(\theta)$, and $\beta > \beta_{c}(\theta)$ are called subcritical (high temperature), critical, and supercritical (low temperature), respectively.

For every fixed $\theta > - 1$ and $\beta > 0$, \citet{AOR2019} established the following first-order asymptotics in probability as $T \to \infty$, unconditionally when $ - 1 < \theta \leq 3$ and under the Riemann hypothesis when $\theta > 3$ (conditionally on RH only for the upper bounds):
\begin{equation}\label{eq:background.maximum}
\frac{1}{\log \log T}\max_{h \in I_{\theta}}\log|\zeta(1/2 + \ii \tau + \ii h)| \xrightarrow{\PP} m(\theta) \leqdef \frac{1 + (\theta \wedge 0)}{2}\,\beta_{c}(\theta)
\end{equation}
and
\begin{equation}\label{eq:background.moments}
\frac{1}{\log \log T}\log\int_{I_{\theta}} |\zeta(1/2 + \ii \tau + \ii h)|^{\beta}\rd{h} \xrightarrow{\PP}
\begin{cases}
\frac{\beta^{2}}{4}(1 + (\theta \wedge 0)) + \theta, & \beta \leq \beta_{c}(\theta), \\
\beta m(\theta) - 1, & \beta > \beta_{c}(\theta).
\end{cases}
\end{equation}
The quadratic dependence on $\beta$ in the first line of \eqref{eq:background.moments} corresponds to contributions from many regions of $I_{\theta}$. Beyond $\beta_{c}(\theta)$, the limit is linear in $\beta$, with slope $m(\theta)$, reflecting the leading contribution of values close to the local maximum. The related problem of quantifying the dependence between the absolute values of zeta at different spatial positions was investigated through shifted moments by \citet{Chandee2011,NgShenWong2024,Curran2024Correlations,Curran2024MomentsOfMoments,Curran2024LowerShifted,Kovaleva2025}, among others.

For shrinking intervals $I_{\theta}$ with $ - 1 < \theta < 0$, \citet{ArguinHamdan2026} isolated an approximately common small-prime contribution and showed that it is asymptotically Gaussian after normalization. After subtracting this contribution from $\log|\zeta|$, they obtained an upper bound of the conjectured sharp order for the positive right tail of the local maximum centered at $(1 + \theta)\log \log T - (3/4)\log \log \log T$, together with an upper bound in probability for the suitably normalized integral of $|\zeta|^{2}$, which extends the case $\theta = 0$ treated by \citet{arXiv:1906.05783}.

\subsection{Gibbs measures, overlaps, and the Euler approximation}\label{sec:gibbs.problem}

\subsubsection{Parameter range and Gibbs measure}

Let $t \leqdef \log \log T$, and fix constants $\phi > 7/2$, $\Phi > \phi + 1$, and $C_{\sig} > 0$. The real part $\sig = \sig(T)$ and the auxiliary prime-power cutoff $X = X(T)$ are deterministic functions of $T$. Write $\eta_{\sig} \leqdef \sig - 1/2$ for the displacement from the critical line and assume throughout the remainder of the paper that, as $T\to \infty$,
\begin{equation}\label{eq:fixed.parameters}
T^{1/t^{\Phi}} \leq X \leq T^{1/t^{\phi}}, \qquad 0 < \eta_{\sig} \leq \frac{C_{\sig}t^{\Phi}}{\log T}, \qquad \frac{\eta_{\sig}\log{X}}{t} \to \infty.
\end{equation}
Subsubsection~\ref{sec:parameter.restrictions} explains these restrictions.

The assumptions imply
\begin{equation}\label{eq:parameter.scales}
t - \Phi\log{t} \leq \log{\log{X}} \leq t - \phi\log{t}, \qquad \log(1/\eta_{\sig}) = t + \OO(\log{t}).
\end{equation}
Indeed, $\eta_{\sig}\log{X} \geq t$ eventually, so $\eta_{\sig} \geq t^{\phi + 1}/\log T$, while the upper bound on $\eta_{\sig}$ in \eqref{eq:fixed.parameters} gives the other side of the second estimate in \eqref{eq:parameter.scales}. For example, the choices $X = T^{1/t^{4}}$ and $\sig = 1/2 + \delta t^{6}/\log T$, for any fixed $\delta > 0$, satisfy \eqref{eq:fixed.parameters} with $\Phi = 6$, $\phi = 4$, and $C_{\sig} \geq \delta$. In this example, $\eta_{\sig}\log{X}/t = \delta t \to \infty$.

For inverse temperature $\beta > 0$, define the Gibbs measure associated with the logarithmic field $\{\log|\zeta(\sig + \ii \tau + \ii h)|\}_{h \in I_{\theta}}$ by
\[
\mathcal{G}_{\beta,T}(\rd{h}) \leqdef \frac{|\zeta(\sig + \ii \tau + \ii h)|^{\beta}}{\int_{I_{\theta}}|\zeta(\sig + \ii \tau + \ii h')|^{\beta}\rd{h'}} \rd {h},
\]
where $\rd{h}$ denotes the Lebesgue measure. The normalization makes $\mathcal{G}_{\beta,T}$ a random probability measure on $I_{\theta}$. For fixed $\tau$, increasing $\beta$ gives greater relative weight to positions where $|\zeta(\sig + \ii \tau + \ii h)|$ is larger. Throughout, the free energy is the normalized log-partition function $\smash{t^{-1}\log\int_{I_{\theta}}|\zeta(\sig + \ii \tau + \ii h)|^{\beta}\rd{h}}$, with the analogous definition for the other fields.

The estimates \eqref{eq:background.maximum} and \eqref{eq:background.moments} describe local maxima and moments, but do not determine the spatial organization of the Gibbs mass. The aim here is to characterize this organization through the joint overlap distributions of positions sampled from $\mathcal{G}_{\beta,T}$ and to determine how these distributions depend on $\beta$ and $\theta$.

\subsubsection{Overlaps and joint sampling}

Conditional on $\tau$, positions $h_{1},\ldots,h_{n}$ sampled independently from $\mathcal{G}_{\beta,T}$ have joint law $\mathcal{G}_{\beta,T}^{\times n}$. More generally, for a probability measure $\mu$ and a measurable function $\Psi$ of $n$ positions, $\mu^{\times n}[\Psi(h_{1},\ldots,h_{n})]$ denotes integration against its $n$-fold product; the same notation applies to a finite measure before normalization. Throughout the paper, the sampled positions in Gibbs expectations and probabilities are displayed explicitly, either individually or through a defined vector such as $\bb{h} \leqdef (h_{1},\ldots,h_{n})$.

To describe the correlation between field values at sampled positions, define the overlap of the logarithmic zeta field at $h,h' \in I_{\theta}$ by
\begin{equation}\label{eq:overlaps.log.zeta}
\rho_{T}(h,h') \leqdef \frac{\EE[\log|\zeta(\sig + \ii \tau + \ii h)|\log|\zeta(\sig + \ii \tau + \ii h')|]}{\sqrt{\EE[(\log|\zeta(\sig + \ii \tau + \ii h)|)^{2}]}\sqrt{\EE[(\log|\zeta(\sig + \ii \tau + \ii h')|)^{2}]}}.
\end{equation}
The expectations in this definition are taken over $\tau$ at fixed positions $h$ and $h'$. The means $\EE[\log|\zeta(\sig + \ii \tau + \ii h)|]$ are $\OO(1)$, whereas the second moments have order $\log \log T$. Centering therefore changes $\rho_{T}(h,h')$ by $\oo(1)$, so it is asymptotically the usual correlation. By definition, $\rho_{T}(h,h) = 1$. The distribution of $\rho_{T}(h_{1},h_{2})$ under $\EE\mathcal{G}_{\beta,T}^{\times 2}$ is called the two-overlap distribution and includes the expectation over $\tau$. An expectation taken over both the field randomness and the sampled positions is called an annealed expectation.

The goal of the present article is to determine, on expanding intervals ($\theta > 0$), the limiting joint distribution of the overlap array $\big(\rho_{T}(h_{j},h_{k})\big)_{1 \leq j,k \leq n}$ under $\EE\mathcal{G}_{\beta,T}^{\times n}$ as $T \to \infty$, for every fixed $n \in \N$. Under this annealed law, $\tau$ is sampled first, followed by positions sampled independently conditional on $\tau$. The logarithmic prime-scale picture in Figure~\ref{fig:prime.scale.regimes} suggests a limiting overlap range from $0$ to $1$, according to how many logarithmic prime scales two positions share. The sampling problem asks which overlap values persist in the limit, how the sampled positions form groups with large mutual overlaps, and whether the critical inverse temperature $\beta_{c}(\theta)$ from the moment asymptotics also governs this grouping. Determining the joint law for every fixed $n$ describes the relative masses of these groups through their sampling probabilities. The main results are stated in Section~\ref{sec:main.results}, together with a proposed extension to bounded or shrinking intervals ($-1 < \theta \leq 0$) in Conjecture~\ref{conj:bounded.overlap.laws}.

\subsubsection{Euler approximation}

The heuristic approximation in \eqref{eq:log.zeta.approximation} motivates the Euler polynomial defined below. Its definition starts from the absolutely convergent expansion
\begin{equation}\label{eq:expand.logarithm.zeta}
\log|\zeta(s)| = \Re\sum_{p}\sum_{j = 1}^{\infty}\frac{1}{j}p^{-js}, \qquad \Re(s) > 1.
\end{equation}
Within the critical strip, the singularities of $\log|\zeta|$ at zeros of $\zeta$ prevent uniform approximation by Euler polynomials on regions containing such zeros. Nevertheless, with a suitable displacement to the right of the critical line and a suitable truncation level, a pointwise approximation holds for a proportion of ordinates tending to one, as shown, for example, by combining Propositions~3 and~4 of \citet{MR3832861}. Under \eqref{eq:fixed.parameters}, Section~\ref{sec:analytic} proves that the expected total variation distance between the zeta and Euler Gibbs measures tends to zero and that the corresponding overlaps differ by $\oo(1)$ uniformly over all pairs of spatial positions in $I_{\theta}$. These comparisons transfer the limiting overlap laws and approximate Ghirlanda--Guerra (GG) identities from the Euler field to the logarithmic zeta field.

With the prime-power cutoff $X$, define the Euler polynomial by
\[
\mathcal{Q}(s) \leqdef \sum_{j = 1}^{\infty}\sum_{p^{j} \leq X}\frac{1}{j}p^{-js},
\]
and its two-term reduction by
\[
\mathcal{Q}^{(2)}(s) \leqdef \sum_{p \leq X}p^{-s} + \frac{1}{2}\sum_{p^{2} \leq X}p^{-2s}.
\]

The random field $\{\Re{\mathcal{Q}(\sig + \ii \tau + \ii h)}\}_{h \in I_{\theta}}$ is called the Euler field or arithmetic field. The Gibbs measure and overlap associated with this field are
\[
G_{\beta,T}^{\mathrm{A}}(\rd{h}) \leqdef \frac{\exp(\beta \Re{\mathcal{Q}(\sig + \ii \tau + \ii h)})}{\int_{I_{\theta}}\exp(\beta \Re{\mathcal{Q}(\sig + \ii \tau + \ii h')})\rd{h'}}\rd{h},
\]
and
\[
R_{T}^{\mathrm{A}}(h,h') \leqdef \frac{\EE[\Re{\mathcal{Q}(\sig + \ii \tau + \ii h)}\Re{\mathcal{Q}(\sig + \ii \tau + \ii h')}]}{\sqrt{\EE[(\Re{\mathcal{Q}(\sig + \ii \tau + \ii h)})^{2}]}\sqrt{\EE[(\Re{\mathcal{Q}(\sig + \ii \tau + \ii h')})^{2}]}}.
\]
Here the superscript $\mathrm{A}$ stands for ``arithmetic''. Mertens' estimate, the mean-value formula, and classical bounds for $\zeta$ near the line of real part $1$ give, uniformly for $h,h' \in I_{\theta}$, the estimate deduced from Lemma~\ref{analytic:lem:covariance}:
\begin{equation}\label{eq:overlaps.Re.Q.asymptotics}
R_{T}^{\mathrm{A}}(h,h') = \frac{\log_{+}((\log T) \wedge |h - h'|^{-1})}{\log \log T} + \oo(1).
\end{equation}
Here $\log_{+}{y} \leqdef \max\{\log{y},0\}$ is the positive part of the logarithm for $y > 0$, and $|h - h'|^{-1} = \infty$ when $h = h'$. As for $\rho_{T}$, centering changes $R_{T}^{\mathrm{A}}$ by $\oo(1)$, so this overlap is also asymptotically the usual correlation. Two positions receive approximately common contributions from primes for which $|h - h'|\log{p}$ is small, while primes on larger logarithmic scales have smaller aggregate covariance. Formula~\eqref{eq:overlaps.Re.Q.asymptotics} makes this interpretation precise and places every possible limiting overlap in $[|\theta \wedge 0|,1]$. An overlap tending to $1$ corresponds to agreement through almost all logarithmic prime scales, without requiring the positions to coincide. The comparison in Section~\ref{sec:analytic} gives the same spatial description for $\rho_{T}$.

Two features make the Gibbs sampling problem more delicate than the first-order analysis of local maxima and free energies of $\log|\zeta|$ and $\Re{\mathcal{Q}}$. First, the contribution from prime powers of order at least three, $\smash{\Re{\mathcal{Q}} - \Re{\mathcal{Q}^{(2)}}}$, satisfies
\begin{equation}\label{eq:higher.prime.powers}
\sup_{\Re(s) \geq 1/2}|\Re{\mathcal{Q}(s)} - \Re{\mathcal{Q}^{(2)}(s)}| \leq \sum_{p}\sum_{j = 3}^{\infty}\frac{1}{j}p^{-j/2} < \infty.
\end{equation}
By \eqref{eq:higher.prime.powers}, replacing $\Re{\mathcal{Q}}$ by $\Re{\mathcal{Q}^{(2)}}$ changes the local maximum, $\max_{h \in I_{\theta}}\Re{\mathcal{Q}(\sig + \ii\tau + \ii h)}$, by $\OO(1)$ and changes the free energy, $t^{-1} \log\int_{I_{\theta}}e^{\beta\Re{\mathcal{Q}(\sig + \ii\tau + \ii h)}}\rd{h}$, by $\OO(t^{-1})$. However, a bounded spatially varying error in the field can change the Gibbs weights by factors that do not tend to one. The Gibbs measures and sampling arguments will therefore retain the full Euler field $\Re{\mathcal{Q}}$, including every prime power $p^{j} \leq X$.

Second, the arithmetic phases $p^{-\ii\tau}$ all depend on the same random variable $\tau$. For every fixed finite collection of primes, the joint moments of these phases converge, as $T \to \infty$, to the corresponding moments of independent random variables uniformly distributed on the unit circle; see, for example, Lemma~1.2 in \citet{arXiv:1904.08204}. This motivates the independent Haar field defined in \eqref{eq:haar.field} (henceforth the Haar field), where ``independent'' refers to phases attached to distinct primes and ``Haar'' to their uniform distribution on the unit circle. This convergence of joint moments for fixed finite collections of primes does not by itself justify replacing the arithmetic phases for all primes $p \leq X$ by independent Haar phases in the Gibbs measure, since $X$ grows with $T$. Any replacement of the Euler field by the Haar field must therefore be justified at the level of the joint laws of positions sampled from the Gibbs measures: this is the content of Theorem~\ref{cap:sampling}, for $\theta > 0$ and $\beta > \beta_{c}(\theta)$.

\subsubsection{Discussion of the parameter restrictions in \eqref{eq:fixed.parameters}}\label{sec:parameter.restrictions}

The restriction $\phi > 7/2$ comes from comparing the joint distributions of the Euler and Haar fields. In the proof of Lemma~\ref{an:smoothing}, their characteristic functions are approximated by Taylor polynomials of degree of order $D_{T} \asymp t^{7/2}(\log{t})^{2}$. The expectations of these polynomials are compared by averaging over $\tau \in [T,2T]$. The error in this averaging contains the factor $X^{2D_{T}}/T$; see \eqref{an:smoothing.detail:polynomial.error}. Since $\phi > 7/2$ gives $D_{T} = \oo(t^{\phi})$, the cutoff bound $X \leq T^{1/t^{\phi}}$ ensures that $X^{2D_{T}} = T^{\oo(1)}$, and the error in comparing the polynomial expectations is bounded by $T^{-1 + \oo(1)}$. The exponent $7/2$ has not been optimized. Lowering $\phi$ to $3$ or below appears to require a refinement of the Gaussian smoothing comparison in Lemma~\ref{an:smoothing} or of the removal of the perturbation in Step~5 of the proof of Proposition~\ref{cap:moment.measures}. Such refinements are beyond the scope of this paper.

The condition $\eta_{\sig}\log{X}/t \to \infty$ is used to approximate the logarithmic zeta field by the Euler field. The proof first replaces the sharp prime-power cutoff $p^{j} \leq X$ in $\mathcal{Q}$ by a smooth weight, producing a smoothed Dirichlet series. On a zero-free event whose probability tends to one, the contour shift in Lemma~\ref{analytic:lem:contour} compares $\log{\zeta}$ with this smoothed series, with error bounded by $\exp( - \eta_{\sig}\log{X}/2 + \OO_{\theta}(t))$. Lemma~\ref{analytic:lem:smoothing} then controls the error in replacing the smoothed series by $\mathcal{Q}$: the expectation over $\tau$ of the square of the largest absolute difference over $h \in I_{\theta}$ is bounded by $\exp( - 2\eta_{\sig}\log{X} + \OO_{\theta}(t))$. Thus, the assumption $\eta_{\sig}\log{X}/t \to \infty$ ensures that the negative terms in both exponents dominate the $\OO_{\theta}(t)$ terms, so both error bounds tend to zero. Together with control of the exceptional event, these estimates yield the mean-square field approximation in \eqref{analytic:eq:uniform-L2} of Proposition~\ref{analytic:prop:approximation}.

The Euler--Haar comparison requires an upper bound on $X$, whereas the approximation of the logarithmic zeta field by the Euler field requires $\eta_{\sig}\log{X}$ to grow faster than $t$. The restriction $\Phi > \phi + 1$ makes these requirements compatible with the prescribed upper bound on $\eta_{\sig}$. Indeed, the upper bounds in \eqref{eq:fixed.parameters} imply $\eta_{\sig}\log{X}/t \leq C_{\sig}t^{\Phi - \phi - 1}$. If $\Phi \leq \phi + 1$, the quantity $C_{\sig}t^{\Phi - \phi - 1}$ remains bounded, so the required divergence of $\eta_{\sig}\log{X}/t$ is impossible. Conversely, if $\Phi > \phi + 1$, the choices $X = T^{1/t^{\phi}}$ and $\eta_{\sig} = C_{\sig}t^{\Phi}/\log T$ satisfy the size bounds and give $\eta_{\sig}\log{X}/t = C_{\sig}t^{\Phi - \phi - 1} \to \infty$.

\subsection{Additional terminology and notation}


Let $\N \leqdef \{1,2,\ldots\}$ and $\N_0 \leqdef \{0,1,\ldots\}$. For $n \in \N$, $\mathrm{Id}_{n}$ is the $n \times n$ identity matrix, and $\bb{1}_{n}$ is the column vector of $n$ ones. For an event $E$, $\ind_{E}$ denotes its indicator, equal to one on $E$ and zero otherwise. The notation $\mathcal{L}(Y)$ denotes the probability law of a random variable or random vector $Y$. For a Borel set $\mathcal{J}$ in a Euclidean space, its Lebesgue measure is denoted by $|\mathcal{J}|$. For a finite collection $\mathcal{S}$ of sets, such as a partition into intervals, $|\mathcal{S}|$ denotes the number of its elements. For a finite signed measure $\mu$, its total variation norm is $\|\mu\|_{\mathrm{TV}} \leqdef \sup_{\|f\|_{\infty} \leq 1}|\int f\rd{\mu}|$, where the supremum is over measurable functions on its underlying space. With this convention, the total variation distance between probability measures is at most $2$. For real $a$ and $b$, write $a \wedge b \leqdef \min\{a,b\}$ and $a \vee b \leqdef \max\{a,b\}$; the positive and negative parts of $a$ are $a_{+} \leqdef \max\{a,0\}$ and $a_{-} \leqdef \max\{ - a,0\}$, respectively. For a positive comparison function $g$, the notation $f = \OO(g)$ or $f \ll g$ means that $|f|$ is bounded by a constant multiple of $g$, while $f = \oo(g)$ means that $f/g \to 0$. A subscript, as in $\OO_{\theta}$ or $\ll_{\theta}$, emphasizes the dependence of the implicit constant on the displayed parameter. For positive $f$ and $g$, the notation $f \asymp g$ means that their ratio is bounded above and below by positive constants, while $f \sim g$ means that $f/g \to 1$. Unless another limit is specified, all limits are taken as $T \to \infty$ along the chosen admissible parameter functions. Positive finite constants, usually denoted by $c$ or $C$, may change from one occurrence to the next. They may depend on the fixed parameters, including $\theta$, $\beta$, $\Phi$, $\phi$, and $C_{\sig}$, and on explicitly fixed cutoffs, but not on $T$ or on the particular values of $(\sig,X)$ in the stated range.

\section{Main results}\label{sec:main.results}

\subsection{Statements}

In this subsection, the limiting two-overlap law (Theorem~\ref{thm:limiting.two.overlap.distribution}), approximate GG identities (Theorem~\ref{thm:approximate.extended.GG.identities}), and limiting joint overlap law (Theorem~\ref{thm:limiting.joint.overlap.distribution}) are stated for expanding intervals ($\theta > 0$). Approximate ultrametricity (Proposition~\ref{cor:ultrametricity}) and the zeta--Euler comparison (Theorem~\ref{thm:gibbs.measure.comparison}) hold throughout $\theta > - 1$. The notation in Theorems~\ref{thm:limiting.two.overlap.distribution} and~\ref{thm:limiting.joint.overlap.distribution} is kept general so that their extension, together with that of Theorem~\ref{thm:approximate.extended.GG.identities}, to bounded or shrinking intervals ($ - 1 < \theta \leq 0$) can be stated with the same formulas in Conjecture~\ref{conj:bounded.overlap.laws} below. For $\theta > 0$, the simplifications are $q_{0} = 0$, $\delta_{|\theta \wedge 0|} = \delta_{0}$, and $\beta_{c}(\theta) = 2\sqrt{1 + \theta}$.

The following theorem establishes, for the logarithmic zeta field on expanding intervals to the right of the critical line, an analogue of the one-step replica symmetry breaking (1-RSB) two-overlap law conjectured by \citet[Conjecture~1]{arXiv:1706.08462}.

\begin{theorem}[Limiting two-overlap distribution]\label{thm:limiting.two.overlap.distribution}
Let $\theta > 0$ and $\beta > 0$. Then, as $T \to \infty$,
\[
\EE\mathcal{G}_{\beta,T}^{\times 2}\big[\rho_{T}(h_{1},h_{2}) \in \,\cdot\,\big]
\stackrel{\mathrm{weakly}}{\longrightarrow}
\begin{cases}
\delta_{|\theta \wedge 0|}, & \beta \leq \beta_{c}(\theta), \\
\displaystyle\frac{\beta_{c}(\theta)}{\beta}\delta_{|\theta \wedge 0|} + \left(1 - \frac{\beta_{c}(\theta)}{\beta}\right)\delta_{1}, & \beta > \beta_{c}(\theta).
\end{cases}
\]
Here $\delta_{y}$ denotes the Dirac probability measure at $y$, and weak convergence means convergence of integrals against every bounded continuous test function. The same conclusion holds for the Euler field $\Re{\mathcal{Q}}$, with $(\mathcal{G}_{\beta,T},\rho_{T})$ replaced by $(G_{\beta,T}^{\mathrm{A}},R_{T}^{\mathrm{A}})$.
\end{theorem}

Theorem~\ref{thm:limiting.two.overlap.distribution} describes the overlap between two positions under the annealed sampling law: first $\tau$ is sampled uniformly on $[T,2T]$, and then two positions are sampled independently from $\mathcal{G}_{\beta,T}$ conditionally on $\tau$. Since the theorem concerns expanding intervals ($\theta > 0$), the smaller limiting overlap is $|\theta \wedge 0| = 0$. If $\beta > \beta_{c}(\theta)$, that is, at low temperature, the overlap converges in distribution under this annealed law to a random variable taking the values $0$ and $1$ with probabilities $\beta_{c}(\theta)/\beta$ and $1 - \beta_{c}(\theta)/\beta$, respectively. In the schematic prime-scale picture of Figure~\ref{fig:prime.scale.regimes}, limiting overlap $0$ corresponds to positions sharing no positive fraction of the logarithmic prime scales, whereas limiting overlap $1$ corresponds to positions sharing asymptotically all of those scales, without requiring the positions to coincide. These alternatives are interpreted as sampling positions from different groups or from the same group, respectively, with the limiting joint grouping law specified in Theorem~\ref{thm:limiting.joint.overlap.distribution}. For every fixed $\e \in (0,1/2)$, the annealed probability that the overlap belongs to $[\e,1 - \e]$ therefore tends to zero. If $0 < \beta \leq \beta_{c}(\theta)$, that is, at critical or high temperature, the overlap between two sampled positions converges to $0$ under the annealed law. For the toy model $F_{h}^{\mathrm{rand}}$ in \eqref{eq:toy.model.zeta} on a bounded interval, \citet[Theorem~1]{arXiv:1706.08462} established the analogous low-temperature two-overlap law, with masses $2/\beta$ at $0$ and $1 - 2/\beta$ at $1$ for $\beta > 2$.

The next result gives approximate GG identities, which express consistency relations among the annealed joint distributions of overlaps of finitely many sampled positions. These identities were introduced by \citet{MR1662161} in the study of the Sherrington--Kirkpatrick model, whose Hamiltonian is a Gaussian random field indexed by spin configurations in $\{ - 1,1\}^{N}$. Analogous identities were established for generalized random energy models (GREM), whose Gaussian Hamiltonians have hierarchical covariance, by \citet{MR2070334,MR2070335}, and for a nonhierarchical log-correlated Gaussian field related to the logarithmic random energy model by \citet{MR3211001}. The GG approach was applied to the two-dimensional discrete Gaussian free field (GFF) by \citet{MR3354619} and generalized to the scale-inhomogeneous setting by \citet{MR3731796}. Further instances include the branching random walk with Gaussian increments on a binary tree, studied by \citet{MR3539644}, and the toy model of $\log|\zeta|$, $F_{h}^{\mathrm{rand}}$, discussed in Subsubsection~\ref{sec:randomized.zeta}; see \citet{MR3841407}. For mixed $p$-spin models, GG identities were established in the Gaussian setting by \citet{MR2600075,MR3052333} and extended to non-Gaussian disorder by \citet{MR3901654}.

\begin{theorem}[Approximate GG identities]\label{thm:approximate.extended.GG.identities}
Let $\theta > 0$, $\beta > 0$, $n \in \N$, and $k \in \{1,\ldots,n\}$. For every measurable function $\Psi_{T}:I_{\theta}^{\times n} \to \R$ (possibly dependent on $T$) with $\|\Psi_{T}\|_{\infty} \leq 1$, we have, as $T \to \infty$,
\[
\begin{aligned}
\EE\mathcal{G}_{\beta,T}^{\times(n+1)}[\rho_{T}(h_{k},h_{n+1})\Psi_{T}(\bb{h})]
&= \frac{1}{n}\EE\mathcal{G}_{\beta,T}^{\times 2}[\rho_{T}(h_{1},h_{2})]\EE\mathcal{G}_{\beta,T}^{\times n}[\Psi_{T}(\bb{h})] \\
&\qquad + \frac{1}{n}\sum_{\substack{1 \leq \ell \leq n \\ \ell \neq k}}\EE\mathcal{G}_{\beta,T}^{\times n}[\rho_{T}(h_{k},h_{\ell})\Psi_{T}(\bb{h})] + \oo(1),
\end{aligned}
\]
where the $\oo(1)$ term is uniform over $\Psi_{T}$. The same identity holds for the Euler field $\Re{\mathcal{Q}}$, with $(\mathcal{G}_{\beta,T},\rho_{T})$ replaced by $(G_{\beta,T}^{\mathrm{A}},R_{T}^{\mathrm{A}})$.
\end{theorem}

\begin{figure}[!ht]
\centering
\includegraphics[width = \textwidth]{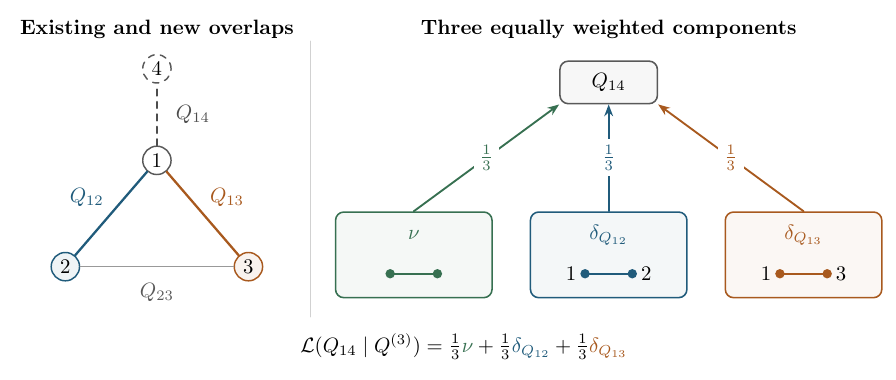}
\caption{Sampling interpretation of the GG identities for $n = 3$ and $k = 1$. Here $(Q_{ij})_{1 \leq i,j \leq 4}$ has the limiting overlap law under $\EE\mathcal{G}_{\beta,T}^{\times 4}$, $Q^{(3)} \leqdef (Q_{ij})_{1 \leq i,j \leq 3}$, and $\nu \leqdef \mathcal{L}(Q_{12})$ is the two-overlap distribution. Conditionally on the overlaps among the first three replicas, $Q_{14}$ is distributed as an equal mixture of a fresh draw from $\nu$, independent of $Q^{(3)}$, and the two already observed values $Q_{12}$ and $Q_{13}$. The colored edges identify these two observed overlaps, while the isolated pair represents a fresh draw from the two-overlap distribution. The point-mass components reproduce overlap values and do not identify sampled positions. For general $n$, the mixture assigns weight $1/n$ to $\nu$ and to each point mass $\delta_{Q_{k\ell}}$, $\ell \in \{1,\ldots,n\} \setminus \{k\}$.}\label{fig:GG.sampling}
\end{figure}

Theorem~\ref{thm:approximate.extended.GG.identities} gives a direct sampling interpretation under the limiting annealed law. Conditional on $\tau$, let $h_{1},\ldots,h_{n + 1}$ be independent with common law $\mathcal{G}_{\beta,T}$, and fix $k \in \{1,\ldots,n\}$. Write $Q_{ij}$ for the limiting overlap corresponding to $\rho_{T}(h_{i},h_{j})$, whose joint law is given by Theorem~\ref{thm:limiting.joint.overlap.distribution}. Conditionally on the overlaps $(Q_{ij})_{1 \leq i,j \leq n}$ among the first $n$ sampled positions, the new overlap $Q_{k,n + 1}$ is distributed as a mixture with weight $1/n$ on a fresh draw from the limiting two-overlap distribution, independent of those observed overlaps, and weight $1/n$ on each of the $n - 1$ already observed values $Q_{k\ell}$, $\ell \neq k$. Figure~\ref{fig:GG.sampling} illustrates this mixture for $n = 3$ and $k = 1$.

In spin glass theory, \citet{MR2999044} proved that the GG identities imply ultrametricity for random probability measures on the unit ball of a separable Hilbert space, without imposing any restriction on the two-overlap distribution; an elementary argument when only finitely many overlap values are possible was given by \citet{MR2825947}. For sequences of random probability measures on the unit ball of a separable Hilbert space, \citet{MR3628881} proved that the corresponding approximate GG identities and weak convergence of the averaged two-overlap distribution yield approximately ultrametric decompositions into nested clusters at any finite collection of admissible overlap levels, with suitably ordered cluster masses converging in law to the weights of Ruelle probability cascades. This gives a universal description of the limiting cluster weights, determined by the limiting two-overlap distribution, with applications at sufficiently low temperature to generic mixed $p$-spin models on the hypercube and the generalized random energy model, among others.

\begin{proposition}[Approximate ultrametricity]\label{cor:ultrametricity}
For every fixed $\theta > - 1$, there exists $C_{\theta} < \infty$ such that, for every sufficiently large $T$ and every $\bb{h} \in I_{\theta}^{\times 3}$,
\begin{equation}\label{eq:main.ultrametricity}
\rho_{T}(h_{1},h_{2}) \geq \min\{\rho_{T}(h_{1},h_{3}),\rho_{T}(h_{2},h_{3})\} - C_{\theta}\frac{\log{t}}{t}.
\end{equation}
Consequently, for every $\beta > 0$ and $\e > 0$ and for all sufficiently large $T$,
\[
\EE\mathcal{G}_{\beta,T}^{\times 3}\big[\rho_{T}(h_{1},h_{2}) \geq \min\{\rho_{T}(h_{1},h_{3}),\rho_{T}(h_{2},h_{3})\} - \e\big] = 1.
\]
Both statements also hold for the Euler field $\Re{\mathcal{Q}}$, with $(\mathcal{G}_{\beta,T},\rho_{T})$ replaced by $(G_{\beta,T}^{\mathrm{A}},R_{T}^{\mathrm{A}})$.
\end{proposition}

\begin{remark}
For $\log|\zeta|$, approximate ultrametricity follows directly from the triangle inequality and the uniform approximation of $\rho_{T}(h,h')$ by $\log_{+}((\log T) \wedge |h - h'|^{-1})/t$ in \eqref{analytic:eq:overlap-geometry} of Theorem~\ref{analytic:geometry}, without using the approximate GG identities.
\end{remark}

\begin{figure}[!t]
\centering
\includegraphics[width = \textwidth]{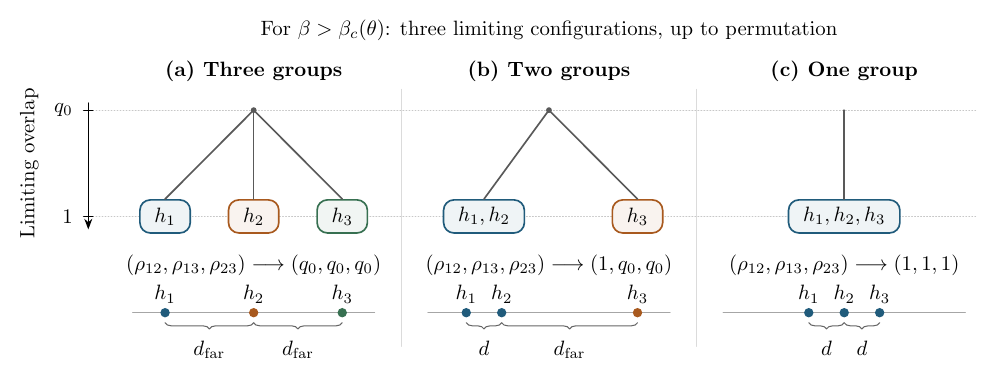}
\caption{Limiting group configurations, up to permutation, of three positions sampled under the annealed law $\EE\mathcal{G}_{\beta,T}^{\times 3}$, when $\theta > 0$ and $\beta > \beta_{c}(\theta)$. Here $\rho_{ij} \leqdef \rho_{T}(h_{i},h_{j})$. Each enclosure represents one group: the limiting overlap is $1$ within a group and $q_{0}$ between groups. The lower sketches are not to scale. The labeled separations $d_{\mathrm{far}} \asymp (\log T)^{\theta}$ and $d \leq (\log T)^{-1}$ give examples yielding overlaps tending to $q_{0}$ and $1$, respectively, as $T \to \infty$. Other separation scales can yield the same limits; in particular, limiting overlap $1$ does not require coincident positions. When $0 < \beta \leq \beta_{c}(\theta)$, only configuration (a) occurs in the limit.}\label{fig:ultrametricity}
\end{figure}

Figure~\ref{fig:ultrametricity} displays the three possible group configurations of three sampled positions, up to permutation, at low temperature ($\beta > \beta_{c}(\theta)$). Together with the two-overlap law in Theorem~\ref{thm:limiting.two.overlap.distribution}, the GG identities determine the limiting joint law of the group assignments for any fixed number of sampled positions: positions in the same group have limiting overlap $1$, and positions in different groups have limiting overlap $q_{0} = |\theta \wedge 0|$ (here, $q_{0} = 0$ since $\theta > 0$). When $\beta \leq \beta_{c}(\theta)$, distinct sampled indices have only the smaller limiting overlap, as in panel (a). At low temperature, the limiting group assignments can be described by random weights representing the groups' shares of the total Gibbs mass: conditional on these weights, sampled positions are assigned independently to groups with probabilities given by the weights. These weights are described by a Poisson--Dirichlet process, defined next.

A Poisson--Dirichlet process $\bb{\xi} \leqdef (\xi_{j})_{j \in \N}$ with parameter $\lambda \in (0,1)$ is a random decreasing sequence of nonnegative weights with sum one, constructed as follows. Let $(w_{j})_{j \in \N}$ be the atoms of a Poisson random measure on $(0,\infty)$ with intensity $\lambda v^{-1-\lambda}\rd{v}$: counts in disjoint measurable sets are independent, and the count in a set $A$ is Poisson with mean $\int_{A}\lambda v^{-1-\lambda}\rd{v}$ whenever this integral is finite. The sum $\sum_{k}w_{k}$ is finite and strictly positive almost surely. Then $\bb{\xi}$ is the decreasing rearrangement of $(w_{j}/\sum_{k}w_{k})_{j \in \N}$. Its weights sum to one almost surely. Write $\EE_{\bb{\xi}}$ for expectation with respect to these weights.

\begin{theorem}[Limiting joint distribution of the overlaps]\label{thm:limiting.joint.overlap.distribution}
Let $\theta > 0$, $\beta > 0$, and $n \in \N$. If $\beta \leq \beta_{c}(\theta)$, define the deterministic limiting overlap matrix
\[
Q^{\mathrm{high}} \leqdef (1 - q_{0})\mathrm{Id}_{n} + q_{0}\bb{1}_{n}\bb{1}_{n}^{\mathsf{T}}.
\]
If $\beta > \beta_{c}(\theta)$, let $\bb{\xi}$ have Poisson--Dirichlet parameter $\beta_{c}(\theta)/\beta$, sample conditionally independent group labels $K_{1},\ldots,K_{n}$ with probabilities $\bb{\xi}$, and define the random limiting overlap matrix
\[
Q^{\mathrm{low}} \leqdef \left(q_{0} + (1 - q_{0})\ind_{\{K_{\ell} = K_{\ell'}\}}\right)_{1 \leq \ell,\ell' \leq n}.
\]
For every continuous function $\psi:[ - 1,1]^{n \times n} \to \R$,
\begin{equation}\label{eq:main.joint.overlap}
\lim_{T \to \infty}\EE\mathcal{G}_{\beta,T}^{\times n}\left[\psi\left((\rho_{T}(h_{\ell},h_{\ell'}))_{1 \leq \ell,\ell' \leq n}\right)\right]
= \begin{cases}
\psi(Q^{\mathrm{high}}), & \beta \leq \beta_{c}(\theta), \\
\EE_{\bb{\xi}}\big[\EE[\psi(Q^{\mathrm{low}}) \mid \bb{\xi}]\big], & \beta > \beta_{c}(\theta).
\end{cases}
\end{equation}
The same limit holds for the Euler field $\Re{\mathcal{Q}}$, with $(\mathcal{G}_{\beta,T},\rho_{T})$ replaced by $(G_{\beta,T}^{\mathrm{A}},R_{T}^{\mathrm{A}})$.
\end{theorem}

Both limiting matrices $Q^{\mathrm{high}}$ and $Q^{\mathrm{low}}$ have diagonal entries equal to one. At critical or high temperature ($\beta \leq \beta_{c}(\theta)$), every off-diagonal entry equals $q_{0} = |\theta \wedge 0|$ (here, $q_{0} = 0$ since $\theta > 0$). At low temperature ($\beta > \beta_{c}(\theta)$), the theorem has a concrete sampling interpretation. Conditional on $\bb{\xi}$, choose $K_{1},\ldots,K_{n}$ independently with $\PP(K_{\ell} = j \mid \bb{\xi}) = \xi_{j}$. The limiting overlap of samples indexed by $\ell$ and $\ell'$ is $1$ when $K_{\ell} = K_{\ell'}$ and $0$ otherwise. Thus Theorem~\ref{thm:limiting.two.overlap.distribution} gives the law for one pair, while Theorem~\ref{thm:limiting.joint.overlap.distribution} gives the joint law for every fixed number of sampled points.

For comparison, \citet{GangulyNam2025} proved Poisson--Dirichlet universality at low temperature for a broad class of discrete logarithmically correlated Gaussian fields, using sharp estimates for the sizes of near-maximum level sets. Their result does not apply to the logarithmic zeta, Euler, or Haar fields considered herein, none of which are Gaussian.

To establish Theorems~\ref{thm:limiting.two.overlap.distribution}, \ref{thm:approximate.extended.GG.identities}, and~\ref{thm:limiting.joint.overlap.distribution} for the logarithmic zeta field $\log|\zeta|$, it suffices to prove the corresponding conclusions for the Euler field $\Re\mathcal{Q}$, as done in Section~\ref{sec:expanding}. At low temperature, the Euler conclusions are obtained by transferring the Haar overlap laws and approximate GG identities through a comparison of annealed sampling laws, despite the dependence of the arithmetic phases $p^{-\ii\tau}$ on the same random variable $\tau$. The following comparison theorem then transfers the Euler conclusions to $\log|\zeta|$.

\begin{theorem}[Zeta--Euler comparison]\label{thm:gibbs.measure.comparison}
For every fixed $\theta > - 1$, $\beta > 0$, and $n \in \N$, the following convergences hold as $T \to \infty$:
\begin{equation}\label{eq:main.comparison}
\sup_{h,h' \in I_{\theta}}|\rho_{T}(h,h') - R_{T}^{\mathrm{A}}(h,h')| \to 0, \qquad
\EE\left[\int_{I_{\theta}^{\times n}}\left|\frac{\rd{\mathcal{G}_{\beta,T}^{\times n}}}{\rd{\bb{h}}}(\bb{h}) - \frac{\rd{(G_{\beta,T}^{\mathrm{A}})^{\times n}}}{\rd{\bb{h}}}(\bb{h})\right|\rd{\bb{h}}\right] \to 0.
\end{equation}
Here $\rd{\mathcal{G}_{\beta,T}^{\times n}}/\rd{\bb{h}}$ and $\rd{(G_{\beta,T}^{\mathrm{A}})^{\times n}}/\rd{\bb{h}}$ denote the densities of the $n$-fold products of the zeta and Euler Gibbs measures, respectively, with respect to the Lebesgue measure on $I_{\theta}^{\times n}$. The second limit in \eqref{eq:main.comparison} implies convergence of the annealed product-measure differences: with the supremum taken over measurable functions $\Psi_{T}:I_{\theta}^{\times n} \to \R$ (possibly dependent on $T$), we have
\begin{equation}\label{eq:annealed.product.diff}
\lim_{T \to \infty} \sup_{\|\Psi_{T}\|_{\infty} \leq 1} \Big|\EE\,\mathcal{G}_{\beta,T}^{\times n}\big[\Psi_{T}(h_{1},\ldots,h_{n})\big] - \EE\,(G_{\beta,T}^{\mathrm{A}})^{\times n}\big[\Psi_{T}(h_{1},\ldots,h_{n})\big]\Big| = 0.
\end{equation}
In particular, taking $\Psi_{T}(h_{1},\ldots,h_{n}) \leqdef \psi\big((\rho_{T}(h_{\ell},h_{\ell'}))_{1 \leq \ell,\ell' \leq n}\big)$ for any continuous function $\psi : [ - 1,1]^{n \times n} \to \R$ with $\|\psi\|_{\infty} \leq 1$, and using the first limit in \eqref{eq:main.comparison}, shows that the joint laws of the respective overlap arrays under $\EE\mathcal{G}_{\beta,T}^{\times n}$ and $\EE(G_{\beta,T}^{\mathrm{A}})^{\times n}$ have the same subsequential limits, for every fixed $n$.
\end{theorem}

The extension of the overlap distributions and approximate GG identities to bounded or shrinking intervals ($ - 1 < \theta \leq 0$) is conjectural, as stated next.

\begin{conjecture}[Overlap laws on bounded or shrinking intervals]\label{conj:bounded.overlap.laws}
For every $ - 1 < \theta \leq 0$ and $\beta > 0$, the conclusions of Theorems~\ref{thm:limiting.two.overlap.distribution}, \ref{thm:approximate.extended.GG.identities}, and~\ref{thm:limiting.joint.overlap.distribution} hold for both the logarithmic zeta field $\log|\zeta|$, with $(\mathcal{G}_{\beta,T},\rho_{T})$, and the Euler field $\Re{\mathcal{Q}}$, with $(G_{\beta,T}^{\mathrm{A}},R_{T}^{\mathrm{A}})$.
\end{conjecture}

\subsection{Proof structure and deduction of the main results}\label{sec:outline}

The main results follow from the analytic comparison and uniform overlap geometry in Theorems~\ref{analytic:comparison} and~\ref{analytic:geometry}, proved in Section~\ref{sec:analytic}, and the Euler overlap conclusions in Theorem~\ref{expand:main}, proved in Section~\ref{sec:expanding}. The deductions are given below. The analytic comparison uses the field approximations \eqref{analytic:eq:uniform-probability} and \eqref{analytic:eq:uniform-L2}: the first controls the ratio of the zeta and Euler Gibbs densities, and the second controls the difference between the overlaps $\rho_{T}$ and $R_{T}^{\mathrm{A}}$.

Section~\ref{sec:expanding} treats expanding intervals ($\theta > 0$) through separate arguments for the two temperature regimes. For $0 < \beta \leq \beta_{c}(\theta)$, estimates showing that the Euler Gibbs mass spreads over many unit intervals imply vanishing off-diagonal overlaps, which yield the overlap laws and approximate GG identities. For $\beta > \beta_{c}(\theta)$, these conclusions are transferred from the Haar field by comparing the annealed sampling laws of the Euler and Haar Gibbs measures. The required Haar results are taken from \citet[Proposition~3.1]{OuimetHaar2026}. The main step in this transfer is to pass from arithmetic comparisons of finite collections of field values to comparisons of Gibbs sampling laws, using mass estimates and a normalization argument. Each remaining section states its required conclusions before developing the estimates used to prove them.

\begin{proof}[Proof of Theorem~\ref{thm:gibbs.measure.comparison}]
Equations \eqref{analytic:eq:TV} and \eqref{analytic:eq:rho-R} in Theorem~\ref{analytic:comparison} give \eqref{eq:main.comparison}. The absolute difference of the annealed expectations of any measurable function $\Psi_{T}$ with $\|\Psi_{T}\|_{\infty} \leq 1$ is bounded by the expected integral in the second limit of \eqref{eq:main.comparison}, which proves \eqref{eq:annealed.product.diff}.
\end{proof}

\begin{proof}[Proofs of Theorems~\ref{thm:limiting.two.overlap.distribution}, \ref{thm:approximate.extended.GG.identities}, and~\ref{thm:limiting.joint.overlap.distribution}]
For $\theta > 0$, Theorem~\ref{expand:main} supplies the analogous Euler field conclusions. For each fixed $n \in \N$, \eqref{eq:annealed.product.diff} gives an $\oo(1)$ difference between annealed expectations under $\mathcal{G}_{\beta,T}^{\times n}$ and $(G_{\beta,T}^{\mathrm{A}})^{\times n}$, uniformly over $\|\Psi_{T}\|_{\infty} \leq 1$. Together with the first limit in \eqref{eq:main.comparison} and the bounds $\|\rho_{T}\|_{\infty},\|R_{T}^{\mathrm{A}}\|_{\infty} \leq 1$, this transfers every GG term, including products of expectations, with the same uniform $\oo(1)$ error. The two-overlap and joint overlap distributions transfer by uniform continuity of continuous functions on $[ - 1,1]^{n \times n}$.
\end{proof}

\begin{proof}[Proof of Proposition~\ref{cor:ultrametricity}]
The result follows from \eqref{analytic:eq:deterministic-ultrametricity} in Theorem~\ref{analytic:geometry}.
\end{proof}

\begin{figure}[H]
\centering
\begin{tikzpicture}[x = 1cm, y = 1cm, every node/.style = {font = \small, align = center}, result/.style = {draw, rounded corners = 2pt, fill = black!15, inner sep = 5pt}, input/.style = {result, fill = white}, dependency/.style = {-{Stealth[length = 2mm]}, line width = 0.5pt}]
\path[fill = black!8, draw = black!20, rounded corners = 4pt] ( - 2.05,1.25) rectangle (2.05, - 1.6);
\path[fill = black!8, draw = black!20, rounded corners = 4pt] (5.9,1.1) rectangle (13.8, - 1.35);
\node[font = \small] at (0, - 1.32) {Section~\ref{sec:expanding}};
\node[font = \small, anchor = west] at (6.1, - 1.07) {Section~\ref{sec:analytic}};

\node[input] (eulerthm) at (0,0) {Theorem~\ref{expand:main} \\Euler overlap laws \\and approximate \\GG identities};
\node[result] (comparison) at (4,0) {Theorem~\ref{thm:gibbs.measure.comparison} \\Zeta--Euler \\comparison};
\node[input] (analyticthm) at (8,0) {Theorem~\ref{analytic:comparison} \\Gibbs measure and \\overlap comparison};
\node[input] (geometry) at (12,0) {Theorem~\ref{analytic:geometry} \\Uniform overlap \\geometry};
\node[result] (pair) at (0, - 3.15) {Theorem~\ref{thm:limiting.two.overlap.distribution} \\Limiting two-overlap \\distribution};
\node[result] (gg) at (4, - 3.15) {Theorem~\ref{thm:approximate.extended.GG.identities} \\Approximate \\GG identities};
\node[result] (joint) at (8, - 3.15) {Theorem~\ref{thm:limiting.joint.overlap.distribution} \\Limiting joint overlap \\distribution};
\node[result] (ultrametricity) at (12, - 3.15) {Proposition~\ref{cor:ultrametricity} \\Approximate \\ultrametricity};

\coordinate (transfer) at (4, - 1.675);
\draw[dependency] (analyticthm.west) to (comparison.east);
\draw[line width = 0.5pt] (eulerthm.south east) to[out = -35, in = 180] (transfer);
\draw[line width = 0.5pt] (comparison.south) -- (transfer);
\fill (transfer) circle[radius = 1.6pt];
\draw[dependency] (transfer) to (pair.north east);
\draw[dependency] (transfer) to (gg.north);
\draw[dependency] (transfer) to (joint.north west);
\draw[dependency] (geometry.south) to (ultrametricity.north);
\end{tikzpicture}
\caption{Dependencies among the main results. Gray boxes identify the main results stated in Section~\ref{sec:main.results}. White boxes identify intermediate results, grouped by section in the shaded regions. The junction combines the Euler conclusions with the zeta--Euler comparison.}\label{fig:main.local.dependencies}
\end{figure}
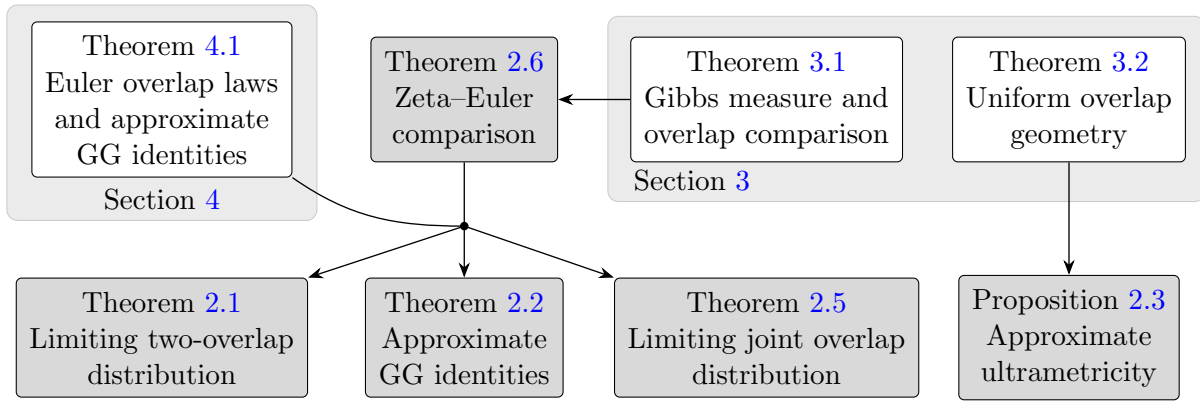

Figure~\ref{fig:main.local.dependencies} summarizes these dependencies. Each arrow points from an input to the result or proof step that uses it.

\subsection{Relations to other random fields}\label{sec:relations.other.models}

Randomized prime sums and log-modulus of characteristic polynomials of random unitary matrices provide two fundamental models for studying the extreme values of $\log|\zeta|$ and the role of its logarithmic correlations.

\subsubsection{Randomized prime sums}\label{sec:randomized.zeta}

Let $(U_{p})_{p}$ be independent random variables uniformly distributed on the complex unit circle. In the logarithmic prime-scale approximation \eqref{eq:log.zeta.approximation}, the arithmetic phases $p^{-\ii\tau}$ depend on the same random shift $\tau$ and are deterministic once it is fixed. Following \citet{arXiv:1304.0677}, replacing them with independent Haar phases gives the basic model
\begin{equation}\label{eq:toy.model.zeta}
F_{h}^{\mathrm{rand}} \leqdef \Re\sum_{p \leq T}\frac{U_{p}}{p^{1/2 + \ii h}}, \qquad h \in I_{\theta},
\end{equation}
for studying large values of $\log|\zeta|$ and its logarithmic correlations. \citet{MR3619786} proved the first two orders of the maximum of this model on $[0,1]$ through a second-moment argument and estimates across logarithmic prime scales. Combining the Gaussian coupling of \citet{arXiv:1604.08378,MR4164452} with the tightness theorem of \citet{MR3272323}, one can deduce tightness of the maximum after recentering by $\log \log T - (3/4)\log \log \log T$. For the Gaussian-prime counterpart $\Re\sum_{p \leq T}g_{p}p^{-1/2 - \ii h}$, where the coefficients $g_{p}$ are independent complex Gaussian variables with independent centered real and imaginary parts of variance $1/2$, \citet{ADH_2021} determined the first two orders of the maximum over $I_{\theta}$ for fixed $\theta > 0$ and for $\theta = (\log \log T)^{-\alpha}$ with a fixed $\alpha \in (0,1)$. In the latter regime, \citet{Chang2024} proved tightness of the recentered maximum for every $\alpha \in (0,1)$, together with a sharp estimate for the critical moment when $\alpha \in (1/2,1)$.

For a centered Gaussian approximation $F_{h}^{\mathrm{G}}$, normalization by the pointwise expectation gives
\[
\frac{e^{\beta F_{h}^{\mathrm{G}}}}{\EE[e^{\beta F_{h}^{\mathrm{G}}}]}\rd{h} = \exp\left\{\beta F_{h}^{\mathrm{G}} - \frac{\beta^{2}}{2}\Var(F_{h}^{\mathrm{G}})\right\}\rd{h}.
\]
Suitable nonzero limits of these measures for logarithmically correlated Gaussian fields are called Gaussian multiplicative chaos. For $F_{h}^{\mathrm{rand}}$ on $[0,1]$ and $0 < \beta < 2$, \citet{arXiv:1604.08378,MR4164452} proved that the measures $e^{\beta F_{h}^{\mathrm{rand}}}\rd{h}/\EE[e^{\beta F_{h}^{\mathrm{rand}}}]$ converge almost surely in the weak topology to a nontrivial non-Gaussian multiplicative-chaos measure. \citet{arXiv:1906.08573} identified its total mass as the limit in probability of the Lebesgue measure of points in $[0,1]$ where $F_{h}^{\mathrm{rand}} > (\beta/2) \log \log T$, divided by its expectation.

For $F_{h}^{\mathrm{rand}}$ on $[0,1]$ at low temperature ($\beta > 2$), \citet{arXiv:1706.08462} identified the limiting two-overlap distribution recalled after Theorem~\ref{thm:limiting.two.overlap.distribution}. In the same regime, \citet{MR3841407} established the GG identities for the limiting overlap array and described the joint overlap distributions through Poisson--Dirichlet weights of parameter $2/\beta$. These proofs used the independence of the phases $U_{p}$ to apply approximate integration by parts one prime at a time while holding the remaining phases fixed.

This argument does not apply directly to the Euler field, since its prime phases depend on the same random shift $\tau$. For expanding intervals ($\theta > 0$) at low temperature ($\beta > \beta_{c}(\theta)$), the proof instead compares the joint laws of positions sampled from the Euler and Haar Gibbs measures after averaging over their respective randomness. Theorem~\ref{cap:sampling} establishes this comparison uniformly over measurable functions of any fixed number of sampled positions with supremum norm at most one, even when these functions depend on $T$. Together with a uniform comparison of the Euler and Haar overlaps, this transfers the Haar overlap laws and approximate GG identities proved by \citet{OuimetHaar2026} to the Euler field.

\subsubsection{Characteristic polynomials of random unitary matrices}\label{sec:CUE.field}

For $N \in \N$, let $U_{N}$ be a random $N \times N$ unitary matrix with Haar distribution on the unitary group. If its eigenangles are $\vartheta_{1},\ldots,\vartheta_{N} \in [0,2\pi)$, so its eigenvalues are $e^{\ii\vartheta_{1}},\ldots,e^{\ii\vartheta_{N}}$, then the log-modulus of its characteristic polynomial is
\[
F_{h}^{\mathrm{CUE}} \leqdef \log|\det(\mathrm{Id}_{N} - U_{N}e^{-\ii h})| = \sum_{j = 1}^{N}\log|1 - e^{\ii(\vartheta_{j} - h)}|, \qquad h \in [0,2\pi].
\]
The Fourier expansion of $\log|1 - e^{\ii x}|$ gives
\[
F_{h}^{\mathrm{CUE}} = - \sum_{k = 1}^{\infty}\frac{1}{k}\Re\!\left(e^{-\ii kh}\operatorname{Tr}(U_{N}^{k})\right),
\]
with pointwise convergence away from the eigenangles. The covariance between $F_{h}^{\mathrm{CUE}}$ and $F_{h'}^{\mathrm{CUE}}$ grows logarithmically as the circular distance between $h$ and $h'$ decreases from order one to the microscopic scale $N^{-1}$. This parallels the covariance of the prime-sum model under the correspondence $N \sim \log T$, which identifies the microscopic scales $N^{-1}$ and $(\log T)^{-1}$.

The leading order of the maximum over the unit circle and the freezing transition of the limiting free energy were proved by \citet{MR3594368}. The second order of the maximum was obtained by \citet{MR3848227}, and tightness after recentering by $\log{N} - (3/4)\log{\log{N}}$ was established by \citet[Theorem~1.2]{MR3848391}. For $\beta_{\mathrm{circ}} > 0$, the circular ensemble with parameter $\beta_{\mathrm{circ}}$ is the law of angles $(\vartheta_{1},\ldots,\vartheta_{N}) \in [0,2\pi)^{N}$ with joint density proportional to $\prod_{1 \leq j < k \leq N}|e^{\ii\vartheta_{j}} - e^{\ii\vartheta_{k}}|^{\beta_{\mathrm{circ}}}$. The case $\beta_{\mathrm{circ}} = 2$ gives the eigenangle law of $U_{N}$. For these circular ensembles, \citet{PaquetteZeitouni2025} gave a detailed description of the extremal landscape and established a randomly shifted Gumbel law for the suitably recentered maximum.

Moments of the partition function of $F_{h}^{\mathrm{CUE}}$ provide another point of comparison. For positive integer moment orders and positive even integer powers of the characteristic-polynomial modulus, \citet{MR4019917} proved that these moments are polynomials in the matrix size. Conditional on the shifted-moment conjecture of \citet{MR2149530}, \citet{MR4231530} derived corresponding asymptotics for zeta moments of moments with the same integer parameters. \citet{arXiv:2104.07403} examined Radziwi{\l\l}'s conjecture that the zeta moment constant predicted by \citet{KeatingSnaith2000} appears as a multiplicative correction to the Gaussian upper-tail approximation for $\PP(\log|\zeta(1/2 + \ii\tau)| > k\log\log T)$, where $k > 0$ is fixed. They also investigated the effect of this correction on the maximum of $\log|\zeta|$ over expanding intervals.

For $0 < \beta_{\mathrm{CUE}} < 2$, consider the finite random measures
\[
\frac{|\det(\mathrm{Id}_{N} - U_{N}e^{-\ii h})|^{\beta_{\mathrm{CUE}}}}{\EE[|\det(\mathrm{Id}_{N} - U_{N}e^{-\ii h})|^{\beta_{\mathrm{CUE}}}]}\frac{\rd{h}}{2\pi},
\]
where the expectation is taken with respect to $U_{N}$. Convergence in distribution of these measures to Gaussian multiplicative chaos was proved for $0 < \beta_{\mathrm{CUE}} < \sqrt{2}$ by \citet{MR3407221} and extended to $\sqrt{2} \leq \beta_{\mathrm{CUE}} < 2$ by \citet{MR4105514}. Throughout $0 < \beta_{\mathrm{CUE}} < 2$, the same limit in distribution can be recovered by restricting $\rd{h}/(2\pi)$ to the thick-point set $\{h:F_{h}^{\mathrm{CUE}} \geq (\beta_{\mathrm{CUE}}/2)\log{N}\}$ and dividing by the corresponding one-point exceedance probability, as shown by \citet{JunnilaLambertWebb2024}. The work of \citet{LambertNajnudel2026} extended multiplicative-chaos convergence to circular $\beta_{\mathrm{circ}}$-ensembles for modulus powers $\gamma$ satisfying $\gamma > - 1$ and $|\gamma| < \sqrt{2\beta_{\mathrm{circ}}}$, including negative powers, and also treated the eigenvalue counting field.

These subcritical Gaussian multiplicative chaos limits are almost surely non-atomic. Dividing the finite measures above by their random total mass gives the corresponding Gibbs probability measures. A low-temperature description of the CUE Gibbs measure through Poisson--Dirichlet overlap arrays and approximate GG identities does not appear to have been established in the literature. Deducing such conclusions from an extremal-process description on bounded height windows requires control of the normalized Gibbs masses and of the contribution from field values below the extremal range.

The analogy between the logarithmic zeta field and the CUE field also depends on whether the spatial interval shrinks, stays fixed, or expands. Under the correspondence $N \sim \log T$ used by \citet{MR3151088}, the regime $ - 1 < \theta < 0$ corresponds to mesoscopic arcs of length comparable to $N^{\theta}$, which shrink while typically containing a diverging number of eigenangles, and $\theta = 0$ corresponds to fixed arcs. Thus the bounded and shrinking regimes $ - 1 < \theta \leq 0$ have direct random-matrix analogues. For $\theta > 0$, the corresponding angular interval would have diverging length, so there is no direct analogue for the characteristic polynomial of a single random unitary matrix: its $2\pi$-periodicity simply repeats the same field values. An expanding zeta interval instead contains a growing number of distinct bounded subintervals, whose competition for Gibbs mass affects both the critical inverse temperature and the limiting overlap law.

\section{Transfer from \texorpdfstring{$\Re \mathcal{Q}$}{Re Q} to \texorpdfstring{$\log|\zeta|$}{log|zeta|} for \texorpdfstring{$\theta > -1$}{theta > -1}}\label{sec:analytic}

Throughout this section, we assume that $\theta > - 1$ and that $\sig$ and $X$ satisfy \eqref{eq:fixed.parameters}. The conclusions for the Euler field $\Re \mathcal{Q}$ imply the corresponding conclusions for the logarithmic zeta field $\log|\zeta|$ once two comparisons are available: the normalized Gibbs densities must be close, and their deterministic overlap functions $R_{T}^{\mathrm{A}}$ and $\rho_{T}$ must be close. Theorem~\ref{analytic:comparison} supplies precisely these inputs, and Theorem~\ref{analytic:geometry} approximates both overlap functions uniformly by a logarithmic distance kernel, which yields Proposition~\ref{cor:ultrametricity}. The deductions of both theorems are given first, so that the approximation and second-moment estimates proved afterward have an explicit purpose.

Recall the displacement $\eta_{\sig}$ from the critical line, and define the spatial half-length $L_{T}$ and the contour truncation height $H_{T}$:
\[
\eta_{\sig} = \sig - \frac{1}{2},\qquad L_{T} \leqdef (\log T)^{\theta},\qquad H_{T} \leqdef \log T.
\]
As assumed in \eqref{eq:fixed.parameters}, $\eta_{\sig}\log{X}/t \to \infty$. Rewrite the Euler polynomial
\[
\mathcal{Q}(s) = \sum_{j = 1}^{\infty}\sum_{p^{j}\leq X}\frac{p^{-js}}{j} = \sum_{2\leq n\leq X}\frac{\Lambda(n)}{\log{n}}n^{-s}
\]
(here $\Lambda(n)$ is the von Mangoldt function, equal to $\log{p}$ if $n = p^{j}$ for a prime $p$ and $j \in \N$, and zero otherwise) and introduce the shorthand notation $U_{h}$ and $V_{h}$ for the logarithmic zeta field and the Euler field, respectively:
\begin{equation}\label{analytic:eq:fields}
U_{h} \leqdef \log|\zeta(\sig + \ii\tau + \ii h)|, \qquad V_{h} \leqdef \Re{\mathcal{Q}(\sig + \ii\tau + \ii h)}.
\end{equation}
Values of $U_{h}$ at zeta zeros can be assigned arbitrarily on the corresponding null set when an expectation at a fixed $h\in I_{\theta}$ is taken. The fourth-moment estimate \eqref{analytic:eq:log-L4} in Lemma~\ref{analytic:lem:L4} proves the required integrability.

For fixed $\beta > 0$, let
\[
\mathcal{G}_{\beta,T}(\rd{h}) = \frac{e^{\beta U_{h}}\rd{h}}{\int_{I_{\theta}}e^{\beta U_{h'}}\rd{h'}},\qquad G_{\beta,T}^{\mathrm{A}}(\rd{h}) = \frac{e^{\beta V_{h}}\rd{h}}{\int_{I_{\theta}}e^{\beta V_{h'}}\rd{h'}}.
\]
The overlaps are
\begin{equation}\label{analytic:eq:overlaps}
\rho_{T}(h,h') = \frac{\EE[U_{h}U_{h'}]}{\sqrt{\EE[U_{h}^{2}]\EE[U_{h'}^{2}]}},\qquad R_{T}^{\mathrm{A}}(h,h') = \frac{\EE[V_{h}V_{h'}]}{\sqrt{\EE[V_{h}^{2}]\EE[V_{h'}^{2}]}}.
\end{equation}

\begin{theorem}[Gibbs measure and overlap comparison]\label{analytic:comparison}
For every fixed $\beta > 0$ and $n \in \N$, as $T \to \infty$,
\begin{equation}\label{analytic:eq:TV}
\EE\left[\int_{I_{\theta}^{\times n}}\left|\frac{\rd{\mathcal{G}_{\beta,T}^{\times n}}}{\rd{\bb{h}}}(\bb{h}) - \frac{\rd{(G_{\beta,T}^{\mathrm{A}})^{\times n}}}{\rd{\bb{h}}}(\bb{h})\right|\rd{\bb{h}}\right] = \oo(1),
\end{equation}
and
\begin{equation}\label{analytic:eq:rho-R}
\sup_{h,h'\in I_{\theta}}|\rho_{T}(h,h') - R_{T}^{\mathrm{A}}(h,h')| \to 0.
\end{equation}
\end{theorem}

Define the deterministic kernel
\begin{equation}\label{analytic:eq:log.kernel}
r_{T}(h,h') \leqdef \frac{\logp\big((\log T)\wedge|h - h'|^{-1}\big)}{t},\qquad r_{T}(h,h) = 1.
\end{equation}

\begin{theorem}[Uniform overlap geometry]\label{analytic:geometry}
As $T \to \infty$, uniformly for $h,h'\in I_{\theta}$,
\begin{equation}\label{analytic:eq:overlap-geometry}
\rho_{T}(h,h') = r_{T}(h,h') + \OO_{\theta}\big((\log{t})/t\big),\qquad R_{T}^{\mathrm{A}}(h,h') = r_{T}(h,h') + \OO_{\theta}\big((\log{t})/t\big).
\end{equation}
There is a constant $C_{\theta} < \infty$ such that, for all sufficiently large $T$ and every $h_{1},h_{2},h_{3}\in I_{\theta}$,
\begin{equation}\label{analytic:eq:deterministic-ultrametricity}
\rho_{T}(h_{1},h_{2}) \geq \min\{\rho_{T}(h_{1},h_{3}),\rho_{T}(h_{2},h_{3})\} - C_{\theta}\frac{\log{t}}{t}.
\end{equation}
The same conclusion holds with $R_{T}^{\mathrm{A}}$ in place of $\rho_{T}$.
\end{theorem}

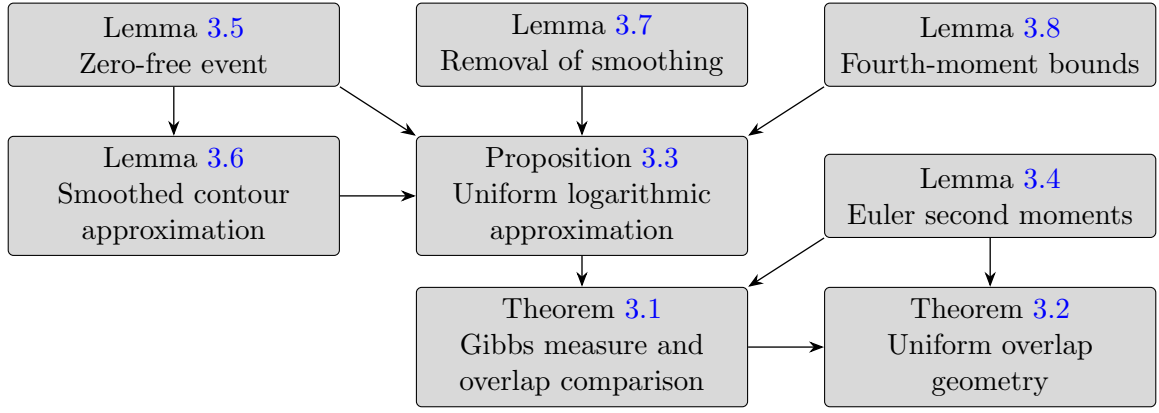
\begin{figure}[!ht]
\centering
\begin{tikzpicture}[x = 1cm, y = 1cm, every node/.style = {font = \small, align = center, execute at begin node = {\hyphenpenalty = 10000\exhyphenpenalty = 10000\relax}}, result/.style = {draw, rounded corners = 2pt, fill = black!15, text width = 4.1cm, minimum height = 1.1cm, inner sep = 4pt}, dependency/.style = {-{Stealth[length = 2mm]}, line width = 0.5pt}]
\node[result] (zero) at ( - 5.4,0) {Lemma~\ref{analytic:lem:zero-free} \\Zero-free event};
\node[result] (smoothing) at (0,0) {Lemma~\ref{analytic:lem:smoothing} \\Removal of smoothing};
\node[result] (moments) at (5.4,0) {Lemma~\ref{analytic:lem:L4} \\Fourth-moment bounds};
\node[result] (contour) at ( - 5.4, - 2) {Lemma~\ref{analytic:lem:contour} \\Smoothed contour approximation};
\node[result] (approximation) at (0, - 2) {Proposition~\ref{analytic:prop:approximation} \\Uniform logarithmic approximation};
\node[result] (covariance) at (5.4, - 2) {Lemma~\ref{analytic:lem:covariance} \\Euler second moments};
\node[result] (comparison) at (0, - 4) {Theorem~\ref{analytic:comparison} \\Gibbs measure and overlap comparison};
\node[result] (geometry) at (5.4, - 4) {Theorem~\ref{analytic:geometry} \\Uniform overlap geometry};
\draw[dependency] (zero) -- (contour);
\draw[dependency] (zero.south east) -- (approximation.north west);
\draw[dependency] (contour) -- (approximation);
\draw[dependency] (smoothing) -- (approximation);
\draw[dependency] (moments.south west) -- (approximation.north east);
\draw[dependency] (approximation) -- (comparison);
\draw[dependency] (covariance.south west) -- (comparison.north east);
\draw[dependency] (comparison) -- (geometry);
\draw[dependency] (covariance) -- (geometry);
\end{tikzpicture}
\caption{Proof dependencies between the results of Section~\ref{sec:analytic}.}\label{fig:analytic.dependencies}
\end{figure}

Figure~\ref{fig:analytic.dependencies} records the dependencies among the results of this section.

\subsection{The two analytic inputs and the comparison proofs}\label{analytic:sec:deduction}

A uniform approximation in probability controls $\rd{\mathcal{G}_{\beta,T}}/\rd{h}$ and $\rd{G_{\beta,T}^{\mathrm{A}}}/\rd{h}$, the normalized Gibbs densities. A pointwise approximation in $L^{2}$, uniform in the spatial argument $h$, controls their overlaps $\rho_{T}$ and $R_{T}^{\mathrm{A}}$. The following proposition separates these two uses of the logarithmic approximation.

\begin{proposition}[Uniform logarithmic approximation]\label{analytic:prop:approximation}
The fields $U_{h}$ and $V_{h}$ in \eqref{analytic:eq:fields} satisfy, as $T \to \infty$,
\begin{equation}\label{analytic:eq:uniform-probability}
\sup_{h\in I_{\theta}}|U_{h} - V_{h}|\xrightarrow{\PP}0.
\end{equation}
Furthermore, for all sufficiently large $T$,
\begin{equation}\label{analytic:eq:uniform-L2}
\sup_{h\in I_{\theta}}\EE[|U_{h} - V_{h}|^{2}] \leq \exp\big( - \eta_{\sig}\log{X} + \OO_{\theta}(t)\big).
\end{equation}
\end{proposition}

The second input identifies the scale of the Euler second moments $\EE[V_{h}V_{h'}]$ and their dependence on spatial separation.

\begin{lemma}[Euler second moments]\label{analytic:lem:covariance}
As $T \to \infty$, uniformly for $h,h'\in I_{\theta}$,
\begin{equation}\label{analytic:eq:covariance}
\EE[V_{h}V_{h'}] = \frac{1}{2}\logp\left(\frac{1}{\eta_{\sig} + |h - h'|}\right) + \OO_{\theta}(\log{t}).
\end{equation}
On the diagonal the sharper estimate
\begin{equation}\label{analytic:eq:variance}
\EE[V_{h}^{2}] = \frac{1}{2}\log(1/\eta_{\sig}) + \OO(1) = \frac{t}{2} + \OO(\log{t})
\end{equation}
holds uniformly in $h$.
\end{lemma}

The proofs of Proposition~\ref{analytic:prop:approximation} and Lemma~\ref{analytic:lem:covariance} are given in Subsections~\ref{analytic:sec:approximation} and \ref{analytic:sec:covariance}, respectively. Together they already imply both theorems at the start of this section.

\begin{proof}[Proof of Theorem~\ref{analytic:comparison}]
Fix $\beta > 0$ and $n \in \N$, let $\Delta_{T}^{\mathrm{fld}} \leqdef \sup_{h \in I_{\theta}}|U_{h} - V_{h}|$, and write $\mathrm{TV}_{n,T}$ for the integral in \eqref{analytic:eq:TV}. On $\{\Delta_{T}^{\mathrm{fld}} \leq \e\}$, the ratio of the normalized Gibbs densities lies in $[e^{-2\beta\e},e^{2\beta\e}]$, so the ratio of their $n$-fold product densities lies in $[e^{-2n\beta\e},e^{2n\beta\e}]$. Since each product density integrates to one, $\mathrm{TV}_{n,T} \leq e^{2n\beta\e} - 1$ on this event, while $\mathrm{TV}_{n,T} \leq 2$ always. Hence
\[
\EE[\mathrm{TV}_{n,T}] \leq e^{2n\beta\e} - 1 + 2\PP(\Delta_{T}^{\mathrm{fld}} > \e).
\]
Using \eqref{analytic:eq:uniform-probability} in Proposition~\ref{analytic:prop:approximation}, first let $T \to \infty$ and then $\e \downarrow 0$ to obtain $\EE[\mathrm{TV}_{n,T}] \to 0$, proving \eqref{analytic:eq:TV}.

Let $\smash{\Delta_{T}^{(2)}}\leqdef \sup_{h \in I_{\theta}}(\EE[|U_{h} - V_{h}|^{2}])^{1/2}$. The bound \eqref{analytic:eq:uniform-L2}, supplied by Proposition~\ref{analytic:prop:approximation}, gives $\smash{\Delta_{T}^{(2)}}\to 0$, while \eqref{analytic:eq:variance} in Lemma~\ref{analytic:lem:covariance} gives $\EE[V_{h}^{2}] \asymp t$ uniformly. The reverse triangle inequality in $L^{2}$ shows that the same estimate holds for $\EE[U_{h}^{2}]$. Moreover, uniformly for $h,h' \in I_{\theta}$,
\begin{equation}\label{analytic:eq:second.moment.comparison}
|\EE[U_{h}U_{h'}] - \EE[V_{h}V_{h'}]| \leq \Delta_{T}^{(2)}\big(\|U_{h'}\|_{2} + \|V_{h}\|_{2}\big) \ll \Delta_{T}^{(2)}\sqrt{t}.
\end{equation}
The bound \eqref{analytic:eq:second.moment.comparison} also applies to the differences of the diagonal second moments. Dividing by the normalizing factors in \eqref{analytic:eq:overlaps} yields
\begin{equation}\label{analytic:eq:quantitative.overlap.comparison}
\sup_{h,h' \in I_{\theta}}|\rho_{T}(h,h') - R_{T}^{\mathrm{A}}(h,h')| \ll \Delta_{T}^{(2)}/\sqrt{t},
\end{equation}
which proves \eqref{analytic:eq:rho-R}.
\end{proof}

\begin{proof}[Proof of Theorem~\ref{analytic:geometry}]
The estimates \eqref{analytic:eq:covariance}--\eqref{analytic:eq:variance} in Lemma~\ref{analytic:lem:covariance}, together with $\log(1/\eta_{\sig}) = t + \OO(\log{t})$ from \eqref{eq:parameter.scales}, give
\begin{equation}\label{analytic:eq:Euler.log.kernel}
R_{T}^{\mathrm{A}}(h,h') = \frac{\logp(1/(\eta_{\sig} + |h - h'|))}{t} + \OO_{\theta}\big((\log{t})/t\big).
\end{equation}
Replacing $\eta_{\sig}$ by $(\log T)^{-1}$ inside the logarithm changes the numerator by $\OO(\log{t})$, uniformly in $h,h'$. Also, for a real separation $d$,
\begin{equation}\label{analytic:eq:log.kernel.difference}
\left|\logp\left(\frac{1}{(\log T)^{-1} + |d|}\right) - \logp\big((\log T)\wedge|d|^{-1}\big)\right| \leq \log{2}.
\end{equation}
Together, \eqref{analytic:eq:Euler.log.kernel}, \eqref{eq:parameter.scales}, and \eqref{analytic:eq:log.kernel.difference} prove the assertion for $R_{T}^{\mathrm{A}}$ in \eqref{analytic:eq:overlap-geometry}. For $\rho_{T}$, the bound \eqref{analytic:eq:uniform-L2} in Proposition~\ref{analytic:prop:approximation} gives $\smash{\Delta_{T}^{(2)}} \leq \exp( - \eta_{\sig}\log{X}/2 + \OO_{\theta}(t))$. Since $\eta_{\sig}\log{X}/t \to \infty$, this bound is smaller than any fixed negative power of $t$ for all sufficiently large $T$. In particular, $\smash{\Delta_{T}^{(2)}}/\sqrt{t} = \OO((\log{t})/t)$, so the quantitative comparison \eqref{analytic:eq:quantitative.overlap.comparison} and the assertion for $R_{T}^{\mathrm{A}}$ prove the assertion for $\rho_{T}$.

The triangle inequality gives $|h_{1} - h_{2}| \leq 2\max\{|h_{1} - h_{3}|,|h_{2} - h_{3}|\}$. Since doubling the argument of the nonincreasing function $r \mapsto \logp((\log T)\wedge r^{-1})$ changes its value by at most $\log{2}$,
\begin{equation}\label{analytic:eq:kernel.ultrametricity}
r_{T}(h_{1},h_{2}) \geq \min\{r_{T}(h_{1},h_{3}),r_{T}(h_{2},h_{3})\} - \frac{\log{2}}{t}.
\end{equation}
Substituting \eqref{analytic:eq:overlap-geometry} into \eqref{analytic:eq:kernel.ultrametricity} proves \eqref{analytic:eq:deterministic-ultrametricity}.
\end{proof}

\subsection{The logarithmic approximation from three error bounds}\label{analytic:sec:approximation}

To prove Proposition~\ref{analytic:prop:approximation}, insert a smoothed Dirichlet series between $\log{\zeta}$ and $\mathcal{Q}$. Define the Gaussian smoothing function $W$ and the smoothed Dirichlet series $\mathcal{S}_{X}$ by
\[
W(v) \leqdef \frac{1}{2\sqrt{\pi}}\int_{-\infty}^{v}e^{-u^{2}/4}\rd{u} = \frac{1 + \operatorname{erf}(v/2)}{2},\qquad \mathcal{S}_{X}(s) \leqdef \sum_{n = 2}^{\infty}\frac{\Lambda(n)}{\log{n}}n^{-s}W(\log(X/n)).
\]
Here $\operatorname{erf}$ is the error function, given by $\operatorname{erf}(z) \leqdef (2/\sqrt{\pi})\int_{0}^{z} e^{-u^{2}}\rd{u}$ for real $z$. Note that the Gaussian decay of $W(\log(X/n))$ as $n\to\infty$ makes the series and all its derivatives absolutely convergent.

Denote the difference between the smoothed and sharply truncated series by
\[
\Delta_{X}^{\mathrm{sm}}(y) \leqdef \mathcal{S}_{X}(\sig + \ii y) - \mathcal{Q}(\sig + \ii y),
\]
and let $\mathcal{E}_{T}$ be the zero-free event that no nontrivial zero of $\zeta$ has real part at least $1/2 + \eta_{\sig}/4$ and ordinate within $L_{T} + H_{T} + 4$ of $\tau$.

On $\mathcal{E}_{T}$, a contour shift controls $\log{\zeta} - \mathcal{S}_{X}$ (Lemma~\ref{analytic:lem:contour}), using the analytic branch of the logarithm supplied by the zero-free event (Lemma~\ref{analytic:lem:zero-free}). A mean-square estimate controls $\Delta_{X}^{\mathrm{sm}}$, including its spatial supremum (Lemma~\ref{analytic:lem:smoothing}). Finally, the fourth moments in \eqref{analytic:eq:log-L4} and \eqref{analytic:eq:Q-L4} in Lemma~\ref{analytic:lem:L4} combine with the probability bound \eqref{analytic:eq:bad-event} in Lemma~\ref{analytic:lem:zero-free} to make the exceptional event $\mathcal{E}_{T}^{c}$ negligible in $L^{2}$ in the proof of Proposition~\ref{analytic:prop:approximation}.

\begin{lemma}[A zero-free event]\label{analytic:lem:zero-free}
The event $\mathcal{E}_{T}$, whose complement is denoted by $\mathcal{E}_{T}^{c}$, satisfies, as $T \to \infty$,
\begin{equation}\label{analytic:eq:bad-event}
\PP(\mathcal{E}_{T}^{c}) \leq \exp\big( - \eta_{\sig}\log T/3 + \OO_{\theta}(t)\big)
\end{equation}
and, on $\mathcal{E}_{T}$, an analytic branch of $\log{\zeta}$ exists on
\begin{equation}\label{analytic:eq:half-strip}
\left\{z:\Re(z) > \frac{1}{2} + \frac{\eta_{\sig}}{4},\quad |\Im(z) - \tau| < L_{T} + H_{T} + 4\right\}.
\end{equation}
The branch agreeing with the absolutely convergent Dirichlet series $\log{\zeta(z)} = \sum_{n=2}^{\infty}\frac{\Lambda(n)}{n^{z}\log{n}}$ on $\Re(z) > 1$ satisfies
\begin{equation}\label{analytic:eq:log-bound}
|\log{\zeta(x + \ii y)}| \ll \frac{\log T}{\eta_{\sig}},\qquad \frac{1}{2} + \frac{\eta_{\sig}}{2} \leq x \leq 2,\quad |y - \tau| \leq L_{T} + H_{T}.
\end{equation}
\end{lemma}

\begin{lemma}[Smoothed contour approximation]\label{analytic:lem:contour}
On $\mathcal{E}_{T}$, as $T \to \infty$,
\[
\sup_{h\in I_{\theta}}|\log{\zeta(\sig + \ii\tau + \ii h)} - \mathcal{S}_{X}(\sig + \ii\tau + \ii h)| \leq \exp\big( - \eta_{\sig}\log{X}/2 + \OO_{\theta}(t)\big).
\]
\end{lemma}

\begin{lemma}[Removal of smoothing]\label{analytic:lem:smoothing}
The smoothing error $\Delta_{X}^{\mathrm{sm}}(\tau + h)$ satisfies, as $T \to \infty$,
\[
\EE\left[\sup_{h\in I_{\theta}}|\Delta_{X}^{\mathrm{sm}}(\tau + h)|^{2}\right] \ll X^{-2\eta_{\sig}}\big(1 + L_{T}(1 + \log{X})\big) = \exp\big( - 2\eta_{\sig}\log{X} + \OO_{\theta}(t)\big).
\]
\end{lemma}

\begin{lemma}[Fourth-moment bounds]\label{analytic:lem:L4}
As $T \to \infty$, uniformly for $1/2 \leq x \leq 2$,
\begin{equation}\label{analytic:eq:log-L4}
\int_{T/2}^{3T}|\log|\zeta(x + \ii y)||^{4}\rd{y} \ll T(\log T)^{4}.
\end{equation}
Also, uniformly for $h \in I_{\theta}$,
\begin{equation}\label{analytic:eq:Q-L4}
\EE[|\mathcal{Q}(\sig + \ii\tau + \ii h)|^{4}] \ll \eta_{\sig}^{-4} \ll (\log T)^{4}.
\end{equation}
\end{lemma}

\begin{proof}[Proof of Proposition~\ref{analytic:prop:approximation}]
On $\mathcal{E}_{T}$, the real part of the logarithm in Lemma~\ref{analytic:lem:contour} equals $U_{h}$, and the identity $U_{h} - V_{h} = \Re\big(\log{\zeta(\sig + \ii\tau + \ii h)} - \mathcal{S}_{X}(\sig + \ii\tau + \ii h)\big) + \Re{\Delta_{X}^{\mathrm{sm}}(\tau + h)}$ splits the approximation error into the contour error and the smoothing error $\Delta_{X}^{\mathrm{sm}}$. The probability bound \eqref{analytic:eq:bad-event} in Lemma~\ref{analytic:lem:zero-free}, together with Lemmas~\ref{analytic:lem:contour} and \ref{analytic:lem:smoothing}, proves \eqref{analytic:eq:uniform-probability}. Moreover, Lemmas~\ref{analytic:lem:contour}--\ref{analytic:lem:smoothing} prove $\sup_{h\in I_{\theta}} \EE[|U_{h} - V_{h}|^{2}\ind_{\mathcal{E}_{T}}] \ll \exp( - \eta_{\sig}\log{X} + \OO_{\theta}(t))$. On the exceptional event $\mathcal{E}_{T}^{c}$, Cauchy--Schwarz and \eqref{analytic:eq:log-L4} and \eqref{analytic:eq:Q-L4} in Lemma~\ref{analytic:lem:L4} give
\begin{equation}\label{analytic:eq:exceptional.L2}
\sup_{h\in I_{\theta}}\EE[|U_{h} - V_{h}|^{2}\ind_{\mathcal{E}_{T}^{c}}] \ll (\log T)^{2}\PP(\mathcal{E}_{T}^{c})^{1/2}.
\end{equation}
By \eqref{analytic:eq:bad-event} in Lemma~\ref{analytic:lem:zero-free}, the right-hand side of \eqref{analytic:eq:exceptional.L2} is at most $\exp( - \eta_{\sig}\log T/6 + \OO_{\theta}(t))$. Since $\log T/\log{X} \geq t^{\phi}$ and $\eta_{\sig}\log{X}/t \to \infty$ by \eqref{eq:fixed.parameters}, this is smaller than the claimed bound \eqref{analytic:eq:uniform-L2}, so that bound follows.
\end{proof}

\subsubsection{Classical analytic estimates}

The analytic comparison uses the following standard estimates. Ingham's zero-density estimate \citep{Ingham1940}, in a uniform form near the critical line, is
\begin{equation}\label{analytic:eq:Ingham}
N(u,Y) \ll Y^{3(1-u)/(2-u)}(\log{Y})^{5},\qquad \frac{1}{2} \leq u \leq \frac{3}{4},
\end{equation}
where $N(u,Y)$ counts, with multiplicity, nontrivial zeros $\varrho = \sig_{\varrho} + \ii\gamma_{\varrho}$, with real part $\sig_{\varrho}$ and ordinate $\gamma_{\varrho}$, satisfying $\sig_{\varrho} \geq u$ and $0 < \gamma_{\varrho} \leq Y$. An explicit version is given in \citet[Corollary~1, Eq.~(3) and Table~1, p.~2]{ChourasiyaSimonic2025}.

The mean-value estimate of \citet[Corollary~3, Eq.~(1.11), p.~75]{MV1974} gives, for any Dirichlet polynomial with complex coefficients $b_{n}$, length $N \in \N$, and real endpoints $A < B$,
\begin{equation}\label{analytic:eq:MV}
\int_{A}^{B}\left|\sum_{n\leq N}b_{n} n^{-\ii y}\right|^{2}\rd{y} = (B - A)\sum_{n\leq N}|b_{n}|^{2} + \OO\left(\sum_{n\leq N}n|b_{n}|^{2}\right).
\end{equation}
The implicit constant is absolute and independent of $A$ and $B$. The estimate extends by truncation to an absolutely convergent Dirichlet series for which $\sum_{n=1}^{\infty} n|b_{n}|^{2}$ is finite.

The prime number theorem with its classical error term, used in Section~\ref{sec:expanding} but not in the analytic comparison, is given by \citet[Theorem~6.9, Eq.~(6.14), p.~179]{MR2378655},
\begin{equation}\label{analytic:eq:PNT}
\pi(y) = \operatorname{Li}(y) + \OO\big(y e^{-c\sqrt{\log{y}}}\big),
\end{equation}
where $\pi(y)$ counts primes at most $y$, $\operatorname{Li}(y) \leqdef \int_{2}^{y}\rd{u}/\log{u}$ is the logarithmic integral, and $c > 0$ is an absolute constant. The analytic comparison also uses Mertens' estimate, as stated in \citet[Eq.~(3.14.5), p.~65]{MR882550},
\begin{equation}\label{analytic:eq:Mertens}
\sum_{p\leq y}p^{-1} = \log{\log{y}} + \OO(1).
\end{equation}

\subsubsection{The zero-free event and the contour error}

The contour estimate requires an analytic logarithm $\log{\zeta}$ throughout a rectangle containing the whole spatial interval and the shifted contour. The proof of Lemma~\ref{analytic:lem:zero-free} first constructs this logarithm with sufficiently small exceptional probability; the proof of Lemma~\ref{analytic:lem:contour} then uses it in the contour shift to bound $|\mathcal{S}_{X}(s) - \log{\zeta(s)}|$ uniformly for $s = \sig + \ii\tau + \ii h$, $h\in I_{\theta}$.

\begin{figure}[!ht]
\centering
\includegraphics[width = 1\textwidth]{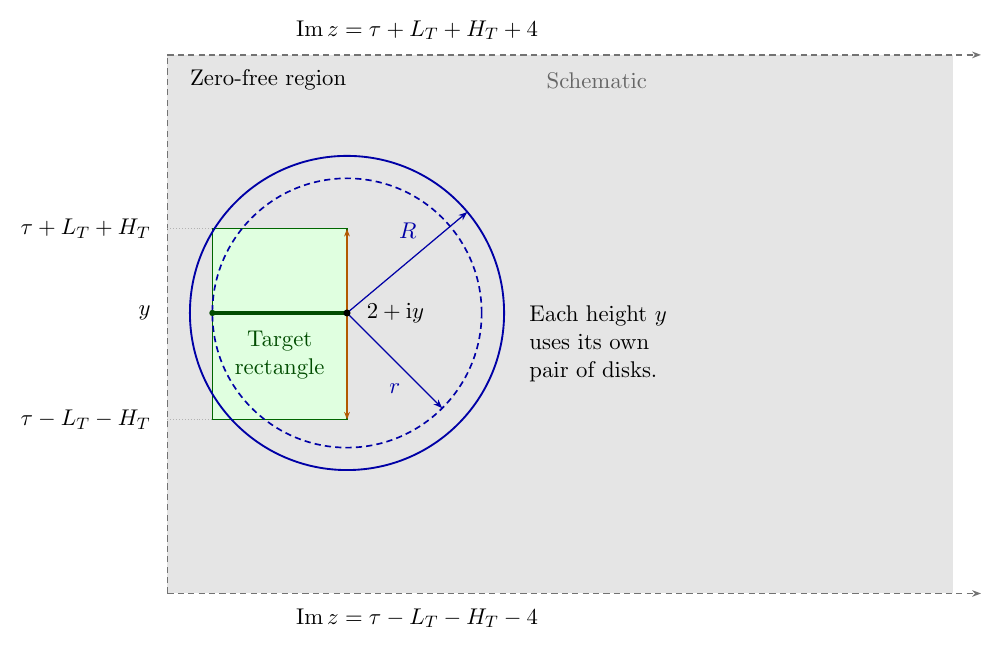}
\caption{Geometry of the analytic logarithm bound.}\label{fig:analytic.zero.free.geometry}
\end{figure}

\begin{proof}[Proof of Lemma~\ref{analytic:lem:zero-free}]
By the definition of $\mathcal{E}_{T}$, its complement requires a zero with real part at least $1/2 + \eta_{\sig}/4$ and ordinate within $L_{T} + H_{T} + 4$ of $\tau$. Since $L_{T} + H_{T} = \oo(T)$, a union bound and \eqref{analytic:eq:Ingham} give
\[
\PP(\mathcal{E}_{T}^{c}) \leq \frac{2(L_{T} + H_{T} + 4)}{T}N(1/2 + \eta_{\sig}/4,3T) \ll (L_{T} + H_{T} + 4)T^{-\eta_{\sig}/(3 - \eta_{\sig}/2)}(\log T)^{5}.
\]
The inequality
\[
- \frac{\eta_{\sig}}{3 - \eta_{\sig}/2} \leq - \frac{\eta_{\sig}}{3}
\]
proves \eqref{analytic:eq:bad-event}. The region in \eqref{analytic:eq:half-strip} is simply connected, and it does not contain the pole of $\zeta$ because $\tau \geq T > L_{T} + H_{T} + 4$ for all sufficiently large $T$. On $\mathcal{E}_{T}$, it also contains no zero of $\zeta$. Consequently, on $\mathcal{E}_{T}$, the function $\zeta$ admits an analytic logarithm on this region, uniquely specified by requiring it to agree with the absolutely convergent Dirichlet series $\sum_{n=2}^{\infty}\frac{\Lambda(n)}{n^{z}\log{n}}$ on the intersection of this region with $\Re(z) > 1$.

Fix $y$ such that $|y - \tau| \leq L_{T} + H_{T}$, as stated. Every point $x + \ii y$ in \eqref{analytic:eq:log-bound} lies in the closed disk centered at $2 + \ii y$ with radius $r \leqdef 3/2 - \eta_{\sig}/2$ (notice that only $x$ varies). The distance from this center to the left boundary of \eqref{analytic:eq:half-strip} is $3/2 - \eta_{\sig}/4$. Choose the outer radius $R \leqdef 3/2 - 3\eta_{\sig}/8$ halfway between this distance and $r$. This gives $R - r = \eta_{\sig}/8$ and leaves the same gap between the outer disk and the left boundary. The distances from the center to the horizontal boundaries are at least $4 > R$, so the closed outer disk lies in \eqref{analytic:eq:half-strip}. Figure~\ref{fig:analytic.zero.free.geometry} shows this geometry schematically: the region of \eqref{analytic:eq:log-bound} is green shaded, and its target segment at ordinate $y$ is drawn in bold and contained in the inner disk.

Polynomial growth of $\zeta$ in fixed vertical strips \citep[see, e.g.,][p.~95]{MR882550} gives $\Re{\log{\zeta}} \leq C\log T$ on the outer disk, while the Dirichlet series gives $|\log{\zeta(2 + \ii y)}| \leq C$. The Borel--Carath\'{e}odory inequality in \citet[Lemma~6.2, p.~169]{MR2378655} therefore yields
\[
|\log{\zeta(x + \ii y)}| \ll \frac{\log T}{R - r} \ll \frac{\log T}{\eta_{\sig}}
\]
throughout the inner disk, proving \eqref{analytic:eq:log-bound}.
\end{proof}

\begin{proof}[Proof of Lemma~\ref{analytic:lem:contour}]
Work on $\mathcal{E}_{T}$ and fix $h \in I_{\theta}$. Set $s \leqdef \sig + \ii\tau + \ii h$ and $c_{T} \leqdef 2 - \sig > 0$. On the line $w = c_{T} + \ii u$, the Dirichlet series for $\log{\zeta(s + w)}$ is absolutely convergent because $\Re(s + w) = 2$, while $|e^{w^{2}}| = e^{c_{T}^{2} - u^{2}}$. This Gaussian decay justifies termwise integration in
\begin{equation}\label{analytic:eq:Perron}
\mathcal{S}_{X}(s) = \frac{1}{2\pi\ii}\int_{c_{T} - \ii\infty}^{c_{T} + \ii\infty}\big(\log{\zeta(s + w)}\big)X^{w}e^{w^{2}}\frac{\rd{w}}{w}.
\end{equation}
To obtain \eqref{analytic:eq:Perron}, take $a = \sqrt{2}$ in the Gaussian inverse Mellin formula of \citet[Section~5.1.1, Exercise~4(c), p.~145]{MR2378655} and integrate it with respect to the logarithmic cutoff from $ - \infty$ to $\log{X}$.

It remains to bound $|\mathcal{S}_{X}(s) - \log{\zeta(s)}|$ uniformly over $h \in I_{\theta}$ on $\mathcal{E}_{T}$, as required in Lemma~\ref{analytic:lem:contour}. To this end, truncate \eqref{analytic:eq:Perron} at imaginary parts $\pm H_{T}$ and shift the contour to $\Re(w) = - \eta_{\sig}/2$. Throughout the rectangle swept out by this shift,
\[
\frac{1}{2} + \frac{\eta_{\sig}}{2} \leq \Re(s + w) \leq 2,\qquad |\Im(s + w) - \tau| \leq |h| + H_{T} \leq L_{T} + H_{T}.
\]
Lemma~\ref{analytic:lem:zero-free} therefore supplies the analytic branch of $\log{\zeta(s + w)}$ and, by \eqref{analytic:eq:log-bound} in that lemma, the uniform bound $|\log{\zeta(s + w)}| \ll (\log T)/\eta_{\sig}$ throughout this rectangle. The only pole of the integrand crossed during the shift is $w = 0$, whose residue is $\log{\zeta(s)}$. By the residue theorem, $|\mathcal{S}_{X}(s) - \log{\zeta(s)}|$ is consequently bounded by the contributions of the new vertical segment, the two horizontal segments, and the discarded tails. Figure~\ref{fig:analytic.contour} shows the truncated original line in blue, the new vertical segment in green, the horizontal segments in orange, and the discarded tails in purple.

\begin{figure}[!ht]
\centering
\includegraphics[width = 1\textwidth]{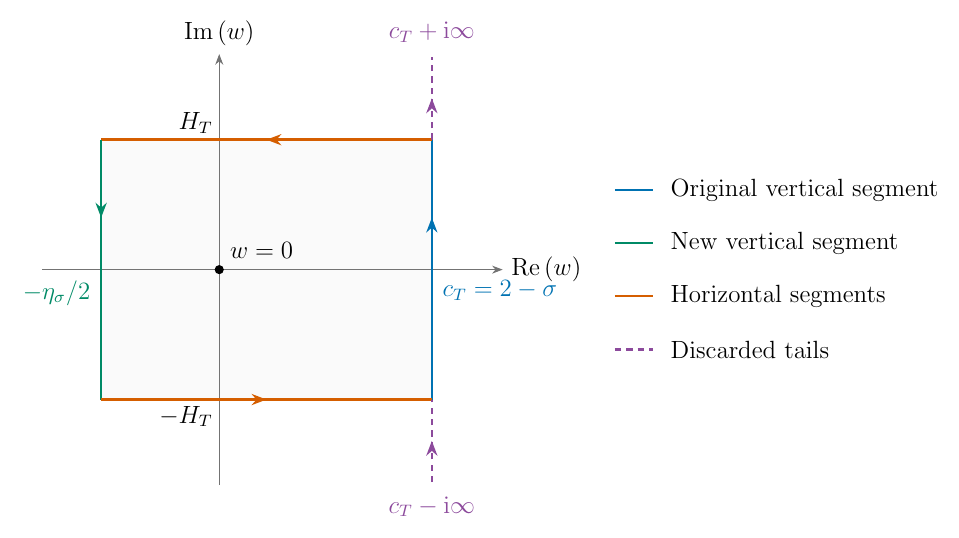}
\caption{Contour shift in the $w$-plane.}\label{fig:analytic.contour}
\end{figure}

On the new vertical segment, write $w = - \eta_{\sig}/2 + \ii v$. Its contribution is at most
\begin{equation}\label{analytic:eq:vertical.contribution}
C\frac{\log T}{\eta_{\sig}}X^{-\eta_{\sig}/2}\int_{-H_{T}}^{H_{T}}\frac{e^{-v^{2} + \eta_{\sig}^{2}/4}}{| - \eta_{\sig}/2 + \ii v|}\rd{v} \leq C\frac{\log T}{\eta_{\sig}^{2}}e^{-\eta_{\sig}(\log{X})/2}.
\end{equation}
The last inequality in \eqref{analytic:eq:vertical.contribution} uses $| - \eta_{\sig}/2 + \ii v| \geq \eta_{\sig}/2$ and the boundedness of $e^{\eta_{\sig}^{2}/4}\int_{\R}e^{-v^{2}}\rd{v}$ for large $T$.

On either horizontal segment, $w = u \pm \ii H_{T}$ with $ - \eta_{\sig}/2 \leq u \leq c_{T}$. Hence $|w| \geq H_{T}$, $|X^{w}| \leq X^{c_{T}}$, and $|e^{w^{2}}| = e^{u^{2} - H_{T}^{2}} \leq e^{-H_{T}^{2} + 4}$ for large $T$. Since each segment has bounded length, their combined contribution is at most
\begin{equation}\label{analytic:eq:horizontal.contribution}
C\frac{\log T}{\eta_{\sig}}\int_{-\eta_{\sig}/2}^{c_{T}}\frac{X^{u}e^{u^{2} - H_{T}^{2}}}{|u + \ii H_{T}|}\rd{u} \leq C\frac{\log T}{\eta_{\sig}H_{T}}X^{c_{T}}e^{-H_{T}^{2} + 4}.
\end{equation}

For the discarded tails on the original line, the absolutely convergent Dirichlet series gives $|\log{\zeta(2 + \ii u)}| \leq \log{\zeta(2)}$ for every real $u$. Their contribution is therefore bounded by
\begin{equation}\label{analytic:eq:tail.contribution}
CX^{c_{T}}e^{c_{T}^{2}}\int_{|v| > H_{T}}\frac{e^{-v^{2}}}{|c_{T} + \ii v|}\rd{v} \leq CX^{c_{T}}e^{-H_{T}^{2} + 4}.
\end{equation}
Because $H_{T} = \log T$ and $\log{X} \leq (\log T)/t^{\phi}$ by \eqref{eq:fixed.parameters}, the horizontal and tail contributions in \eqref{analytic:eq:horizontal.contribution} and \eqref{analytic:eq:tail.contribution} are both at most $\exp( - (\log T)^{2}/2)$ for large $T$.

Finally, \eqref{eq:parameter.scales} gives $\log((\log T)/\eta_{\sig}^{2}) = 3t + \OO(\log{t}) = \OO(t)$. The vertical contribution in \eqref{analytic:eq:vertical.contribution} is thus bounded by $\exp( - \eta_{\sig}\log{X}/2 + \OO(t))$, and the other contributions in \eqref{analytic:eq:horizontal.contribution} and \eqref{analytic:eq:tail.contribution} are smaller. All these estimates are uniform in $h \in I_{\theta}$, proving the lemma.
\end{proof}

\subsubsection{Removing the smoothing uniformly in the spatial shift}

The contour argument in the proof of Lemma~\ref{analytic:lem:contour} has replaced $\log \zeta$ by $\mathcal{S}_{X}$. To prove Lemma~\ref{analytic:lem:smoothing}, the mean-square bounds for the smoothing error $\Delta_{X}^{\mathrm{sm}}$ and its derivative $\big(\Delta_{X}^{\mathrm{sm}}\big)'$ are combined to control the supremum over $I_{\theta}$.

\begin{proof}[Proof of Lemma~\ref{analytic:lem:smoothing}]
Define the Dirichlet coefficient of the smoothing error by
\[
b_{X}^{\mathrm{sm}}(n) \leqdef \frac{\Lambda(n)}{\log{n}}\left(W(\log(X/n)) - \ind_{\{n\leq X\}}\right).
\]
An elementary Gaussian tail bound gives
\[
|W(v) - \ind_{\{v\geq0\}}| = \frac{1}{2\sqrt{\pi}}\int_{|v|}^{\infty}e^{-u^{2}/4}\rd{u} \leq \frac{1}{2}e^{-v^{2}/4},\qquad v \in \R.
\]
Since $0 \leq \Lambda(n)/\log{n} \leq 1$, taking $v = \log(X/n)$ yields
\begin{equation}\label{analytic:eq:dX-bound}
|b_{X}^{\mathrm{sm}}(n)| \leq \frac{1}{2}\exp\left( - \frac{(\log(n/X))^{2}}{4}\right).
\end{equation}

To estimate the weighted sums of these coefficients, partition the integers $n \geq 2$ into the logarithmic blocks
\[
\mathcal{B}_{k} \leqdef \{n \in \N:n \geq 2,\ Xe^{k} \leq n < Xe^{k + 1}\},\qquad k \in \Z.
\]
Each block contains at most $Xe^{k + 1}$ integers. For $n \in \mathcal{B}_{k}$, the quantity $u \leqdef \log(n/X)$ belongs to $[k,k + 1)$, so $u^{2} \geq k^{2}/2 - 1$. Consequently, \eqref{analytic:eq:dX-bound} gives $|b_{X}^{\mathrm{sm}}(n)|^{2} \ll e^{-k^{2}/4}$. Also,
\begin{equation}\label{analytic:eq:block.logarithm}
\log{n} \leq \log{X} + |k| + 1 \leq (1 + \log{X})(1 + |k|).
\end{equation}
For $m,j \in \{0,1\}$, the exponent $m - 2\sig$ is negative, so $n^{m - 2\sig} \leq (Xe^{k})^{m - 2\sig}$ on $\mathcal{B}_{k}$. Combining \eqref{analytic:eq:dX-bound} and \eqref{analytic:eq:block.logarithm} with $n^{m - 2\sig} \leq (Xe^{k})^{m - 2\sig}$ and the number of integers in the block gives
\begin{equation}\label{analytic:eq:block.weighted.sum}
\begin{aligned}
\sum_{n\in\mathcal{B}_{k}}|b_{X}^{\mathrm{sm}}(n)|^{2}n^{m - 2\sig}(\log{n})^{2j}
&\ll Xe^{k}(Xe^{k})^{m - 2\sig}e^{-k^{2}/4}(1 + \log{X})^{2j}(1 + |k|)^{2j} \\
&= X^{m - 2\eta_{\sig}}(1 + \log{X})^{2j}e^{(m - 2\eta_{\sig})k - k^{2}/4}(1 + |k|)^{2j}.
\end{aligned}
\end{equation}
For large $T$, the condition $0 < \eta_{\sig} < 1/4$ implies $|m - 2\eta_{\sig}| \leq 1$. Thus
\begin{equation}\label{analytic:eq:block.summability}
\sum_{k\in\Z}e^{(m - 2\eta_{\sig})k - k^{2}/4}(1 + |k|)^{2j} \leq \sum_{k\in\Z}e^{|k| - k^{2}/4}(1 + |k|)^{2} < \infty,
\end{equation}
with a bound independent of $T$. Summing \eqref{analytic:eq:block.weighted.sum} over the blocks, using \eqref{analytic:eq:block.summability}, and taking $m = 0$ and $m = 1$, respectively, proves, for $j = 0,1$,
\begin{equation}\label{analytic:eq:weighted-sums}
\begin{aligned}
\sum_{n = 2}^{\infty}|b_{X}^{\mathrm{sm}}(n)|^{2}n^{-2\sig}(\log{n})^{2j}&\ll X^{-2\eta_{\sig}}(1 + \log{X})^{2j}, \\
\sum_{n = 2}^{\infty}n|b_{X}^{\mathrm{sm}}(n)|^{2}n^{-2\sig}(\log{n})^{2j}&\ll X^{1 - 2\eta_{\sig}}(1 + \log{X})^{2j}.
\end{aligned}
\end{equation}

Write $\big(\Delta_{X}^{\mathrm{sm}}\big)^{(j)}$ for the $j$th derivative of $\Delta_{X}^{\mathrm{sm}}$ with respect to the real argument, with $\big(\Delta_{X}^{\mathrm{sm}}\big)^{(0)} \leqdef \Delta_{X}^{\mathrm{sm}}$. The infinite-series version of the mean-value estimate in \eqref{analytic:eq:MV}, combined with \eqref{analytic:eq:weighted-sums}, yields, uniformly for every real $v$,
\begin{equation}\label{analytic:eq:D-moments}
\frac{1}{T}\int_{T + v}^{2T + v}|\big(\Delta_{X}^{\mathrm{sm}}\big)^{(j)}(y)|^{2}\rd{y} \ll X^{-2\eta_{\sig}}(1 + \log{X})^{2j},\qquad j = 0,1,
\end{equation}
because $X/T \leq 1$. No interval-dependent constant is present. The fundamental theorem of calculus gives
\begin{equation}\label{analytic:eq:smoothing.calculus}
\sup_{|h|\leq L_{T}}|\Delta_{X}^{\mathrm{sm}}(\tau + h)|^{2} \leq |\Delta_{X}^{\mathrm{sm}}(\tau)|^{2} + 2\int_{-L_{T}}^{L_{T}}|\Delta_{X}^{\mathrm{sm}}(\tau + h)\big(\Delta_{X}^{\mathrm{sm}}\big)'(\tau + h)|\rd{h}.
\end{equation}
Average \eqref{analytic:eq:smoothing.calculus} in $\tau$, apply Cauchy--Schwarz for each $h$, and use \eqref{analytic:eq:D-moments}. This proves the lemma.
\end{proof}

\subsubsection{Fourth-moment bounds}

Selberg's mean-value approximation and the Dirichlet mean-value estimate \eqref{analytic:eq:MV} give the required fourth-moment bounds.

\begin{proof}[Proof of Lemma~\ref{analytic:lem:L4}]
By Selberg's mean-value approximation in \citet[p.~491]{MR2678896}, applied with $k = 2$ and $a = 1/2$, the fourth moment of $\log|\zeta(x + \ii y)| - \Re\sum_{p\leq (3T)^{1/4}}p^{-x - \ii y}$ over $[T/2,3T]$ is $\OO(T)$, uniformly for $1/2 \leq x \leq 1$. Applying \eqref{analytic:eq:MV} to the square of this prime sum and using Mertens' estimate \eqref{analytic:eq:Mertens} gives
\[
\int_{T/2}^{3T}|\log|\zeta(x + \ii y)||^{4}\rd{y} \ll T\left(1 + \left(\sum_{p\leq (3T)^{1/4}}p^{-2x}\right)^{2}\right) \ll T(\log{\log T})^{2},
\]
which proves \eqref{analytic:eq:log-L4} on this range.

For $1 \leq x \leq 2$, the classical estimate $\smash{\frac{\zeta'}{\zeta}(x + \ii y)} \ll \log T$ holds uniformly for $T/2 \leq y \leq 3T$ by \citet[Eq.~(3.11.7), p.~60]{MR882550}. Integrating its real part from $2$ to $x$, and using $|\log|\zeta(2 + \ii y)|| \ll 1$, gives $|\log|\zeta(x + \ii y)|| \ll \log T$. This proves \eqref{analytic:eq:log-L4} on the remaining range.

Write $\mathcal{Q}(s)^{2} = \sum_{n\leq X^{2}}b_{X}^{(2)}(n)n^{-s}$. Since $0 \leq \Lambda(n)/\log{n} \leq 1$, its coefficients satisfy $0 \leq b_{X}^{(2)}(n) \leq d(n)$, where $d(n)$ counts the positive divisors of $n$. The divisor-series identity in \citet[Section~1.3.1, Exercise~5, p.~27]{MR2378655} and \eqref{analytic:eq:MV} therefore give
\[
\EE[|\mathcal{Q}(\sig + \ii\tau + \ii h)|^{4}] \ll \left(1 + \frac{X^{2}}{T}\right)\sum_{n = 1}^{\infty}\frac{d(n)^{2}}{n^{1 + 2\eta_{\sig}}} = \left(1 + \frac{X^{2}}{T}\right)\frac{\zeta(1 + 2\eta_{\sig})^{4}}{\zeta(2 + 4\eta_{\sig})} \ll \eta_{\sig}^{-4},
\]
uniformly in $h$, because $X^{2}/T \to 0$ and $\zeta$ has a simple pole at $1$. As noted after \eqref{eq:parameter.scales}, the assumptions in \eqref{eq:fixed.parameters} give $\eta_{\sig} \geq t^{\phi + 1}/\log T \geq (\log T)^{-1}$ for all sufficiently large $T$, so this proves \eqref{analytic:eq:Q-L4}.
\end{proof}

\subsection{Euler second moments and overlap normalization}\label{analytic:sec:covariance}

The second-moment estimate follows by adapting the prime-sum comparison in the proof of \citet[Lemma~3.3, pp.~495--496]{MR2678896}. Expressing the prime sum through $\log|\zeta|$ near the line of real part $1$ gives both the logarithmic dependence on spatial separation and the sharper diagonal estimate.

\begin{proof}[Proof of Lemma~\ref{analytic:lem:covariance}]
Let $d \leqdef h - h'$ and write
\[
\mathcal{Q}_{h} \leqdef \mathcal{Q}(\sig + \ii\tau + \ii h) = \sum_{2\leq n\leq X}\frac{\Lambda(n)}{\log{n}}n^{-\sig - \ii\tau - \ii h}.
\]
Since $V_{h} = \Re\mathcal{Q}_{h}$, the identity $2\Re(z)\Re(w) = \Re(z\overline{w} + zw)$ gives
\[
\EE[V_{h}V_{h'}] = \frac{1}{2}\Re\big(\EE[\mathcal{Q}_{h}\overline{\mathcal{Q}_{h'}}]\big) + \frac{1}{2}\Re\big(\EE[\mathcal{Q}_{h}\mathcal{Q}_{h'}]\big).
\]
To estimate the first term, apply \eqref{analytic:eq:MV} to $\mathcal{Q}_{h} + \mathcal{Q}_{h'}$ and $\mathcal{Q}_{h} - \mathcal{Q}_{h'}$, then subtract the resulting identities. Using $|z + w|^{2} - |z - w|^{2} = 4\Re(z\overline{w})$ yields
\begin{equation}\label{analytic:eq:conjugate.moment}
\Re\big(\EE[\mathcal{Q}_{h}\overline{\mathcal{Q}_{h'}}]\big) = \sum_{2\leq n\leq X}\frac{\Lambda(n)^{2}\cos(d\log{n})}{n^{2\sig}(\log{n})^{2}} + \OO\left(\frac{X}{T}\sum_{2\leq n\leq X}\frac{\Lambda(n)^{2}}{n^{2\sig}(\log{n})^{2}}\right).
\end{equation}
The error in \eqref{analytic:eq:conjugate.moment} follows because each coefficient of either polynomial has modulus at most $2\Lambda(n)/(n^{\sig}\log{n})$ and $n \leq X$. Mertens' estimate \eqref{analytic:eq:Mertens} and the bounded contribution from higher prime powers give $\sum_{2\leq n\leq X}\Lambda(n)^{2}/(n^{2\sig}(\log{n})^{2}) \ll t$, so this error is $\OO(Xt/T)$, uniformly in $h,h'$.

For the second term, expand the product and use $|\EE[(nm)^{-\ii\tau}]| \leq 2/(T\log(nm))$ to obtain
\begin{equation}\label{analytic:eq:nonconjugate.moment}
|\EE[\mathcal{Q}_{h}\mathcal{Q}_{h'}]| \leq \frac{2}{T}\sum_{2\leq n,m\leq X}\frac{\Lambda(n)\Lambda(m)}{(nm)^{\sig}\log{n}\log{m}\log(nm)} \ll \frac{X^{2 - 2\sig}}{T} \leq \frac{X}{T}.
\end{equation}
Here $0 \leq \Lambda(n)/\log{n} \leq 1$, $\log(nm) \geq \log{4}$, and $\sum_{2\leq n\leq X}n^{-\sig} \ll X^{1 - \sig}$ for large $T$, since $\sig \to 1/2$.

Finally, $\Lambda(p^{j})/\log(p^{j}) = 1/j$, and $\Lambda(n) = 0$ when $n$ is not a prime power. The terms with $j \geq 2$ in the conjugated sum are bounded in absolute value by
\begin{equation}\label{analytic:eq:higher.power.covariance}
\sum_{p}\sum_{j = 2}^{\infty}\frac{1}{j^{2}p^{2j\sig}} \leq \sum_{p}\sum_{j = 2}^{\infty}\frac{1}{p^{j}} = \sum_{p}\frac{1}{p(p - 1)} < \infty.
\end{equation}
Keeping the terms with $j = 1$ in \eqref{analytic:eq:conjugate.moment}, using \eqref{analytic:eq:nonconjugate.moment} and \eqref{analytic:eq:higher.power.covariance}, and noting that $Xt/T \to 0$ therefore gives
\begin{equation}\label{analytic:eq:prime-covariance}
\EE[V_{h}V_{h'}] = \frac{1}{2}\sum_{p\leq X}\frac{\cos(d\log{p})}{p^{1 + 2\eta_{\sig}}} + \OO(1),
\end{equation}
uniformly in $h,h'$.

The prime sum can be extended to all primes because
\[
\sum_{p>X}\frac{1}{p^{1 + 2\eta_{\sig}}} \leq \sum_{n>X}\frac{1}{n^{1 + 2\eta_{\sig}}} \leq \int_{\lfloor X\rfloor}^{\infty}\frac{\rd{u}}{u^{1 + 2\eta_{\sig}}} \ll \frac{X^{-2\eta_{\sig}}}{\eta_{\sig}} = \exp( - 2\eta_{\sig}\log{X} + \OO(t)) = \oo(1).
\]
The absolutely convergent Euler product for $\zeta(1 + 2\eta_{\sig} + \ii d)$ then identifies this infinite prime sum with $\log|\zeta(1 + 2\eta_{\sig} + \ii d)|$, up to the uniformly bounded contribution from higher prime powers in \eqref{analytic:eq:higher.power.covariance}. Consequently, \eqref{analytic:eq:prime-covariance} becomes
\begin{equation}\label{analytic:eq:zeta.covariance}
\EE[V_{h}V_{h'}] = \frac{1}{2}\Re\sum_{p}\frac{e^{-\ii d\log{p}}}{p^{1 + 2\eta_{\sig}}} + \OO(1) = \frac{1}{2}\log|\zeta(1 + 2\eta_{\sig} + \ii d)| + \OO(1).
\end{equation}

For $|d| \leq 1$, the simple pole at $1$ and the nonvanishing of $\zeta$ on the line of real part $1$ give
\begin{equation}\label{analytic:eq:near.zeta.covariance}
\log|\zeta(1 + 2\eta_{\sig} + \ii d)| = - \log|2\eta_{\sig} + \ii d| + \OO(1) = \logp\left(\frac{1}{\eta_{\sig} + |d|}\right) + \OO(1).
\end{equation}
Uniformity follows because $(z - 1)\zeta(z)$, continued at $z = 1$, is bounded above and away from zero on the compact rectangle $1 \leq \Re(z) \leq 3/2$, $|\Im(z)| \leq 1$.

For $|d| \geq 1$, the bounds in \citet[Theorem~3.5, p.~49, and Eq.~(3.11.8), p.~60]{MR882550} give
\begin{equation}\label{analytic:eq:distant.zeta.bound}
|\zeta(1 + 2\eta_{\sig} + \ii d)| + |\zeta(1 + 2\eta_{\sig} + \ii d)|^{-1} \ll \log(3 + |d|).
\end{equation}
By \eqref{analytic:eq:distant.zeta.bound}, $|\log|\zeta(1 + 2\eta_{\sig} + \ii d)|| \ll \log{\log(3 + |d|)} \ll_{\theta} \log{t}$ for $1 \leq |d| \leq 2(\log T)^{\theta}$. Since $\logp(1/(\eta_{\sig} + |d|)) = 0$ in this range, substituting \eqref{analytic:eq:near.zeta.covariance} for $|d| \leq 1$ and the logarithmic bound obtained from \eqref{analytic:eq:distant.zeta.bound} for $|d| > 1$ into \eqref{analytic:eq:zeta.covariance} proves \eqref{analytic:eq:covariance}.

Finally, taking $d = 0$ in \eqref{analytic:eq:zeta.covariance} gives
\[
\EE[V_{h}^{2}] = \frac{1}{2}\log{\zeta(1 + 2\eta_{\sig})} + \OO(1) = \frac{1}{2}\log(1/\eta_{\sig}) + \OO(1) = \frac{t}{2} + \OO(\log{t}),
\]
using the residue $1$ of $\zeta$ at its pole and \eqref{eq:parameter.scales}. This proves \eqref{analytic:eq:variance}.
\end{proof}

\section{Overlap distributions of \texorpdfstring{$\Re \mathcal{Q}$}{Re Q} for \texorpdfstring{$\theta > 0$}{theta > 0}}\label{sec:expanding}

Throughout this section, assume that $\theta > 0$ and that $\sig$ and $X$ satisfy \eqref{eq:fixed.parameters}. The goal is to prove the overlap conclusions for the Euler Gibbs measure; see Theorem~\ref{expand:main}.

\subsection{Fields, scales, and the Euler overlap theorem}\label{expand:scales}

Let $(U_{p})_{p}$ be independent random variables uniformly distributed on the complex unit circle. Replacing each arithmetic phase $p^{-\ii\tau}$ in the Euler field by $U_{p}$ gives the Haar field $H$ and its Gibbs measure $G_{\beta,T}^{\mathrm{H}}$:
\begin{equation}\label{eq:haar.field}
H(h) \leqdef \Re\sum_{p^{j} \leq X}\frac{U_{p}^{j} p^{-\ii jh}}{j p^{j\sig}}, \qquad G_{\beta,T}^{\mathrm{H}}(\rd{h}) \leqdef \frac{\exp(\beta H(h))\rd{h}}{\int_{I_{\theta}}\exp(\beta H(h'))\rd{h'}}.
\end{equation}
Write $\PP_{\mathrm{H}}$ and $\EE_{\mathrm{H}}$ for probability and expectation with respect to $(U_{p})_{p}$. The probability and expectation over the common random shift $\tau$ in the Euler field are denoted by $\PP_{\mathrm{A}} \leqdef \PP$ and $\EE_{\mathrm{A}} \leqdef \EE$ when both fields occur together.

The Haar field is centered and stationary, since multiplying each $U_{p}$ by a deterministic phase preserves their joint law. Its variance $V_{T}$ and overlap $R_{T}^{\mathrm{H}}$ are therefore
\begin{equation}\label{eq:variance.Haar}
V_{T} \leqdef \EE_{\mathrm{H}}[H(h)^{2}] = \frac{1}{2}\sum_{p^{j} \leq X}\frac{1}{j^{2}p^{2j\sig}}, \qquad R_{T}^{\mathrm{H}}(h,h') \leqdef \frac{\EE_{\mathrm{H}}[H(h)H(h')]}{V_{T}}.
\end{equation}
The Euler field, its Gibbs measure, and its overlap retain the notation of Section~\ref{sec:intro} and Section~\ref{sec:analytic}. Use the aliases
\[
F_{h}^{\mathrm{A}} \leqdef V_{h}, \qquad F_{h}^{\mathrm{H}} \leqdef H(h).
\]
Throughout this section, $\star \in \{\mathrm{A},\mathrm{H}\}$. Unsubscripted probability and expectation also apply to auxiliary random objects when their law is clear from the context.

Since $\eta_{\sig}\log{X}/t \to \infty$ by \eqref{eq:fixed.parameters}, the truncation level $\log{X}$ lies far beyond $\eta_{\sig}^{-1}$. The prime terms oscillate in $h$ with frequencies $\log{p}$, and the leading covariance sums of the Euler and Haar fields contain the weights $p^{-2\sig} = p^{-1} e^{-2\eta_{\sig}\log{p}}$. The exponential factor suppresses frequencies $\log{p} \gg \eta_{\sig}^{-1}$, so extending the sums beyond this range contributes little additional variance. This is a gradual damping rather than a sharp truncation. In this sense, up to a constant factor, the effective logarithmic cutoff is $\eta_{\sig}^{-1}$ for all sufficiently large $T$. The reciprocal of the natural damping frequency $(2\eta_{\sig})^{-1}$ gives the microscopic spatial scale. Define this scale, its logarithmic inverse, and the rescaled interval length by
\begin{equation}\label{eq:N.T}
\ell_{T} \leqdef 2\eta_{\sig}, \qquad s_{T} \leqdef \log\left(\frac{1}{\ell_{T}}\right), \qquad N_{T} \leqdef \frac{|I_{\theta}|}{\ell_{T}}.
\end{equation}
These scales satisfy
\begin{equation}\label{eq:effective.cutoff.bounds}
s_{T} = \log\left(\frac{1}{2\eta_{\sig}}\right), \qquad \log{N_{T}} = s_{T} + \theta t + \log{2}.
\end{equation}

\newpage
\noindent
By \eqref{eq:parameter.scales}, $s_{T} = t + \OO(\log{t})$ and $\log{\log{X}} = t + \OO(\log{t})$. Since $\phi > 7/2$, the polynomial degree $t^{7/2}(\log{t})^{2}$ used in Lemma~\ref{an:smoothing} is $\oo(t^{\phi})$, as required for the arithmetic comparison at this cutoff.

Recall that $m(\theta) = \sqrt{1 + \theta}$ and $\beta_{c}(\theta) = 2m(\theta)$ for $\theta > 0$. The variance $V_{T}$ of the Haar field is defined in \eqref{eq:variance.Haar}. Define its height quantile $u_{T}$ by
\begin{equation}\label{eq:height.quantile}
N_{T}\PP_{\mathrm{H}}(H(0) > u_{T}) = 1.
\end{equation}
Equation \eqref{an:covariance.comparison} in Lemma~\ref{an:covariance} and \eqref{an:height.scales} in Lemma~\ref{an:height} show that $V_{T} = s_{T}/2 + \OO(1)$, $u_{T}/s_{T} \to m(\theta)$, and $u_{T}/V_{T} \to \beta_{c}(\theta)$. The exact quantile in \eqref{eq:height.quantile} is used before taking limits.

\begin{theorem}[Overlap distributions of the Euler field on expanding intervals]\label{expand:main}
Let $\theta > 0$ and $\beta > 0$. All limits are taken as $T \to \infty$. The Euler Gibbs measure $G_{\beta,T}^{\mathrm{A}}$ has limiting two-overlap distribution
\[
\EE (G_{\beta,T}^{\mathrm{A}})^{\times 2}\big[R_{T}^{\mathrm{A}}(h_{1},h_{2}) \in \, \cdot \,\big] \stackrel{\mathrm{weakly}}{\longrightarrow}
\begin{cases}
\delta_{0}, & \beta \leq \beta_{c}(\theta), \\
\dfrac{\beta_{c}(\theta)}{\beta}\delta_{0} + \left(1 - \dfrac{\beta_{c}(\theta)}{\beta}\right)\delta_{1}, & \beta > \beta_{c}(\theta).
\end{cases}
\]
For every fixed $n \in \N$, $1 \leq k \leq n$, and every measurable function $\Psi_{T} : I_{\theta}^{\times n} \to \R$ of the sampled positions with uniformly bounded supremum norm,
\begin{equation}\label{expand:GG}
\begin{aligned}
\EE (G_{\beta,T}^{\mathrm{A}})^{\times(n + 1)}[R_{T}^{\mathrm{A}}(h_{k},h_{n + 1})\Psi_{T}(\bb{h})]
&= \frac{1}{n}\EE (G_{\beta,T}^{\mathrm{A}})^{\times 2}[R_{T}^{\mathrm{A}}(h_{1},h_{2})]\EE (G_{\beta,T}^{\mathrm{A}})^{\times n}[\Psi_{T}(\bb{h})] \\
&\quad + \frac{1}{n}\sum_{\substack{1 \leq \ell \leq n \\ \ell \neq k}}\EE (G_{\beta,T}^{\mathrm{A}})^{\times n}[R_{T}^{\mathrm{A}}(h_{k},h_{\ell})\Psi_{T}(\bb{h})] + \oo(1).
\end{aligned}
\end{equation}
The $\oo(1)$ term is uniform over $\|\Psi_{T}\|_{\infty} \leq 1$. If $\beta > \beta_{c}(\theta)$, let $\bb{\xi}$ have Poisson--Dirichlet parameter $\beta_{c}(\theta)/\beta$ and, conditional on $\bb{\xi}$, let the group labels $K_{1},\ldots,K_{n}$ be independent with $\PP(K_{i} = j \mid \bb{\xi}) = \xi_{j}$. For every continuous function $\psi:[ - 1,1]^{n \times n} \to \R$,
\[
\lim_{T \to \infty}\EE (G_{\beta,T}^{\mathrm{A}})^{\times n}\left[\psi\left((R_{T}^{\mathrm{A}}(h_{i},h_{j}))_{1 \leq i,j \leq n}\right)\right]
= \begin{cases}
\psi(\mathrm{Id}_n), & \beta \leq \beta_{c}(\theta), \\
\EE_{\bb{\xi}}\big[\EE[\psi((\ind_{\{K_{i} = K_{j}\}})_{1 \leq i,j \leq n}) \mid \bb{\xi}]\big], & \beta > \beta_{c}(\theta).
\end{cases}
\]
\end{theorem}

At high and critical temperatures ($\beta \leq \beta_c(\theta)$), the proof establishes vanishing normalized masses on unit intervals. At low temperature ($\beta > \beta_c(\theta)$), a comparison of the annealed distributions of sampled positions transfers the Haar conclusions of \citet[Proposition~3.1]{OuimetHaar2026}. This comparison controls the random Gibbs normalization through capped measures and does not require a coupling of the Euler and Haar fields.

Subsection~\ref{expand:arithmetic} establishes the finite-dimensional arithmetic estimates. Subsection~\ref{sec:temperature} treats the unit-interval masses and proves a lower bound for the rescaled partition function. The normalization principle in Subsection~\ref{expand:normalization} is then applied in Subsection~\ref{sec:capped} to compare the annealed sampling laws above the critical inverse temperature. Subsection~\ref{sec:main.proof} combines these results with the Haar overlap laws of \citet[Proposition~3.1]{OuimetHaar2026}.

\begin{figure}[H]
\centering
\begin{tikzpicture}[x = 1cm, y = 1cm, every node/.style = {font = \small, align = center, execute at begin node = {\hyphenpenalty = 10000\exhyphenpenalty = 10000\relax}}, result/.style = {draw, rounded corners = 2pt, fill = black!15, text width = 4.1cm, minimum height = 1.05cm, inner sep = 4pt}, input/.style = {result, fill = white}, dependency/.style = {-{Stealth[length = 2mm]}, line width = 0.5pt}]
\path[fill = black!8, draw = black!20, rounded corners = 4pt] ( - 2.45,0.9) rectangle (7.65, - 3.05);
\path[fill = black!8, draw = black!20, rounded corners = 4pt] (7.95, - 0.7) rectangle (12.85, - 5.7);
\node[font = \small, anchor = west] at ( - 2.4,0.62) {Subsection~\ref{expand:arithmetic}};
\node[font = \small, anchor = west] at (8, - 0.95) {Subsection~\ref{sec:capped}};

\node[result] (polynomial) at (2.6,0) {Lemma~\ref{an:polynomial} \\Polynomial comparison};
\node[result] (normalization) at (10.4,0.2) {Proposition~\ref{norm:transfer} \\Normalization principle};

\node[result] (covariance) at (0, - 2.1) {Lemma~\ref{an:covariance} \\Uniform overlap estimates};
\node[result] (arithmetic) at (5.2, - 2.1) {Lemmas~\ref{an:transform}--\ref{an:smoothing} \\Height and smoothing \\estimates};
\node[result] (retained) at (10.4, - 2.1) {Lemmas~\ref{cap:moments}--\ref{cap:inequalities}, \ref{cap:loss} \\Propositions~\ref{cap:moment.measures}, \ref{cap:normalized.capture} \\Capped measures};

\node[result] (vanishing) at (0, - 4.5) {Corollary~\ref{temp:overlap} \\Vanishing off-diagonal overlaps \\$0 < \beta \leq \beta_{c}(\theta)$};
\node[result] (temperature) at (5.2, - 4.5) {Lemma~\ref{temp:critical.moment} \\Propositions~\ref{temp:critical}, \ref{temp:high} \\Corollary~\ref{temp:lower.mass} \\Mass estimates};
\node[result] (sampling) at (10.4, - 4.5) {Theorem~\ref{cap:sampling} \\Annealed sampling comparison \\$\beta > \beta_{c}(\theta)$};

\node[result] (main) at (5.2, - 6.6) {Theorem~\ref{expand:main} \\Euler overlap laws};
\node[input, fill = green!15] (haar) at (10.4, - 6.6) {\mbox{\citet{OuimetHaar2026}} \\Proposition~3.1 \\Haar overlap laws};

\draw[dependency] (polynomial.south west) to (covariance.north);
\draw[dependency] (polynomial.south east) to (arithmetic.north);
\draw[dependency] (covariance.east) to (arithmetic.west);
\draw[dependency] (arithmetic.south) to (temperature.north);
\draw[dependency] (arithmetic.east) to (retained.west);
\draw[dependency] (temperature.north east) to (retained.south west);
\draw[dependency] (retained.south) to (sampling.north);
\draw[dependency] (normalization.east) -- (13.15,0.2) -- (13.15, - 4.5) -- (sampling.east);
\draw[dependency] (temperature.west) to (vanishing.east);
\draw[dependency] (covariance.south) to (vanishing.north);
\draw[dependency] (vanishing.south east) to (main.north west);
\draw[dependency] (sampling.south west) to (main.north east);
\draw[dependency] (covariance.west) -- ( - 2.7, - 2.1) -- ( - 2.7, - 6.6) -- (main.west);
\draw[dependency] (haar.west) to (main.east);
\end{tikzpicture}
\caption{Dependencies among the results of Section~\ref{sec:expanding}. Shaded regions group selected results from the indicated subsections. The light green box identifies the established external Haar input. The dependencies within Subsections~\ref{expand:arithmetic} and~\ref{sec:capped} are detailed in Figures~\ref{fig:arithmetic.dependencies} and~\ref{fig:capped.dependencies}, respectively.}\label{fig:expanding.dependencies}
\end{figure}
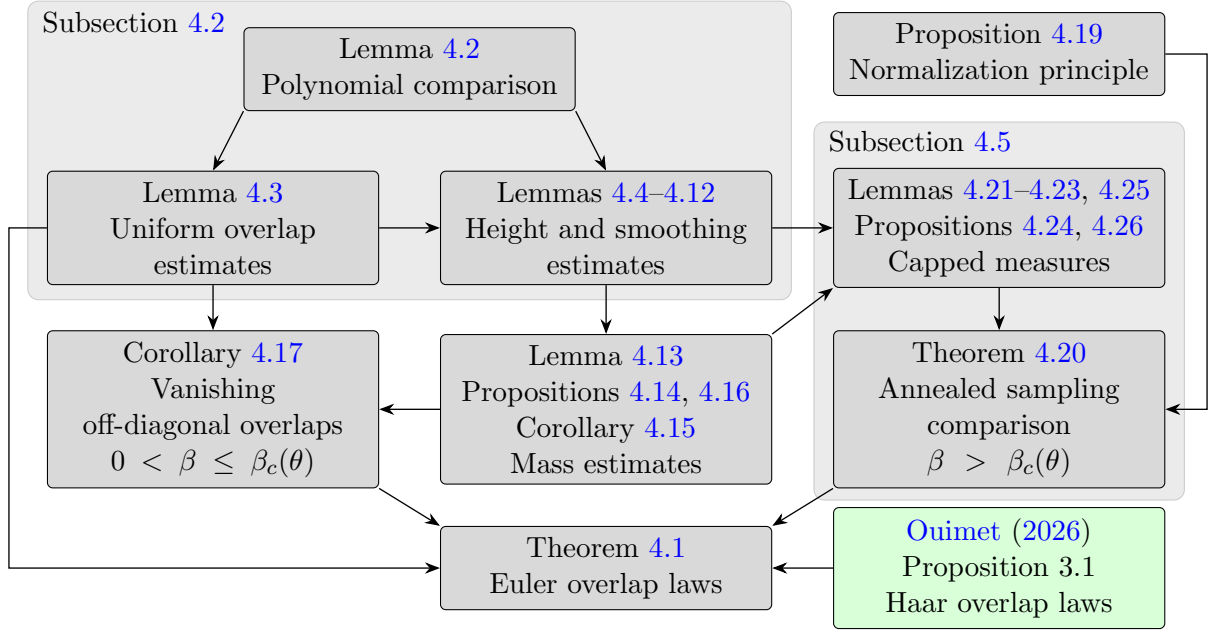

Figure~\ref{fig:expanding.dependencies} summarizes the proof dependencies in this section.

\subsection{Arithmetic estimates and height laws}\label{expand:arithmetic}

This subsection compares finite collections of field values before any Gibbs normalization. Constants may also depend on the number of selected coordinates. Every power of a given prime is grouped into the same independent summand under the Haar law.

\begin{figure}[H]
\centering
\begin{tikzpicture}[x = 1cm, y = 1cm, every node/.style = {font = \small, align = center, execute at begin node = {\hyphenpenalty = 10000\exhyphenpenalty = 10000\relax}}, result/.style = {draw, rounded corners = 2pt, fill = black!15, text width = 4.1cm, minimum height = 1.05cm, inner sep = 4pt}, dependency/.style = {-{Stealth[length = 2mm]}, line width = 0.5pt}]
\node[result] (polynomial) at (0,0) {Lemma~\ref{an:polynomial} \\Polynomial comparison \\and moment bounds};
\node[result] (transform) at (5.2,0) {Lemma~\ref{an:transform} \\Finite-dimensional \\truncated transforms};
\node[result] (spatial) at (10.4,0) {Lemma~\ref{an:spatial.transform} \\Simultaneous spatial \\truncated transforms};
\node[result] (smoothing) at (0, - 2.1) {Lemma~\ref{an:smoothing} \\Comparison after \\Gaussian smoothing};
\node[result] (covariance) at (5.2, - 2.1) {Lemma~\ref{an:covariance} \\Uniform overlap \\estimates};
\node[result] (conditional) at (10.4, - 2.1) {Lemma~\ref{an:conditional} \\Conditional Laplace \\estimates};
\node[result] (height) at (5.2, - 4.2) {Lemma~\ref{an:height} \\One-point heights \\and tails};
\node[result] (annular) at (10.4, - 4.2) {Lemma~\ref{an:annular} \\Annular fractional \\moments};
\node[result] (capped) at (0, - 6.3) {Lemma~\ref{an:capped.subcritical} \\Capped exponential moments \\at high temperature};
\node[result] (remote) at (5.2, - 6.3) {Lemma~\ref{an:remote} \\Height kernels at separated positions};
\node[result] (joint) at (10.4, - 6.3) {Lemma~\ref{an:joint.annular} \\Localization under joint \\height conditioning};
\draw[dependency] (polynomial.east) to (transform.west);
\draw[dependency] (transform.east) to (spatial.west);
\draw[dependency] (polynomial.south) to (smoothing.north);
\draw[dependency] (polynomial.south east) to (covariance.north west);
\draw[dependency] (spatial.south) to (conditional.north);
\draw[dependency] (covariance.south) to (height.north);
\draw[dependency] (conditional.south) to (annular.north);
\draw[dependency] (covariance.south east) to (annular.north west);
\draw[dependency, preaction = {draw = white, line width = 2pt}] (conditional.south west) to (height.north east);
\draw[dependency] (height.south) to (remote.north);
\draw[dependency] (remote.west) to (capped.east);
\draw[dependency] (annular.south) to (joint.north);
\draw[dependency] (remote.east) to (joint.west);
\end{tikzpicture}
\caption{Dependencies among the lemmas in Subsection~\ref{expand:arithmetic}. An arrow indicates that the target proof uses the source lemma or an estimate established in its proof. Arrows already represented by directed paths are omitted.}\label{fig:arithmetic.dependencies}
\end{figure}
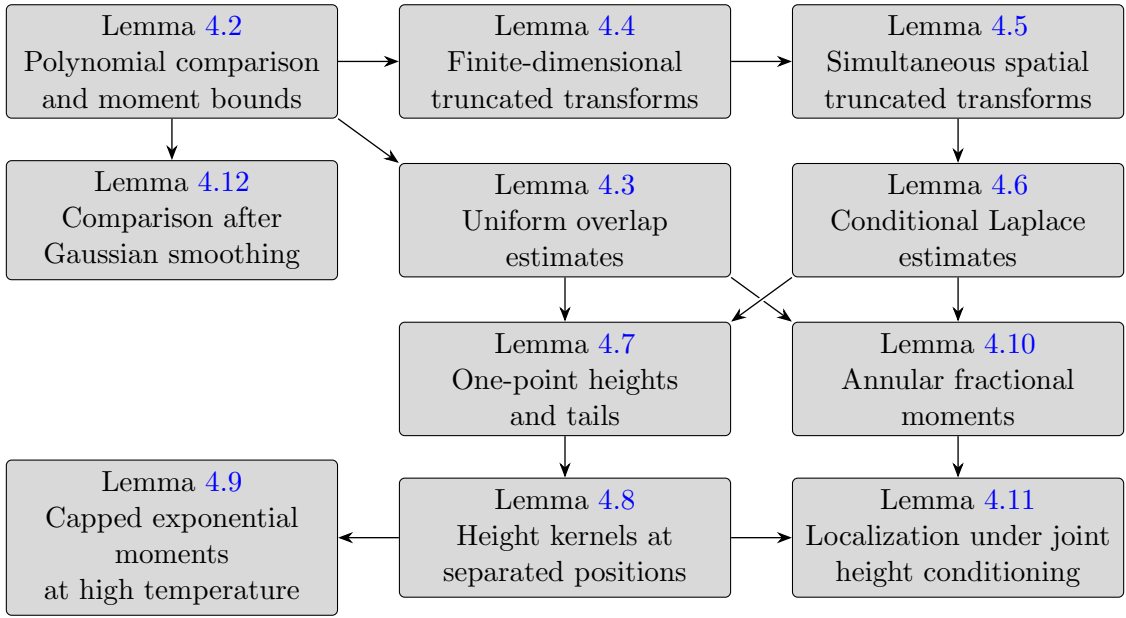

Figure~\ref{fig:arithmetic.dependencies} summarizes the proof dependencies among the lemmas in this subsection.

\subsubsection{Polynomial comparison and overlap estimates}

\begin{lemma}[Polynomial expectations and exponential moments]\label{an:polynomial}
Let $\bb{Y}_{T}^{\mathrm{A}}$ and $\bb{Y}_{T}^{\mathrm{H}}$ have a fixed number $n_{\mathrm{dim}}$ of coordinates of the form
\[
Y_{T,r}^{\mathrm{A}} \leqdef \Re\sum_{p^{j} \leq X}a_{r,p,j}p^{-\ii j\tau}, \qquad Y_{T,r}^{\mathrm{H}} \leqdef \Re\sum_{p^{j} \leq X}a_{r,p,j}U_{p}^{j}, \qquad 1 \leq r \leq n_{\mathrm{dim}},
\]
where the deterministic complex coefficients are the same under both laws and satisfy $|a_{r,p,j}| \leq C p^{-j/2}$ uniformly in $T$. For each prime $p$ and coordinate $r$, define
\[
A_{p,r} \leqdef \sum_{j:p^{j}\leq X}|a_{r,p,j}|.
\]
If $P$ is a polynomial of degree at most $D$ in these coordinates, then
\begin{equation}\label{an:polynomial.bound}
\left|\EE[P(\bb{Y}_{T}^{\mathrm{A}})] - \EE_{\mathrm{H}}[P(\bb{Y}_{T}^{\mathrm{H}})]\right| \leq \frac{CX^{2D}}{T}\left(\EE_{\mathrm{H}}[|P(\bb{Y}_{T}^{\mathrm{H}})|^{2}]\right)^{1/2}.
\end{equation}
For every $\bb{\lambda} \in \R^{n_{\mathrm{dim}}}$ and $k \in \N$,
\begin{equation}\label{an:haar.moments}
\EE_{\mathrm{H}}[e^{\bb{\lambda}\cdot\bb{Y}_{T}^{\mathrm{H}}}] \leq \exp(C_{n_{\mathrm{dim}}}t\|\bb{\lambda}\|_{2}^{2}), \qquad \EE_{\mathrm{H}}[|\bb{\lambda}\cdot\bb{Y}_{T}^{\mathrm{H}}|^{2k}] \leq (C_{n_{\mathrm{dim}}}kt\|\bb{\lambda}\|_{2}^{2})^{k}.
\end{equation}
The moment bound in \eqref{an:haar.moments}, with a changed constant, holds under the arithmetic law whenever $k = \oo(t^{\phi})$ as $T \to \infty$.

If a coordinate $r$ satisfies the stronger bound $\sum_{p \leq X}A_{p,r}^{2} \leq C$ uniformly in $T$, then, for every $\lambda \in \R$ and $k \in \N$,
\begin{equation}\label{an:bounded.moments}
\EE_{\mathrm{H}}[e^{\lambda Y_{T,r}^{\mathrm{H}}}] \leq e^{C\lambda^{2}}, \qquad \EE_{\mathrm{H}}[|Y_{T,r}^{\mathrm{H}}|^{2k}] \leq (Ck)^{k}.
\end{equation}
The moment bound in \eqref{an:bounded.moments} also holds under the arithmetic law whenever $k = \oo(t^{\phi})$ as $T \to \infty$. If every coordinate $r$ satisfies this stronger bound, then $t$ can be replaced by $1$ throughout \eqref{an:haar.moments} and its arithmetic moment counterpart.
\end{lemma}

\begin{proof}[Proof of Lemma~\ref{an:polynomial}]
The coefficient assumption gives $A_{p,r} \leq C p^{-1/2}$. Mertens' estimate \eqref{analytic:eq:Mertens} and the cutoff bound in \eqref{eq:fixed.parameters} give $\sum_{p \leq X}p^{-1} = t + \OO(\log{t})$, so $\sum_{p \leq X}A_{p,r}^{2} \leq Ct$. The centered contribution of a prime to $\bb{\lambda}\cdot\bb{Y}_{T}^{\mathrm{H}}$ is bounded in absolute value by $A_{p} \leqdef \sum_{r = 1}^{n_{\mathrm{dim}}}|\lambda_{r}|A_{p,r}$. Cauchy--Schwarz gives $\sum_{p}A_{p}^{2} \leq C_{n_{\mathrm{dim}}}t\|\bb{\lambda}\|_{2}^{2}$. Hoeffding's lemma for bounded variables and independence over primes give the first bound in \eqref{an:haar.moments}. For $\bb{\lambda} \neq \bb{0}$ and $x \geq 0$, Markov's inequality gives the sub-Gaussian tail bound
\begin{equation}\label{an:detail:subgaussian.tail}
\PP_{\mathrm{H}}(|\bb{\lambda}\cdot\bb{Y}_{T}^{\mathrm{H}}| \geq x) \leq 2\exp( - cx^{2}/(t\|\bb{\lambda}\|_{2}^{2})).
\end{equation}
Integrating the tail probability in \eqref{an:detail:subgaussian.tail} via the identity $\EE_{\mathrm{H}}[|Z|^{2k}] = \int_{0}^{\infty}2kx^{2k - 1}\PP_{\mathrm{H}}(|Z| \geq x)\rd{x}$ with $Z = \bb{\lambda}\cdot\bb{Y}_{T}^{\mathrm{H}}$ gives the second bound in \eqref{an:haar.moments}. The case $\bb{\lambda} = \bb{0}$ is immediate.

After equal rational characters are collected, a polynomial of degree $D$ has the form
\[
P(\bb{Y}_{T}^{\mathrm{A}}) = \sum_{\substack{n,m \leq X^{D} \\ (n,m)=1}}c_{n,m}(n/m)^{-\ii\tau}, \qquad
P(\bb{Y}_{T}^{\mathrm{H}}) = \sum_{\substack{n,m \leq X^{D} \\ (n,m)=1}}c_{n,m}U_{n}\overline{U}_{m},
\]
where $\smash{U_{n} \leqdef \prod_{p^{j} \parallel n} U_{p}^{j}}$. Unique factorization identifies its Haar expectation with $c_{1,1}$ (since $\EE_{\mathrm{H}}[U_{n}\overline{U}_{m}] = 0$ for $n \neq m$) and its squared Haar $L^{2}$ norm with $\sum_{n,m}|c_{n,m}|^{2}$ because the orthogonality of the independent random phases annihilates all cross terms, yielding
\begin{equation}\label{an:detail:character.orthogonality}
\EE_{\mathrm{H}}\left[|P(\bb{Y}_{T}^{\mathrm{H}})|^{2}\right] = \sum_{\substack{n,m \leq X^{D} \\ (n,m)=1}} \sum_{\substack{n',m' \leq X^{D} \\ (n',m')=1}} c_{n,m}\overline{c}_{n',m'} \EE_{\mathrm{H}}\left[U_{n}\overline{U}_{n'}\overline{U}_{m}U_{m'}\right] = \sum_{\substack{n,m \leq X^{D} \\ (n,m)=1}} |c_{n,m}|^{2}.
\end{equation}
For $n \neq m$, the mean value theorem applied to distinct positive integers bounded by $X^{D}$ implies $|\log(n/m)| \geq |n - m|/\max\{n,m\} \geq X^{-D}$, and hence
\begin{equation}\label{an:detail:character.mean}
\left|\frac{1}{T}\int_{T}^{2T}(n/m)^{-\ii v}\rd{v}\right| \leq \frac{1}{T}\left|\frac{(n/m)^{-2\ii T} - (n/m)^{-\ii T}}{-\ii\log(n/m)}\right| \leq \frac{2}{T|\log(n/m)|} \leq \frac{CX^{D}}{T}.
\end{equation}
Subtracting the Haar expectation isolates the oscillating terms. Using \eqref{an:detail:character.mean} and \eqref{an:detail:character.orthogonality} and applying the Cauchy--Schwarz inequality over the at most $CX^{2D}$ distinct characters yields
\[
\left|\EE[P(\bb{Y}_{T}^{\mathrm{A}})] - \EE_{\mathrm{H}}[P(\bb{Y}_{T}^{\mathrm{H}})]\right|
\leq \frac{CX^{D}}{T}\sum_{\substack{n,m \leq X^{D} \\ (n,m)=1,\,n\neq m}}|c_{n,m}|
\leq \frac{CX^{2D}}{T}\bigg(\sum_{\substack{n,m \leq X^{D} \\ (n,m)=1}}|c_{n,m}|^{2}\bigg)^{1/2},
\]
which proves \eqref{an:polynomial.bound}. Apply \eqref{an:polynomial.bound} to the polynomial $P(\bb{y}) \leqdef |\bb{\lambda}\cdot\bb{y}|^{2k}$, which has degree $D = 2k$, and use the second inequality in \eqref{an:haar.moments} with $k$ replaced by $2k$ to evaluate its Haar $L^{2}$ norm as $\smash{\left(\EE_{\mathrm{H}}[|\bb{\lambda}\cdot\bb{Y}_{T}^{\mathrm{H}}|^{4k}]\right)^{1/2} \leq (C_{n_{\mathrm{dim}}}2kt\|\bb{\lambda}\|_{2}^{2})^{k}}$. Because $\log{X} \leq (\log T)/t^{\phi}$ by \eqref{eq:fixed.parameters}, we have $X^{4k} \leq \exp(4k(\log T)/t^{\phi})$, meaning the error prefactor in \eqref{an:polynomial.bound} becomes $X^{4k}/T = T^{-1 + \oo(1)}$ whenever $k = \oo(t^{\phi})$. This vanishing error term ensures the arithmetic moment $\EE[|\bb{\lambda}\cdot\bb{Y}_{T}^{\mathrm{A}}|^{2k}]$ is bounded by the same leading term as the Haar moment $\EE_{\mathrm{H}}[|\bb{\lambda}\cdot\bb{Y}_{T}^{\mathrm{H}}|^{2k}]$, up to a negligible error, establishing the arithmetic counterpart of the second inequality in \eqref{an:haar.moments}. If a coordinate satisfies the stronger bound $\sum_{p \leq X}A_{p,r}^{2} \leq C$, the variance proxy in Hoeffding's lemma loses its dependence on $t$, and \eqref{an:polynomial.bound} yields the improved $t$-independent moment bounds for that coordinate, as well as for the full vector if this stricter condition holds uniformly across all coordinates.
\end{proof}

\begin{lemma}[Uniform overlap estimates]\label{an:covariance}
Let $\theta > 0$. As $T \to \infty$, uniformly for $h,h' \in I_{\theta}$,
\begin{equation}\label{an:overlap.geometry}
R_{T}^{\mathrm{A}}(h,h') = r_{T}(h,h') + \OO\left(\frac{\log{t}}{t}\right), \qquad R_{T}^{\mathrm{H}}(h,h') = r_{T}(h,h') + \OO\left(\frac{\log{t}}{t}\right),
\end{equation}
where $r_{T}$ is the logarithmic distance kernel defined in \eqref{analytic:eq:log.kernel}. Moreover,
\begin{equation}\label{an:covariance.comparison}
V_{T} = \frac{s_{T}}{2} + \OO(1), \qquad \sup_{h,h' \in I_{\theta}}|R_{T}^{\mathrm{A}}(h,h') - R_{T}^{\mathrm{H}}(h,h')| = \oo(1).
\end{equation}
In particular, $\sup_{|h - h'|\geq 1}|R_{T}^{\mathrm{A}}(h,h')| = \OO((\log{t})/t)$.
\end{lemma}

\begin{proof}[Proof of Lemma~\ref{an:covariance}]
Let $d \leqdef h - h'$. Independence of the Haar phases and the identities $\smash{\EE_{\mathrm{H}}[U_{p}^{j}\overline{U}_{q}^{k}] = \ind_{\{p = q,\,j = k\}}}$ and $\EE_{\mathrm{H}}[U_{p}^{j}U_{q}^{k}] = 0$, for $j,k \in \N$, give
\begin{equation}\label{an:detail:cov.haar.sum}
\EE_{\mathrm{H}}[H(h)H(h')] = \frac{1}{2}\sum_{p^{j} \leq X}\frac{\cos(jd\log{p})}{j^{2}p^{2j\sig}} = \frac{1}{2}\sum_{p \leq X}\frac{e^{-2\eta_{\sig}\log{p}}\cos(d\log{p})}{p} + \OO(1).
\end{equation}
The error in \eqref{an:detail:cov.haar.sum} is uniform in $d$, because $2\sig \geq 1$ and
\[
\sum_{p}\sum_{j = 2}^{\infty}\frac{1}{j^{2}p^{2j\sig}} \leq \sum_{p}\sum_{j = 2}^{\infty}p^{-j} \leq 2\sum_{p}p^{-2} < \infty.
\]
The assumptions in \eqref{eq:fixed.parameters} imply $\eta_{\sig}\log{X} \to \infty$. They therefore give
\begin{equation}\label{an:effective.cutoff.estimates}
\int_{\log{2}}^{\log{X}}e^{-2\eta_{\sig}v}\frac{\rd{v}}{v} = \log(1/\eta_{\sig}) + \OO(1), \qquad
\ell_{T}^{2j}\int_{0}^{\log{X}}v^{2j - 1}e^{-2\eta_{\sig}v}\rd{v} \leq C_{j}, \quad j \in \N.
\end{equation}
For the first estimate in \eqref{an:effective.cutoff.estimates}, split the integral at $\eta_{\sig}^{-1}$. On its lower part, $|e^{-2\eta_{\sig}v} - 1| \leq 2\eta_{\sig}v$, so replacing the damping factor by $1$ costs at most $\smash{2\eta_{\sig}\int_{0}^{\eta_{\sig}^{-1}}\rd{v} = 2}$. On its upper part, the substitution $w \leqdef 2\eta_{\sig}v$ bounds the integral by $\int_{2}^{\infty}e^{-w}\rd{w}/w$. The same substitution, together with $\ell_{T} = 2\eta_{\sig}$, bounds the second integral in \eqref{an:effective.cutoff.estimates} by $\int_{0}^{\infty}w^{2j - 1}e^{-w}\rd{w} = (2j - 1)!$.

Partial summation using the prime number theorem in \eqref{analytic:eq:PNT} gives, at $d = 0$,
\begin{equation}\label{an:detail:cov.variance}
V_{T} = \frac{1}{2}\sum_{p^{j} \leq X}\frac{1}{j^{2}p^{2j\sig}} = \frac{1}{2}\int_{\log{2}}^{\log{X}}e^{-2\eta_{\sig}v}\frac{\rd{v}}{v} + \OO(1) = \frac{1}{2}\log(1/\eta_{\sig}) + \OO(1) = \frac{s_{T}}{2} + \OO(1).
\end{equation}
Here the error in \eqref{analytic:eq:PNT} is integrable uniformly in $\eta_{\sig}$.

To obtain uniformity for nonzero $d$, apply \eqref{analytic:eq:PNT} above $e^{t^{3}}$. Partial summation yields
\begin{equation}\label{an:detail:cov.prime.integral}
\sum_{\exp(t^{3}) < p \leq X}\frac{e^{-2\eta_{\sig}\log{p}}\cos(d\log{p})}{p} = \int_{t^{3}}^{\log{X}}e^{-2\eta_{\sig}v}\cos(dv)\frac{\rd{v}}{v} + \OO\left((1 + |d|)t^{C}e^{-ct^{3/2}}\right).
\end{equation}
The derivative of $\cos(dv)$ accounts for the factor $|d|$. Since $|d| \leq 2|I_{\theta}| = \exp(\OO_{\theta}(t))$, the error in \eqref{an:detail:cov.prime.integral} is $\oo(1)$ uniformly in the required range. The primes and the integral below $e^{t^{3}}$ each have absolute contribution at most $C\log{t}$. Consequently, \eqref{an:detail:cov.haar.sum} and \eqref{an:detail:cov.prime.integral} imply
\begin{equation}\label{an:detail:cov.full.integral}
\EE_{\mathrm{H}}[H(h)H(h')] = \frac{1}{2}\int_{\log{2}}^{\log{X}}e^{-\ell_{T}v}\cos(dv)\frac{\rd{v}}{v} + \OO(\log{t}).
\end{equation}
Let $a \leqdef \ell_{T} \vee |d|$. If $\log{2} < a^{-1}$, then on $\log{2} \leq v \leq a^{-1}$ the inequalities $|e^{-\ell_{T}v} - 1| \leq \ell_{T}v$ and $|\cos(dv) - 1| \leq d^{2}v^{2}/2$ show that replacing the integrand by $1/v$ costs $\OO(1)$. Above $a^{-1}$, exponential decay gives a bounded integral when $a = \ell_{T}$; when $a = |d|$, integration by parts gives the same conclusion because $e^{-\ell_{T}v}/v$ is decreasing and $|d|a^{-1} = 1$. If $\log{2} \geq a^{-1}$, the same integration by parts starts at $\log{2}$ and bounds the entire integral. Extending the upper limit from $\log{X}$ to infinity costs $\oo(1)$, since $\ell_{T}\log{X} \to \infty$. Finally, $a \leq \ell_{T} + |d| \leq 2a$, so \eqref{an:detail:cov.full.integral} gives
\begin{equation}\label{an:detail:cov.logarithmic}
\EE_{\mathrm{H}}[H(h)H(h')] = \frac{1}{2}\logp\left(\frac{1}{\ell_{T} + |d|}\right) + \OO(\log{t}).
\end{equation}

For $d \neq 0$, the quantity $\logp(1/(\ell_{T} + |d|))$ differs from $\min\{s_{T},\logp(1/|d|)\}$ by at most $\log{2}$; for $d = 0$ it equals $s_{T}$. Since $s_{T} = t + \OO(\log{t})$ by \eqref{eq:parameter.scales} and \eqref{eq:effective.cutoff.bounds}, this proves
\begin{equation}\label{an:detail:logarithmic.kernel}
\logp\left(\frac{1}{\ell_{T} + |d|}\right) = t r_{T}(h,h') + \OO(\log{t}).
\end{equation}
Divide \eqref{an:detail:cov.logarithmic} by \eqref{an:detail:cov.variance}, using \eqref{an:detail:logarithmic.kernel} and $0 \leq r_{T} \leq 1$, to obtain the Haar assertion in \eqref{an:overlap.geometry}.

For the arithmetic comparison, apply \eqref{an:polynomial.bound} in Lemma~\ref{an:polynomial} to the degree-two polynomial given by the product of the two field coordinates. Cauchy--Schwarz and the fourth-moment case of \eqref{an:haar.moments} in Lemma~\ref{an:polynomial} bound its Haar $L^{2}$ norm by $Ct$, uniformly in $h,h'$. Hence
\begin{equation}\label{an:detail:cov.arithmetic.error}
\sup_{h,h' \in I_{\theta}}\left|\EE[V_{h}V_{h'}] - \EE_{\mathrm{H}}[H(h)H(h')]\right| \leq \frac{CX^{4}t}{T} = T^{-1 + \oo(1)}.
\end{equation}
The last equality follows from $\log{X} \leq (\log T)/t^{\phi}$ in \eqref{eq:fixed.parameters}. Estimate \eqref{an:detail:cov.arithmetic.error} with $h' = h$ gives $\EE[V_{h}^{2}] = V_{T} + \OO(X^{4}t/T)$ uniformly in $h$. Thus the arithmetic overlap denominator $\sqrt{\EE[V_{h}^{2}]\EE[V_{h'}^{2}]}$ equals $V_{T} + \OO(X^{4}t/T)$, and is bounded below by $ct$ for all sufficiently large $T$ by \eqref{an:detail:cov.variance} and \eqref{eq:parameter.scales}. Since $|\EE_{\mathrm{H}}[H(h)H(h')]| \leq V_{T}$, subtracting the two normalized ratios in \eqref{an:detail:cov.arithmetic.error} gives
\begin{equation}\label{an:detail:normalized.covariance.error}
\sup_{h,h' \in I_{\theta}}|R_{T}^{\mathrm{A}}(h,h') - R_{T}^{\mathrm{H}}(h,h')| \leq \frac{CX^{4}}{T} = \oo(1).
\end{equation}
The error in \eqref{an:detail:normalized.covariance.error} is also $\oo((\log{t})/t)$, proving the arithmetic assertion in \eqref{an:overlap.geometry}. For $|h - h'| \geq 1$, the definition gives $r_{T}(h,h') = 0$, so the final assertion follows directly from \eqref{an:overlap.geometry}.
\end{proof}

\subsubsection{Truncated transforms and rare-height conditioning}

\begin{lemma}[Finite-dimensional truncated transforms]\label{an:transform}
Fix $0 < \e_{0} < 1/4$, $n_{\mathrm{dim}} \in \N$, and $C_{0} > 0$. Let $\bb{Y}_{T}^{\mathrm{A}}$ and $\bb{Y}_{T}^{\mathrm{H}}$ have coordinates
\[
Y_{T,r}^{\mathrm{A}} \leqdef \Re\sum_{p^{j} \leq X}a_{r,p,j}p^{-\ii j\tau}, \qquad Y_{T,r}^{\mathrm{H}} \leqdef \Re\sum_{p^{j} \leq X}a_{r,p,j}U_{p}^{j}, \qquad 1 \leq r \leq n_{\mathrm{dim}},
\]
where the deterministic complex coefficients are the same under both laws and satisfy $|a_{r,p,j}| \leq C_{0}p^{-j/2}$ uniformly in $T$. For each $\star \in \{\mathrm{A},\mathrm{H}\}$, define
\[
\mathcal{A}_{T}^{\star} \leqdef \left\{\max_{1 \leq r \leq n_{\mathrm{dim}}}|Y_{T,r}^{\star}| \leq t^{1 + \e_{0}}\right\}.
\]
For every fixed $M > 0$, there are constants $c,C > 0$ such that, for all sufficiently large $T$,
\begin{equation}\label{an:exceptional}
\PP_{\star}((\mathcal{A}_{T}^{\star})^{c}) \leq \exp( - ct^{1 + 2\e_{0}}),
\end{equation}
and, uniformly for $\bb{z} \in \C^{n_{\mathrm{dim}}}$ with $|\Re(z_{r})| \leq M$ and $|\Im(z_{r})| \leq t^{\e_{0}/2}$ for every $1 \leq r \leq n_{\mathrm{dim}}$,
\begin{equation}\label{an:detail:conditional.transform}
\left|\EE_{\star}\left[e^{\bb{z}\cdot\bb{Y}_{T}^{\star}}\ind_{\mathcal{A}_{T}^{\star}}\right] - \EE_{\mathrm{H}}\left[e^{\bb{z}\cdot\bb{Y}_{T}^{\mathrm{H}}}\right]\right| \leq \exp\{ - ct^{1 + 2\e_{0}} + Ct\}.
\end{equation}
The constants are uniform over all coefficients satisfying the stated bound and may depend on $\e_{0}$, $n_{\mathrm{dim}}$, $C_{0}$, $M$, and the fixed parameters of the paper.
\end{lemma}

\begin{proof}[Proof of Lemma~\ref{an:transform}]
Let $k \leqdef \lfloor c_{0}t^{1 + 2\e_{0}}\rfloor$, where $c_{0} > 0$ is sufficiently small. Since $1 + 2\e_{0} < 3/2 < \phi$, the arithmetic counterpart of the moment bound \eqref{an:haar.moments} in Lemma~\ref{an:polynomial} is available at order $2k$. Markov's inequality and \eqref{an:haar.moments} in Lemma~\ref{an:polynomial}, or its arithmetic counterpart, give for each coordinate $Y_{T,r}^{\star}$ satisfying the stated coefficient bound
\begin{equation}\label{an:detail:conditional.coordinate.tail}
\PP_{\star}(|Y_{T,r}^{\star}| > t^{1 + \e_{0}}) \leq \frac{(Ckt)^{k}}{t^{2k(1 + \e_{0})}} = \left(\frac{Ck}{t^{1 + 2\e_{0}}}\right)^{k} \leq \exp( - ct^{1 + 2\e_{0}}).
\end{equation}
The last inequality in \eqref{an:detail:conditional.coordinate.tail} follows by choosing $Cc_{0} < 1$. A union bound over the $n_{\mathrm{dim}}$ coordinates contributes only the fixed factor $n_{\mathrm{dim}}$, which is absorbed by reducing $c$ for all sufficiently large $T$. This proves \eqref{an:exceptional} for $\mathcal{A}_{T}^{\star}$.

Let $\bb{z}$ be a complex parameter vector with bounded real parts and $|\Im(z_{r})| \leq t^{\e_{0}/2}$, and let $n_{T} \leqdef \lfloor t^{2}\rfloor$. On $\mathcal{A}_{T}^{\star}$, $w \leqdef \bb{z}\cdot\bb{Y}_{T}^{\star}$ satisfies $|w| \leq Ct^{1 + 3\e_{0}/2} = \oo(n_{T})$. The exponential Taylor remainder therefore satisfies, on $\mathcal{A}_{T}^{\star}$,
\begin{equation}\label{an:detail:conditional.taylor}
\left|e^{w} - \sum_{j = 0}^{n_{T}}\frac{w^{j}}{j!}\right| \leq \sum_{j = n_{T} + 1}^{\infty}\frac{|w|^{j}}{j!} \leq e^{|w|}\frac{|w|^{n_{T} + 1}}{(n_{T} + 1)!} \leq \exp( - ct^{2}\log{t}).
\end{equation}
Indeed, the logarithm of $e^{|w|}|w|^{n_{T} + 1}/(n_{T} + 1)!$ in \eqref{an:detail:conditional.taylor} is at most $|w| + (n_{T} + 1)\log(e|w|/(n_{T} + 1))$, and $1 + 3\e_{0}/2 < 2$ makes its second term negative of order $t^{2}\log{t}$.

Write $P_{T,\bb{z}}(\bb{y}) \leqdef \sum_{j = 0}^{n_{T}}(\bb{z}\cdot\bb{y})^{j}/j!$. Minkowski's inequality, \eqref{an:haar.moments} in Lemma~\ref{an:polynomial} applied separately to $\Re(\bb{z})$ and $\Im(\bb{z})$, and the factorial bound $j! \geq (j/e)^{j}$ for $j \in \N$ give, after changing $C$,
\begin{equation}\label{an:detail:conditional.polynomial.norm}
\begin{aligned}
\left(\EE_{\mathrm{H}}[|P_{T,\bb{z}}(\bb{Y}_{T}^{\mathrm{H}})|^{2}]\right)^{1/2}
&\leq 1 + \sum_{j = 1}^{n_{T}}\frac{1}{j!}\left(\EE_{\mathrm{H}}[|\bb{z}\cdot\bb{Y}_{T}^{\mathrm{H}}|^{2j}]\right)^{1/2}
\leq 1 + \sum_{j = 1}^{n_{T}}\frac{j^{j/2}}{j!}(Ct\|\bb{z}\|_{2}^{2})^{j/2} \\
&\leq \sum_{j = 0}^{\infty}\frac{(Ct\|\bb{z}\|_{2}^{2})^{j/2}}{\sqrt{j!}}
\leq \left(\sum_{j = 0}^{\infty}2^{-j}\right)^{1/2}\left(\sum_{j = 0}^{\infty}\frac{(2Ct\|\bb{z}\|_{2}^{2})^{j}}{j!}\right)^{1/2} \\
&= \sqrt{2}\exp(Ct\|\bb{z}\|_{2}^{2}) \leq \exp\{C(1 + \|\bb{z}\|_{2}^{2})t\} \leq \exp(Ct^{1 + \e_{0}}).
\end{aligned}
\end{equation}
The fourth inequality is Cauchy--Schwarz. To obtain the corresponding Haar $L^{4}$ bound, apply Minkowski's inequality in $L^{4}$ and use \eqref{an:haar.moments} in Lemma~\ref{an:polynomial} at moment order $4j$, separating real and imaginary parts. The same summation estimate as in \eqref{an:detail:conditional.polynomial.norm} gives
\begin{equation}\label{an:detail:conditional.polynomial.fourth.norm}
\begin{aligned}
\left(\EE_{\mathrm{H}}[|P_{T,\bb{z}}(\bb{Y}_{T}^{\mathrm{H}})|^{4}]\right)^{1/4}
&\leq 1 + \sum_{j = 1}^{n_{T}}\frac{1}{j!}\left(\EE_{\mathrm{H}}[|\bb{z}\cdot\bb{Y}_{T}^{\mathrm{H}}|^{4j}]\right)^{1/4} \\
&\leq 1 + \sum_{j = 1}^{n_{T}}\frac{j^{j/2}}{j!}(Ct\|\bb{z}\|_{2}^{2})^{j/2}
\leq \exp(Ct^{1 + \e_{0}}).
\end{aligned}
\end{equation}
This fourth-moment bound is needed to transfer the Haar $L^{2}$ bound to the arithmetic law. Indeed, $|P_{T,\bb{z}}(\bb{y})|^{2}$ is a polynomial in the real coordinates of $\bb{y}$ of degree at most $2n_{T}$, and its Haar $L^{2}$ norm is $(\EE_{\mathrm{H}}[|P_{T,\bb{z}}(\bb{Y}_{T}^{\mathrm{H}})|^{4}])^{1/2}$. Applying \eqref{an:polynomial.bound} in Lemma~\ref{an:polynomial} to this squared modulus, taking square roots, and using $\sqrt{a + b} \leq \sqrt{a} + \sqrt{b}$ for $a,b \geq 0$ therefore give
\begin{equation}\label{an:detail:conditional.polynomial.arithmetic.norm}
\begin{aligned}
\left(\EE_{\mathrm{A}}[|P_{T,\bb{z}}(\bb{Y}_{T}^{\mathrm{A}})|^{2}]\right)^{1/2}
&\leq \left(\EE_{\mathrm{H}}[|P_{T,\bb{z}}(\bb{Y}_{T}^{\mathrm{H}})|^{2}]\right)^{1/2} + \left(\frac{CX^{4n_{T}}}{T}\right)^{1/2}\left(\EE_{\mathrm{H}}[|P_{T,\bb{z}}(\bb{Y}_{T}^{\mathrm{H}})|^{4}]\right)^{1/4} \\
&\leq \exp(Ct^{1 + \e_{0}}).
\end{aligned}
\end{equation}
The last inequality uses \eqref{an:detail:conditional.polynomial.norm}, \eqref{an:detail:conditional.polynomial.fourth.norm}, and $X^{4n_{T}}/T = T^{-1 + \oo(1)}$, which follows from \eqref{eq:fixed.parameters} because $n_{T} = \oo(t^{\phi})$.

Applying \eqref{an:polynomial.bound} in Lemma~\ref{an:polynomial} directly to $P_{T,\bb{z}}$, whose degree is at most $n_{T}$, and using \eqref{an:detail:conditional.polynomial.norm} give, for each $\star \in \{\mathrm{A},\mathrm{H}\}$,
\begin{equation}\label{an:detail:conditional.polynomial.comparison}
\begin{aligned}
\left|\EE_{\star}\left[P_{T,\bb{z}}(\bb{Y}_{T}^{\star})\right] - \EE_{\mathrm{H}}\left[P_{T,\bb{z}}(\bb{Y}_{T}^{\mathrm{H}})\right]\right|
&\leq \frac{CX^{2n_{T}}}{T}\left(\EE_{\mathrm{H}}[|P_{T,\bb{z}}(\bb{Y}_{T}^{\mathrm{H}})|^{2}]\right)^{1/2} \\
&\leq T^{-1 + \oo(1)}\exp(Ct^{1 + \e_{0}}).
\end{aligned}
\end{equation}
For $\star = \mathrm{H}$, the left-hand side is zero.

For each $\star \in \{\mathrm{A},\mathrm{H}\}$, Cauchy--Schwarz, \eqref{an:exceptional}, and the $L^{2}$ bounds \eqref{an:detail:conditional.polynomial.norm} and \eqref{an:detail:conditional.polynomial.arithmetic.norm} give
\begin{equation}\label{an:detail:polynomial.box.error}
\left|\EE_{\star}\left[P_{T,\bb{z}}(\bb{Y}_{T}^{\star})\ind_{(\mathcal{A}_{T}^{\star})^{c}}\right]\right|
\leq \exp\{ - ct^{1 + 2\e_{0}} + Ct^{1 + \e_{0}}\} \leq \exp( - c't^{1 + 2\e_{0}})
\end{equation}
for some $c' > 0$ and all sufficiently large $T$. To compare $\smash{\EE_{\star}[e^{\bb{z}\cdot\bb{Y}_{T}^{\star}}\ind_{\mathcal{A}_{T}^{\star}}]}$ with $\smash{\EE_{\mathrm{H}}[e^{\bb{z}\cdot\bb{Y}_{T}^{\mathrm{H}}}]}$, the contribution from $(\mathcal{A}_{T}^{\mathrm{H}})^{c}$ must also be controlled. Since $\bb{Y}_{T}^{\mathrm{H}}$ is real-valued, the modulus of its complex exponential depends only on the real parts of $\bb{z}$. Applying \eqref{an:haar.moments} in Lemma~\ref{an:polynomial} to the real vector $2\Re(\bb{z})$ gives
\begin{equation}\label{an:detail:haar.exponential.norm}
\left(\EE_{\mathrm{H}}\left[\big|e^{\bb{z}\cdot\bb{Y}_{T}^{\mathrm{H}}}\big|^{2}\right]\right)^{1/2}
= \left(\EE_{\mathrm{H}}\left[e^{2\Re(\bb{z})\cdot\bb{Y}_{T}^{\mathrm{H}}}\right]\right)^{1/2}
\leq \exp(Ct\|\Re(\bb{z})\|_{2}^{2}) \leq e^{Ct},
\end{equation}
where the last inequality uses $|\Re(z_{r})| \leq M$ and the fixed dimension $n_{\mathrm{dim}}$. Cauchy--Schwarz, \eqref{an:exceptional}, and \eqref{an:detail:haar.exponential.norm} therefore yield
\begin{equation}\label{an:detail:haar.exponential.tail}
\left|\EE_{\mathrm{H}}\left[e^{\bb{z}\cdot\bb{Y}_{T}^{\mathrm{H}}}\ind_{(\mathcal{A}_{T}^{\mathrm{H}})^{c}}\right]\right|
\leq e^{Ct}\PP_{\mathrm{H}}((\mathcal{A}_{T}^{\mathrm{H}})^{c})^{1/2}
\leq \exp\{ - ct^{1 + 2\e_{0}} + Ct\},
\end{equation}
after reducing $c$.

Finally, by the triangle inequality,
\begin{equation}\label{an:detail:conditional.decomposition}
\begin{aligned}
&\left|\EE_{\star}\left[e^{\bb{z}\cdot\bb{Y}_{T}^{\star}}\ind_{\mathcal{A}_{T}^{\star}}\right] - \EE_{\mathrm{H}}\left[e^{\bb{z}\cdot\bb{Y}_{T}^{\mathrm{H}}}\right]\right| \\
&\qquad \leq \left|\EE_{\star}\left[\left(e^{\bb{z}\cdot\bb{Y}_{T}^{\star}} - P_{T,\bb{z}}(\bb{Y}_{T}^{\star})\right)\ind_{\mathcal{A}_{T}^{\star}}\right]\right|
+ \left|\EE_{\star}\left[P_{T,\bb{z}}(\bb{Y}_{T}^{\star})\ind_{(\mathcal{A}_{T}^{\star})^{c}}\right]\right| \\
&\qquad + \left|\EE_{\star}\left[P_{T,\bb{z}}(\bb{Y}_{T}^{\star})\right] - \EE_{\mathrm{H}}\left[P_{T,\bb{z}}(\bb{Y}_{T}^{\mathrm{H}})\right]\right|
+ \left|\EE_{\mathrm{H}}\left[P_{T,\bb{z}}(\bb{Y}_{T}^{\mathrm{H}})\ind_{(\mathcal{A}_{T}^{\mathrm{H}})^{c}}\right]\right| \\
&\qquad + \left|\EE_{\mathrm{H}}\left[\left(P_{T,\bb{z}}(\bb{Y}_{T}^{\mathrm{H}}) - e^{\bb{z}\cdot\bb{Y}_{T}^{\mathrm{H}}}\right)\ind_{\mathcal{A}_{T}^{\mathrm{H}}}\right]\right|
+ \left|\EE_{\mathrm{H}}\left[e^{\bb{z}\cdot\bb{Y}_{T}^{\mathrm{H}}}\ind_{(\mathcal{A}_{T}^{\mathrm{H}})^{c}}\right]\right|.
\end{aligned}
\end{equation}
The first and fifth terms on the right-hand side of \eqref{an:detail:conditional.decomposition} are controlled by the exponential approximation in \eqref{an:detail:conditional.taylor}, the second and fourth by \eqref{an:detail:polynomial.box.error}, the third by \eqref{an:detail:conditional.polynomial.comparison}, and the sixth by \eqref{an:detail:haar.exponential.tail}, proving \eqref{an:detail:conditional.transform} after adjusting $c$ and $C$.
\end{proof}

\begin{lemma}[Simultaneous spatial truncated transforms]\label{an:spatial.transform}
Let $\theta > 0$. Fix $0 < \e_{0} < 1/4$ and $B_{0},C_{0} > 0$. Let $m_{T} \in \N$ satisfy $m_{T} \leq t^{B_{0}}$. For $1 \leq m \leq m_{T}$, let
\[
Y_{T,m}^{\mathrm{A}}(h) \leqdef \Re\sum_{p^{j} \leq X}a_{m,p,j}(h)p^{-\ii j\tau}, \qquad
Y_{T,m}^{\mathrm{H}}(h) \leqdef \Re\sum_{p^{j} \leq X}a_{m,p,j}(h)U_{p}^{j},
\]
where the deterministic complex coefficients are the same under both laws, are continuously differentiable on $I_{\theta}$, and satisfy $\sup_{h \in I_{\theta}}\left\{|a_{m,p,j}(h)| + \ell_{T}|a_{m,p,j}'(h)|\right\} \leq C_{0}p^{-j/2}$. For each $\star \in \{\mathrm{A},\mathrm{H}\}$, define the common spatial event
\[
\mathcal{A}_{T,\mathrm{sp}}^{\star} \leqdef \left\{\max_{1 \leq m \leq m_{T}}\sup_{h \in I_{\theta}}|Y_{T,m}^{\star}(h)| \leq t^{1 + \e_{0}}\right\}.
\]
There are constants $c,C > 0$ such that, for all sufficiently large $T$,
\begin{equation}\label{an:spatial.exceptional}
\PP_{\star}((\mathcal{A}_{T,\mathrm{sp}}^{\star})^{c}) \leq \exp\{ - ct^{1 + 2\e_{0}} + Ct\}.
\end{equation}

Fix $n_{\mathrm{dim}},L_{0} \in \N$, $C_{1} \geq 1$, and $M > 0$, independently of $T$. For each $1 \leq r \leq n_{\mathrm{dim}}$, choose an integer $1 \leq L_{r} \leq L_{0}$, indices $m_{r,q} \in \{1,\ldots,m_{T}\}$, positions $h_{r,q} \in I_{\theta}$, and real coefficients $c_{r,q}$, for $1 \leq q \leq L_{r}$, such that $\smash{\sum_{q = 1}^{L_{r}}|c_{r,q}| \leq C_{1}}$. All these choices are deterministic, may depend on $T$, and are the same under both laws. Define $\bb{Z}_{T}^{\star} \leqdef (Z_{T,1}^{\star},\ldots,Z_{T,n_{\mathrm{dim}}}^{\star})$ by
\[
Z_{T,r}^{\star} \leqdef \sum_{q = 1}^{L_{r}}c_{r,q}Y_{T,m_{r,q}}^{\star}(h_{r,q}), \qquad 1 \leq r \leq n_{\mathrm{dim}}.
\]
A single evaluation is included by taking $L_{r} = 1$ and $c_{r,1} = 1$. Uniformly over these choices and over $\bb{z} \in \C^{n_{\mathrm{dim}}}$ with $|\Re(z_{r})| \leq M$ and $|\Im(z_{r})| \leq t^{\e_{0}/2}$ for every $r$,
\begin{equation}\label{an:spatial.transform.bound}
\left|\EE_{\star}\left[e^{\bb{z}\cdot\bb{Z}_{T}^{\star}}\ind_{\mathcal{A}_{T,\mathrm{sp}}^{\star}}\right] - \EE_{\mathrm{H}}\left[e^{\bb{z}\cdot\bb{Z}_{T}^{\mathrm{H}}}\right]\right| \leq \exp\{ - ct^{1 + 2\e_{0}} + Ct\}.
\end{equation}
The constants may depend on $\e_{0}$, $B_{0}$, $C_{0}$, $n_{\mathrm{dim}}$, $L_{0}$, $C_{1}$, $M$, and the fixed parameters of the paper. Once the family $(Y_{T,m}^{\star})_{1 \leq m \leq m_{T}}$ is fixed, the same event $\mathcal{A}_{T,\mathrm{sp}}^{\star}$ is used for every permitted choice of the indices $m_{r,q}$, positions $h_{r,q}$, and coefficients $c_{r,q}$.
\end{lemma}

\begin{proof}[Proof of Lemma~\ref{an:spatial.transform}]
By the assumed coefficient bounds, Lemma~\ref{an:polynomial} applies to both $Y_{T,m}^{\star}(h)$ and $\ell_{T}(Y_{T,m}^{\star})'(h)$, with constants uniform in $m$ and $h$.

Recall $N_{T} = |I_{\theta}|/\ell_{T}$ from \eqref{eq:N.T}, and cover $I_{\theta}$ by at most $1 + N_{T}$ intervals of positive length at most $\ell_{T}$. Fix one such interval $B$, let $h_{B}$ be its left endpoint, and write $Y^{\star}(h) \leqdef Y_{T,m}^{\star}(h)$ for one member of the family. The fundamental theorem of calculus gives
\begin{equation}\label{an:detail:cell.calculus}
\sup_{h \in B}|Y^{\star}(h)| \leq |Y^{\star}(h_{B})| + \ell_{T}^{-1}\int_{B}|\ell_{T}(Y^{\star})'(x)|\rd{x}.
\end{equation}
For every $k \in \N$, Jensen's inequality with respect to the probability measure $\rd{x}/|B|$ on $B$ yields
\begin{equation}\label{an:spatial.detail:cell.jensen}
\begin{aligned}
\left(\ell_{T}^{-1}\int_{B}|\ell_{T}(Y^{\star})'(x)|\rd{x}\right)^{2k}
&\leq \left(\frac{|B|}{\ell_{T}}\right)^{2k - 1}\ell_{T}^{-1}\int_{B}|\ell_{T}(Y^{\star})'(x)|^{2k}\rd{x} \\
&\leq \ell_{T}^{-1}\int_{B}|\ell_{T}(Y^{\star})'(x)|^{2k}\rd{x},
\end{aligned}
\end{equation}
where the last inequality uses $|B| \leq \ell_{T}$. Combining \eqref{an:detail:cell.calculus}, \eqref{an:spatial.detail:cell.jensen}, and $(a + b)^{2k} \leq 2^{2k - 1}(a^{2k} + b^{2k})$, and then applying \eqref{an:haar.moments} in Lemma~\ref{an:polynomial} or its arithmetic counterpart, gives
\begin{equation}\label{an:detail:conditional.cell.moments}
\begin{aligned}
\EE_{\star}\left[\sup_{h \in B}|Y^{\star}(h)|^{2k}\right]
&\leq 2^{2k - 1}\left(\EE_{\star}[|Y^{\star}(h_{B})|^{2k}] + \ell_{T}^{-1}\int_{B}\EE_{\star}[|\ell_{T}(Y^{\star})'(x)|^{2k}]\rd{x}\right) \\
&\leq (Ckt)^{k}.
\end{aligned}
\end{equation}
The bound is uniform in $m$ and $B$, and holds under the arithmetic law whenever $k = \oo(t^{\phi})$. Let $k \leqdef \lfloor c_{0}t^{1 + 2\e_{0}}\rfloor$, where $c_{0} > 0$ is sufficiently small. This choice is admissible because $1 + 2\e_{0} < 3/2 < \phi$. Markov's inequality and \eqref{an:detail:conditional.cell.moments} imply
\begin{equation}\label{an:spatial.detail:cell.tail}
\PP_{\star}\left(\sup_{h \in B}|Y_{T,m}^{\star}(h)| > t^{1 + \e_{0}}\right)
\leq \frac{(Ckt)^{k}}{t^{2k(1 + \e_{0})}}
= \left(\frac{Ck}{t^{1 + 2\e_{0}}}\right)^{k}
\leq \exp( - ct^{1 + 2\e_{0}}).
\end{equation}
A union bound using \eqref{an:spatial.detail:cell.tail} over the intervals and the at most $t^{B_{0}}$ family members gives
\[
\PP_{\star}((\mathcal{A}_{T,\mathrm{sp}}^{\star})^{c})
\leq (1 + N_{T})t^{B_{0}}\exp( - ct^{1 + 2\e_{0}})
\leq \exp\{ - ct^{1 + 2\e_{0}} + Ct\},
\]
which proves \eqref{an:spatial.exceptional}. The last inequality uses $\log(1 + N_{T}) = \OO(t)$, a consequence of \eqref{eq:effective.cutoff.bounds} and \eqref{eq:parameter.scales}. Thus the term $Ct$ accounts for the number of spatial intervals and family members.

Fix the evaluations and linear combinations in the statement. Since the coefficients $c_{r,q}$ are real, each $Z_{T,r}^{\star}$ has the required prime-power representation with coefficients
\begin{equation}\label{an:spatial.detail:linear.coefficients}
\widetilde{a}_{r,p,j} \leqdef \sum_{q = 1}^{L_{r}}c_{r,q}a_{m_{r,q},p,j}(h_{r,q}), \qquad
|\widetilde{a}_{r,p,j}| \leq C_{0}p^{-j/2}\sum_{q = 1}^{L_{r}}|c_{r,q}| \leq C_{0}C_{1}p^{-j/2}.
\end{equation}
Consequently, the polynomial and moment estimates in Lemma~\ref{an:polynomial} apply to $\bb{Z}_{T}^{\star}$ uniformly over all permitted choices.

Let $n_{T} \leqdef \lfloor t^{2}\rfloor$ and define $P_{T,\bb{z}}(\bb{y}) \leqdef \sum_{j = 0}^{n_{T}}(\bb{z}\cdot\bb{y})^{j}/j!$. On $\mathcal{A}_{T,\mathrm{sp}}^{\star}$, the bound $\sum_{q = 1}^{L_{r}}|c_{r,q}| \leq C_{1}$ gives $|Z_{T,r}^{\star}| \leq C_{1}t^{1 + \e_{0}}$. Hence $|\bb{z}\cdot\bb{Z}_{T}^{\star}| \leq Ct^{1 + 3\e_{0}/2} = \oo(n_{T})$ on this event. The exponential Taylor remainder estimate \eqref{an:detail:conditional.taylor} therefore gives
\begin{equation}\label{an:spatial.detail:taylor}
\left|e^{\bb{z}\cdot\bb{Z}_{T}^{\star}} - P_{T,\bb{z}}(\bb{Z}_{T}^{\star})\right|\ind_{\mathcal{A}_{T,\mathrm{sp}}^{\star}}
\leq \exp( - ct^{2}\log{t}).
\end{equation}

The global Haar polynomial estimates \eqref{an:detail:conditional.polynomial.norm} and \eqref{an:detail:conditional.polynomial.fourth.norm} depend only on the coefficient bound and therefore apply with $\bb{Y}_{T}^{\mathrm{H}}$ replaced by $\bb{Z}_{T}^{\mathrm{H}}$, using \eqref{an:spatial.detail:linear.coefficients}. After changing $C$, they give
\begin{equation}\label{an:spatial.detail:haar.polynomial.norms}
\left(\EE_{\mathrm{H}}[|P_{T,\bb{z}}(\bb{Z}_{T}^{\mathrm{H}})|^{2}]\right)^{1/2}
+ \left(\EE_{\mathrm{H}}[|P_{T,\bb{z}}(\bb{Z}_{T}^{\mathrm{H}})|^{4}]\right)^{1/4}
\leq \exp(Ct^{1 + \e_{0}}).
\end{equation}
To obtain the arithmetic $L^{2}$ bound, apply \eqref{an:polynomial.bound} in Lemma~\ref{an:polynomial} to $|P_{T,\bb{z}}|^{2}$, whose degree is at most $2n_{T}$. Taking square roots, using $\sqrt{a + b} \leq \sqrt{a} + \sqrt{b}$ for $a,b \geq 0$, and applying \eqref{an:spatial.detail:haar.polynomial.norms} yield
\begin{equation}\label{an:spatial.detail:polynomial.norm}
\begin{aligned}
\left(\EE_{\star}[|P_{T,\bb{z}}(\bb{Z}_{T}^{\star})|^{2}]\right)^{1/2}
&\leq \left(\EE_{\mathrm{H}}[|P_{T,\bb{z}}(\bb{Z}_{T}^{\mathrm{H}})|^{2}]\right)^{1/2}
+ \left(\frac{CX^{4n_{T}}}{T}\right)^{1/2}\left(\EE_{\mathrm{H}}[|P_{T,\bb{z}}(\bb{Z}_{T}^{\mathrm{H}})|^{4}]\right)^{1/4} \\
&\leq \exp(Ct^{1 + \e_{0}})
\end{aligned}
\end{equation}
for each $\star \in \{\mathrm{A},\mathrm{H}\}$. For $\star = \mathrm{H}$, the first inequality is trivial. The last inequality uses $X^{4n_{T}}/T = T^{-1 + \oo(1)}$, which follows from \eqref{eq:fixed.parameters} because $n_{T} = \oo(t^{\phi})$. The Haar $L^{4}$ estimate in \eqref{an:spatial.detail:haar.polynomial.norms} thus supplies the error control needed to obtain the arithmetic $L^{2}$ bound.

Applying \eqref{an:polynomial.bound} in Lemma~\ref{an:polynomial} directly to $P_{T,\bb{z}}$ and using \eqref{an:spatial.detail:haar.polynomial.norms} give
\begin{equation}\label{an:spatial.detail:polynomial.comparison}
\begin{aligned}
\left|\EE_{\star}\left[P_{T,\bb{z}}(\bb{Z}_{T}^{\star})\right] - \EE_{\mathrm{H}}\left[P_{T,\bb{z}}(\bb{Z}_{T}^{\mathrm{H}})\right]\right|
&\leq \frac{CX^{2n_{T}}}{T}\left(\EE_{\mathrm{H}}[|P_{T,\bb{z}}(\bb{Z}_{T}^{\mathrm{H}})|^{2}]\right)^{1/2} \\
&\leq T^{-1 + \oo(1)}\exp(Ct^{1 + \e_{0}}).
\end{aligned}
\end{equation}
For $\star = \mathrm{H}$, the left-hand side is zero.

The restriction to the common spatial event can now be removed from the polynomial expectations. By Cauchy--Schwarz, \eqref{an:spatial.exceptional}, and \eqref{an:spatial.detail:polynomial.norm},
\begin{equation}\label{an:spatial.detail:polynomial.tail}
\begin{aligned}
\left|\EE_{\star}\left[P_{T,\bb{z}}(\bb{Z}_{T}^{\star})\ind_{(\mathcal{A}_{T,\mathrm{sp}}^{\star})^{c}}\right]\right|
&\leq \left(\EE_{\star}[|P_{T,\bb{z}}(\bb{Z}_{T}^{\star})|^{2}]\right)^{1/2}\PP_{\star}((\mathcal{A}_{T,\mathrm{sp}}^{\star})^{c})^{1/2} \\
&\leq \exp\{ - ct^{1 + 2\e_{0}} + Ct^{1 + \e_{0}} + Ct\}
\leq \exp( - c't^{1 + 2\e_{0}})
\end{aligned}
\end{equation}
for some $c' > 0$ and all sufficiently large $T$.

To compare $\smash{\EE_{\star}[e^{\bb{z}\cdot\bb{Z}_{T}^{\star}}\ind_{\mathcal{A}_{T,\mathrm{sp}}^{\star}}]}$ with $\smash{\EE_{\mathrm{H}}[e^{\bb{z}\cdot\bb{Z}_{T}^{\mathrm{H}}}]}$, the contribution from $(\mathcal{A}_{T,\mathrm{sp}}^{\mathrm{H}})^{c}$ must also be controlled. Applying \eqref{an:haar.moments} in Lemma~\ref{an:polynomial} to the real vector $2\Re(\bb{z})$, and using $|\Re(z_{r})| \leq M$ and the fixed dimension $n_{\mathrm{dim}}$, gives
\begin{equation}\label{an:spatial.detail:haar.exponential.norm}
\left(\EE_{\mathrm{H}}\left[\big|e^{\bb{z}\cdot\bb{Z}_{T}^{\mathrm{H}}}\big|^{2}\right]\right)^{1/2}
= \left(\EE_{\mathrm{H}}\left[e^{2\Re(\bb{z})\cdot\bb{Z}_{T}^{\mathrm{H}}}\right]\right)^{1/2}
\leq \exp(Ct\|\Re(\bb{z})\|_{2}^{2}) \leq e^{Ct}.
\end{equation}
Cauchy--Schwarz, \eqref{an:spatial.exceptional}, and \eqref{an:spatial.detail:haar.exponential.norm} therefore imply
\begin{equation}\label{an:spatial.detail:haar.exponential.tail}
\left|\EE_{\mathrm{H}}\left[e^{\bb{z}\cdot\bb{Z}_{T}^{\mathrm{H}}}\ind_{(\mathcal{A}_{T,\mathrm{sp}}^{\mathrm{H}})^{c}}\right]\right|
\leq e^{Ct}\PP_{\mathrm{H}}((\mathcal{A}_{T,\mathrm{sp}}^{\mathrm{H}})^{c})^{1/2}
\leq \exp\{ - ct^{1 + 2\e_{0}} + Ct\},
\end{equation}
after adjusting $c$ and $C$.

Finally, the triangle inequality gives
\begin{equation}\label{an:spatial.detail:decomposition}
\begin{aligned}
&\left|\EE_{\star}\left[e^{\bb{z}\cdot\bb{Z}_{T}^{\star}}\ind_{\mathcal{A}_{T,\mathrm{sp}}^{\star}}\right] - \EE_{\mathrm{H}}\left[e^{\bb{z}\cdot\bb{Z}_{T}^{\mathrm{H}}}\right]\right| \\
&\qquad \leq \left|\EE_{\star}\left[\left(e^{\bb{z}\cdot\bb{Z}_{T}^{\star}} - P_{T,\bb{z}}(\bb{Z}_{T}^{\star})\right)\ind_{\mathcal{A}_{T,\mathrm{sp}}^{\star}}\right]\right|
+ \left|\EE_{\star}\left[P_{T,\bb{z}}(\bb{Z}_{T}^{\star})\ind_{(\mathcal{A}_{T,\mathrm{sp}}^{\star})^{c}}\right]\right| \\
&\qquad + \left|\EE_{\star}\left[P_{T,\bb{z}}(\bb{Z}_{T}^{\star})\right] - \EE_{\mathrm{H}}\left[P_{T,\bb{z}}(\bb{Z}_{T}^{\mathrm{H}})\right]\right|
+ \left|\EE_{\mathrm{H}}\left[P_{T,\bb{z}}(\bb{Z}_{T}^{\mathrm{H}})\ind_{(\mathcal{A}_{T,\mathrm{sp}}^{\mathrm{H}})^{c}}\right]\right| \\
&\qquad + \left|\EE_{\mathrm{H}}\left[\left(P_{T,\bb{z}}(\bb{Z}_{T}^{\mathrm{H}}) - e^{\bb{z}\cdot\bb{Z}_{T}^{\mathrm{H}}}\right)\ind_{\mathcal{A}_{T,\mathrm{sp}}^{\mathrm{H}}}\right]\right|
+ \left|\EE_{\mathrm{H}}\left[e^{\bb{z}\cdot\bb{Z}_{T}^{\mathrm{H}}}\ind_{(\mathcal{A}_{T,\mathrm{sp}}^{\mathrm{H}})^{c}}\right]\right|.
\end{aligned}
\end{equation}
The first and fifth terms on the right-hand side of \eqref{an:spatial.detail:decomposition} are controlled by the exponential approximation in \eqref{an:spatial.detail:taylor}, the second and fourth by \eqref{an:spatial.detail:polynomial.tail}, the third by \eqref{an:spatial.detail:polynomial.comparison}, and the sixth by \eqref{an:spatial.detail:haar.exponential.tail}, proving \eqref{an:spatial.transform.bound} after adjusting $c$ and $C$.
\end{proof}

\begin{lemma}[Conditional Laplace estimates]\label{an:conditional}
Let $\theta > 0$. Fix $0 < \e_{0} < 1/4$, $C_{0} > 0$, and $k \in \N$. Let $\bb{Y}_{T}^{\mathrm{A}}$ and $\bb{Y}_{T}^{\mathrm{H}}$ be as in Lemma~\ref{an:transform} with $n_{\mathrm{dim}} = k + 1$, with $Y_{T,r}^{\star} = F_{h_{r}}^{\star}$ for $1 \leq r \leq k$ and $\star \in \{\mathrm{A},\mathrm{H}\}$. Let $\mathcal{A}_{T}^{\star}$ be the event in that lemma for these vectors. Define
\[
\Sigma_{T} \leqdef \left(\Cov_{\mathrm{H}}(Y_{T,r}^{\mathrm{H}},Y_{T,s}^{\mathrm{H}})\right)_{1 \leq r,s \leq k}, \quad
\bb{c}_{T} \leqdef \left(\Cov_{\mathrm{H}}(Y_{T,r}^{\mathrm{H}},Y_{T,k + 1}^{\mathrm{H}})\right)_{1 \leq r \leq k}, \quad
\sig_{T}^{2} \leqdef \Var_{\mathrm{H}}(Y_{T,k + 1}^{\mathrm{H}}).
\]
Assume every eigenvalue of $\Sigma_{T}$ lies in $[c_{\Sigma}t,C_{\Sigma}t]$ for fixed $0 < c_{\Sigma} \leq C_{\Sigma} < \infty$. Fix a compact set $K \subseteq \R^{k}$ and $M > 0$. Then there is a constant $C > 0$ such that, for all sufficiently large $T$, uniformly over $\bb{v}/s_{T} \in K$ and real $|\lambda_{k + 1}| \leq M$,
\begin{equation}\label{an:conditional.bound}
\begin{aligned}
&\EE_{\star}\left[e^{\lambda_{k + 1}Y_{T,k + 1}^{\star}}\ind_{\mathcal{A}_{T}^{\star}}\,\middle|\,(Y_{T,1}^{\star},\ldots,Y_{T,k}^{\star}) \in \bb{v} + [0,1]^{k}\right] \\
&\qquad \leq C\exp\left\{\lambda_{k + 1}\bb{c}_{T}^{\mathsf{T}}\Sigma_{T}^{-1}\bb{v} + \frac{\lambda_{k + 1}^{2}}{2}\left(\sig_{T}^{2} - \bb{c}_{T}^{\mathsf{T}}\Sigma_{T}^{-1}\bb{c}_{T}\right)\right\}.
\end{aligned}
\end{equation}
If the vector $\bb{Y}_{T}^{\star}$ also has the representation of $\bb{Z}_{T}^{\star}$ in Lemma~\ref{an:spatial.transform} with $n_{\mathrm{dim}} = k + 1$, so that each coordinate is a permitted finite real linear combination of evaluations $Y_{T,m}^{\star}(h)$, then \eqref{an:conditional.bound} holds with $\mathcal{A}_{T}^{\star}$ replaced by $\mathcal{A}_{T,\mathrm{sp}}^{\star}$. The conditioning events have positive probability for all sufficiently large $T$. The constant $C$ in \eqref{an:conditional.bound} is uniform over all positions and coefficients satisfying the hypotheses and may depend on $\e_{0}$, $C_{0}$, $k$, $K$, $M$, $c_{\Sigma}$, $C_{\Sigma}$, and the fixed parameters of the paper. For the version with $\mathcal{A}_{T,\mathrm{sp}}^{\star}$ in place of $\mathcal{A}_{T}^{\star}$, it may additionally depend on $B_{0}$, $L_{0}$, and $C_{1}$ from Lemma~\ref{an:spatial.transform}. Estimate \eqref{an:conditional.bound} holds with $[0,1]^{k}$ replaced by $\prod_{j = 1}^{k}[0,b_{j}]$, uniformly for $0 < b_{-} \leq b_{j} \leq b_{+} < \infty$, with additional dependence only on the fixed bounds $b_{-},b_{+}$. In the Haar case, \eqref{an:conditional.bound} also holds with $\ind_{\mathcal{A}_{T}^{\mathrm{H}}}$ replaced by $1$ in the finite-dimensional version and with $\smash{\ind_{\mathcal{A}_{T,\mathrm{sp}}^{\mathrm{H}}}}$ replaced by $1$ in the spatial version.
\end{lemma}

\begin{proof}[Proof of Lemma~\ref{an:conditional}]
For each prime $p \leq X$ and $1 \leq r \leq k + 1$, define
\[
B_{T,p,r} \leqdef \Re\sum_{j:p^{j}\leq X}a_{r,p,j}U_{p}^{j}, \qquad \bb{B}_{T,p} \leqdef (B_{T,p,1},\ldots,B_{T,p,k + 1}).
\]
Thus $\bb{Y}_{T}^{\mathrm{H}} = \sum_{p \leq X}\bb{B}_{T,p}$, and the vectors $\bb{B}_{T,p}$ are independent because each depends only on $U_{p}$. The coefficient assumption inherited from Lemma~\ref{an:transform} gives
\begin{equation}\label{an:conditional.detail:prime.bound}
|B_{T,p,r}| \leq \sum_{j:p^{j}\leq X}|a_{r,p,j}| \leq C_{0}\sum_{j = 1}^{\infty}p^{-j/2} = \frac{C_{0}p^{-1/2}}{1 - p^{-1/2}} \leq Cp^{-1/2}, \qquad \|\bb{B}_{T,p}\|_{2} \leq Cp^{-1/2}.
\end{equation}
Moreover, $\EE_{\mathrm{H}}[B_{T,p,r}] = 0$, since $\EE_{\mathrm{H}}[U_{p}^{j}] = 0$ for every $j \in \N$.

Let $\bb{\eta} \leqdef (\lambda_{1},\ldots,\lambda_{k},\lambda_{k + 1}) \in \R^{k + 1}$, and define the prime and joint log moment generating functions by
\[
\mathcal{L}_{T,p}(\bb{\eta}) \leqdef \log\EE_{\mathrm{H}}[e^{\bb{\eta}\cdot\bb{B}_{T,p}}], \qquad
\mathcal{L}_{T}(\bb{\lambda},\lambda_{k + 1}) \leqdef \log\EE_{\mathrm{H}}[e^{\sum_{r = 1}^{k}\lambda_{r}Y_{T,r}^{\mathrm{H}} + \lambda_{k + 1}Y_{T,k + 1}^{\mathrm{H}}}] = \sum_{p \leq X}\mathcal{L}_{T,p}(\bb{\eta}).
\]
For a function $g$ of the prime vector, write
\[
\EE_{p}^{\bb{\eta}}[g(\bb{B}_{T,p})] \leqdef \frac{\EE_{\mathrm{H}}[g(\bb{B}_{T,p})e^{\bb{\eta}\cdot\bb{B}_{T,p}}]}{\EE_{\mathrm{H}}[e^{\bb{\eta}\cdot\bb{B}_{T,p}}]}.
\]
The denominator is positive for every real $\bb{\eta}$. Differentiation under the expectation is justified by \eqref{an:conditional.detail:prime.bound}. If $D^{j}$ denotes the $j$th derivative as a multilinear form, then, for $\bb{u}_{1},\bb{u}_{2},\bb{u}_{3} \in \R^{k + 1}$,
\begin{equation}\label{an:conditional.detail:prime.derivatives}
\begin{aligned}
D\mathcal{L}_{T,p}(\bb{\eta})[\bb{u}_{1}]
&= \EE_{p}^{\bb{\eta}}[\bb{u}_{1}\cdot\bb{B}_{T,p}], \\
D^{2}\mathcal{L}_{T,p}(\bb{\eta})[\bb{u}_{1},\bb{u}_{2}]
&= \Cov_{p}^{\bb{\eta}}(\bb{u}_{1}\cdot\bb{B}_{T,p},\bb{u}_{2}\cdot\bb{B}_{T,p}), \\
D^{3}\mathcal{L}_{T,p}(\bb{\eta})[\bb{u}_{1},\bb{u}_{2},\bb{u}_{3}]
&= \EE_{p}^{\bb{\eta}}\left[\prod_{a = 1}^{3}\left\{\bb{u}_{a}\cdot\left(\bb{B}_{T,p} - \EE_{p}^{\bb{\eta}}[\bb{B}_{T,p}]\right)\right\}\right].
\end{aligned}
\end{equation}
Here $\Cov_{p}^{\bb{\eta}}$ is covariance under the prime tilt just defined. The last identity and \eqref{an:conditional.detail:prime.bound} imply
\[
\left|D^{3}\mathcal{L}_{T,p}(\bb{\eta})[\bb{u}_{1},\bb{u}_{2},\bb{u}_{3}]\right| \leq Cp^{-3/2}\prod_{a = 1}^{3}\|\bb{u}_{a}\|_{2}, \qquad \sum_{p \leq X}p^{-3/2} \leq \sum_{n = 2}^{\infty}n^{-3/2} < \infty.
\]
Taylor's formula at $\bb{\eta} = \bb{0}$, the vanishing first derivatives there, and independence over primes now give
\begin{equation}\label{an:joint.cumulants}
\mathcal{L}_{T}(\bb{\lambda},\lambda_{k + 1})
= \frac{1}{2}\left(\bb{\lambda}^{\mathsf{T}}\Sigma_{T}\bb{\lambda} + 2\lambda_{k + 1}\bb{\lambda}^{\mathsf{T}}\bb{c}_{T} + \lambda_{k + 1}^{2}\sig_{T}^{2}\right) + \mathcal{R}_{T}(\bb{\lambda},\lambda_{k + 1}).
\end{equation}
More explicitly, viewing $\mathcal{R}_{T}$ as a function of $\bb{\eta}$, its integral remainder is
\[
\mathcal{R}_{T}(\bb{\eta}) \leqdef \frac{1}{2}\sum_{p \leq X}\int_{0}^{1}(1 - u)^{2}D^{3}\mathcal{L}_{T,p}(u\bb{\eta})[\bb{\eta},\bb{\eta},\bb{\eta}]\rd{u}.
\]
Taylor's formula applied also to the first and second derivatives in \eqref{an:conditional.detail:prime.derivatives} therefore yields, for every fixed $M_{1} > 0$,
\begin{equation}\label{an:conditional.detail:remainder}
\sup_{\|\bb{\eta}\|_{2}\leq M_{1}}\left\{|\mathcal{R}_{T}(\bb{\eta})| + \|\nabla \mathcal{R}_{T}(\bb{\eta})\|_{2} + \|\nabla^{2}\mathcal{R}_{T}(\bb{\eta})\|_{\mathrm{op}}\right\} \leq C_{M_{1}}\sum_{p \leq X}p^{-3/2} \leq C_{M_{1}}.
\end{equation}

Let $\PP_{\mathrm{H}}^{\bb{\lambda},\lambda_{k + 1}}$ denote the normalized Haar tilt with density
\[
\frac{\rd{\PP}_{\mathrm{H}}^{\bb{\lambda},\lambda_{k + 1}}}{\rd{\PP}_{\mathrm{H}}}
\leqdef \exp\left\{\sum_{r = 1}^{k}\lambda_{r}Y_{T,r}^{\mathrm{H}} + \lambda_{k + 1}Y_{T,k + 1}^{\mathrm{H}} - \mathcal{L}_{T}(\bb{\lambda},\lambda_{k + 1})\right\}.
\]
Since this density is the product of the normalized prime densities, the prime vectors remain independent under the tilt, with expectations $\EE_{p}^{\bb{\eta}}$. Differentiating \eqref{an:joint.cumulants} and using \eqref{an:conditional.detail:remainder} gives the following mean vector and covariance matrix for the first $k$ coordinates:
\begin{equation}\label{an:conditional.detail:tilted.mean.covariance}
\begin{aligned}
\bb{\mu}_{T} \leqdef \left(\EE_{\mathrm{H}}^{\bb{\lambda},\lambda_{k + 1}}[Y_{T,r}^{\mathrm{H}}]\right)_{1 \leq r \leq k}
&= \nabla_{\bb{\lambda}}\mathcal{L}_{T}(\bb{\lambda},\lambda_{k + 1}) = \Sigma_{T}\bb{\lambda} + \lambda_{k + 1}\bb{c}_{T} + \OO(1), \\
\left(\Cov_{\mathrm{H}}^{\bb{\lambda},\lambda_{k + 1}}(Y_{T,r}^{\mathrm{H}},Y_{T,s}^{\mathrm{H}})\right)_{1 \leq r,s \leq k}
&= \nabla_{\bb{\lambda}}^{2}\mathcal{L}_{T}(\bb{\lambda},\lambda_{k + 1}) = \Sigma_{T} + \OO(1).
\end{aligned}
\end{equation}
All the error bounds are uniform when $\bb{\lambda}$ and $\lambda_{k + 1}$ range over fixed compact sets.

To estimate probabilities of fixed bins, a bound on the tilted characteristic function is needed throughout the Fourier cube in Lemma~\ref{an:transform}. Define
\begin{equation}\label{eq:varphi}
\varphi_{T}^{\bb{\lambda},\lambda_{k + 1}}(\bb{\omega}) \leqdef \EE_{\mathrm{H}}^{\bb{\lambda},\lambda_{k + 1}}\left[e^{\ii\sum_{r = 1}^{k}\omega_{r}Y_{T,r}^{\mathrm{H}}}\right]
= \prod_{p \leq X}\EE_{p}^{\bb{\eta}}\left[e^{\ii\sum_{r = 1}^{k}\omega_{r}B_{T,p,r}}\right].
\end{equation}
The factors corresponding to primes $p \leq t^{2}$ in \eqref{eq:varphi} have modulus at most one and can therefore be discarded when bounding the modulus of the product from above. For $p > t^{2}$ and $\|\bb{\omega}\|_{\infty} \leq t^{\e_{0}/2}$, \eqref{an:conditional.detail:prime.bound} gives $|\sum_{r = 1}^{k}\omega_{r}B_{T,p,r}| \leq Ct^{\e_{0}/2 - 1} = \oo(1)$. Expanding the characteristic function of this bounded variable after centering it, and then using $\log(1 + z) = z + \OO(|z|^{2})$ near zero, gives
\[
\log\left|\EE_{p}^{\bb{\eta}}\left[e^{\ii\sum_{r = 1}^{k}\omega_{r}B_{T,p,r}}\right]\right|
= - \frac{1}{2}\Var_{p}^{\bb{\eta}}\left(\sum_{r = 1}^{k}\omega_{r}B_{T,p,r}\right) + \OO(\|\bb{\omega}\|_{2}^{3}p^{-3/2}).
\]
The total variance is $\bb{\omega}^{\mathsf{T}}\Sigma_{T}\bb{\omega} + \OO(\|\bb{\omega}\|_{2}^{2})$ by \eqref{an:conditional.detail:tilted.mean.covariance}. The contribution of primes $p \leq t^{2}$ is at most $C\|\bb{\omega}\|_{2}^{2}\sum_{n \leq t^{2}}n^{-1} \leq C\log{t}\|\bb{\omega}\|_{2}^{2}$. Consequently, for all sufficiently large $T$,
\begin{equation}\label{an:conditional.detail:large.prime.variance}
\begin{aligned}
\sum_{t^{2}<p\leq X}\Var_{p}^{\bb{\eta}}\left(\sum_{r = 1}^{k}\omega_{r}B_{T,p,r}\right)
&\geq (c_{\Sigma}t - C\log{t} - C)\|\bb{\omega}\|_{2}^{2} \geq \frac{c_{\Sigma}}{2}t\|\bb{\omega}\|_{2}^{2}, \\
\|\bb{\omega}\|_{2}^{3}\sum_{p>t^{2}}p^{-3/2}
&\leq \frac{C}{t}\|\bb{\omega}\|_{2}^{3} \leq Ct^{\e_{0}/2 - 1}\|\bb{\omega}\|_{2}^{2} = \oo(t\|\bb{\omega}\|_{2}^{2}).
\end{aligned}
\end{equation}
Multiplying the bounds for the factors with $t^{2} < p \leq X$ in \eqref{eq:varphi} and using \eqref{an:conditional.detail:large.prime.variance} therefore proves
\begin{equation}\label{an:detail:conditional.characteristic}
\left|\varphi_{T}^{\bb{\lambda},\lambda_{k + 1}}(\bb{\omega})\right| \leq \exp( - ct\|\bb{\omega}\|_{2}^{2}), \qquad \|\bb{\omega}\|_{\infty} \leq t^{\e_{0}/2}.
\end{equation}
For $\|\bb{\omega}\|_{2}$ sufficiently small, \eqref{an:conditional.detail:prime.bound} gives $|\sum_{r = 1}^{k}\omega_{r}B_{T,p,r}| \leq C\|\bb{\omega}\|_{2}p^{-1/2}$, so Taylor expansion of each prime characteristic function after subtracting its tilted mean, followed by expansion of its complex logarithm near one, applies to every prime, uniformly in $T$. Its cubic errors sum to at most $C\|\bb{\omega}\|_{2}^{3}\sum_{p \leq X}p^{-3/2} \leq C\|\bb{\omega}\|_{2}^{3}$. Using \eqref{an:conditional.detail:tilted.mean.covariance} to replace the tilted covariance matrix by $\Sigma_{T}$ introduces an additional error $\OO(\|\bb{\omega}\|_{2}^{2})$. Hence
\begin{equation}\label{an:conditional.detail:small.frequency}
e^{-\ii\bb{\omega}\cdot\bb{v}}\varphi_{T}^{\bb{\lambda},\lambda_{k + 1}}(\bb{\omega})
= \exp\left\{\ii\bb{\omega}\cdot(\bb{\mu}_{T} - \bb{v}) - \frac{1}{2}\bb{\omega}^{\mathsf{T}}\Sigma_{T}\bb{\omega} + \OO(\|\bb{\omega}\|_{2}^{2} + \|\bb{\omega}\|_{2}^{3})\right\}.
\end{equation}
In particular, if $\|\bb{\mu}_{T} - \bb{v}\|_{2} = \oo(\sqrt{t})$, then the phase and error in \eqref{an:conditional.detail:small.frequency} tend uniformly to zero on $\|\bb{\omega}\|_{2} \leq R_{0}/\sqrt{t}$ for every fixed $R_{0} > 0$.

The next step converts the characteristic-function estimates \eqref{an:detail:conditional.characteristic}--\eqref{an:conditional.detail:small.frequency} into weighted bin estimates. Let $\mathcal{J}_{1},\ldots,\mathcal{J}_{k}$ be intervals with endpoints in a fixed bounded set and lengths bounded below by a fixed positive constant. The choices $\mathcal{J}_{j} = [0,1]$ and $\mathcal{J}_{j} = [0,b_{j}]$ in the statement are included. For a bounded real $\lambda$ and one such interval $\mathcal{J}$, let $q_{\mathcal{J},\lambda}(x) \leqdef e^{-\lambda x}\ind_{\mathcal{J}}(x)$. Let $\Omega > 0$ (to be chosen later). The Fourier convention is $\widehat{g}(\bb{\omega}) \leqdef \int_{\R^{k}}e^{-\ii\bb{\omega}\cdot\bb{x}}g(\bb{x})\rd{\bb{x}}$, with the corresponding one-dimensional convention below.

For $0 \leq u \leq \|q_{\mathcal{J},\lambda}\|_{\infty}$, let $\mathcal{I}_{u} \leqdef \{x:q_{\mathcal{J},\lambda}(x) > u\}$. The Selberg--Beurling construction in \citet[Corollary~12, equations~(4.8)--(4.9), and the subsequent discussion, p.~201]{Vaaler1985}, with bandwidth parameter $2\Omega/\pi$, gives real-valued continuous $L^{1}(\R)$ functions $S_{\mathcal{I}_{u},\Omega}^{\pm}$ that extend to entire functions of exponential type at most $4\Omega$ and satisfy
\[
\begin{aligned}
&S_{\mathcal{I}_{u},\Omega}^{-} \leq \ind_{\mathcal{I}_{u}} \leq S_{\mathcal{I}_{u},\Omega}^{+}, \qquad
&&\|S_{\mathcal{I}_{u},\Omega}^{+} - S_{\mathcal{I}_{u},\Omega}^{-}\|_{1} \leq C\Omega^{-1}, \\
&\|S_{\mathcal{I}_{u},\Omega}^{\pm}\|_{1} \leq |\mathcal{I}_{u}| + C\Omega^{-1}, \qquad
&&\operatorname{supp}(\widehat{S}_{\mathcal{I}_{u},\Omega}^{\pm}) \subseteq [ - 4\Omega,4\Omega].
\end{aligned}
\]
The lower approximant is obtained by applying the majorant construction to the negative of the interval indicator and changing the sign. These bounds are uniform over the interval endpoints and lengths, without any positive lower bound on $|\mathcal{I}_{u}|$. By continuity, the majorant is at least one at both endpoints and the minorant is at most zero there, so the inequalities hold pointwise for open, closed, and half-open intervals. For an empty interval, let both approximants be zero. Their standard endpoint formulas give a measurable family in $u$.

To obtain approximants for $q_{\mathcal{J},\lambda}$, use the identity and definition
\[
q_{\mathcal{J},\lambda}(x) = \int_{0}^{\|q_{\mathcal{J},\lambda}\|_{\infty}}\ind_{\mathcal{I}_{u}}(x)\rd{u}, \qquad
q_{\mathcal{J},\lambda,\Omega}^{\pm}(x) \leqdef \int_{0}^{\|q_{\mathcal{J},\lambda}\|_{\infty}}S_{\mathcal{I}_{u},\Omega}^{\pm}(x)\rd{u}.
\]
For fixed $\Omega$, Fourier inversion bounds the interval approximants uniformly, so dominated convergence shows that $q_{\mathcal{J},\lambda,\Omega}^{\pm}$ are continuous. Since $\|q_{\mathcal{J},\lambda}\|_{\infty}$ and $|\mathcal{J}|$ are uniformly bounded, integration of the preceding estimates and Fubini's theorem give
\begin{equation}\label{an:detail:interval.approximants}
\begin{aligned}
&q_{\mathcal{J},\lambda,\Omega}^{-} \leq q_{\mathcal{J},\lambda} \leq q_{\mathcal{J},\lambda,\Omega}^{+}, \qquad
&&\|q_{\mathcal{J},\lambda,\Omega}^{+} - q_{\mathcal{J},\lambda,\Omega}^{-}\|_{1} \leq C\Omega^{-1}, \\
&\|q_{\mathcal{J},\lambda,\Omega}^{\pm}\|_{1} \leq C(1 + \Omega^{-1}), \qquad
&&\operatorname{supp}(\widehat{q}_{\mathcal{J},\lambda,\Omega}^{\pm}) \subseteq [ - 4\Omega,4\Omega].
\end{aligned}
\end{equation}
Figure~\ref{fig:conditional.weighted.bin.setup} illustrates the interval setup and the associated weighted approximants.

\begin{figure}[!ht]
\centering
\includegraphics[width = \textwidth]{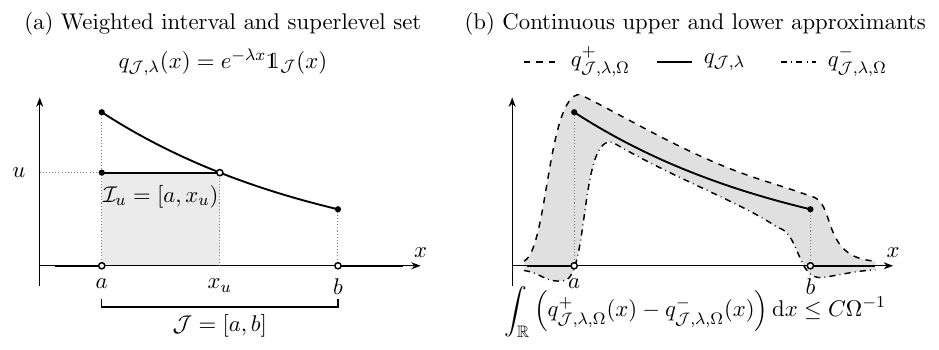}
\caption{Weighted interval approximation for $\mathcal{J} = [a,b]$ and $\lambda > 0$. In (a), a level $e^{-\lambda b} < u < e^{-\lambda a}$ gives $x_{u} \leqdef - (\log{u})/\lambda$ and $\mathcal{I}_{u} = [a,x_{u})$. In (b), the schematic continuous curves satisfy $\smash{q_{\mathcal{J},\lambda,\Omega}^{-} \leq q_{\mathcal{J},\lambda} \leq q_{\mathcal{J},\lambda,\Omega}^{+}}$, and the shaded gap has total area at most $C\Omega^{-1}$. The approximants extend beyond the displayed window and have Fourier support in $[ - 4\Omega,4\Omega]$.}\label{fig:conditional.weighted.bin.setup}
\end{figure}

For the remaining estimates, choose $\Omega \leqdef t^{\e_{0}/2}/4$. Then $\Omega \to \infty$ and $[ - 4\Omega,4\Omega] = [ - t^{\e_{0}/2},t^{\e_{0}/2}]$, so \eqref{an:detail:interval.approximants} gives $\smash{\|q_{\mathcal{J},\lambda,\Omega}^{\pm}\|_{1} \leq C}$ for all sufficiently large $T$. For each $j \in \{1,\ldots,k\}$, let $q_{\mathcal{J}_{j},\lambda_{j}}(x) \leqdef e^{-\lambda_{j}x}\ind_{\mathcal{J}_{j}}(x)$ and abbreviate its approximants by $\smash{q_{j}^{\pm} \leqdef q_{\mathcal{J}_{j},\lambda_{j},\Omega}^{\pm}}$. Define, for $\bb{x} \in \R^{k}$,
\[
\begin{aligned}
f(\bb{x}) &\leqdef \prod_{j = 1}^{k}q_{\mathcal{J}_{j},\lambda_{j}}(x_{j}), \qquad f^{+}(\bb{x}) \leqdef \prod_{j = 1}^{k}q_{j}^{+}(x_{j}), \\
f^{-}(\bb{x}) &\leqdef f^{+}(\bb{x}) - \sum_{j = 1}^{k}(q_{j}^{+}(x_{j}) - q_{j}^{-}(x_{j}))\prod_{i \neq j}q_{i}^{+}(x_{i}).
\end{aligned}
\]
Each $q_{j}^{+}$ is nonnegative. Telescoping the product and using $0 \leq q_{\mathcal{J}_{j},\lambda_{j}} \leq q_{j}^{+}$ gives
\[
0 \leq f^{+} - f \leq \sum_{j = 1}^{k}(q_{j}^{+} - q_{\mathcal{J}_{j},\lambda_{j}})\prod_{i \neq j}q_{i}^{+} \leq \sum_{j = 1}^{k}(q_{j}^{+} - q_{j}^{-})\prod_{i \neq j}q_{i}^{+} = f^{+} - f^{-}.
\]
Thus $f^{-} \leq f \leq f^{+}$, and Fubini's theorem together with \eqref{an:detail:interval.approximants} yields
\begin{equation}\label{an:conditional.detail:product.approximants}
\begin{aligned}
\|f^{+} - f^{-}\|_{1}
&= \sum_{j = 1}^{k}\|q_{j}^{+} - q_{j}^{-}\|_{1}\prod_{i \neq j}\|q_{i}^{+}\|_{1} \leq C\Omega^{-1}, \\
\|f^{\pm} - f\|_{1} &\leq C\Omega^{-1}, \qquad \|f^{\pm}\|_{1} \leq C, \qquad
\operatorname{supp}(\widehat{f}^{\pm}) \subseteq [ - 4\Omega,4\Omega]^{k}.
\end{aligned}
\end{equation}
Here and below the constants may depend on $k$ and the fixed bounds on the intervals and real tilts. Since $f$ has uniformly bounded support and integral, \eqref{an:conditional.detail:product.approximants} also gives
\begin{equation}\label{an:conditional.detail:fourier.continuity}
\begin{aligned}
|\widehat{f}^{\pm}(\bb{\omega}) - \widehat{f}^{\pm}(\bb{0})|
&\leq \int_{\R^{k}}|e^{-\ii\bb{\omega}\cdot\bb{x}} - 1||f^{\pm}(\bb{x}) - f(\bb{x})|\rd{\bb{x}} + \int_{\R^{k}}|e^{-\ii\bb{\omega}\cdot\bb{x}} - 1|f(\bb{x})\rd{\bb{x}} \\
&\leq 2\|f^{\pm} - f\|_{1} + \int_{\R^{k}}|e^{-\ii\bb{\omega}\cdot\bb{x}} - 1|f(\bb{x})\rd{\bb{x}}
\leq C(\Omega^{-1} + \|\bb{\omega}\|_{2}), \\
\widehat{f}^{\pm}(\bb{0})
&= \int_{\R^{k}}f^{\pm}(\bb{x})\rd{\bb{x}}
= \int_{\R^{k}}f(\bb{x})\rd{\bb{x}} + \OO(\Omega^{-1}).
\end{aligned}
\end{equation}

Suppose now that $\|\bb{\mu}_{T} - \bb{v}\|_{2} = \oo(\sqrt{t})$, uniformly over the parameters under consideration. The functions $f^{\pm}$ are continuous, and their Fourier transforms in \eqref{an:conditional.detail:product.approximants} are bounded and compactly supported, so Fourier inversion and Fubini's theorem give
\begin{equation}\label{an:conditional.detail:fourier.inversion}
\EE_{\mathrm{H}}^{\bb{\lambda},\lambda_{k + 1}}\left[f^{\pm}\left((Y_{T,r}^{\mathrm{H}} - v_{r})_{1 \leq r \leq k}\right)\right] = \frac{1}{(2\pi)^{k}}\int_{[ - 4\Omega,4\Omega]^{k}}\widehat{f}^{\pm}(\bb{\omega})e^{-\ii\bb{\omega}\cdot\bb{v}}\varphi_{T}^{\bb{\lambda},\lambda_{k + 1}}(\bb{\omega})\rd{\bb{\omega}}.
\end{equation}
The leading contribution comes from frequencies of order $t^{-1/2}$. More precisely, for every fixed $R_{0} \geq 1$, equations \eqref{an:conditional.detail:small.frequency} and \eqref{an:conditional.detail:fourier.continuity} replace the integrand on $\|\bb{\omega}\|_{2} \leq R_{0}/\sqrt{t}$ by $\left(\int_{\R^{k}}f(\bb{x})\rd{\bb{x}}\right)e^{-\bb{\omega}^{\mathsf{T}}\Sigma_{T}\bb{\omega}/2}$ with integral error $\oo(t^{-k/2})$, since the phase and error in \eqref{an:conditional.detail:small.frequency} tend uniformly to zero on this ball.

To show that the remaining frequencies contribute negligibly, use $|\widehat{f}^{\pm}(\bb{\omega})| \leq \|f^{\pm}\|_{1} \leq C$ from \eqref{an:conditional.detail:product.approximants} and the characteristic-function estimate \eqref{an:detail:conditional.characteristic}. These bounds control the absolute contribution from the part of the Fourier cube $[ - 4\Omega,4\Omega]^{k}$ outside the ball by
\[
C\int_{\|\bb{\omega}\|_{2}>R_{0}/\sqrt{t}}e^{-ct\|\bb{\omega}\|_{2}^{2}}\rd{\bb{\omega}} = Ct^{-k/2}\int_{\|\bb{z}\|_{2}>R_{0}}e^{-c\|\bb{z}\|_{2}^{2}}\rd{\bb{z}} \leq Ct^{-k/2}e^{-c'R_{0}^{2}}
\]
for some $c' > 0$. Thus restricting the original Fourier integral in \eqref{an:conditional.detail:fourier.inversion} to the ball introduces an error that becomes negligible relative to $t^{-k/2}$ as $R_{0} \to \infty$.

The Gaussian approximation on this ball can then be extended to all of $\R^{k}$ at the same cost. Indeed, the inequality $\Sigma_{T} \geq c_{\Sigma}t\mathrm{Id}_{k}$ in the Loewner order gives
\[
\frac{1}{(2\pi)^{k}}\int_{\|\bb{\omega}\|_{2}>R_{0}/\sqrt{t}}e^{-\bb{\omega}^{\mathsf{T}}\Sigma_{T}\bb{\omega}/2}\rd{\bb{\omega}} \leq Ct^{-k/2}e^{-c'R_{0}^{2}},
\]
and multiplication by the uniformly bounded factor $\int_{\R^{k}}f(\bb{x})\rd{\bb{x}}$ preserves this bound. The resulting integral over $\R^{k}$ can be evaluated explicitly and supplies the leading coefficient:
\[
\frac{1}{(2\pi)^{k}}\int_{\R^{k}}e^{-\bb{\omega}^{\mathsf{T}}\Sigma_{T}\bb{\omega}/2}\rd{\bb{\omega}} = (2\pi)^{-k/2}(\det(\Sigma_{T}))^{-1/2}.
\]
Combining the approximation inside the ball $\{\bb{\omega} \in \R^{k}:\|\bb{\omega}\|_{2} \leq R_{0}/\sqrt{t}\}$ with the bounds for the discarded Fourier contribution and the added Gaussian tail therefore gives, for every fixed $R_{0} \geq 1$,
\[
\begin{aligned}
&\EE_{\mathrm{H}}^{\bb{\lambda},\lambda_{k + 1}}\left[f^{\pm}\left((Y_{T,r}^{\mathrm{H}} - v_{r})_{1 \leq r \leq k}\right)\right] \\
&\qquad = (2\pi)^{-k/2}(\det(\Sigma_{T}))^{-1/2}\int_{\R^{k}}f(\bb{x})\rd{\bb{x}} + \oo(t^{-k/2}) + \OO(t^{-k/2}e^{-c'R_{0}^{2}}).
\end{aligned}
\]

It remains to pass from the approximants to the original weighted bin function $f$. The pointwise bounds $f^{-} \leq f \leq f^{+}$ imply $0 \leq f^{+} - f \leq f^{+} - f^{-}$, so the expected gap between the approximants controls the error in replacing $f^{+}$ by $f$. Subtracting the identities \eqref{an:conditional.detail:fourier.inversion} for $f^{+}$ and $f^{-}$, taking absolute values inside the resulting integral, and using \eqref{an:detail:conditional.characteristic} yields
\[
\begin{aligned}
0 &\leq \EE_{\mathrm{H}}^{\bb{\lambda},\lambda_{k + 1}}\left[(f^{+} - f)\left((Y_{T,r}^{\mathrm{H}} - v_{r})_{1 \leq r \leq k}\right)\right] \\
&\leq \EE_{\mathrm{H}}^{\bb{\lambda},\lambda_{k + 1}}\left[(f^{+} - f^{-})\left((Y_{T,r}^{\mathrm{H}} - v_{r})_{1 \leq r \leq k}\right)\right] \\
&\leq \frac{\|f^{+} - f^{-}\|_{1}}{(2\pi)^{k}}\int_{[ - 4\Omega,4\Omega]^{k}}e^{-ct\|\bb{\omega}\|_{2}^{2}}\rd{\bb{\omega}}
\leq C\Omega^{-1}t^{-k/2}.
\end{aligned}
\]
Here $|\widehat{f}^{+}(\bb{\omega}) - \widehat{f}^{-}(\bb{\omega})| \leq \|f^{+} - f^{-}\|_{1}$. The last inequality uses $\|f^{+} - f^{-}\|_{1} \leq C\Omega^{-1}$ from \eqref{an:conditional.detail:product.approximants} and bounds the integral by its extension to $\R^{k}$, which is a constant multiple of $t^{-k/2}$. Since $\Omega \to \infty$, this replacement contributes only $\oo(t^{-k/2})$.

First letting $T \to \infty$ and then $R_{0} \to \infty$ shows that the total error is $\oo(t^{-k/2})$, viz.,
\[
\EE_{\mathrm{H}}^{\bb{\lambda},\lambda_{k + 1}}\left[f\left((Y_{T,r}^{\mathrm{H}} - v_{r})_{1 \leq r \leq k}\right)\right]
= (2\pi)^{-k/2}(\det(\Sigma_{T}))^{-1/2}\int_{\R^{k}}f(\bb{x})\rd{\bb{x}} + \oo(t^{-k/2}).
\]
The same expansion holds with $f^{+}$ or $f^{-}$ in place of $f$ inside the tilted expectation, with the integral $\int_{\R^{k}}f(\bb{x})\rd{\bb{x}}$ on the right-hand side unchanged. The eigenvalue bounds give $(\det(\Sigma_{T}))^{-1/2} \asymp t^{-k/2}$, so this error becomes $\oo(1)$ after factoring out $(2\pi)^{-k/2}(\det(\Sigma_{T}))^{-1/2}$. Finally, integrating the product that defines $f$ gives the product of the one-dimensional weighted bin integrals. These integrals are uniformly bounded away from zero because the tilts and endpoints are bounded and the bin lengths have a fixed positive lower bound. Consequently,
\begin{equation}\label{an:local.factor}
\EE_{\mathrm{H}}^{\bb{\lambda},\lambda_{k + 1}}\left[\prod_{j = 1}^{k}q_{\mathcal{J}_{j},\lambda_{j}}(Y_{T,j}^{\mathrm{H}} - v_{j})\right]
= (2\pi)^{-k/2}(\det(\Sigma_{T}))^{-1/2} \prod_{j = 1}^{k}\int_{\mathcal{J}_{j}}e^{-\lambda_{j}x}\rd{x} \times \left\{1 + \oo(1)\right\}.
\end{equation}
The $\oo(1)$ term is uniform over the permitted intervals and tilts.

To obtain the arithmetic counterpart of the tilted Haar local estimate \eqref{an:local.factor}, let $\mathscr{Q} \leqdef \prod_{j = 1}^{k}\mathcal{J}_{j}$. The transform comparison \eqref{an:detail:conditional.transform} in Lemma~\ref{an:transform}, applied with parameters $(\lambda_{1} + \ii\omega_{1},\ldots,\lambda_{k} + \ii\omega_{k},\lambda_{k + 1})$, gives
\begin{equation}\label{an:conditional.detail:tilted.transform.comparison}
\begin{aligned}
&\left|\EE_{\mathrm{A}}\left[e^{\sum_{r = 1}^{k}(\lambda_{r} + \ii\omega_{r})Y_{T,r}^{\mathrm{A}} + \lambda_{k + 1}Y_{T,k + 1}^{\mathrm{A}}}\ind_{\mathcal{A}_{T}^{\mathrm{A}}}\right]
- e^{\mathcal{L}_{T}(\bb{\lambda},\lambda_{k + 1})}\varphi_{T}^{\bb{\lambda},\lambda_{k + 1}}(\bb{\omega})\right| \\
&\qquad \leq \exp\{ - ct^{1 + 2\e_{0}} + Ct\}, \qquad \bb{\omega} \in [ - 4\Omega,4\Omega]^{k}.
\end{aligned}
\end{equation}
Fourier inversion with $f^{+}$ or $f^{-}$, followed by \eqref{an:conditional.detail:tilted.transform.comparison}, therefore yields
\begin{equation}\label{an:conditional.detail:inversion.comparison}
\begin{aligned}
&\left|\EE_{\mathrm{A}}\left[e^{\sum_{r = 1}^{k}\lambda_{r}Y_{T,r}^{\mathrm{A}} + \lambda_{k + 1}Y_{T,k + 1}^{\mathrm{A}}}f^{\pm}\left((Y_{T,r}^{\mathrm{A}} - v_{r})_{1 \leq r \leq k}\right)\ind_{\mathcal{A}_{T}^{\mathrm{A}}}\right]\right. \\
&\qquad \left. - e^{\mathcal{L}_{T}(\bb{\lambda},\lambda_{k + 1})}\EE_{\mathrm{H}}^{\bb{\lambda},\lambda_{k + 1}}\left[f^{\pm}\left((Y_{T,r}^{\mathrm{H}} - v_{r})_{1 \leq r \leq k}\right)\right]\right| \\
&\qquad = \frac{1}{(2\pi)^{k}}\Bigg|\int_{[ - 4\Omega,4\Omega]^{k}}\widehat{f}^{\pm}(\bb{\omega})e^{-\ii\bb{\omega}\cdot\bb{v}} \\
&\qquad\qquad\qquad \times \Bigg\{\EE_{\mathrm{A}}\left[e^{\sum_{r = 1}^{k}(\lambda_{r} + \ii\omega_{r})Y_{T,r}^{\mathrm{A}} + \lambda_{k + 1}Y_{T,k + 1}^{\mathrm{A}}}\ind_{\mathcal{A}_{T}^{\mathrm{A}}}\right] - e^{\mathcal{L}_{T}(\bb{\lambda},\lambda_{k + 1})}\varphi_{T}^{\bb{\lambda},\lambda_{k + 1}}(\bb{\omega})\Bigg\}\rd{\bb{\omega}}\Bigg| \\
&\qquad \leq \frac{\|\widehat{f}^{\pm}\|_{\infty}(8\Omega)^{k}}{(2\pi)^{k}}\exp\{ - ct^{1 + 2\e_{0}} + Ct\}
\leq \exp\{ - ct^{1 + 2\e_{0}} + C't\}.
\end{aligned}
\end{equation}
The last inequality uses $\|\widehat{f}^{\pm}\|_{\infty} \leq \|f^{\pm}\|_{1} \leq C$ and $\Omega^{k} = \OO(t^{k\e_{0}/2})$. Multiplying the pointwise inequalities $f^{-} \leq f \leq f^{+}$ by the nonnegative factor $e^{\sum_{r = 1}^{k}\lambda_{r}Y_{T,r}^{\mathrm{A}} + \lambda_{k + 1}Y_{T,k + 1}^{\mathrm{A}}}\ind_{\mathcal{A}_{T}^{\mathrm{A}}}$, taking arithmetic expectations, and applying \eqref{an:conditional.detail:inversion.comparison} give
\begin{equation}\label{an:conditional.detail:weighted.sandwich}
\begin{aligned}
&e^{\mathcal{L}_{T}(\bb{\lambda},\lambda_{k + 1})}\EE_{\mathrm{H}}^{\bb{\lambda},\lambda_{k + 1}}\left[f^{-}\left((Y_{T,r}^{\mathrm{H}} - v_{r})_{1 \leq r \leq k}\right)\right]
- \exp\{ - ct^{1 + 2\e_{0}} + C't\} \\
&\qquad \leq \EE_{\mathrm{A}}\left[e^{\sum_{r = 1}^{k}\lambda_{r}Y_{T,r}^{\mathrm{A}} + \lambda_{k + 1}Y_{T,k + 1}^{\mathrm{A}}}f\left((Y_{T,r}^{\mathrm{A}} - v_{r})_{1 \leq r \leq k}\right)\ind_{\mathcal{A}_{T}^{\mathrm{A}}}\right] \\
&\qquad \leq e^{\mathcal{L}_{T}(\bb{\lambda},\lambda_{k + 1})}\EE_{\mathrm{H}}^{\bb{\lambda},\lambda_{k + 1}}\left[f^{+}\left((Y_{T,r}^{\mathrm{H}} - v_{r})_{1 \leq r \leq k}\right)\right]
+ \exp\{ - ct^{1 + 2\e_{0}} + C't\}.
\end{aligned}
\end{equation}
The Haar local estimate \eqref{an:local.factor}, including its versions with $f^{+}$ or $f^{-}$ inside the tilted expectation, gives the same leading term for both outer expectations in \eqref{an:conditional.detail:weighted.sandwich}. Finally, the identity
\begin{equation}\label{an:conditional.detail:undo.tilt}
e^{\lambda_{k + 1}Y_{T,k + 1}^{\star}}\ind_{\{(Y_{T,r}^{\star})_{1\leq r\leq k}\in\bb{v} + \mathscr{Q}\}}
= e^{-\bb{\lambda}\cdot\bb{v}}e^{\sum_{r = 1}^{k}\lambda_{r}Y_{T,r}^{\star} + \lambda_{k + 1}Y_{T,k + 1}^{\star}}f\left((Y_{T,r}^{\star} - v_{r})_{1 \leq r \leq k}\right)
\end{equation}
undoes the tilt. Combining this identity, \eqref{an:local.factor}, and \eqref{an:conditional.detail:weighted.sandwich} gives
\begin{equation}\label{an:conditional.detail:weighted.bin}
\EE_{\mathrm{A}}\left[e^{\lambda_{k + 1}Y_{T,k + 1}^{\mathrm{A}}}\ind_{\mathcal{A}_{T}^{\mathrm{A}}}\ind_{\{(Y_{T,r}^{\mathrm{A}})_{1\leq r\leq k}\in\bb{v} + \mathscr{Q}\}}\right]
= \frac{e^{\mathcal{L}_{T}(\bb{\lambda},\lambda_{k + 1}) - \bb{\lambda}\cdot\bb{v}}}{(2\pi)^{k/2}\sqrt{\det(\Sigma_{T})}}\prod_{j = 1}^{k}\int_{\mathcal{J}_{j}}e^{-\lambda_{j}x}\rd{x} \times \left\{1 + \oo(1)\right\} + \mathcal{E}_{T}.
\end{equation}
This formula holds whenever the real tilt vector $(\bb{\lambda},\lambda_{k + 1})$ is uniformly bounded and $\|\bb{\mu}_{T} - \bb{v}\|_{2} = \oo(\sqrt{t})$. If also $\|\bb{v}\|_{2} = \OO(t)$, multiplication by $e^{-\bb{\lambda}\cdot\bb{v}}$ in undoing the tilt costs at most $e^{Ct}$, so the additional transform comparison error satisfies
\begin{equation}\label{an:detail:conditional.bin.error}
|\mathcal{E}_{T}| \leq \exp\{ - ct^{1 + 2\e_{0}} + Ct\}.
\end{equation}
The approximation of $f$ by $f^{-}$ and $f^{+}$ contributes to the $\oo(1)$ term in the multiplicative factor in \eqref{an:conditional.detail:weighted.bin}, whereas \eqref{an:detail:conditional.bin.error} bounds the additional absolute error from comparing the truncated arithmetic transform with the Haar transform. Under the Haar law, \eqref{an:local.factor} and \eqref{an:conditional.detail:undo.tilt} give \eqref{an:conditional.detail:weighted.bin} without the event indicator $\ind_{\mathcal{A}_{T}^{\mathrm{H}}}$ and with $\mathcal{E}_{T} = 0$.

It remains to apply \eqref{an:conditional.detail:weighted.bin} twice, choosing the tilt parameters separately to bound the conditioning-bin probability $\PP_{\star}((Y_{T,r}^{\star})_{1\leq r\leq k} \in \bb{v} + \mathscr{Q})$ from below and the exponentially weighted bin expectation $\smash{\EE_{\star}[e^{\lambda_{k + 1}Y_{T,k + 1}^{\star}}\ind_{\mathcal{A}_{T}^{\star}}\ind_{\{(Y_{T,r}^{\star})_{1\leq r\leq k}\in\bb{v} + \mathscr{Q}\}}]}$ from above. Dividing the latter upper bound by the former lower bound will give \eqref{an:conditional.bound}.

Set $\mathscr{Q} = [0,1]^{k}$ initially. Since $s_{T} \sim t$ and $\bb{v}/s_{T} \in K$, the vectors $\bb{v}$ satisfy $\|\bb{v}\|_{2} \leq Ct$. The moment estimate \eqref{an:haar.moments} in Lemma~\ref{an:polynomial}, Cauchy--Schwarz, and the eigenvalue assumption on $\Sigma_{T}$ also give
\[
\sig_{T}^{2} \leq Ct, \qquad |c_{T,r}| \leq \sqrt{(\Sigma_{T})_{rr}\sig_{T}^{2}} \leq Ct, \qquad \|\Sigma_{T}^{-1}\|_{\mathrm{op}} \leq \frac{1}{c_{\Sigma}t}.
\]
For every fixed $M_{1} > 0$ and all $|\lambda_{j}| \leq M_{1}$, the weighted unit-bin integrals satisfy
\begin{equation}\label{an:conditional.detail:unit.bin.integrals}
e^{-kM_{1}} \leq \prod_{j = 1}^{k}\int_{0}^{1}e^{-\lambda_{j}x}\rd{x} \leq e^{kM_{1}}.
\end{equation}
For the denominator of the conditional expectation in \eqref{an:conditional.bound}, choose $\lambda_{k + 1} = 0$ and $\bb{\lambda}^{(0)} \leqdef \Sigma_{T}^{-1}\bb{v}$. This tilt is bounded, and its mean is $\bb{v} + \OO(1)$ by \eqref{an:conditional.detail:tilted.mean.covariance}. Formula \eqref{an:joint.cumulants} gives
\[
\mathcal{L}_{T}(\bb{\lambda}^{(0)},0) - (\bb{\lambda}^{(0)})^{\mathsf{T}}\bb{v}
= - \frac{1}{2}\bb{v}^{\mathsf{T}}\Sigma_{T}^{-1}\bb{v} + \mathcal{R}_{T}(\bb{\lambda}^{(0)},0).
\]
Apply \eqref{an:conditional.detail:weighted.bin} with this tilt. For the arithmetic law, the event indicator may be removed because \eqref{an:exceptional} in Lemma~\ref{an:transform} implies
\[
0 \leq \PP_{\mathrm{A}}((Y_{T,r}^{\mathrm{A}})_{1\leq r\leq k} \in \bb{v} + \mathscr{Q}) - \PP_{\mathrm{A}}((Y_{T,r}^{\mathrm{A}})_{1\leq r\leq k} \in \bb{v} + \mathscr{Q},\mathcal{A}_{T}^{\mathrm{A}})
\leq \PP_{\mathrm{A}}((\mathcal{A}_{T}^{\mathrm{A}})^{c}) \leq \exp\{-ct^{1 + 2\e_{0}}\}.
\]
As noted after \eqref{an:detail:conditional.bin.error}, the Haar version of \eqref{an:conditional.detail:weighted.bin} has no event indicator $\smash{\ind_{\mathcal{A}_{T}^{\mathrm{H}}}}$. Applying \eqref{an:conditional.detail:weighted.bin} with $(\bb{\lambda},\lambda_{k + 1}) = (\Sigma_{T}^{-1}\bb{v},0)$ and using the remainder bound \eqref{an:conditional.detail:remainder}, the lower bound \eqref{an:conditional.detail:unit.bin.integrals} on the bin integrals, and the error estimate \eqref{an:detail:conditional.bin.error} gives, for each $\star \in \{\mathrm{A},\mathrm{H}\}$,
\begin{equation}\label{an:detail:conditional.bin.lower}
\PP_{\star}((Y_{T,r}^{\star})_{1\leq r\leq k} \in \bb{v} + \mathscr{Q})
\geq c(\det(\Sigma_{T}))^{-1/2}\exp\left( - \frac{1}{2}\bb{v}^{\mathsf{T}}\Sigma_{T}^{-1}\bb{v}\right) \geq e^{-Ct}.
\end{equation}
Indeed, $\smash{0 \leq \bb{v}^{\mathsf{T}}\Sigma_{T}^{-1}\bb{v} \leq Ct}$, so the main term is at least $ct^{-k/2}e^{-Ct}$, whereas the errors just displayed are smaller than this by a factor tending to zero. This also proves that the conditioning event has positive probability for all sufficiently large $T$.

To bound the numerator obtained by writing the conditional expectation in \eqref{an:conditional.bound} as a ratio, choose
\[
\bb{\lambda} \leqdef \Sigma_{T}^{-1}(\bb{v} - \lambda_{k + 1}\bb{c}_{T}), \qquad
\|\bb{\lambda}\|_{2} \leq \frac{\|\bb{v}\|_{2} + M\|\bb{c}_{T}\|_{2}}{c_{\Sigma}t} \leq C.
\]
Its tilted mean is again $\bb{v} + \OO(1)$ by \eqref{an:conditional.detail:tilted.mean.covariance}, so \eqref{an:conditional.detail:weighted.bin} applies. Define its quadratic exponent after undoing the tilt by
\begin{equation}\label{an:conditional.detail:quadratic.exponent}
E_{T}(\lambda_{k + 1}) \leqdef - \frac{1}{2}(\bb{v} - \lambda_{k + 1}\bb{c}_{T})^{\mathsf{T}}\Sigma_{T}^{-1}(\bb{v} - \lambda_{k + 1}\bb{c}_{T}) + \frac{\lambda_{k + 1}^{2}}{2}\sig_{T}^{2}.
\end{equation}
Substituting $\Sigma_{T}\bb{\lambda} = \bb{v} - \lambda_{k + 1}\bb{c}_{T}$ into \eqref{an:joint.cumulants} gives
\[
\begin{aligned}
\mathcal{L}_{T}(\bb{\lambda},\lambda_{k + 1}) - \bb{\lambda}^{\mathsf{T}}\bb{v}
&= \frac{1}{2}\bb{\lambda}^{\mathsf{T}}\Sigma_{T}\bb{\lambda} + \lambda_{k + 1}\bb{\lambda}^{\mathsf{T}}\bb{c}_{T} + \frac{\lambda_{k + 1}^{2}}{2}\sig_{T}^{2} - \bb{\lambda}^{\mathsf{T}}\bb{v} + \mathcal{R}_{T}(\bb{\lambda},\lambda_{k + 1}) \\
&= - \frac{1}{2}(\bb{v} - \lambda_{k + 1}\bb{c}_{T})^{\mathsf{T}}\Sigma_{T}^{-1}(\bb{v} - \lambda_{k + 1}\bb{c}_{T}) + \frac{\lambda_{k + 1}^{2}}{2}\sig_{T}^{2} + \mathcal{R}_{T}(\bb{\lambda},\lambda_{k + 1}) \\
&= E_{T}(\lambda_{k + 1}) + \mathcal{R}_{T}(\bb{\lambda},\lambda_{k + 1}).
\end{aligned}
\]
Let $W_{T}^{\mathrm{A}}(\lambda_{k + 1})$ denote the left-hand side of \eqref{an:conditional.detail:weighted.bin}, and let $W_{T}^{\mathrm{H}}(\lambda_{k + 1})$ denote the corresponding Haar expectation without the event indicator $\ind_{\mathcal{A}_{T}^{\mathrm{H}}}$. The uniform upper bound on the bin integrals in \eqref{an:conditional.detail:unit.bin.integrals}, the remainder estimate \eqref{an:conditional.detail:remainder}, and the error bound \eqref{an:detail:conditional.bin.error}, applied to \eqref{an:conditional.detail:weighted.bin}, yield
\begin{equation}\label{an:conditional.detail:numerator.bound}
W_{T}^{\star}(\lambda_{k + 1}) \leq C(\det(\Sigma_{T}))^{-1/2}e^{E_{T}(\lambda_{k + 1})} + \exp\{ - ct^{1 + 2\e_{0}} + Ct\}.
\end{equation}
Dividing by the sharper first lower bound in \eqref{an:detail:conditional.bin.lower}, with $E_{T}(0) = - \bb{v}^{\mathsf{T}}\Sigma_{T}^{-1}\bb{v}/2$, cancels the determinant factors and gives
\begin{equation}\label{an:conditional.detail:ratio}
\frac{W_{T}^{\star}(\lambda_{k + 1})}{\PP_{\star}((Y_{T,r}^{\star})_{1\leq r\leq k} \in \bb{v} + \mathscr{Q})}
\leq Ce^{E_{T}(\lambda_{k + 1}) - E_{T}(0)} + \exp\{ - ct^{1 + 2\e_{0}} + C't\}.
\end{equation}
Expanding the quadratic form $(\bb{v} - \lambda_{k + 1}\bb{c}_{T})^{\mathsf{T}}\Sigma_{T}^{-1}(\bb{v} - \lambda_{k + 1}\bb{c}_{T})$ in \eqref{an:conditional.detail:quadratic.exponent} gives
\[
\begin{aligned}
E_{T}(\lambda_{k + 1}) - E_{T}(0)
&= - \frac{1}{2}\left(\bb{v}^{\mathsf{T}}\Sigma_{T}^{-1}\bb{v} - 2\lambda_{k + 1}\bb{c}_{T}^{\mathsf{T}}\Sigma_{T}^{-1}\bb{v} + \lambda_{k + 1}^{2}\bb{c}_{T}^{\mathsf{T}}\Sigma_{T}^{-1}\bb{c}_{T}\right)
+ \frac{\lambda_{k + 1}^{2}}{2}\sig_{T}^{2} + \frac{1}{2}\bb{v}^{\mathsf{T}}\Sigma_{T}^{-1}\bb{v} \\
&= \lambda_{k + 1}\bb{c}_{T}^{\mathsf{T}}\Sigma_{T}^{-1}\bb{v} + \frac{\lambda_{k + 1}^{2}}{2}\left(\sig_{T}^{2} - \bb{c}_{T}^{\mathsf{T}}\Sigma_{T}^{-1}\bb{c}_{T}\right).
\end{aligned}
\]
The absolute value of this exponent is at most $Ct$. Hence the second term on the right-hand side of \eqref{an:conditional.detail:ratio}, divided by $e^{E_{T}(\lambda_{k + 1}) - E_{T}(0)}$, is at most $\exp\{ - ct^{1 + 2\e_{0}} + C''t\} = \oo(1)$ and can be absorbed into the constant $C$. This proves \eqref{an:conditional.bound} for the arithmetic law. By the definition of $W_{T}^{\mathrm{H}}$, the Haar ratio gives explicitly
\begin{equation}\label{an:conditional.detail:haar.untruncated}
\begin{aligned}
&\EE_{\mathrm{H}}\left[e^{\lambda_{k + 1}Y_{T,k + 1}^{\mathrm{H}}}\,\middle|\,(Y_{T,1}^{\mathrm{H}},\ldots,Y_{T,k}^{\mathrm{H}}) \in \bb{v} + [0,1]^{k}\right] \\
&\qquad \leq C\exp\left\{\lambda_{k + 1}\bb{c}_{T}^{\mathsf{T}}\Sigma_{T}^{-1}\bb{v} + \frac{\lambda_{k + 1}^{2}}{2}\left(\sig_{T}^{2} - \bb{c}_{T}^{\mathsf{T}}\Sigma_{T}^{-1}\bb{c}_{T}\right)\right\}.
\end{aligned}
\end{equation}
Multiplying the nonnegative integrand of the conditional expectation in \eqref{an:conditional.detail:haar.untruncated} by $\ind_{\mathcal{A}_{T}^{\mathrm{H}}}$ cannot increase the expectation, which proves \eqref{an:conditional.bound} for the Haar law as well.

For $\mathscr{Q} = \prod_{j = 1}^{k}[0,b_{j}]$ with $b_{-} \leq b_{j} \leq b_{+}$, the approximation and inversion estimates above are unchanged, and the bin integrals satisfy
\begin{equation}\label{an:conditional.detail:variable.bin.integrals}
(b_{-}e^{-M_{1}b_{+}})^{k} \leq \prod_{j = 1}^{k}\int_{0}^{b_{j}}e^{-\lambda_{j}x}\rd{x} \leq (b_{+}e^{M_{1}b_{+}})^{k}
\end{equation}
for $|\lambda_{j}| \leq M_{1}$. The tilts $\bb{\lambda}^{(0)} = \Sigma_{T}^{-1}\bb{v}$ and $\bb{\lambda} = \Sigma_{T}^{-1}(\bb{v} - \lambda_{k + 1}\bb{c}_{T})$ do not depend on $b_{j}$, so their bounds and the estimate $\bb{\mu}_{T} = \bb{v} + \OO(1)$ remain uniform. Applying \eqref{an:conditional.detail:weighted.bin} with $\mathcal{J}_{j} = [0,b_{j}]$, the lower bound in \eqref{an:conditional.detail:variable.bin.integrals} gives the denominator estimate \eqref{an:detail:conditional.bin.lower} for this $\mathscr{Q}$, while its upper bound gives the same numerator estimate for $W_{T}^{\star}(\lambda_{k + 1})$ in \eqref{an:conditional.detail:numerator.bound}. The exponents $E_{T}(\lambda_{k + 1})$ and $E_{T}(0)$ in \eqref{an:conditional.detail:quadratic.exponent} are unchanged. Division therefore cancels the determinant factors and yields \eqref{an:conditional.detail:ratio} for this $\mathscr{Q}$. The transform comparison error $\mathcal{E}_{T}$ in \eqref{an:conditional.detail:weighted.bin} remains negligible relative to the main term by \eqref{an:detail:conditional.bin.error} and the uniform lower bound \eqref{an:detail:conditional.bin.lower}. Thus \eqref{an:conditional.bound} holds with $[0,1]^{k}$ replaced by $\mathscr{Q}$, with additional constant dependence only on $b_{-}$ and $b_{+}$.

Finally, for the common arithmetic spatial event $\smash{\mathcal{A}_{T,\mathrm{sp}}^{\mathrm{A}}}$, replace \eqref{an:detail:conditional.transform} by \eqref{an:spatial.transform.bound} from Lemma~\ref{an:spatial.transform} when deriving \eqref{an:conditional.detail:tilted.transform.comparison}. With this substitution, the Fourier comparison \eqref{an:conditional.detail:inversion.comparison} and the bin estimate \eqref{an:conditional.detail:weighted.bin} hold with $\smash{\ind_{\mathcal{A}_{T,\mathrm{sp}}^{\mathrm{A}}}}$ in place of $\smash{\ind_{\mathcal{A}_{T}^{\mathrm{A}}}}$. Removing the restriction to $\smash{\mathcal{A}_{T,\mathrm{sp}}^{\mathrm{A}}}$ from an unweighted arithmetic bin probability now costs at most
\[
\PP_{\mathrm{A}}((\mathcal{A}_{T,\mathrm{sp}}^{\mathrm{A}})^{c}) \leq \exp\{ - ct^{1 + 2\e_{0}} + Ct\}
\]
by \eqref{an:spatial.exceptional} in Lemma~\ref{an:spatial.transform}, which is still negligible relative to \eqref{an:detail:conditional.bin.lower}. The same choice of tilts and division therefore prove \eqref{an:conditional.bound} under the arithmetic law with $\smash{\mathcal{A}_{T,\mathrm{sp}}^{\mathrm{A}}}$ in place of $\mathcal{A}_{T}^{\mathrm{A}}$, both for $\mathscr{Q} = [0,1]^{k}$ and for $\mathscr{Q} = \prod_{j = 1}^{k}[0,b_{j}]$ with $b_{-} \leq b_{j} \leq b_{+}$, with the additional constant dependencies as stated.

For the Haar law, \eqref{an:local.factor} and \eqref{an:conditional.detail:undo.tilt} give the Haar version of \eqref{an:conditional.detail:weighted.bin} without $\smash{\ind_{\mathcal{A}_{T}^{\mathrm{H}}}}$ and with $\mathcal{E}_{T} = 0$. Consequently, $W_{T}^{\mathrm{H}}(\lambda_{k + 1})$ in \eqref{an:conditional.detail:ratio} is an untruncated expectation, and \eqref{an:conditional.detail:ratio} gives \eqref{an:conditional.detail:haar.untruncated} after expanding the exponent in \eqref{an:conditional.detail:quadratic.exponent}. The bin bounds \eqref{an:conditional.detail:variable.bin.integrals} give the same conclusion with $[0,1]^{k}$ replaced by $\smash{\mathscr{Q} = \prod_{j = 1}^{k}[0,b_{j}]}$, uniformly for $b_{-} \leq b_{j} \leq b_{+}$. The additional spatial representation changes neither the vector $\bb{Y}_{T}^{\mathrm{H}}$ nor the conditioning event $\{(Y_{T,r}^{\mathrm{H}})_{1 \leq r \leq k} \in \bb{v} + \mathscr{Q}\}$, so \eqref{an:conditional.detail:haar.untruncated} applies to the spatial version as well. Thus both Haar versions of \eqref{an:conditional.bound} hold on the whole probability space, with $\smash{\ind_{\mathcal{A}_{T}^{\mathrm{H}}}}$ and $\smash{\ind_{\mathcal{A}_{T,\mathrm{sp}}^{\mathrm{H}}}}$, respectively, replaced by $1$. Inserting either $\smash{\ind_{\mathcal{A}_{T}^{\mathrm{H}}}}$ or $\smash{\ind_{\mathcal{A}_{T,\mathrm{sp}}^{\mathrm{H}}}}$ cannot increase the nonnegative conditional expectation, which also proves both truncated Haar versions of \eqref{an:conditional.bound}.
\end{proof}

\begin{lemma}[One-point heights and tails]\label{an:height}
Let $\theta > 0$. The Haar quantile $u_{T}$ defined in \eqref{eq:height.quantile} satisfies, as $T \to \infty$,
\begin{equation}\label{an:height.scales}
\frac{u_{T}}{s_{T}} \to m(\theta), \qquad \frac{u_{T}}{V_{T}} \to \beta_{c}(\theta) = 2m(\theta).
\end{equation}
For each $\star \in \{\mathrm{A},\mathrm{H}\}$, as $T \to \infty$, uniformly in $h \in I_{\theta}$,
\begin{equation}\label{an:height.measure}
N_{T}\PP_{\star}(F_{h}^{\star} - u_{T} \in \rd{x}) \longrightarrow \beta_{c}(\theta) e^{-\beta_{c}(\theta) x}\rd{x}
\end{equation}
weakly on every bounded interval, and $N_{T}\PP_{\star}(F_{h}^{\star} - u_{T} > x) \to e^{-\beta_{c}(\theta) x}$ for every fixed $x$. Uniformly over integers $0 \leq j \leq u_{T}$,
\begin{equation}\label{an:depth.bound}
N_{T}\PP_{\star}(F_{h}^{\star} \in [u_{T} - j - 1,u_{T} - j]) \leq C\exp\left\{\frac{u_{T}j}{V_{T}} - \frac{j^{2}}{2V_{T}}\right\}.
\end{equation}
For every $\beta > \beta_{c}(\theta)$, every fixed $c_{\beta} \in (0,\beta - \beta_{c}(\theta))$, and every fixed $L,K \geq 0$,
\begin{equation}\label{an:height.tails}
\begin{aligned}
\limsup_{T \to \infty}\EE_{\star}\left[\ell_{T}^{-1}\int_{I_{\theta}}e^{\beta(F_{h}^{\star} - u_{T})}\ind_{\{F_{h}^{\star} - u_{T}< - L\}}\rd{h}\right] &\leq C e^{-c_{\beta}L}, \\
\limsup_{T \to \infty}\PP_{\star}\left(\sup_{h \in I_{\theta}}(F_{h}^{\star} - u_{T}) > K\right) &\leq C e^{-\beta_{c}(\theta) K}.
\end{aligned}
\end{equation}
The constants $C$ in the first and second estimates in \eqref{an:height.tails} are independent of $L$ and $K$, respectively.
\end{lemma}

\begin{proof}[Proof of Lemma~\ref{an:height}]
The proof first obtains a local estimate at heights of order $s_{T}$, then uses the exact quantile $u_{T}$ to normalize it. The two tail estimates require additional summation and spatial oscillation bounds.

For the local estimate, consider the single Haar coordinate $H(0)$ and take the additional coordinate in Lemma~\ref{an:conditional} to be zero. The scalar decomposition underlying \eqref{an:joint.cumulants} gives
\[
\log\EE_{\mathrm{H}}[e^{\lambda H(0)}] - \frac{\lambda^{2}V_{T}}{2}
= \sum_{p \leq X}\left\{\log\EE_{\mathrm{H}}\left[\exp\{\lambda\Re\sum_{j:p^{j}\leq X}U_{p}^{j}/(jp^{j\sig})\}\right] - \frac{\lambda^{2}}{4}\sum_{j:p^{j}\leq X}\frac{1}{j^{2}p^{2j\sig}}\right\}.
\]
For each prime, the random variable inside the exponential is centered and has variance $\frac{1}{2}\sum_{j:p^{j}\leq X}(j^{2}p^{2j\sig})^{-1}$. Moreover, the scalar bound \eqref{an:conditional.detail:prime.bound} gives
\begin{equation}\label{an:detail:height.prime.bound}
\left|\Re\sum_{j:p^{j}\leq X}\frac{U_{p}^{j}}{jp^{j\sig}}\right|
\leq \sum_{j = 1}^{\infty}\frac{1}{jp^{j/2}}
\leq Cp^{-1/2}.
\end{equation}
Taylor expansion of the exponential and then of the logarithm therefore yields
\begin{equation}\label{an:detail:height.prime.expansion}
\log\EE_{\mathrm{H}}\left[\exp\{\lambda\Re\sum_{j:p^{j}\leq X}U_{p}^{j}/(jp^{j\sig})\}\right]
= \frac{\lambda^{2}}{4}\sum_{j:p^{j}\leq X}\frac{1}{j^{2}p^{2j\sig}} + \OO(p^{-3/2}),
\end{equation}
uniformly for $\lambda$ in compact real sets. Thus the subtracted term cancels the entire quadratic contribution in each prime summand, leaving only the third- and higher-order remainder, bounded by $Cp^{-3/2}$.

For each fixed prime $p$, \eqref{eq:fixed.parameters} gives $\sig \to 1/2$ and $X \to \infty$. Domination by the convergent series $\sum_{j = 1}^{\infty}(jp^{j/2})^{-1}$ then gives uniform convergence of the prime variable in $U_{p}$ and locally uniform convergence of the prime remainder in $\lambda$. Since $\sum_{p}p^{-3/2} < \infty$, the remainder bound in \eqref{an:detail:height.prime.expansion} and dominated convergence give locally uniform convergence of the series of remainders. Exponentiating yields the locally uniform limit
\begin{equation}\label{an:detail:height.mod.gaussian}
\mathfrak{F}(\lambda) \leqdef \lim_{T \to \infty}e^{-\lambda^{2}V_{T}/2}\EE_{\mathrm{H}}[e^{\lambda H(0)}].
\end{equation}
The same argument applies to derivatives. Indeed, Taylor's formula applied to the first two identities in \eqref{an:conditional.detail:prime.derivatives}, followed by repeated differentiation of the tilted expectations for higher orders, bounds every fixed $\lambda$-derivative of each prime remainder by $Cp^{-3/2}$ on compact real sets. These derivatives also converge for each fixed prime, by the uniform convergence just established and the bound \eqref{an:detail:height.prime.bound}. Dominated convergence therefore gives the convergence in \eqref{an:detail:height.mod.gaussian} after any fixed number of differentiations with respect to $\lambda$.

For real $\lambda$, the limit $\mathfrak{F}(\lambda)$ is the exponential of an absolutely convergent real series. Thus it is positive and smooth, its derivatives are bounded on compact sets, and its minimum on every compact real interval is positive. Haar stationarity gives the same limit, with the same derivative convergence, for $H(h)$ at every $h \in I_{\theta}$.

For a target height $v$ satisfying $m_{-} \leq v/s_{T} \leq m_{+}$, where $m_{-} < m_{+}$ are fixed real numbers, choose the tilt $\lambda_{v} \leqdef v/V_{T}$. This tilt is bounded by \eqref{an:detail:cov.variance} in Lemma~\ref{an:covariance}. Under the normalized Haar tilt with density $e^{\lambda_{v}H(h)}/\EE_{\mathrm{H}}[e^{\lambda_{v}H(h)}]$, \eqref{an:conditional.detail:tilted.mean.covariance} gives mean $v + \OO(1)$ and variance $V_{T} + \OO(1)$. For fixed $a < b$, apply \eqref{an:local.factor} to the weighted bin function $q_{[a,b],\lambda_{v}}(x) \leqdef e^{-\lambda_{v}x}\ind_{[a,b]}(x)$. Undoing the tilt gives
\[
\begin{aligned}
\PP_{\mathrm{H}}(H(h) \in v + [a,b])
&= \EE_{\mathrm{H}}[e^{\lambda_{v}H(h)}]\EE_{\mathrm{H}}^{\lambda_{v},0}\left[e^{-\lambda_{v}H(h)}\ind_{\{H(h) \in v + [a,b]\}}\right] \\
&= e^{-\lambda_{v}v}\EE_{\mathrm{H}}[e^{\lambda_{v}H(h)}]\EE_{\mathrm{H}}^{\lambda_{v},0}[q_{[a,b],\lambda_{v}}(H(h) - v)] \\
&= \frac{\mathfrak{F}(v/V_{T})e^{-v^{2}/(2V_{T})}}{\sqrt{2\pi V_{T}}} \int_{a}^{b}e^{-(v/V_{T})x}\rd{x} \times \{1 + \oo(1)\},
\end{aligned}
\]
where \eqref{an:detail:height.mod.gaussian} supplies the exponential moment. The lower bound for $\mathfrak{F}$ on compact sets makes the error from this limit relative, and the weighted bin integral is uniformly bounded away from zero because $a < b$ is fixed and $v/V_{T}$ is bounded. For the arithmetic law, \eqref{an:conditional.detail:weighted.bin} adds an absolute error at most $\exp\{ - ct^{1 + 2\e_{0}} + Ct\}$ by \eqref{an:detail:conditional.bin.error}. The main term above is at least $e^{-Ct}$, since $|v| \leq \max\{|m_{-}|,|m_{+}|\}s_{T} = \OO(t)$ and $V_{T} \asymp t$, so this error is negligible relative to it. Removing the exceptional event costs at most its probability in \eqref{an:exceptional} of Lemma~\ref{an:transform}, which is also negligible relative to the main term. Consequently, for each $\star \in \{\mathrm{A},\mathrm{H}\}$,
\begin{equation}\label{an:detail:height.local.bin}
\PP_{\star}(F_{h}^{\star} \in v + [a,b]) = \frac{\mathfrak{F}(v/V_{T})e^{-v^{2}/(2V_{T})}}{\sqrt{2\pi V_{T}}} \int_{a}^{b}e^{-(v/V_{T})x}\rd{x} \times \{1 + \oo(1)\},
\end{equation}
uniformly for $m_{-} \leq v/s_{T} \leq m_{+}$ and $h \in I_{\theta}$. The same argument applies to open and half-open bins, by the pointwise approximants in \eqref{an:detail:interval.approximants}.

To pass from bins to upper tails, restrict $m_{-} \leq v/s_{T} \leq m_{+}$ with fixed $0 < m_{-} < m_{+}$, and partition $(v,\infty)$ into the disjoint bins $v + (j,j + 1]$, $j \in \N_0$. Choose a fixed $C > m_{+}$. For $j \in \N_0$ satisfying $v + j \leq Cs_{T}$, the shifted slope satisfies $m_{-} \leq (v + j)/s_{T} \leq C$, so \eqref{an:detail:height.local.bin} applies with lower and upper slope bounds $m_{-}$ and $C$. On this range, the ratio of its prefactor at $v + j$ to that at $v$ satisfies
\begin{equation}\label{an:detail:height.bin.ratio}
\begin{aligned}
\frac{\mathfrak{F}((v + j)/V_{T})e^{-(v + j)^{2}/(2V_{T})}}{\mathfrak{F}(v/V_{T})e^{-v^{2}/(2V_{T})}}
&\leq C_{1}\exp\left\{ - \frac{vj}{V_{T}} - \frac{j^{2}}{2V_{T}}\right\} \leq C_{1}e^{-m_{-}j}.
\end{aligned}
\end{equation}
The last inequality uses $v/V_{T} \geq m_{-}$ for all sufficiently large $T$, by \eqref{an:detail:cov.variance}. The weighted unit-bin integrals in \eqref{an:detail:height.local.bin} are also uniformly bounded for $m_{-} \leq (v + j)/s_{T} \leq C$. Thus \eqref{an:detail:height.bin.ratio} supplies a summable bound independent of $T$.

The remaining heights above $Cs_{T}$ are negligible for this fixed choice $C > m_{+}$. Indeed, fix $0 < \e_{0} < 1/4$. For each $h \in I_{\theta}$ and $\star \in \{\mathrm{A},\mathrm{H}\}$, let $\mathcal{A}_{T}^{\star} \leqdef \{|F_{h}^{\star}| \leq t^{1 + \e_{0}}\}$ be the truncation event in Lemma~\ref{an:transform} for the single coordinate $F_{h}^{\star}$, with its dependence on $h$ suppressed. Exponential Markov on $\mathcal{A}_{T}^{\star}$, the transform estimate \eqref{an:detail:conditional.transform} and exceptional-probability bound \eqref{an:exceptional} in Lemma~\ref{an:transform}, and the cumulant estimate \eqref{an:joint.cumulants} give, for fixed $\lambda > 0$,
\[
\begin{aligned}
\PP_{\star}(F_{h}^{\star} > Cs_{T})
&\leq e^{-\lambda Cs_{T}}\EE_{\star}[e^{\lambda F_{h}^{\star}}\ind_{\mathcal{A}_{T}^{\star}}] + \PP_{\star}((\mathcal{A}_{T}^{\star})^{c}) \\
&\leq \exp\left\{ - \lambda Cs_{T} + \frac{\lambda^{2}V_{T}}{2} + \OO(1)\right\} + \exp\{ - ct^{1 + 2\e_{0}} + C_{1}t\}.
\end{aligned}
\]
At $\lambda = 2C$, the first exponent is $ - C^{2}s_{T} + \OO(1)$ by \eqref{an:detail:cov.variance}. Since $C > m_{+}$, both exponential terms in the preceding upper bound are negligible relative to $\smash{V_{T}^{-1/2}e^{-v^{2}/(2V_{T})}}$, uniformly for $m_{-} \leq v/s_{T} \leq m_{+}$. For each fixed $j$, the ratio in \eqref{an:detail:height.bin.ratio} equals $e^{-(v/V_{T})j}\{1 + \oo(1)\}$. Moreover,
\[
\sum_{j = 0}^{\infty}e^{-(v/V_{T})j}\int_{0}^{1}e^{-(v/V_{T})x}\rd{x} = \int_{0}^{\infty}e^{-(v/V_{T})x}\rd{x} = \frac{V_{T}}{v}.
\]
The local bin estimate \eqref{an:detail:height.local.bin}, the summable bound \eqref{an:detail:height.bin.ratio}, and the negligible contribution above $Cs_{T}$ justify dominated summation. The preceding identity then gives
\begin{equation}\label{an:detail:height.scalar.tail}
\begin{aligned}
\PP_{\star}(F_{h}^{\star} > v)
&= \sum_{j = 0}^{\infty}\PP_{\star}(F_{h}^{\star} \in v + (j,j + 1]) \\
&= \frac{\mathfrak{F}(v/V_{T})e^{-v^{2}/(2V_{T})}}{\sqrt{2\pi V_{T}}}\sum_{j = 0}^{\infty}e^{-(v/V_{T})j}\int_{0}^{1}e^{-(v/V_{T})x}\rd{x} \times \{1 + \oo(1)\} \\
&= \frac{\mathfrak{F}(v/V_{T})}{(v/V_{T})\sqrt{2\pi V_{T}}}e^{-v^{2}/(2V_{T})} \times \{1 + \oo(1)\},
\end{aligned}
\end{equation}
uniformly for $m_{-} \leq v/s_{T} \leq m_{+}$ and $h \in I_{\theta}$.

For fixed $\gamma > 0$, apply \eqref{an:detail:height.scalar.tail} at $v = \gamma s_{T}$ with fixed bounds $0 < m_{-} < \gamma < m_{+}$. Together with the variance estimate \eqref{an:detail:cov.variance} and the scale identities \eqref{eq:effective.cutoff.bounds}, this gives
\begin{equation}\label{an:detail:height.slope}
\begin{aligned}
\log\left\{N_{T}\PP_{\mathrm{H}}(H(0) > \gamma s_{T})\right\}
&= \log{N_{T}} - \frac{\gamma^{2}s_{T}^{2}}{2V_{T}} - \frac{1}{2}\log{V_{T}} + \OO(1) \\
&= (1 + \theta - \gamma^{2})s_{T} + \oo(s_{T}).
\end{aligned}
\end{equation}
Taking $\gamma$ successively below and above $m(\theta) = \sqrt{1 + \theta}$ makes the last expression in \eqref{an:detail:height.slope} tend to $ + \infty$ and $ - \infty$, respectively. The quantile identity \eqref{eq:height.quantile} therefore places $u_{T}$ between the corresponding levels. It follows that $u_{T}/s_{T} \to m(\theta)$, and then $u_{T}/V_{T} \to 2m(\theta) = \beta_{c}(\theta)$ by \eqref{an:detail:cov.variance}. This proves \eqref{an:height.scales}.

Applying \eqref{an:detail:height.scalar.tail} at $v = u_{T}$, with fixed bounds $0 < m_{-} < m(\theta) < m_{+}$, gives the normalization
\begin{equation}\label{an:detail:height.quantile.prefactor}
N_{T}\frac{\mathfrak{F}(u_{T}/V_{T})e^{-u_{T}^{2}/(2V_{T})}}{\sqrt{2\pi V_{T}}} = \frac{u_{T}}{V_{T}}\{1 + \oo(1)\}.
\end{equation}
For each fixed $x$, both $u_{T}/s_{T}$ and $(u_{T} + x)/s_{T}$ lie in $[m_{-},m_{+}]$ for all sufficiently large $T$. The changes in the factors $e^{-u_{T}^{2}/(2V_{T})}$ and $\mathfrak{F}(u_{T}/V_{T})$ under the shift $u_{T} \mapsto u_{T} + x$ are described by
\begin{equation}\label{an:detail:height.shift}
\frac{(u_{T} + x)^{2} - u_{T}^{2}}{2V_{T}} = \frac{u_{T}x}{V_{T}} + \frac{x^{2}}{2V_{T}} \to \beta_{c}(\theta)x, \qquad \frac{\mathfrak{F}((u_{T} + x)/V_{T})}{\mathfrak{F}(u_{T}/V_{T})} \to 1.
\end{equation}
Divide \eqref{an:detail:height.scalar.tail} at $u_{T} + x$ by its Haar version at $u_{T}$, using \eqref{eq:height.quantile}, \eqref{an:height.scales}, and \eqref{an:detail:height.shift}. This gives $N_{T}\PP_{\star}(F_{h}^{\star} > u_{T} + x) \to e^{-\beta_{c}(\theta)x}$. Similarly, \eqref{an:detail:height.local.bin} at $v = u_{T}$, the normalization \eqref{an:detail:height.quantile.prefactor}, and the limit $u_{T}/V_{T} \to \beta_{c}(\theta)$ in \eqref{an:height.scales} yield
\[
N_{T}\PP_{\star}(F_{h}^{\star} - u_{T} \in [a,b]) \to \int_{a}^{b}\beta_{c}(\theta)e^{-\beta_{c}(\theta)x}\rd{x}.
\]
On each fixed bounded interval $[a,b]$, the rescaled masses of all subintervals are uniformly bounded for large $T$ and all $h \in I_{\theta}$. Step-function approximation using convergence of interval masses then proves \eqref{an:height.measure} restricted to $[a,b]$, uniformly in $h\in I_{\theta}$.

For the depth bound \eqref{an:depth.bound}, let $0 \leq j \leq u_{T}$ be an integer and take $v = u_{T} - j$ and $[a,b] = [ - 1,0]$ in \eqref{an:detail:height.local.bin}, using the lower slope bound $m_{-} = 0$ and any fixed upper slope bound $m_{+} > m(\theta)$. By \eqref{an:height.scales}, $0 \leq v \leq u_{T}$ then implies $m_{-} \leq v/s_{T} \leq m_{+}$ for all sufficiently large $T$. The factors $\mathfrak{F}(v/V_{T})$, $\mathfrak{F}(u_{T}/V_{T})^{-1}$, and $\int_{ - 1}^{0}e^{-(v/V_{T})x}\rd{x}$ remain bounded on this range. Multiplying \eqref{an:detail:height.local.bin} by $N_{T}$ and using \eqref{an:detail:height.quantile.prefactor} therefore gives
\[
\begin{aligned}
&N_{T}\PP_{\star}(F_{h}^{\star} \in [u_{T} - j - 1,u_{T} - j]) \\
&\qquad= \frac{N_{T}\mathfrak{F}((u_{T} - j)/V_{T})e^{-(u_{T} - j)^{2}/(2V_{T})}}{\sqrt{2\pi V_{T}}}\int_{ - 1}^{0}e^{-((u_{T} - j)/V_{T})x}\rd{x} \times \{1 + \oo(1)\} \\
&\qquad= \frac{u_{T}}{V_{T}}\frac{\mathfrak{F}((u_{T} - j)/V_{T})}{\mathfrak{F}(u_{T}/V_{T})}\exp\left\{\frac{u_{T}^{2} - (u_{T} - j)^{2}}{2V_{T}}\right\} \int_{ - 1}^{0}e^{-((u_{T} - j)/V_{T})x}\rd{x} \times \{1 + \oo(1)\} \\
&\qquad\leq C\exp\left\{\frac{u_{T}^{2} - (u_{T} - j)^{2}}{2V_{T}}\right\}
= C\exp\left\{\frac{u_{T}j}{V_{T}} - \frac{j^{2}}{2V_{T}}\right\}.
\end{aligned}
\]
This proves \eqref{an:depth.bound}.

To bound the low-height mass, observe that the weight $e^{\beta(F_{h}^{\star} - u_{T})}$ is at most $e^{-\beta j}$ on a bin at depth $j$. Estimate \eqref{an:depth.bound} therefore gives
\begin{equation}\label{an:detail:height.low.bin}
N_{T}\EE_{\star}\left[e^{\beta(F_{h}^{\star} - u_{T})}\ind_{\{F_{h}^{\star}\in[u_{T} - j - 1,u_{T} - j]\}}\right] \leq C\exp\left\{ - \left(\beta - \frac{u_{T}}{V_{T}}\right)j - \frac{j^{2}}{2V_{T}}\right\}
\end{equation}
for integers $0 \leq j \leq u_{T}$. Since $c_{\beta} < \beta - \beta_{c}(\theta)$, $u_{T}/V_{T} \to \beta_{c}(\theta)$ in \eqref{an:height.scales} makes the right side of \eqref{an:detail:height.low.bin} at most $Ce^{-c_{\beta}j}$ for all sufficiently large $T$. The nonnegative heights satisfying $F_{h}^{\star} - u_{T} < - L$ are covered by the bins $[u_{T} - j - 1,u_{T} - j]$ with integer $\max\{0,\lceil L - 1\rceil\} \leq j \leq u_{T}$. The part below zero satisfies
\begin{equation}\label{an:detail:height.negative.mass}
N_{T}\EE_{\star}\left[e^{\beta(F_{h}^{\star} - u_{T})}\ind_{\{F_{h}^{\star}<0\}}\right]
\leq N_{T}e^{-\beta u_{T}}
= \exp\{(m(\theta)^{2} - \beta m(\theta) + \oo(1))s_{T}\} = \oo(1),
\end{equation}
by \eqref{eq:effective.cutoff.bounds} and \eqref{an:height.scales}; its exponent is negative because $\beta > \beta_{c}(\theta) = 2m(\theta)$. Tonelli's theorem, $\ell_{T}^{-1}|I_{\theta}| = N_{T}$, \eqref{an:detail:height.low.bin} for the nonnegative heights, and \eqref{an:detail:height.negative.mass} for the negative heights now give
\[
\begin{aligned}
\EE_{\star}\left[\ell_{T}^{-1}\int_{I_{\theta}}e^{\beta(F_{h}^{\star} - u_{T})}\ind_{\{F_{h}^{\star} - u_{T}< - L\}}\rd{h}\right]
&\leq \sup_{h \in I_{\theta}}N_{T}\EE_{\star}\left[e^{\beta(F_{h}^{\star} - u_{T})}\ind_{\{F_{h}^{\star} - u_{T}< - L\}}\right] \\
&\leq C\sum_{\substack{j\in\Z\\\max\{0,\lceil L - 1\rceil\}\leq j\leq u_{T}}}e^{-c_{\beta}j} + N_{T}e^{-\beta u_{T}} \\
&\leq Ce^{-c_{\beta}L} + \oo(1).
\end{aligned}
\]
Taking the limit superior proves the first estimate in \eqref{an:height.tails}.

For the supremum estimate in \eqref{an:height.tails}, the field must be controlled between a discrete set of endpoints. Apply Lemma~\ref{an:spatial.transform} to $F^{\star}$ and its first scaled derivative $\ell_{T}(F^{\star})'$, where the prime denotes differentiation with respect to $h$. Its coefficient hypotheses hold because, for $r = 0,1,2$,
\[
\frac{\ell_{T}^{r}(j\log{p})^{r}}{jp^{j\sig}}
= p^{-j/2}\frac{(2\eta_{\sig}j\log{p})^{r}e^{-\eta_{\sig}j\log{p}}}{j}
\leq C_{r}p^{-j/2}.
\]
The case $r = 2$ verifies the derivative hypothesis for the family member $\ell_{T}(F^{\star})'$. Let $\mathcal{A}_{T,\mathrm{sp}}^{\star}$ be the resulting common event. Partition $I_{\theta}$ into at most $CN_{T}$ intervals $B$ of length at most $\ell_{T}$, and let $h_{B}$ be the left endpoint of $B$.

The primewise absolute bound for the derivative coordinate $\ell_{T}H'(x)$ satisfies
\[
\ell_{T}\log{p}\sum_{j:p^{j}\leq X}p^{-j\sig} \leq C\ell_{T}\log{p}\,p^{-\sig}.
\]
Partial summation using \eqref{analytic:eq:PNT}, followed by \eqref{an:effective.cutoff.estimates} in Lemma~\ref{an:covariance}, gives
\begin{equation}\label{an:detail:height.derivative.prime.bounds}
\sum_{p \leq X}\left(\ell_{T}\log{p}\sum_{j:p^{j}\leq X}p^{-j\sig}\right)^{2}
\leq C\ell_{T}^{2}\sum_{p \leq X}\frac{(\log{p})^{2}}{p^{1 + \ell_{T}}}
\leq C\ell_{T}^{2}\left(1 + \int_{0}^{\log{X}}ve^{-\ell_{T}v}\rd{v}\right) \leq C.
\end{equation}
The exact covariance sum \eqref{an:detail:cov.haar.sum} also gives, after differentiating its finite sum and taking absolute values,
\[
|\Cov_{\mathrm{H}}(H(h_{B}),\ell_{T}H'(x))|
\leq \frac{\ell_{T}}{2}\sum_{p^{j} \leq X}\frac{\log{p}}{jp^{2j\sig}}
\leq C\ell_{T}\left(1 + \int_{0}^{\log{X}}e^{-\ell_{T}v}\rd{v}\right) \leq C.
\]
The terms with $j \geq 2$ in this estimate are bounded by $C\ell_{T}\sum_{p}(\log{p})p^{-2}$. In particular, uniformly in $x,h_{B} \in B$,
\begin{equation}\label{an:detail:height.derivative.covariance}
\Var_{\mathrm{H}}(\ell_{T}H'(x)) \leq C, \qquad |\Cov_{\mathrm{H}}(H(h_{B}),\ell_{T}H'(x))| \leq C.
\end{equation}

Define the cell oscillation by $\operatorname{osc}_{B}(F^{\star}) \leqdef \sup_{h \in B}|F_{h}^{\star} - F_{h_{B}}^{\star}|$. The fundamental theorem of calculus gives \begin{equation}\label{eq:osc.bound}
\operatorname{osc}_{B}(F^{\star}) \leq \ell_{T}^{-1}\int_{B}|\ell_{T}(F^{\star})'(x)|\rd{x}.
\end{equation}
For fixed $\lambda > 0$, exponentiating this bound and applying Jensen's inequality yield
\[
\begin{aligned}
e^{\lambda\operatorname{osc}_{B}(F^{\star})}
&\leq \exp\left\{\frac{\lambda}{\ell_{T}}\int_{B}|\ell_{T}(F^{\star})'(x)|\rd{x}\right\} \\
&= \exp\left\{\left(1 - \frac{|B|}{\ell_{T}}\right)\cdot 0 + \frac{|B|}{\ell_{T}}\left(\frac{\lambda}{|B|}\int_{B}|\ell_{T}(F^{\star})'(x)|\rd{x}\right)\right\} \\
&\leq 1 - \frac{|B|}{\ell_{T}} + \ell_{T}^{-1}\int_{B}e^{\lambda|\ell_{T}(F^{\star})'(x)|}\rd{x}
\leq 1 + \ell_{T}^{-1}\int_{B}\left(e^{\lambda\ell_{T}(F^{\star})'(x)} + e^{-\lambda\ell_{T}(F^{\star})'(x)}\right)\rd{x}.
\end{aligned}
\]
For a short cell, the term $1 - |B|/\ell_{T}$ adds the missing mass at zero in Jensen's inequality. Condition now on $F_{h_{B}}^{\star} \in v + [ - 1,0]$, with $m_{-} \leq v/s_{T} \leq m_{+}$ for fixed real bounds $m_{-} < m_{+}$. In \eqref{an:conditional.bound} of Lemma~\ref{an:conditional}, take $k = 1$, the first coordinate $Y_{T,1}^{\star} = F_{h_{B}}^{\star}$, the last coordinate $Y_{T,2}^{\star} = \ell_{T}(F^{\star})'(x)$, and the last tilt $\lambda_{2} = \pm\lambda$. Then $\Sigma_{T} = (V_{T})$, $\bb{c}_{T} = (\Cov_{\mathrm{H}}(H(h_{B}),\ell_{T}H'(x)))$, and $\sig_{T}^{2} = \Var_{\mathrm{H}}(\ell_{T}H'(x))$. The conditioning interval is $(v - 1) + [0,1]$, so $\bb{v} = (v - 1)$, and its lower-endpoint slope satisfies $m_{-} - 1 \leq (v - 1)/s_{T} \leq m_{+}$ for all sufficiently large $T$. Thus the mean term $\bb{c}_{T}^{\mathsf{T}}\Sigma_{T}^{-1}\bb{v} = \Cov_{\mathrm{H}}(H(h_{B}),\ell_{T}H'(x))(v - 1)/V_{T}$ in \eqref{an:conditional.bound} is bounded in absolute value by $C|v - 1|/V_{T} \leq C$, while the variance term $\sig_{T}^{2} - \bb{c}_{T}^{\mathsf{T}}\Sigma_{T}^{-1}\bb{c}_{T}$ is at most $C$, by \eqref{an:detail:height.derivative.covariance}. Applying \eqref{an:conditional.bound} in Lemma~\ref{an:conditional} with $\lambda_{2} = \pm\lambda$ and the common spatial event $\mathcal{A}_{T,\mathrm{sp}}^{\star}$ bounds both conditional exponential moments by $C_{\lambda}$. Integrating these two conditional bounds in the preceding display and then applying Markov's inequality proves
\begin{equation}\label{an:detail:height.cell.oscillation}
\begin{aligned}
\EE_{\star}\left[e^{\lambda\operatorname{osc}_{B}(F^{\star})}\ind_{\mathcal{A}_{T,\mathrm{sp}}^{\star}}\,\middle|\,F_{h_{B}}^{\star} \in v + [ - 1,0]\right] &\leq C_{\lambda}, \\
\PP_{\star}\left(\{\operatorname{osc}_{B}(F^{\star}) \geq j\}\cap\mathcal{A}_{T,\mathrm{sp}}^{\star}\,\middle|\,F_{h_{B}}^{\star} \in v + [ - 1,0]\right) &\leq C_{\lambda}e^{-\lambda j}.
\end{aligned}
\end{equation}
By the untruncated Haar version of \eqref{an:conditional.bound} in Lemma~\ref{an:conditional}, both estimates in \eqref{an:detail:height.cell.oscillation} also hold with the spatial truncation event $\mathcal{A}_{T,\mathrm{sp}}^{\mathrm{H}}$ replaced by the whole probability space. The same replacement of $\mathcal{A}_{T,\mathrm{sp}}^{\mathrm{H}}$ is valid in \eqref{an:detail:maximum.positive.endpoints} and \eqref{an:detail:maximum.union.bound}, with $\smash{\ind_{\mathcal{A}_{T,\mathrm{sp}}^{\mathrm{H}}}}$ replaced by $1$ and the exceptional-probability term $\PP_{\mathrm{H}}((\mathcal{A}_{T,\mathrm{sp}}^{\mathrm{H}})^{c})$ replaced by $0$.

Fix $K \geq 0$ and let $w \leqdef u_{T} + K$. First, the endpoints already above $w$ contribute at most
\begin{equation}\label{an:detail:maximum.endpoint.tail}
\begin{aligned}
\sum_{B}\PP_{\star}(F_{h_{B}}^{\star} > w)
&\leq CN_{T}\sup_{h \in I_{\theta}}\PP_{\star}(F_{h}^{\star} > u_{T} + K) \\
&= C\sup_{h \in I_{\theta}}\frac{\PP_{\star}(F_{h}^{\star} > u_{T} + K)}{\PP_{\mathrm{H}}(H(0) > u_{T})}
= Ce^{-\beta_{c}(\theta)K} + \oo(1),
\end{aligned}
\end{equation}
by the shifted-tail limit proved from \eqref{an:detail:height.scalar.tail} and \eqref{eq:height.quantile}. Every endpoint in $[0,w]$ belongs to a bin $w - j + [ - 1,0]$ with $j \in \N_0$ and $w - j \geq 0$. For an endpoint in this bin, $F_{h_{B}}^{\star} \leq w - j$, so an exceedance $F_{h}^{\star} > w$ at some $h \in B$ implies $\operatorname{osc}_{B}(F^{\star}) \geq F_{h}^{\star} - F_{h_{B}}^{\star} > j$. The local bin formula \eqref{an:detail:height.local.bin} and the quantile normalization \eqref{an:detail:height.quantile.prefactor} give
\begin{equation}\label{an:detail:height.maximum.bin}
N_{T}\PP_{\star}(F_{h_{B}}^{\star} \in w - j + [ - 1,0])
\leq C\exp\left\{\frac{u_{T}^{2} - (u_{T} + K - j)^{2}}{2V_{T}}\right\}
\leq Ce^{-u_{T}K/V_{T}}e^{u_{T}j/V_{T}}.
\end{equation}
The constant in \eqref{an:detail:height.maximum.bin} is independent of fixed $K$: take $m_{-} = 0$ and $m_{+} = m(\theta) + 1$. For sufficiently large $T$, depending on $K$, all slopes $(w - j)/s_{T}$ under consideration lie in $[m_{-},m_{+}]$.

Choose $\lambda > \beta_{c}(\theta) + 1$. A union bound over the endpoint bins $w - j + [ - 1,0]$, followed by conditioning on each bin, allows the use of \eqref{an:detail:height.maximum.bin} and the second estimate in \eqref{an:detail:height.cell.oscillation} with $v = w - j$. Summing over at most $CN_{T}$ cells cancels the factor $N_{T}^{-1}$ from the bin probabilities and gives
\begin{equation}\label{an:detail:maximum.positive.endpoints}
\begin{aligned}
&\sum_{B}\PP_{\star}\left(\left\{\sup_{h \in B}F_{h}^{\star} > w,\ 0 \leq F_{h_{B}}^{\star} \leq w\right\}\cap\mathcal{A}_{T,\mathrm{sp}}^{\star}\right) \\
&\qquad \leq \sum_{B}\sum_{j = 0}^{\lfloor w\rfloor}\PP_{\star}(F_{h_{B}}^{\star} \in w - j + [ - 1,0]) \\
&\hspace{30mm} \times \PP_{\star}\left(\{\operatorname{osc}_{B}(F^{\star}) \geq j\}\cap\mathcal{A}_{T,\mathrm{sp}}^{\star}\,\middle|\,F_{h_{B}}^{\star} \in w - j + [ - 1,0]\right) \\
&\qquad \leq \frac{C}{N_{T}}e^{-u_{T}K/V_{T}}\sum_{B}\sum_{j = 0}^{\lfloor w\rfloor}e^{-(\lambda - u_{T}/V_{T})j} \\
&\qquad \leq Ce^{-u_{T}K/V_{T}}\sum_{j = 0}^{\infty}e^{-(\lambda - u_{T}/V_{T})j}
= \frac{Ce^{-u_{T}K/V_{T}}}{1 - e^{-(\lambda - u_{T}/V_{T})}}
\leq C'e^{-u_{T}K/V_{T}}.
\end{aligned}
\end{equation}
The denominator $1 - e^{-(\lambda - u_{T}/V_{T})}$ is bounded away from zero by \eqref{an:height.scales}. Hence the limit superior in \eqref{an:detail:maximum.positive.endpoints} is at most $C'e^{-\beta_{c}(\theta)K}$, with $C'$ independent of fixed $K$.

Finally, if the endpoint is negative, an exceedance of $w$ requires $\operatorname{osc}_{B}(F^{\star}) \geq w$. The stronger moment estimate \eqref{an:bounded.moments} in Lemma~\ref{an:polynomial} applies to the derivative by \eqref{an:detail:height.derivative.prime.bounds}. The integral bound \eqref{eq:osc.bound} for the oscillation and Jensen's inequality therefore give
\begin{equation}\label{an:detail:oscillation.moments}
\EE_{\star}[\operatorname{osc}_{B}(F^{\star})^{2D}]
\leq \left(\frac{|B|}{\ell_{T}}\right)^{2D - 1}\ell_{T}^{-1}\int_{B}\EE_{\star}[|\ell_{T}(F^{\star})'(x)|^{2D}]\rd{x}
\leq (CD)^{D}.
\end{equation}
Take $D \leqdef \lfloor c_{1}t^{2}\rfloor$, with fixed $c_{1} > 0$ sufficiently small. This choice satisfies $D = \oo(t^{\phi})$, as required for the arithmetic assertion of \eqref{an:bounded.moments}, because $\phi > 7/2$. Also $w \geq ct$ for all sufficiently large $T$ by \eqref{an:height.scales}. Markov's inequality and \eqref{an:detail:oscillation.moments} imply
\begin{equation}\label{an:detail:maximum.negative.endpoints}
\sum_{B}\PP_{\star}(\operatorname{osc}_{B}(F^{\star}) \geq w)
\leq CN_{T}\left(\frac{CD}{w^{2}}\right)^{D}
\leq \exp\{ - ct^{2} + Ct\} = \oo(1),
\end{equation}
where \eqref{eq:effective.cutoff.bounds} gives $\log{N_{T}} = \OO(t)$. On the common event $\mathcal{A}_{T,\mathrm{sp}}^{\star}$, split the cells according to whether their endpoint heights exceed $w$, are negative, or lie in $[0,w]$. For negative endpoints, an exceedance of $w$ requires $\operatorname{osc}_{B}(F^{\star}) \geq w$. The union bound, the exceptional-probability estimate $\PP_{\star}((\mathcal{A}_{T,\mathrm{sp}}^{\star})^{c}) = \oo(1)$ from \eqref{an:spatial.exceptional} in Lemma~\ref{an:spatial.transform}, the endpoint-tail estimate \eqref{an:detail:maximum.endpoint.tail}, the oscillation bound \eqref{an:detail:maximum.negative.endpoints}, and the nonnegative-endpoint estimate \eqref{an:detail:maximum.positive.endpoints} therefore give
\begin{equation}\label{an:detail:maximum.union.bound}
\begin{aligned}
\PP_{\star}\left(\sup_{h \in I_{\theta}}(F_{h}^{\star} - u_{T}) > K\right)
&\leq \PP_{\star}((\mathcal{A}_{T,\mathrm{sp}}^{\star})^{c}) + \sum_{B}\PP_{\star}(F_{h_{B}}^{\star} > w)
+ \sum_{B}\PP_{\star}(\operatorname{osc}_{B}(F^{\star}) \geq w) \\
&\qquad + \sum_{B}\PP_{\star}\left(\left\{\sup_{h \in B}F_{h}^{\star} > w,\ 0 \leq F_{h_{B}}^{\star} \leq w\right\}\cap\mathcal{A}_{T,\mathrm{sp}}^{\star}\right) \\
&\leq Ce^{-\beta_{c}(\theta)K} + C'e^{-u_{T}K/V_{T}} + \oo(1).
\end{aligned}
\end{equation}
Taking the limit superior and using $u_{T}/V_{T} \to \beta_{c}(\theta)$ from \eqref{an:height.scales} yields
\[
\limsup_{T \to \infty}\PP_{\star}\left(\sup_{h \in I_{\theta}}(F_{h}^{\star} - u_{T}) > K\right)
\leq (C + C')e^{-\beta_{c}(\theta)K},
\]
with $C + C'$ independent of fixed $K$. This proves the second estimate in \eqref{an:height.tails}.
\end{proof}

\subsubsection{Spatially integrated estimates at separated positions}

For $d_{0} > 0$ and $k \in \N$, define the separated position set
\begin{equation}\label{eq:D.k.T}
\mathscr{D}_{k,T}(d_{0}) \leqdef \{\bb{h} \in I_{\theta}^{\times k}: |h_{i} - h_{j}| \geq d_{0} \text{ for every }1 \leq i < j \leq k\},
\end{equation}
with $\mathscr{D}_{1,T}(d_{0}) = I_{\theta}$.

\begin{lemma}[Height kernels at separated positions]\label{an:remote}
Let $\theta > 0$. Fix $d_{0} > 0$ and $k \in \N$, independently of $T$. Let $\Theta(\rd{x}) \leqdef \beta_{c}(\theta)e^{-\beta_{c}(\theta)x}\rd{x}$. For a fixed bounded nonnegative Borel function $f:\R \to \R$ with compact support and a discontinuity set of Lebesgue measure zero, define
\[
\mu_{\star,f}(\rd{h}) \leqdef \ell_{T}^{-1}f(F_{h}^{\star} - u_{T})\rd{h}, \quad h \in I_{\theta}, \qquad \Theta(f) \leqdef \int_{\R}f(x)\Theta(\rd{x}).
\]
Then, for each $\star \in \{\mathrm{A},\mathrm{H}\}$,
\begin{equation}\label{an:remote.bound}
\limsup_{T \to \infty}\EE_{\star}\left[\int_{\mathscr{D}_{k,T}(d_{0})}\prod_{j = 1}^{k}\mu_{\star,f}(\rd{h_{j}})\right] \leq \Theta(f)^{k}.
\end{equation}

More precisely, let $\lambda_{T} \leqdef u_{T}/V_{T}$ and define the deterministic kernel
\begin{equation}\label{eq:kernel.T.h}
\mathcal{K}_{T}(\bb{h}) \leqdef \frac{\EE_{\mathrm{H}}\left[\exp\left\{\lambda_{T}\sum_{j = 1}^{k}H(h_{j})\right\}\right]}{\left(\EE_{\mathrm{H}}[e^{\lambda_{T}H(0)}]\right)^{k}}, \qquad \bb{h} \in \mathscr{D}_{k,T}(d_{0}).
\end{equation}
Fix $m_{0} \in \N$, $L,M > 0$, and $0 < \delta \leq 2L$. Consider real functions of the form
\begin{equation}\label{eq:def.g.finite.sum.ind}
g(\bb{x}) \leqdef \sum_{r = 1}^{m}c_{r}\prod_{j = 1}^{k}\ind_{[a_{r,j},b_{r,j})}(x_{j}), \qquad \bb{x} \in \R^{k},
\end{equation}
where $1 \leq m \leq m_{0}$, $|c_{r}| \leq M$, $ - L \leq a_{r,j} < b_{r,j} \leq L$, and $b_{r,j} - a_{r,j} \geq \delta$. Write $\Theta^{\otimes k}(g) \leqdef \int_{\R^{k}}g(\bb{x})\prod_{j = 1}^{k}\Theta(\rd{x_{j}})$. Then, as $T\to \infty$,
\begin{equation}\label{an:remote.local.approximation}
\begin{aligned}
&\ell_{T}^{-k}\int_{\mathscr{D}_{k,T}(d_{0})}\bigg|\EE_{\star}\left[g(F_{h_{1}}^{\star} - u_{T},\ldots,F_{h_{k}}^{\star} - u_{T})\right] - N_{T}^{-k}\mathcal{K}_{T}(\bb{h})\Theta^{\otimes k}(g)\bigg|\rd{\bb{h}} \to 0,
\end{aligned}
\end{equation}
uniformly over this class of functions and over $\star \in \{\mathrm{A},\mathrm{H}\}$. The coefficients and bin endpoints may depend on $T$ within the stated bounds; $d_{0}$, $k$, $m_{0}$, $L$, $M$, and $\delta$ remain fixed.
\end{lemma}

\begin{proof}[Proof of Lemma~\ref{an:remote}]
The proof first controls the spatial average of the exact Haar kernel for exponential moments. A local estimate under a common tilt then relates this kernel to the centered height bins. The same averaged kernel bound, with tilt $\beta$, will also be used in Lemma~\ref{an:capped.subcritical}.

Fix $0 < \e_{0} < 1/4$ and $M_{0} > \beta_{c}(\theta) + 1$. All estimates involving a real tilt $\lambda$ below are uniform for $|\lambda| \leq M_{0}$. Let
\[
M_{T}(\lambda) \leqdef \EE_{\mathrm{H}}[e^{\lambda H(0)}], \qquad H_{\leq P}(h) \leqdef \Re\sum_{\substack{p \leq P \\ p^{j}\leq X}}\frac{U_{p}^{j}}{jp^{j\sig + \ii jh}}, \qquad P \leqdef t^{B},
\]
where $B > 0$ is fixed and will be chosen sufficiently large in terms of the fixed parameters. In particular, $P < \exp(t^{3}) < X$ for all sufficiently large $T$. Every power of a prime $p \leq P$ permitted by $p^{j} \leq X$ is retained in $H_{\leq P}$. Define
\[
\Pi_{T,P,\lambda}(h) \leqdef \frac{e^{\lambda H_{\leq P}(h)}}{\EE_{\mathrm{H}}[e^{\lambda H_{\leq P}(0)}]}, \qquad \overline{\Pi}_{T,P,\lambda} \leqdef \frac{1}{|I_{\theta}|}\int_{I_{\theta}}\Pi_{T,P,\lambda}(h)\rd{h}.
\]
The denominator in the definition of $\Pi_{T,P,\lambda}(h)$ does not depend on the spatial position $h$, because multiplication of each $U_{p}$ by $p^{-\ii h}$ preserves the joint Haar law. Both this denominator and $M_{T}(\lambda)$ are at least one by Jensen's inequality.

For every fixed $\kappa \in \N$, the following two estimates will be proved below:
\begin{equation}\label{an:low.average}
\EE_{\mathrm{H}}[|\overline{\Pi}_{T,P,\lambda} - 1|^{\kappa}] \to 0, \qquad
\EE_{\mathrm{H}}\left[\prod_{i = 1}^{\kappa}\Pi_{T,P,\lambda}(h_{i})\right] \leq C_{\kappa}\exp(C_{\kappa}\log\log{t}).
\end{equation}
The second estimate in \eqref{an:low.average} bounds the low-prime kernel before averaging, even when the positions coincide. To prove it, the geometric prime bound \eqref{an:conditional.detail:prime.bound} and Mertens' estimate \eqref{analytic:eq:Mertens} give
\[
\sum_{p \leq P}\left(\sum_{j:p^{j}\leq X}\frac{1}{jp^{j\sig}}\right)^{2}
\leq C\sum_{p \leq P}\frac{1}{p}
\leq C\log\log{P} + C
= C\log\log{t} + C_{B}.
\]
This controls every low-prime variance, and Cauchy--Schwarz gives the same bound for the absolute value of every low-prime covariance. Apply \eqref{an:joint.cumulants} and \eqref{an:conditional.detail:remainder} to the coordinates $H_{\leq P}(h_{1}),\ldots,H_{\leq P}(h_{\kappa})$ with common tilt $\lambda$. It follows that
\[
\begin{aligned}
\log\EE_{\mathrm{H}}\left[e^{\lambda\sum_{i = 1}^{\kappa}H_{\leq P}(h_{i})}\right]
&= \frac{\lambda^{2}}{2}\sum_{i,j = 1}^{\kappa}\Cov_{\mathrm{H}}(H_{\leq P}(h_{i}),H_{\leq P}(h_{j})) + \OO_{\kappa}(1) \\
&\leq C_{\kappa}\log\log{t} + C_{\kappa,B}.
\end{aligned}
\]
The normalizing denominator $(\EE_{\mathrm{H}}[e^{\lambda H_{\leq P}(0)}])^{\kappa}$ is at least one, so its logarithm is nonnegative. Subtracting this logarithm and exponentiating proves the second estimate in \eqref{an:low.average}, uniformly over the positions, with a power of $\log{t}$ independent of fixed $B$.

For the first assertion of \eqref{an:low.average}, fix $A > 4$ and let $\widetilde{H}_{\leq P}$ be obtained from $H_{\leq P}$ by also imposing $p^{j} \leq P^{A}$. The omitted terms have the deterministic bound
\begin{equation}\label{an:remote.detail:power.tail}
\sup_{h}|H_{\leq P}(h) - \widetilde{H}_{\leq P}(h)| \leq \sum_{p \leq P}\sum_{j:p^{j}>P^{A}}p^{-j/2} \leq CP^{1 - A/2}.
\end{equation}
For each prime, the inner sum in \eqref{an:remote.detail:power.tail} is a geometric tail bounded by $CP^{-A/2}$, and there are at most $P$ primes. Exponentiating \eqref{an:remote.detail:power.tail} and taking Haar expectations gives the following bounds on the ratio of the normalizing constants:
\begin{equation}\label{an:remote.detail:power.tail.next}
e^{-C|\lambda|P^{1 - A/2}}
\leq \frac{\EE_{\mathrm{H}}[e^{\lambda H_{\leq P}(0)}]}{\EE_{\mathrm{H}}[e^{\lambda\widetilde{H}_{\leq P}(0)}]}
\leq e^{C|\lambda|P^{1 - A/2}}.
\end{equation}
The ratio of the normalized exponentials is therefore between $e^{-C M_{0}P^{1 - A/2}}$ and $e^{C M_{0}P^{1 - A/2}}$, after changing $C$. In particular, \eqref{an:remote.detail:power.tail}, \eqref{an:remote.detail:power.tail.next}, and the moment bound in \eqref{an:low.average} imply
\begin{equation}\label{an:remote.detail:normalized.exponential}
\begin{aligned}
\left\|\Pi_{T,P,\lambda}(h) - \frac{e^{\lambda\widetilde{H}_{\leq P}(h)}}{\EE_{\mathrm{H}}[e^{\lambda\widetilde{H}_{\leq P}(0)}]}\right\|_{L^{\kappa}(\mathrm{H})}
&= \left\|\Pi_{T,P,\lambda}(h)\left(1 - \frac{e^{\lambda(\widetilde{H}_{\leq P}(h) - H_{\leq P}(h))}}{\EE_{\mathrm{H}}[e^{\lambda\widetilde{H}_{\leq P}(0)}]/\EE_{\mathrm{H}}[e^{\lambda H_{\leq P}(0)}]}\right)\right\|_{L^{\kappa}(\mathrm{H})} \\
&\leq \big(e^{CP^{1 - A/2}} - 1\big)\|\Pi_{T,P,\lambda}(h)\|_{L^{\kappa}(\mathrm{H})} \\
&\leq CP^{1 - A/2}\|\Pi_{T,P,\lambda}(h)\|_{L^{\kappa}(\mathrm{H})} \\
&\leq CP^{1 - A/2}(\log{t})^{C_{\kappa}} = \oo(1),
\end{aligned}
\end{equation}
where $\|Y\|_{L^{\kappa}(\mathrm{H})} \leqdef (\EE_{\mathrm{H}}[|Y|^{\kappa}])^{1/\kappa}$ denotes the norm with respect to the Haar probability law $\PP_{\mathrm{H}}$. Thus the power cutoff $p^{j} \leq P^{A}$ can be inserted before proving the spatial-average limit in \eqref{an:low.average}. This cutoff bounds the frequencies of a polynomial approximation to the exponential.

Let $D \leqdef \lfloor c\, t/\log{t}\rfloor$, with $c > 0$ chosen later, and write
\[
Q_{D}(h) \leqdef \sum_{r = 0}^{D}\frac{(\lambda\widetilde{H}_{\leq P}(h))^{r}}{r!}.
\]
The exponential remainder inequality and Cauchy--Schwarz give, for $Z \leqdef \lambda\widetilde{H}_{\leq P}(h)$,
\[
\left\|e^{Z} - \sum_{r = 0}^{D}\frac{Z^{r}}{r!}\right\|_{L^{\kappa}(\mathrm{H})} \leq \frac{\|e^{|Z|}\|_{L^{2\kappa}(\mathrm{H})}\|Z\|_{L^{2\kappa(D + 1)}(\mathrm{H})}^{D + 1}}{(D + 1)!}.
\]
Hoeffding's argument underlying \eqref{an:haar.moments} in Lemma~\ref{an:polynomial}, with squared primewise bounds summing to $\OO(\log\log{t})$, gives
\[
\|e^{|Z|}\|_{L^{2\kappa}(\mathrm{H})} \leq \exp(C_{\kappa}\log\log{t}), \qquad
\|Z\|_{L^{2\kappa(D + 1)}(\mathrm{H})}^{D + 1} \leq (C_{\kappa}(D + 1)\log\log{t})^{(D + 1)/2}.
\]
Consequently,
\begin{equation}\label{an:remote.detail:taylor}
\begin{aligned}
\|e^{\lambda\widetilde{H}_{\leq P}(h)} - Q_{D}(h)\|_{L^{\kappa}(\mathrm{H})}
&\leq \exp(C_{\kappa}\log\log{t})\left(\frac{C_{\kappa}\log\log{t}}{D + 1}\right)^{(D + 1)/2} \\
&\leq \exp\left\{C_{\kappa}\log\log{t} - \frac{D + 1}{2}\log\left(\frac{D + 1}{C_{\kappa}\log\log{t}}\right)\right\}
\leq e^{-c't}.
\end{aligned}
\end{equation}

Every nonzero spatial frequency of $Q_{D}$ is $\log(n_{1}/n_{2})$ with $n_{1},n_{2} \leq P^{AD}$. For distinct such integers, $|\log(n_{1}/n_{2})| \geq |n_{1} - n_{2}|/\max\{n_{1},n_{2}\} \geq P^{-AD}$, so
\[
\left|\frac{1}{|I_{\theta}|}\int_{I_{\theta}}(n_{1}/n_{2})^{-\ii h}\rd{h}\right| \leq \frac{2P^{AD}}{|I_{\theta}|}.
\]
The sum of the absolute coefficients in $\widetilde{H}_{\leq P}$ is at most $\smash{C\sum_{p \leq P}p^{-1/2} \leq C\sqrt{P}}$, and hence the corresponding coefficient sum for $Q_{D}$ is at most $(C\sqrt{P})^{D}$, after increasing $C$. Unique factorization implies that its zero-frequency coefficient is deterministic and equals $\EE_{\mathrm{H}}[Q_{D}(0)]$. It follows, for every realization of the phases, that
\begin{equation}\label{an:remote.detail:spatial.polynomial}
\begin{aligned}
\left|\frac{1}{|I_{\theta}|}\int_{I_{\theta}}Q_{D}(h)\rd{h} - \EE_{\mathrm{H}}[Q_{D}(0)]\right|
&\leq \frac{CP^{AD}(C\sqrt{P})^{D}}{|I_{\theta}|} = \frac{C^{D + 1}}{2}t^{B(A + 1/2)D}e^{-\theta t} \\
&= \exp\{(cB(A + 1/2) - \theta + \oo(1))t\} \\
&= \oo(1),
\end{aligned}
\end{equation}
provided $cB(A + 1/2) < \theta$. The constant $B$ is fixed first, and $c$ is then chosen to satisfy this inequality. To transfer \eqref{an:remote.detail:spatial.polynomial} to the exponential, Minkowski's inequality and the uniform remainder estimate \eqref{an:remote.detail:taylor} give
\[
\begin{aligned}
\left\|\frac{1}{|I_{\theta}|}\int_{I_{\theta}}e^{\lambda\widetilde{H}_{\leq P}(h)}\rd{h} - \EE_{\mathrm{H}}[e^{\lambda\widetilde{H}_{\leq P}(0)}]\right\|_{L^{\kappa}(\mathrm{H})}
&\leq \frac{1}{|I_{\theta}|}\int_{I_{\theta}}\|e^{\lambda\widetilde{H}_{\leq P}(h)} - Q_{D}(h)\|_{L^{\kappa}(\mathrm{H})}\rd{h} \\
&\qquad + \left\|\frac{1}{|I_{\theta}|}\int_{I_{\theta}}Q_{D}(h)\rd{h} - \EE_{\mathrm{H}}[Q_{D}(0)]\right\|_{L^{\kappa}(\mathrm{H})} \\
&\qquad + \|Q_{D}(0) - e^{\lambda\widetilde{H}_{\leq P}(0)}\|_{L^{\kappa}(\mathrm{H})} \\
&\leq 2e^{-c't} + \oo(1) = \oo(1).
\end{aligned}
\]
Dividing by $\EE_{\mathrm{H}}[e^{\lambda\widetilde{H}_{\leq P}(0)}] \geq 1$ and using the normalized-exponential comparison \eqref{an:remote.detail:normalized.exponential} proves the first assertion of \eqref{an:low.average}. Its lower-moment cases also imply
\begin{equation}\label{an:remote.detail:low.moments}
\left|\EE_{\mathrm{H}}[\overline{\Pi}_{T,P,\lambda}^{\kappa}] - 1\right|
\leq \sum_{r = 1}^{\kappa}\binom{\kappa}{r}\EE_{\mathrm{H}}[|\overline{\Pi}_{T,P,\lambda} - 1|^{r}] = \oo(1).
\end{equation}
All limits remain uniform for $|\lambda| \leq M_{0}$.

The next goal is to extend the low-prime average estimate \eqref{an:remote.detail:low.moments} to the full normalized joint exponential moment, yielding the bound \eqref{an:remote.detail:normalized.kernel} below. Independence over primes factors this moment into the low-prime kernel and the contribution from $p > P$. It remains to show that this latter factor can be replaced by one after spatial averaging against the low-prime kernel, with an error tending to zero. The second estimate in \eqref{an:low.average} will control the low-prime weight multiplying this error.

To control the contribution from $p > P$, expand the logarithm of each prime factor to second order and subtract the corresponding contributions from $k\log{M_{T}(\lambda)}$, arising from the normalizing denominator $M_{T}(\lambda)^{k}$. The diagonal variance terms then cancel, leaving the cross-covariance terms and a uniformly vanishing remainder. Split these primes at $\exp(t^{3})$. Above this threshold, the separation $|d| \geq d_{0}$ makes the cross covariances uniformly small. Between $P$ and $\exp(t^{3})$, a mean-square estimate controls their exponential contribution after spatial averaging. For the high primes, the prime summation estimate \eqref{an:detail:cov.prime.integral}, obtained from \eqref{analytic:eq:PNT}, together with the integration-by-parts bound below, gives, uniformly for $d_{0} \leq |d| \leq 2|I_{\theta}|$,
\begin{equation}\label{an:high.remote}
\left|\sum_{\exp(t^{3})<p\leq X}\frac{e^{-2\eta_{\sig}\log{p}}\cos(d\log{p})}{p}\right| \leq \frac{C_{d_{0}}}{t^{3}} + C(1 + |I_{\theta}|)t^{C}e^{-ct^{3/2}} = \oo(1).
\end{equation}
To see the first term, let $w(v) \leqdef e^{-2\eta_{\sig}v}/v$, which is positive and decreasing because $w'(v) = - 2\eta_{\sig}e^{-2\eta_{\sig}v}/v - e^{-2\eta_{\sig}v}/v^{2} = - e^{-2\eta_{\sig}v}(2\eta_{\sig}/v + 1/v^{2}) < 0$ for $v > 0$. Integration by parts gives
\[
\left|\int_{t^{3}}^{\log{X}}w(v)\cos(dv)\rd{v}\right|
\leq \frac{w(t^{3}) + w(\log{X}) - \int_{t^{3}}^{\log{X}}w'(v)\rd{v}}{|d|}
= \frac{2w(t^{3})}{|d|}
\leq \frac{2}{d_{0}t^{3}}.
\]
The contribution to the covariance from prime powers with $j \geq 2$ and $p > P$ is bounded by $C\sum_{p>P}p^{-2} = \OO(P^{-1})$, uniformly over the positions.

For the middle-prime polynomial
\[
\Upsilon_{T}(d) \leqdef \sum_{P<p\leq\exp(t^{3})}p^{-1 - 2\eta_{\sig}}e^{\ii d\log{p}},
\]
the mean-value estimate \eqref{analytic:eq:MV}, the triangle inequality, and Mertens' estimate \eqref{analytic:eq:Mertens} give
\begin{equation}\label{an:middle.square}
\frac{1}{4|I_{\theta}|}\int_{ - 2|I_{\theta}|}^{2|I_{\theta}|}|\Upsilon_{T}(d)|^{2}\rd{d} \leq \frac{C}{P} + \frac{C\log{t}}{|I_{\theta}|}, \qquad |\Upsilon_{T}(d)| \leq C\log{t},
\end{equation}
because
\[
\sum_{P<p\leq\exp(t^{3})}|p^{-1 - 2\eta_{\sig}}|^{2} \leq \sum_{n>P}n^{-2} \leq \frac{C}{P}, \qquad
\sum_{P<p\leq\exp(t^{3})}p|p^{-1 - 2\eta_{\sig}}|^{2} \leq \sum_{p\leq\exp(t^{3})}\frac{1}{p} \leq C\log{t}.
\]
The second estimate in \eqref{an:middle.square} follows directly from $\smash{|\Upsilon_{T}(d)| \leq \sum_{p\leq\exp(t^{3})}p^{-1}}$ and Mertens' estimate \eqref{analytic:eq:Mertens}. For any nonnegative measurable function $\psi$, the change of variables $d \leqdef h - h'$ and Tonelli's theorem give
\begin{equation}\label{an:remote.detail:pair.differences}
\frac{1}{|I_{\theta}|^{2}}\int_{I_{\theta}^{\times 2}}\psi(h - h')\rd{h}\rd{h'} = \frac{1}{|I_{\theta}|^{2}}\int_{ - |I_{\theta}|}^{|I_{\theta}|}(|I_{\theta}| - |d|)\psi(d)\rd{d} \leq \frac{1}{|I_{\theta}|}\int_{ - |I_{\theta}|}^{|I_{\theta}|}\psi(d)\rd{d}.
\end{equation}
For $|\lambda| \leq M_{0}$ and $\bb{h} = (h_{1},\ldots,h_{k}) \in I_{\theta}^{\times k}$, define
\[
S_{T}(\bb{h},\lambda) \leqdef \frac{\lambda^{2}}{2} \sum_{1 \leq i < j \leq k}\Re(\Upsilon_{T}(h_{i} - h_{j})).
\]
The pointwise estimate in \eqref{an:middle.square} gives $|S_{T}(\bb{h},\lambda)| \leq C_{k}\log{t}$. To control its spatial mean square, Cauchy--Schwarz for the finite sum over pairs and the difference identity \eqref{an:remote.detail:pair.differences} give
\begin{equation}\label{an:remote.detail:middle.spatial.square}
\begin{aligned}
\frac{1}{|I_{\theta}|^{k}}\int_{I_{\theta}^{\times k}}|S_{T}(\bb{h},\lambda)|^{2}\rd{\bb{h}}
&\leq \frac{C_{k}}{|I_{\theta}|^{2}}\sum_{i<j}\int_{I_{\theta}^{\times 2}}|\Upsilon_{T}(h - h')|^{2}\rd{h}\rd{h'} \\
&\leq \frac{C_{k}}{|I_{\theta}|}\int_{ - |I_{\theta}|}^{|I_{\theta}|}|\Upsilon_{T}(d)|^{2}\rd{d}
\leq C_{k}\left(P^{-1} + \frac{\log{t}}{|I_{\theta}|}\right).
\end{aligned}
\end{equation}
The last inequality is \eqref{an:middle.square}. Since the normalized joint exponential-moment kernel in \eqref{an:remote.detail:kernel.factorization} below contains the middle-prime factor $e^{S_{T}(\bb{h},\lambda)}$, combine \eqref{an:remote.detail:middle.spatial.square} with $|e^{x} - 1| \leq e^{|x|}|x|$ and Cauchy--Schwarz to obtain
\begin{equation}\label{an:remote.detail:middle.exponential}
\begin{aligned}
\frac{1}{|I_{\theta}|^{k}}\int_{I_{\theta}^{\times k}}|e^{S_{T}(\bb{h},\lambda)} - 1|\rd{\bb{h}}
&\leq t^{C_{k}}\left(\frac{1}{|I_{\theta}|^{k}}\int_{I_{\theta}^{\times k}}|S_{T}(\bb{h},\lambda)|^{2}\rd{\bb{h}}\right)^{1/2} \\
&\leq C_{k}t^{C_{k}}\left(P^{-1} + \frac{\log{t}}{|I_{\theta}|}\right)^{1/2}.
\end{aligned}
\end{equation}
The power $C_{k}$ of $t$ in \eqref{an:remote.detail:middle.exponential} is independent of $B$. Choose $B > 2C_{k} + 2$. The right side of \eqref{an:remote.detail:middle.exponential} then tends to zero after multiplication by the fixed powers of $\log{t}$ in \eqref{an:low.average}: its first term is bounded by $Ct^{-1}(\log{t})^{C}$, and its second term decays exponentially in $t$ because $\theta > 0$. This fixes $B$ before the choice of $c$ in \eqref{an:remote.detail:spatial.polynomial}. When $k = 1$, all pair sums are empty and the estimates \eqref{an:remote.detail:middle.spatial.square} and \eqref{an:remote.detail:middle.exponential} hold with $S_{T} = 0$.

For the prime decomposition of the exponential moments, define
\[
\mathcal{K}_{T}^{\mathrm{low}}(\bb{h},\lambda) \leqdef \EE_{\mathrm{H}}\left[\prod_{j = 1}^{k}\Pi_{T,P,\lambda}(h_{j})\right], \qquad \mathcal{K}_{T}^{\mathrm{mid}}(\bb{h},\lambda) \leqdef e^{S_{T}(\bb{h},\lambda)}.
\]
Here $\mathcal{K}_{T}^{\mathrm{low}}$ retains the exact Haar expectation over the low-prime phases, whereas $\mathcal{K}_{T}^{\mathrm{mid}}$ is the deterministic approximation obtained after averaging over the middle-prime phases and expanding the logarithm of their normalized exponential moment to second order, as justified below. Independence over primes factors both the joint exponential moment and each normalizing factor $M_{T}(\lambda)$ into contributions from $p \leq P$ and $p > P$. The former contributions give $\mathcal{K}_{T}^{\mathrm{low}}$ exactly. For $p > P$, expand the logarithm of each prime factor to second order as in \eqref{an:joint.cumulants}. The linear terms vanish, the diagonal quadratic terms cancel upon division by $M_{T}(\lambda)^{k}$, and the cross terms sum to $\lambda^{2}\sum_{i<j}\Cov_{\mathrm{H}}(H_{>P}(h_{i}),H_{>P}(h_{j}))$, where $H_{>P} \leqdef H - H_{\leq P}$. The prime remainder estimate \eqref{an:conditional.detail:remainder}, restricted to $p > P$, bounds their sum by $C_{k}\sum_{p>P}p^{-3/2} \leq C_{k}P^{-1/2}$. The first-prime covariance between $P$ and $\exp(t^{3})$ equals $\Re(\Upsilon_{T}(h_{i} - h_{j}))/2$, by \eqref{an:detail:cov.haar.sum}. Thus its contribution is exactly $S_{T}(\bb{h},\lambda)$. The higher powers contribute $\OO_{k}(P^{-1})$, and \eqref{an:high.remote} bounds the remaining first-prime terms. More explicitly,
\[
\begin{aligned}
\log\left(\frac{\EE_{\mathrm{H}}\big[e^{\lambda\sum_{j = 1}^{k}H(h_{j})}\big]}{M_{T}(\lambda)^{k}\mathcal{K}_{T}^{\mathrm{low}}(\bb{h},\lambda)}\right)
&= \lambda^{2}\sum_{i<j}\Cov_{\mathrm{H}}(H_{>P}(h_{i}),H_{>P}(h_{j})) + \OO_{k}(P^{-1/2}) \\
&= S_{T}(\bb{h},\lambda) + \OO_{k,d_{0}}\left(P^{-1/2} + t^{-3} + (1 + |I_{\theta}|)t^{C}e^{-ct^{3/2}}\right).
\end{aligned}
\]
The error tends uniformly to zero on $\mathscr{D}_{k,T}(d_{0})$, defined in \eqref{eq:D.k.T}. Exponentiating proves
\begin{equation}\label{an:remote.detail:kernel.factorization}
\frac{\EE_{\mathrm{H}}\big[e^{\lambda\sum_{j = 1}^{k}H(h_{j})}\big]}{M_{T}(\lambda)^{k}} = \mathcal{K}_{T}^{\mathrm{low}}(\bb{h},\lambda)\mathcal{K}_{T}^{\mathrm{mid}}(\bb{h},\lambda)\{1 + \oo(1)\},
\end{equation}
uniformly for $\bb{h} \in \mathscr{D}_{k,T}(d_{0})$ and $|\lambda| \leq M_{0}$. At $\lambda = \lambda_{T} \equiv u_{T}/V_{T}$, the ratio on the left in \eqref{an:remote.detail:kernel.factorization} is exactly $\mathcal{K}_{T}(\bb{h})$ from \eqref{eq:kernel.T.h}. When $k = 1$, this ratio and both factors on the right equal one by stationarity.

The second estimate in \eqref{an:low.average} bounds $\mathcal{K}_{T}^{\mathrm{low}}$ pointwise by a fixed power of $\log{t}$. Multiplying \eqref{an:remote.detail:middle.exponential} by the second bound in \eqref{an:low.average} and using the choice of $B$ gives
\begin{equation}\label{an:remote.detail:kernel.average}
\begin{aligned}
\frac{1}{|I_{\theta}|^{k}}\int_{\mathscr{D}_{k,T}(d_{0})}\mathcal{K}_{T}^{\mathrm{low}}(\bb{h},\lambda)|\mathcal{K}_{T}^{\mathrm{mid}}(\bb{h},\lambda) - 1|\rd{\bb{h}} &= \oo(1), \\
\frac{1}{|I_{\theta}|^{k}}\int_{\mathscr{D}_{k,T}(d_{0})}\mathcal{K}_{T}^{\mathrm{low}}(\bb{h},\lambda)\mathcal{K}_{T}^{\mathrm{mid}}(\bb{h},\lambda)\rd{\bb{h}} &\leq \frac{1}{|I_{\theta}|^{k}}\int_{I_{\theta}^{\times k}}\mathcal{K}_{T}^{\mathrm{low}}(\bb{h},\lambda)\rd{\bb{h}} + \oo(1) \\
&= \EE_{\mathrm{H}}[\overline{\Pi}_{T,P,\lambda}^{k}] + \oo(1) = 1 + \oo(1).
\end{aligned}
\end{equation}
Tonelli's theorem gives the equality, and \eqref{an:remote.detail:low.moments} with $\kappa = k$ gives the last limit. All errors are uniform for $|\lambda| \leq M_{0}$. Combining \eqref{an:remote.detail:kernel.factorization} and \eqref{an:remote.detail:kernel.average} proves
\begin{equation}\label{an:remote.detail:normalized.kernel}
\frac{1}{|I_{\theta}|^{k}M_{T}(\lambda)^{k}}\int_{\mathscr{D}_{k,T}(d_{0})}\EE_{\mathrm{H}}\big[e^{\lambda\sum_{j = 1}^{k}H(h_{j})}\big]\rd{\bb{h}} \leq 1 + \oo(1).
\end{equation}
Estimate \eqref{an:remote.detail:normalized.kernel} holds for every fixed $k \in \N$ and $d_{0} > 0$, uniformly for $|\lambda| \leq M_{0}$. In particular, substituting $\lambda = \lambda_{T}$ into \eqref{an:remote.detail:normalized.kernel} gives
\begin{equation}\label{an:remote.detail:kernel.spatial.average}
\frac{1}{|I_{\theta}|^{k}}\int_{\mathscr{D}_{k,T}(d_{0})}\mathcal{K}_{T}(\bb{h})\rd{\bb{h}} \leq 1 + \oo(1).
\end{equation}

It remains to relate the exact kernel $\mathcal{K}_{T}(\bb{h})$ defined in \eqref{eq:kernel.T.h} to the joint probabilities that the centered heights satisfy $F_{h_{j}}^{\star} - u_{T} \in [a_{j},b_{j})$ for every $1 \leq j \leq k$, where the bin endpoints lie in $[ - L,L]$ and each bin has width at least $\delta$. For $\bb{h} \in \mathscr{D}_{k,T}(d_{0})$, let $\bb{F}^{\star} \leqdef (F_{h_{1}}^{\star},\ldots,F_{h_{k}}^{\star})$ and $\Sigma_{T} \leqdef \Cov_{\mathrm{H}}(\bb{F}^{\mathrm{H}})$. Equation~\eqref{an:overlap.geometry} in Lemma~\ref{an:covariance} gives $\Sigma_{T} = V_{T}\mathrm{Id}_{k} + \OO_{k,d_{0}}(\log{t})$ entrywise. By \eqref{an:covariance.comparison} in Lemma~\ref{an:covariance}, \eqref{eq:effective.cutoff.bounds}, and \eqref{eq:parameter.scales}, $V_{T} = s_{T}/2 + \OO(1) \asymp t$. Since $k$ is fixed, the entrywise covariance error also bounds the operator norm, so
\begin{equation}\label{eq:det.Sigma.T.approx}
\begin{aligned}
\|V_{T}^{-1}\Sigma_{T} - \mathrm{Id}_{k}\|_{\mathrm{op}}
&\leq C_{k,d_{0}}\frac{\log{t}}{t} = \oo(1), \\
\det(\Sigma_{T})
&= V_{T}^{k}\det(V_{T}^{-1}\Sigma_{T})
= V_{T}^{k}\left(1 + \OO_{k,d_{0}}\left(\frac{\log{t}}{t}\right)\right).
\end{aligned}
\end{equation}
In particular, the eigenvalues of $\Sigma_{T}$ are bounded above and below by fixed positive multiples of $t$, as required for the local estimate \eqref{an:local.factor}. By \eqref{an:height.scales} in Lemma~\ref{an:height}, $\lambda_{T} \to \beta_{c}(\theta)$, so $|\lambda_{T}| \leq M_{0}$ for all sufficiently large $T$. The joint Haar tilt with the same parameter $\lambda_{T}$ at all $k$ positions has density
\[
\frac{\rd{\PP}_{\mathrm{H}}^{\lambda_{T}\bb{1}_{k},0}}{\rd{\PP}_{\mathrm{H}}} = \frac{e^{\lambda_{T}\sum_{j = 1}^{k}H(h_{j})}}{M_{T}(\lambda_{T})^{k}\mathcal{K}_{T}(\bb{h})}.
\]
The tilted-mean identity \eqref{an:conditional.detail:tilted.mean.covariance}, applied with an additional zero coordinate, gives
\[
\begin{aligned}
\EE_{\mathrm{H}}^{\lambda_{T}\bb{1}_{k},0}[\bb{F}^{\mathrm{H}}] - u_{T}\bb{1}_{k}
&= \Sigma_{T}\lambda_{T}\bb{1}_{k} - V_{T}\lambda_{T}\bb{1}_{k} + \OO_{k}(1) \\
&= (\Sigma_{T} - V_{T}\mathrm{Id}_{k})\lambda_{T}\bb{1}_{k} + \OO_{k}(1)
= \OO_{k,d_{0}}(\log{t}) = \oo(\sqrt{t}).
\end{aligned}
\]
Thus the Haar tilt with parameter $\lambda_{T}$ at all $k$ positions places the mean within the range required in \eqref{an:local.factor}, namely $\|\bb{\mu}_{T} - \bb{v}\|_{2} = \oo(\sqrt{t})$ with $\bb{v} = u_{T}\bb{1}_{k}$, uniformly over the separated positions. Using this common tilt preserves the exact kernel $\mathcal{K}_{T}(\bb{h})$ in the change of measure.

Let $\mathscr{Q} \leqdef \prod_{j = 1}^{k}[a_{j},b_{j})$ be one product of bins with endpoints in $[ - L,L]$ and widths at least $\delta$. For $q_{[a,b),\lambda_{T}}(x) \leqdef e^{-\lambda_{T}x}\ind_{[a,b)}(x)$, the following bound holds:
\[
|\widehat{q}_{[a,b),\lambda_{T}}(\omega) - \widehat{q}_{[a,b),\lambda_{T}}(0)| = \left|\int_{\R}(e^{-\ii\omega x} - 1)q_{[a,b),\lambda_{T}}(x)\rd{x}\right| \leq |\omega|\int_{ - L}^{L}|x|e^{M_{0}|x|}\rd{x} \leq C_{L,M_{0}}|\omega|.
\]
Thus $\widehat{q}_{[a,b),\lambda_{T}}$ is uniformly continuous at zero over the permitted bins. Together with the approximant bounds in \eqref{an:detail:interval.approximants}, this supplies the Fourier continuity estimate \eqref{an:conditional.detail:fourier.continuity} used in the proof of \eqref{an:local.factor}. The interval majorants and minorants satisfy the uniform $L^{1}$ gap bound in \eqref{an:detail:interval.approximants}, and the product-bin integrals satisfy
\[
(\delta e^{-M_{0}L})^{k} \leq \prod_{j = 1}^{k}\int_{a_{j}}^{b_{j}}e^{-\lambda_{T}x}\rd{x} \leq (2Le^{M_{0}L})^{k}.
\]
Thus the $\oo(1)$ term in the multiplicative factor in \eqref{an:local.factor} is uniform over these bins, including endpoints that depend on $T$. The same interval majorants and minorants apply to the half-open bins $[a,b)$. Undoing the common tilt and applying the local estimate \eqref{an:local.factor}, together with $\det(\Sigma_{T}) = V_{T}^{k}\{1 + \oo(1)\}$ from \eqref{eq:det.Sigma.T.approx}, gives
\begin{equation}\label{an:remote.detail:tilted.bin}
\begin{aligned}
\PP_{\mathrm{H}}(\bb{F}^{\mathrm{H}} - u_{T}\bb{1}_{k} \in \mathscr{Q})
&= \EE_{\mathrm{H}}\left[\prod_{j = 1}^{k}\ind_{[a_{j},b_{j})}(H(h_{j}) - u_{T})\right] \\
&= M_{T}(\lambda_{T})^{k}\mathcal{K}_{T}(\bb{h})\EE_{\mathrm{H}}^{\lambda_{T}\bb{1}_{k},0}\left[e^{-\lambda_{T}\sum_{j = 1}^{k}H(h_{j})}\prod_{j = 1}^{k}\ind_{[a_{j},b_{j})}(H(h_{j}) - u_{T})\right] \\
&= M_{T}(\lambda_{T})^{k}\mathcal{K}_{T}(\bb{h})e^{-k\lambda_{T}u_{T}}\EE_{\mathrm{H}}^{\lambda_{T}\bb{1}_{k},0}\left[\prod_{j = 1}^{k}q_{[a_{j},b_{j}),\lambda_{T}}(H(h_{j}) - u_{T})\right] \\
&= \mathcal{K}_{T}(\bb{h})\left(\frac{M_{T}(\lambda_{T})e^{-\lambda_{T}u_{T}}}{\sqrt{2\pi V_{T}}}\right)^{k}\left\{\prod_{j = 1}^{k}\int_{a_{j}}^{b_{j}}e^{-\lambda_{T}x}\rd{x}\right\}\{1 + \oo(1)\}.
\end{aligned}
\end{equation}
The locally uniform limit in \eqref{an:detail:height.mod.gaussian} gives $M_{T}(\lambda_{T}) = e^{\lambda_{T}^{2}V_{T}/2}\mathfrak{F}(\lambda_{T})\{1 + \oo(1)\}$ with $\lambda_{T} = u_{T}/V_{T}$. Since $\lambda_{T}u_{T} = u_{T}^{2}/V_{T}$, \eqref{an:detail:height.quantile.prefactor} yields
\begin{equation}\label{an:remote.detail:quantile.normalization}
\begin{aligned}
N_{T}\frac{M_{T}(\lambda_{T})e^{-\lambda_{T}u_{T}}}{\sqrt{2\pi V_{T}}}
= N_{T}\frac{\mathfrak{F}(\lambda_{T})e^{-u_{T}^{2}/(2V_{T})}}{\sqrt{2\pi V_{T}}}\{1 + \oo(1)\}
= \lambda_{T}\{1 + \oo(1)\} \to \beta_{c}(\theta).
\end{aligned}
\end{equation}
Moreover, $\lambda_{T}e^{-\lambda_{T}x} \to \beta_{c}(\theta)e^{-\beta_{c}(\theta)x}$ uniformly on $[ - L,L]$. Substituting \eqref{an:remote.detail:quantile.normalization} into \eqref{an:remote.detail:tilted.bin} therefore gives
\begin{equation}\label{an:remote.detail:haar.bin.approximation}
\PP_{\mathrm{H}}(\bb{F}^{\mathrm{H}} - u_{T}\bb{1}_{k} \in \mathscr{Q}) = N_{T}^{-k}\mathcal{K}_{T}(\bb{h})\Theta^{\otimes k}(\mathscr{Q})\{1 + \oo(1)\},
\end{equation}
where the $\oo(1)$ term is uniform over the permitted bins and positions. Here $\Theta(\rd{x}) = \beta_{c}(\theta)e^{-\beta_{c}(\theta)x}\rd{x}$ and $\Theta^{\otimes k}$ denotes its $k$-fold product measure, so
\[
\Theta^{\otimes k}(\mathscr{Q}) = \prod_{j = 1}^{k}\Theta([a_{j},b_{j})) = \prod_{j = 1}^{k}\int_{a_{j}}^{b_{j}}\beta_{c}(\theta)e^{-\beta_{c}(\theta)x}\rd{x}.
\]

The uniform $\oo(1)$ term in \eqref{an:remote.detail:haar.bin.approximation} is bounded in absolute value by a deterministic $\e_{T} \geq 0$ tending to zero. Using $\ell_{T}^{-k}N_{T}^{-k} = |I_{\theta}|^{-k}$ and \eqref{an:remote.detail:kernel.spatial.average} therefore gives
\[
\begin{aligned}
&\ell_{T}^{-k}\int_{\mathscr{D}_{k,T}(d_{0})}\left|\PP_{\mathrm{H}}(\bb{F}^{\mathrm{H}} - u_{T}\bb{1}_{k} \in \mathscr{Q}) - N_{T}^{-k}\mathcal{K}_{T}(\bb{h})\Theta^{\otimes k}(\mathscr{Q})\right|\rd{\bb{h}} \\
&\qquad \leq \e_{T}\ell_{T}^{-k}N_{T}^{-k}\Theta^{\otimes k}(\mathscr{Q})\int_{\mathscr{D}_{k,T}(d_{0})}\mathcal{K}_{T}(\bb{h})\rd{\bb{h}}
\leq \e_{T}\Theta([ - L,L])^{k}\{1 + \oo(1)\} = \oo(1),
\end{aligned}
\]
uniformly over the permitted bins. By linearity of expectation and of integration against $\Theta^{\otimes k}$, together with the triangle inequality, the integrated absolute error for the finite sum defining $g$ in \eqref{eq:def.g.finite.sum.ind} is at most $m_{0}M$ times this uniform bound, since $\sum_{r = 1}^{m}|c_{r}| \leq m_{0}M$. This proves \eqref{an:remote.local.approximation} for $\star = \mathrm{H}$.

To prove \eqref{an:remote.local.approximation} also for $\star = \mathrm{A}$, choose the common spatial event $\smash{\mathcal{A}_{T,\mathrm{sp}}^{\mathrm{A}}}$ of Lemma~\ref{an:spatial.transform} for the Euler field $\{F_{h}^{\mathrm{A}}\}_{h \in I_{\theta}}$, including all prime powers $p^{j} \leq X$. For each prime power, let $a_{p,j}(h) \leqdef p^{-j(\sig + \ii h)}/j$. These coefficients satisfy the hypotheses of Lemma~\ref{an:spatial.transform}, since $\sig = 1/2 + \eta_{\sig}$ and $\ell_{T} = 2\eta_{\sig}$ give
\[
|a_{p,j}(h)| + \ell_{T}|a_{p,j}'(h)| = p^{-j/2}\frac{e^{-\eta_{\sig}j\log{p}}(1 + 2\eta_{\sig}j\log{p})}{j} \leq Cp^{-j/2}.
\]
Here the family consists of a single field, so $m_{T} = 1$ and
\[
\mathcal{A}_{T,\mathrm{sp}}^{\mathrm{A}} = \left\{\sup_{h \in I_{\theta}}|F_{h}^{\mathrm{A}}| \leq t^{1 + \e_{0}}\right\}.
\]
In particular, the spatial event $\smash{\mathcal{A}_{T,\mathrm{sp}}^{\mathrm{A}}}$ is the same for every choice of the positions $h_{1},\ldots,h_{k}$ and bin endpoints. The Haar bin estimate \eqref{an:remote.detail:haar.bin.approximation} holds without restriction to the spatial truncation event $\mathcal{A}_{T,\mathrm{sp}}^{\mathrm{H}}$, and the Haar expectation in \eqref{an:remote.detail:arithmetic.transform} is also untruncated. Thus the Haar quantities in \eqref{an:remote.detail:haar.bin.approximation} and \eqref{an:remote.detail:arithmetic.transform} are evaluated on the whole probability space, without inserting $\smash{\ind_{\mathcal{A}_{T,\mathrm{sp}}^{\mathrm{H}}}}$, whereas the arithmetic expectation in \eqref{an:remote.detail:arithmetic.transform} is restricted to $\smash{\mathcal{A}_{T,\mathrm{sp}}^{\mathrm{A}}}$.

Let $\Omega \leqdef t^{\e_{0}/2}/4$, so $4\Omega = t^{\e_{0}/2}$, and define the weighted product-bin function
\[
q_{T,\mathscr{Q}}(\bb{x}) \leqdef \prod_{j = 1}^{k}q_{[a_{j},b_{j}),\lambda_{T}}(x_{j}) = e^{-\lambda_{T}\sum_{j = 1}^{k}x_{j}}\ind_{\mathscr{Q}}(\bb{x}), \qquad \bb{x} \in \R^{k}.
\]
Let $q_{T,\mathscr{Q}}^{-}$ and $q_{T,\mathscr{Q}}^{+}$ be its lower and upper approximants, obtained by the construction preceding \eqref{an:conditional.detail:product.approximants} with $\lambda_{j} = \lambda_{T}$ and $\mathcal{J}_{j} = [a_{j},b_{j})$. Thus $q_{T,\mathscr{Q}}^{-} \leq q_{T,\mathscr{Q}} \leq q_{T,\mathscr{Q}}^{+}$ pointwise, and \eqref{an:conditional.detail:product.approximants} gives
\begin{equation}\label{an:remote.detail:product.approximants}
\|q_{T,\mathscr{Q}}^{+} - q_{T,\mathscr{Q}}^{-}\|_{1} \leq C\Omega^{-1}, \quad
\|\widehat{q}_{T,\mathscr{Q}}^{\pm}\|_{\infty} \leq \|q_{T,\mathscr{Q}}^{\pm}\|_{1} \leq C, \quad
\operatorname{supp}(\widehat{q}_{T,\mathscr{Q}}^{\pm}) \subseteq [ - 4\Omega,4\Omega]^{k}.
\end{equation}
These bounds are uniform over the permitted bins, since their endpoints lie in $[ - L,L]$ and $|\lambda_{T}| \leq M_{0}$ for all sufficiently large $T$.

Apply \eqref{an:spatial.transform.bound} in Lemma~\ref{an:spatial.transform} to the $k$ evaluations $F_{h_{1}}^{\star},\ldots,F_{h_{k}}^{\star}$, with complex parameters $z_{j} = \lambda_{T} + \ii\omega_{j}$. The real parts are bounded by $M_{0}$ and the imaginary parts are bounded by $t^{\e_{0}/2}$ throughout the Fourier cube $[ - 4\Omega,4\Omega]^{k}$. Consequently,
\begin{equation}\label{an:remote.detail:arithmetic.transform}
\left|\EE_{\mathrm{A}}\left[\exp\{\sum_{j = 1}^{k}(\lambda_{T} + \ii\omega_{j})F_{h_{j}}^{\mathrm{A}}\}\ind_{\mathcal{A}_{T,\mathrm{sp}}^{\mathrm{A}}}\right] - \EE_{\mathrm{H}}\left[\exp\{\sum_{j = 1}^{k}(\lambda_{T} + \ii\omega_{j})F_{h_{j}}^{\mathrm{H}}\}\right]\right|
\leq \exp\{ - ct^{1 + 2\e_{0}} + Ct\},
\end{equation}
uniformly over the positions and Fourier variables. Fourier inversion, as in \eqref{an:conditional.detail:fourier.inversion} and \eqref{an:conditional.detail:inversion.comparison}, expresses the difference of the weighted approximant expectations as an integral of the transform difference in \eqref{an:remote.detail:arithmetic.transform}. Taking absolute values removes the factor $\smash{e^{-\ii u_{T}\sum_{j = 1}^{k}\omega_{j}}}$, whose modulus is one, and gives
\begin{equation}\label{an:remote.detail:arithmetic.inversion}
\begin{aligned}
&\left|\EE_{\mathrm{A}}\left[\exp\{\lambda_{T}\sum_{j = 1}^{k}F_{h_{j}}^{\mathrm{A}}\}q_{T,\mathscr{Q}}^{\pm}(\bb{F}^{\mathrm{A}} - u_{T}\bb{1}_{k})\ind_{\mathcal{A}_{T,\mathrm{sp}}^{\mathrm{A}}}\right] - \EE_{\mathrm{H}}\left[\exp\{\lambda_{T}\sum_{j = 1}^{k}F_{h_{j}}^{\mathrm{H}}\}q_{T,\mathscr{Q}}^{\pm}(\bb{F}^{\mathrm{H}} - u_{T}\bb{1}_{k})\right]\right| \\
&\qquad \leq \frac{1}{(2\pi)^{k}}\int_{[ - 4\Omega,4\Omega]^{k}}|\widehat{q}_{T,\mathscr{Q}}^{\pm}(\bb{\omega})|\exp\{ - ct^{1 + 2\e_{0}} + Ct\}\rd{\bb{\omega}} \\
&\qquad \leq C(8\Omega)^{k}\exp\{ - ct^{1 + 2\e_{0}} + Ct\}.
\end{aligned}
\end{equation}
The last inequality uses the uniform Fourier bound from \eqref{an:remote.detail:product.approximants} and the volume $(8\Omega)^{k}$ of the Fourier cube.

The definition of $q_{T,\mathscr{Q}}$ gives the exact identity
\[
e^{-k\lambda_{T}u_{T}}\exp\{\lambda_{T}\sum_{j = 1}^{k}F_{h_{j}}^{\mathrm{A}}\}q_{T,\mathscr{Q}}(\bb{F}^{\mathrm{A}} - u_{T}\bb{1}_{k}) = \ind_{\{\bb{F}^{\mathrm{A}} - u_{T}\bb{1}_{k} \in \mathscr{Q}\}}.
\]
Multiplying $q_{T,\mathscr{Q}}^{-} \leq q_{T,\mathscr{Q}} \leq q_{T,\mathscr{Q}}^{+}$ by the nonnegative factor $\smash{e^{-k\lambda_{T}u_{T}}\exp\{\lambda_{T}\sum_{j = 1}^{k}F_{h_{j}}^{\mathrm{A}}\}\ind_{\mathcal{A}_{T,\mathrm{sp}}^{\mathrm{A}}}}$ and taking expectations therefore bounds the bin probability on $\smash{\mathcal{A}_{T,\mathrm{sp}}^{\mathrm{A}}}$ between the two arithmetic approximant expectations in \eqref{an:remote.detail:arithmetic.inversion}, multiplied by $e^{-k\lambda_{T}u_{T}}$. For each of the corresponding Haar expectations, undoing the common tilt as in \eqref{an:remote.detail:tilted.bin} gives
\[
\begin{aligned}
&e^{-k\lambda_{T}u_{T}}\EE_{\mathrm{H}}\left[\exp\{\lambda_{T}\sum_{j = 1}^{k}F_{h_{j}}^{\mathrm{H}}\}q_{T,\mathscr{Q}}^{\pm}(\bb{F}^{\mathrm{H}} - u_{T}\bb{1}_{k})\right] \\
&\qquad= e^{-k\lambda_{T}u_{T}}M_{T}(\lambda_{T})^{k}\mathcal{K}_{T}(\bb{h})
\EE_{\mathrm{H}}^{\lambda_{T}\bb{1}_{k},0}\left[q_{T,\mathscr{Q}}^{\pm}(\bb{F}^{\mathrm{H}} - u_{T}\bb{1}_{k})\right]
= N_{T}^{-k}\mathcal{K}_{T}(\bb{h})\Theta^{\otimes k}(\mathscr{Q})\{1 + \oo(1)\}.
\end{aligned}
\]
For the last equality, use the versions of the local estimate \eqref{an:local.factor} for the lower and upper approximants, proved immediately before that equation, together with the determinant estimate \eqref{eq:det.Sigma.T.approx} and the normalization \eqref{an:remote.detail:quantile.normalization}. Both approximants give the same leading bin integral $\smash{\int_{\R^{k}}q_{T,\mathscr{Q}}(\bb{x})\rd{\bb{x}} = \prod_{j = 1}^{k}\int_{a_{j}}^{b_{j}}e^{-\lambda_{T}x}\rd{x}}$, because their $L^{1}$ distances from $q_{T,\mathscr{Q}}$ are $\OO(\Omega^{-1})$ by \eqref{an:remote.detail:product.approximants}. The lower bound on the weighted bin integrals established before \eqref{an:remote.detail:tilted.bin} makes the resulting $\oo(1)$ relative error uniform over the permitted bins.

After multiplication by $e^{-k\lambda_{T}u_{T}}$, the absolute error in \eqref{an:remote.detail:arithmetic.inversion} is at most
\[
Ce^{-k\lambda_{T}u_{T}}(8\Omega)^{k}\exp\{ - ct^{1 + 2\e_{0}} + Ct\} \leq \exp\{ - ct^{1 + 2\e_{0}} + C_{k,L}t\}.
\]
Indeed, $\lambda_{T}$ is bounded, $u_{T} = \OO(s_{T}) = \OO(t)$ by \eqref{an:height.scales} in Lemma~\ref{an:height} and the scale identities \eqref{eq:effective.cutoff.bounds} and \eqref{eq:parameter.scales}, and $(8\Omega)^{k} = 2^{k}t^{k\e_{0}/2}$. Thus the factors preceding the transform comparison error can be absorbed into $\exp(C_{k,L}t)$. Combining the lower and upper bounds from \eqref{an:remote.detail:arithmetic.inversion} with the common Haar leading term proves
\begin{equation}\label{an:remote.detail:arithmetic.bin}
\begin{aligned}
&\PP_{\mathrm{A}}(\{\bb{F}^{\mathrm{A}} - u_{T}\bb{1}_{k} \in \mathscr{Q}\} \cap \mathcal{A}_{T,\mathrm{sp}}^{\mathrm{A}}) \\
&\qquad = N_{T}^{-k}\mathcal{K}_{T}(\bb{h})\Theta^{\otimes k}(\mathscr{Q})\{1 + \oo(1)\} + \OO\left(\exp\{ - ct^{1 + 2\e_{0}} + C_{k,L}t\}\right).
\end{aligned}
\end{equation}
Here the relative error $\oo(1)$ comes from the Haar local estimate and the bin approximations, while the additional absolute error comes from the transform comparison.

To remove the event $\smash{\mathcal{A}_{T,\mathrm{sp}}^{\mathrm{A}}}$ from \eqref{an:remote.detail:arithmetic.bin}, use \eqref{an:spatial.exceptional} in Lemma~\ref{an:spatial.transform}. It gives
\[
\begin{aligned}
0 &\leq \PP_{\mathrm{A}}(\bb{F}^{\mathrm{A}} - u_{T}\bb{1}_{k} \in \mathscr{Q}) - \PP_{\mathrm{A}}(\{\bb{F}^{\mathrm{A}} - u_{T}\bb{1}_{k} \in \mathscr{Q}\} \cap \mathcal{A}_{T,\mathrm{sp}}^{\mathrm{A}}) \\
&\leq \PP_{\mathrm{A}}((\mathcal{A}_{T,\mathrm{sp}}^{\mathrm{A}})^{c}) \leq \exp\{ - ct^{1 + 2\e_{0}} + Ct\}.
\end{aligned}
\]
Together with \eqref{an:remote.detail:arithmetic.bin} and the Haar estimate \eqref{an:remote.detail:haar.bin.approximation}, this gives a deterministic sequence $\e_{T} \geq 0$ tending to zero such that
\begin{equation}\label{an:rare.kernel}
\begin{aligned}
&\left|\PP_{\star}(\bb{F}^{\star} - u_{T}\bb{1}_{k} \in \mathscr{Q}) - N_{T}^{-k}\mathcal{K}_{T}(\bb{h})\Theta^{\otimes k}(\mathscr{Q})\right| \\
&\qquad \leq \e_{T}N_{T}^{-k}\mathcal{K}_{T}(\bb{h}) + \exp\{ - ct^{1 + 2\e_{0}} + C_{k,L}t\},
\end{aligned}
\end{equation}
uniformly over $\bb{h} \in \mathscr{D}_{k,T}(d_{0})$, all the permitted bins, and $\star \in \{\mathrm{A},\mathrm{H}\}$. The bound $\Theta^{\otimes k}(\mathscr{Q}) \leq \Theta([ - L,L])^{k}$ has been absorbed into $\e_{T}$, enlarging the sequence used for the Haar estimate if necessary. The sequence may depend on the fixed $k,d_{0},L,\delta$ and the fixed parameters of the paper, but not on the positions or bin endpoints. The last exponential term can be omitted under the Haar law.

Integrate \eqref{an:rare.kernel} over $\mathscr{D}_{k,T}(d_{0})$ and multiply by $\ell_{T}^{-k}$. The first term is controlled by $\ell_{T}^{-k}N_{T}^{-k} = |I_{\theta}|^{-k}$ and the spatial kernel bound \eqref{an:remote.detail:kernel.spatial.average}. For the second term, use $|\mathscr{D}_{k,T}(d_{0})| \leq |I_{\theta}|^{k}$ and $\ell_{T}^{-k}|I_{\theta}|^{k} = N_{T}^{k}$. It follows that
\begin{equation}\label{an:remote.detail:integrated.bin}
\begin{aligned}
&\ell_{T}^{-k}\int_{\mathscr{D}_{k,T}(d_{0})}\left|\PP_{\star}(\bb{F}^{\star} - u_{T}\bb{1}_{k} \in \mathscr{Q}) - N_{T}^{-k}\mathcal{K}_{T}(\bb{h})\Theta^{\otimes k}(\mathscr{Q})\right|\rd{\bb{h}} \\
&\qquad \leq \frac{\e_{T}}{|I_{\theta}|^{k}}\int_{\mathscr{D}_{k,T}(d_{0})}\mathcal{K}_{T}(\bb{h})\rd{\bb{h}} + N_{T}^{k}\exp\{ - ct^{1 + 2\e_{0}} + C_{k,L}t\} \\
&\qquad \leq \e_{T}\{1 + \oo(1)\} + \exp\{ - ct^{1 + 2\e_{0}} + C_{k,L}t + k\log{N_{T}}\} = \oo(1).
\end{aligned}
\end{equation}
The last equality uses $\log{N_{T}} = (1 + \theta)t + \OO(\log{t})$ from \eqref{eq:effective.cutoff.bounds} and \eqref{eq:parameter.scales}. Since $\e_{0} > 0$, the negative term $ - ct^{1 + 2\e_{0}}$ dominates all terms of order $t$ in this exponent.

For the finite sum defining $g$ in \eqref{eq:def.g.finite.sum.ind}, let $\mathscr{Q}_{r} \leqdef \prod_{j = 1}^{k}[a_{r,j},b_{r,j})$. Linearity of expectation and of integration against $\Theta^{\otimes k}$, followed by the triangle inequality, gives
\[
\begin{aligned}
&\ell_{T}^{-k}\int_{\mathscr{D}_{k,T}(d_{0})}\left|\EE_{\star}[g(\bb{F}^{\star} - u_{T}\bb{1}_{k})] - N_{T}^{-k}\mathcal{K}_{T}(\bb{h})\Theta^{\otimes k}(g)\right|\rd{\bb{h}} \\
&\qquad \leq \sum_{r = 1}^{m}|c_{r}|\ell_{T}^{-k}\int_{\mathscr{D}_{k,T}(d_{0})}\left|\PP_{\star}(\bb{F}^{\star} - u_{T}\bb{1}_{k} \in \mathscr{Q}_{r}) - N_{T}^{-k}\mathcal{K}_{T}(\bb{h})\Theta^{\otimes k}(\mathscr{Q}_{r})\right|\rd{\bb{h}} \\
&\qquad \leq m_{0}M\,\oo(1) = \oo(1).
\end{aligned}
\]
The error in the last line is the uniform single-bin error in \eqref{an:remote.detail:integrated.bin}, and $\sum_{r = 1}^{m}|c_{r}| \leq m_{0}M$. This proves \eqref{an:remote.local.approximation} for each $\star \in \{\mathrm{A},\mathrm{H}\}$, uniformly over the functions $g$ in \eqref{eq:def.g.finite.sum.ind}, including when the coefficients and bin endpoints depend on $T$ within the stated bounds.

To deduce \eqref{an:remote.bound}, first suppose that $f$ is a nonnegative compactly supported step function with finitely many half-open interval bins. Write $f(x) = \sum_{r = 1}^{n}\alpha_{r}\ind_{[a_{r},b_{r})}(x)$ with $a_{r} < b_{r}$ and coefficients $\alpha_{r} \geq 0$. The product $g(\bb{x}) \leqdef \prod_{j = 1}^{k}f(x_{j})$ expands as
\[
g(\bb{x}) = \sum_{r_{1},\ldots,r_{k} = 1}^{n}\left(\prod_{j = 1}^{k}\alpha_{r_{j}}\right)\prod_{j = 1}^{k}\ind_{[a_{r_{j}},b_{r_{j}})}(x_{j}).
\]
Thus $g$ has the form \eqref{eq:def.g.finite.sum.ind}, with $m_{0},L,M,\delta$ chosen independently of $T$ since $f$ is fixed. Tonelli's theorem also gives $\Theta^{\otimes k}(g) = \Theta(f)^{k}$.

By the definition of $\mu_{\star,f}$ in the statement, Tonelli's theorem, and \eqref{an:remote.local.approximation},
\begin{equation}\label{an:remote.detail:step.function}
\begin{aligned}
\EE_{\star}\left[\int_{\mathscr{D}_{k,T}(d_{0})}\prod_{j = 1}^{k}\mu_{\star,f}(\rd{h_{j}})\right]
&= \ell_{T}^{-k}\int_{\mathscr{D}_{k,T}(d_{0})}\EE_{\star}\left[\prod_{j = 1}^{k}f(F_{h_{j}}^{\star} - u_{T})\right]\rd{\bb{h}} \\
&= \Theta(f)^{k}\frac{1}{|I_{\theta}|^{k}}\int_{\mathscr{D}_{k,T}(d_{0})}\mathcal{K}_{T}(\bb{h})\rd{\bb{h}} + \oo(1)
\leq \Theta(f)^{k} + \oo(1).
\end{aligned}
\end{equation}
The second equality uses $\ell_{T}^{-k}N_{T}^{-k} = |I_{\theta}|^{-k}$, and the last inequality is \eqref{an:remote.detail:kernel.spatial.average}.

For a general $f$ as in the statement of the lemma, the bound follows by approximating $f$ from above by nonnegative compactly supported step functions whose integrals against $\Theta$ converge to $\Theta(f)$. Choose $R > 0$ such that its support is contained in $( - R,R)$. Since $f$ is bounded and its discontinuity set has Lebesgue measure zero, the Lebesgue criterion implies that $f$ is Riemann integrable on $[ - R,R]$. The density $\beta_{c}(\theta)e^{-\beta_{c}(\theta)x}$ of $\Theta$ is bounded there by $C_{R} \leqdef \beta_{c}(\theta)e^{\beta_{c}(\theta)R}$. Given $\e > 0$, choose a partition $ - R = x_{0} < x_{1} < \cdots < x_{n} = R$ whose upper Riemann sum exceeds $\int_{ - R}^{R}f(x)\rd{x}$ by at most $\e/C_{R}$, and define
\[
f_{\e}(x) \leqdef \sum_{r = 1}^{n}\left(\sup_{y \in [x_{r - 1},x_{r}]}f(y)\right)\ind_{[x_{r - 1},x_{r})}(x).
\]
This is a nonnegative compactly supported step function. The suprema are taken over closed subintervals, so $f_{\e} \geq f$ also at the interior partition points; at $R$ and outside $[ - R,R]$, $f$ vanishes. Therefore
\begin{equation}\label{an:remote.detail:upper.step}
0 \leq \Theta(f_{\e}) - \Theta(f)
= \int_{ - R}^{R}(f_{\e}(x) - f(x))\beta_{c}(\theta)e^{-\beta_{c}(\theta)x}\rd{x}
\leq C_{R}\int_{ - R}^{R}(f_{\e}(x) - f(x))\rd{x} \leq \e.
\end{equation}
For fixed $\e$, the function $f_{\e}$ and its partition are independent of $T$, so \eqref{an:remote.detail:step.function} applies to $f_{\e}$. The pointwise inequality $0 \leq f \leq f_{\e}$ implies $\prod_{j = 1}^{k}f(F_{h_{j}}^{\star} - u_{T}) \leq \prod_{j = 1}^{k}f_{\e}(F_{h_{j}}^{\star} - u_{T})$. Using this inequality in the definition of $\mu_{\star,f}$, followed by \eqref{an:remote.detail:step.function} and \eqref{an:remote.detail:upper.step}, gives
\[
\begin{aligned}
\limsup_{T \to \infty}\EE_{\star}\left[\int_{\mathscr{D}_{k,T}(d_{0})}\prod_{j = 1}^{k}\mu_{\star,f}(\rd{h_{j}})\right]
&\leq \limsup_{T \to \infty}\EE_{\star}\left[\int_{\mathscr{D}_{k,T}(d_{0})}\prod_{j = 1}^{k}\mu_{\star,f_{\e}}(\rd{h_{j}})\right] \\
&\leq \Theta(f_{\e})^{k} \leq (\Theta(f) + \e)^{k}.
\end{aligned}
\]
Letting $\e \downarrow 0$ proves \eqref{an:remote.bound}.
\end{proof}

\begin{lemma}[Capped exponential moments at high temperature]\label{an:capped.subcritical}
Let $\theta > 0$. Fix $d_{0} > 0$, $0 < \beta < \beta_{c}(\theta)$, and $K \in \R$, independently of $T$. For each $\star \in \{\mathrm{A},\mathrm{H}\}$,
\begin{equation}\label{an:remote.tilt.mean}
\EE_{\star}\left[\int_{I_{\theta}}e^{\beta F_{h}^{\star}}\ind_{\{F_{h}^{\star}\leq u_{T} + K\}}\rd{h}\right] = |I_{\theta}|\EE_{\mathrm{H}}[e^{\beta H(0)}]\{1 + \oo(1)\}, \qquad \text{as } T\to \infty,
\end{equation}
and
\begin{equation}\label{an:remote.tilt.pairs}
\limsup_{T \to \infty}\frac{1}{|I_{\theta}|^{2}\EE_{\mathrm{H}}[e^{\beta H(0)}]^{2}}\EE_{\star}\left[\int_{\mathscr{D}_{2,T}(d_{0})}e^{\beta(F_{h}^{\star} + F_{h'}^{\star})}\ind_{\{F_{h}^{\star}\leq u_{T} + K,\,F_{h'}^{\star}\leq u_{T} + K\}}\rd{h}\rd{h'}\right] \leq 1.
\end{equation}
The error term may depend on $d_{0}$, $\beta$, $K$, and the fixed parameters of the paper.
\end{lemma}

\begin{proof}[Proof of Lemma~\ref{an:capped.subcritical}]
Fix $0 < \e_{0} < 1/4$ and let $M_{T}(\lambda) \leqdef \EE_{\mathrm{H}}[e^{\lambda H(0)}]$. Let $\smash{\mathcal{A}_{T,\mathrm{sp}}^{\mathrm{A}}}$ be the common spatial event of Lemma~\ref{an:spatial.transform} for the full fields, with family size $m_{T} = 1$. The estimates \eqref{an:spatial.exceptional} and \eqref{an:spatial.transform.bound} in Lemma~\ref{an:spatial.transform} then apply to both single evaluations and pairs of evaluations with the same spatial event $\smash{\mathcal{A}_{T,\mathrm{sp}}^{\mathrm{A}}}$. The Haar estimates \eqref{an:remote.detail:tilted.tail} and \eqref{an:remote.detail:haar.pair.moment} hold on the whole probability space, without inserting $\smash{\ind_{\mathcal{A}_{T,\mathrm{sp}}^{\mathrm{H}}}}$ or otherwise restricting to the spatial truncation event $\mathcal{A}_{T,\mathrm{sp}}^{\mathrm{H}}$.

To prove \eqref{an:remote.tilt.mean}, let $b \leqdef \beta_{c}(\theta) - \beta > 0$, so $\beta + b = \beta_{c}(\theta)$. The cumulant estimate \eqref{an:joint.cumulants} gives
\[
\log\left(\frac{M_{T}(\beta + b)}{M_{T}(\beta)}\right) = \left(\beta b + \frac{b^{2}}{2}\right)V_{T} + \OO_{\beta,b}(1).
\]
Exponential Markov under the normalized Haar tilt consequently yields
\begin{equation}\label{an:remote.detail:tilted.tail}
\begin{aligned}
\frac{\EE_{\mathrm{H}}[e^{\beta H(h)}\ind_{\{H(h)>u_{T} + K\}}]}{M_{T}(\beta)}
&\leq e^{-b(u_{T} + K)}\frac{M_{T}(\beta + b)}{M_{T}(\beta)} \\
&= \exp\left\{\left[ - b\left(m(\theta) - \frac{\beta}{2}\right) + \frac{b^{2}}{4} + \oo(1)\right]s_{T}\right\} \\
&= \exp\left\{ - \frac{b^{2}}{4}s_{T} + \oo(s_{T})\right\} = \oo(1),
\end{aligned}
\end{equation}
uniformly in $h$. Here $u_{T}/s_{T} \to m(\theta)$ by \eqref{an:height.scales} in Lemma~\ref{an:height}, $V_{T}/s_{T} \to 1/2$ by \eqref{an:covariance.comparison} in Lemma~\ref{an:covariance}, $K$ is fixed, and $m(\theta) - \beta/2 = b/2$. Subtracting the tail in \eqref{an:remote.detail:tilted.tail} from the full Haar moment $\EE_{\mathrm{H}}[e^{\beta H(h)}] = M_{T}(\beta)$ gives the required capped mean under the Haar law.

Under the arithmetic law, the common event $\smash{\mathcal{A}_{T,\mathrm{sp}}^{\mathrm{A}}}$ and \eqref{an:spatial.transform.bound} in Lemma~\ref{an:spatial.transform} give, for $\lambda \in \{\beta,\beta + b\}$,
\begin{equation}\label{an:remote.detail:capped.transform}
\EE\left[e^{\lambda F_{h}^{\mathrm{A}}}\ind_{\mathcal{A}_{T,\mathrm{sp}}^{\mathrm{A}}}\right] = M_{T}(\lambda) + \OO\left(\exp\{ - ct^{1 + 2\e_{0}} + Ct\}\right),
\end{equation}
uniformly in $h$. For the contribution of $\{F_{h}^{\mathrm{A}} > u_{T} + K\}$ on the common event, \eqref{an:remote.detail:capped.transform} at $\lambda = \beta + b$ gives
\begin{equation}\label{an:remote.detail:arithmetic.tail}
\begin{aligned}
\EE\left[e^{\beta F_{h}^{\mathrm{A}}}\ind_{\{F_{h}^{\mathrm{A}}>u_{T} + K\}}\ind_{\mathcal{A}_{T,\mathrm{sp}}^{\mathrm{A}}}\right]
&\leq e^{-b(u_{T} + K)}\EE\left[e^{(\beta + b)F_{h}^{\mathrm{A}}}\ind_{\mathcal{A}_{T,\mathrm{sp}}^{\mathrm{A}}}\right] \\
&= e^{-b(u_{T} + K)}M_{T}(\beta + b) + \OO\left(\exp\{ - ct^{1 + 2\e_{0}} + C_{\beta,K}t\}\right) \\
&= \oo(M_{T}(\beta)).
\end{aligned}
\end{equation}
The first term is $\oo(M_{T}(\beta))$ by \eqref{an:remote.detail:tilted.tail}. The $\OO$ term is also $\oo(M_{T}(\beta))$, since $\e_{0} > 0$ and $M_{T}(\beta) \geq 1$. By \eqref{an:spatial.exceptional} in Lemma~\ref{an:spatial.transform}, the capped contribution outside the event is at most
\begin{equation}\label{an:remote.detail:capped.complement}
e^{\beta(u_{T} + K)}\PP((\mathcal{A}_{T,\mathrm{sp}}^{\mathrm{A}})^{c}) \leq \exp\{ - ct^{1 + 2\e_{0}} + C_{\beta,K}t\} = \oo(M_{T}(\beta)).
\end{equation}
For the arithmetic law, the capped expectation decomposes as
\begin{equation}\label{an:remote.detail:capped.decomposition}
\begin{aligned}
\EE[e^{\beta F_{h}^{\mathrm{A}}}\ind_{\{F_{h}^{\mathrm{A}}\leq u_{T} + K\}}]
&= \EE[e^{\beta F_{h}^{\mathrm{A}}}\ind_{\mathcal{A}_{T,\mathrm{sp}}^{\mathrm{A}}}]
- \EE[e^{\beta F_{h}^{\mathrm{A}}}\ind_{\{F_{h}^{\mathrm{A}}>u_{T} + K\}}\ind_{\mathcal{A}_{T,\mathrm{sp}}^{\mathrm{A}}}] \\
&\qquad + \EE[e^{\beta F_{h}^{\mathrm{A}}}\ind_{\{F_{h}^{\mathrm{A}}\leq u_{T} + K\}}\ind_{(\mathcal{A}_{T,\mathrm{sp}}^{\mathrm{A}})^{c}}].
\end{aligned}
\end{equation}
In \eqref{an:remote.detail:capped.decomposition}, the first term is $M_{T}(\beta) + \oo(M_{T}(\beta))$ by \eqref{an:remote.detail:capped.transform} at $\lambda = \beta$, the second is $\oo(M_{T}(\beta))$ by \eqref{an:remote.detail:arithmetic.tail}, and the third is $\oo(M_{T}(\beta))$ by \eqref{an:remote.detail:capped.complement}. Together with the Haar conclusion obtained from \eqref{an:remote.detail:tilted.tail}, this proves, for each $\star \in \{\mathrm{A},\mathrm{H}\}$,
\begin{equation}\label{an:remote.detail:capped.mean}
\sup_{h \in I_{\theta}}\left|\frac{\EE_{\star}[e^{\beta F_{h}^{\star}}\ind_{\{F_{h}^{\star}\leq u_{T} + K\}}]}{M_{T}(\beta)} - 1\right| \to 0.
\end{equation}
Integrating \eqref{an:remote.detail:capped.mean} over $I_{\theta}$ proves \eqref{an:remote.tilt.mean}.

To prove \eqref{an:remote.tilt.pairs}, apply \eqref{an:remote.detail:normalized.kernel} from the proof of Lemma~\ref{an:remote} with $k = 2$ and $\lambda = \beta$. This tilt lies in the interval $|\lambda| \leq M_{0}$ used there, since $\beta < \beta_{c}(\theta) < M_{0}$. Positivity permits removal of the height caps under the Haar law and gives
\begin{equation}\label{an:remote.detail:haar.pair.moment}
\frac{1}{|I_{\theta}|^{2}\EE_{\mathrm{H}}[e^{\beta H(0)}]^{2}}\int_{\mathscr{D}_{2,T}(d_{0})}\EE_{\mathrm{H}}[e^{\beta(H(h) + H(h'))}]\rd{h}\rd{h'} \leq 1 + \oo(1).
\end{equation}
For the arithmetic law, split the expectation according to $\smash{\mathcal{A}_{T,\mathrm{sp}}^{\mathrm{A}}}$. Remove the caps on this event and bound the capped integrand by $e^{2\beta(u_{T} + K)}$ on its complement. This gives, uniformly for $(h,h') \in \mathscr{D}_{2,T}(d_{0})$,
\[
\begin{aligned}
&\EE\left[e^{\beta(F_{h}^{\mathrm{A}} + F_{h'}^{\mathrm{A}})}\ind_{\{F_{h}^{\mathrm{A}}\leq u_{T} + K,\,F_{h'}^{\mathrm{A}}\leq u_{T} + K\}}\right] \\
&\qquad \leq \EE\left[e^{\beta(F_{h}^{\mathrm{A}} + F_{h'}^{\mathrm{A}})}\ind_{\mathcal{A}_{T,\mathrm{sp}}^{\mathrm{A}}}\right] + e^{2\beta(u_{T} + K)}\PP((\mathcal{A}_{T,\mathrm{sp}}^{\mathrm{A}})^{c}) \\
&\qquad \leq \EE_{\mathrm{H}}[e^{\beta(H(h) + H(h'))}] + (1 + e^{2\beta(u_{T} + K)})\exp\{ - ct^{1 + 2\e_{0}} + Ct\}.
\end{aligned}
\]
The second inequality uses \eqref{an:spatial.transform.bound} in Lemma~\ref{an:spatial.transform} for the two coordinates with real parameters $(\beta,\beta)$, together with the exceptional-event bound \eqref{an:spatial.exceptional} in the same lemma. Integrating over $\mathscr{D}_{2,T}(d_{0})$, using Tonelli's theorem and $|\mathscr{D}_{2,T}(d_{0})| \leq |I_{\theta}|^{2}$, and dividing by $|I_{\theta}|^{2}M_{T}(\beta)^{2}$ gives
\[
\begin{aligned}
&\frac{1}{|I_{\theta}|^{2}\EE_{\mathrm{H}}[e^{\beta H(0)}]^{2}}\EE\left[\int_{\mathscr{D}_{2,T}(d_{0})}e^{\beta(F_{h}^{\mathrm{A}} + F_{h'}^{\mathrm{A}})}\ind_{\{F_{h}^{\mathrm{A}}\leq u_{T} + K,\,F_{h'}^{\mathrm{A}}\leq u_{T} + K\}}\rd{h}\rd{h'}\right] \\
&\qquad \leq \frac{1}{|I_{\theta}|^{2}M_{T}(\beta)^{2}}\int_{\mathscr{D}_{2,T}(d_{0})}\EE_{\mathrm{H}}[e^{\beta(H(h) + H(h'))}]\rd{h}\rd{h'}
+ \frac{1 + e^{2\beta(u_{T} + K)}}{M_{T}(\beta)^{2}}\exp\{ - ct^{1 + 2\e_{0}} + Ct\} \\
&\qquad \leq 1 + \oo(1).
\end{aligned}
\]
The last error tends to zero because $u_{T} = \OO(t)$ and $M_{T}(\beta) \geq 1$. Taking the limsup proves \eqref{an:remote.tilt.pairs}.
\end{proof}

\subsubsection{Conditional fractional moments on annuli}

For $\beta > 0$, define
\begin{equation}\label{an:annular.parameters}
q \leqdef 1 \wedge \frac{2}{\beta}, \qquad
v(\beta) \leqdef \frac{1}{q} + \frac{q\beta^{2}}{4}
=
\begin{cases}
1 + \beta^{2}/4, & \beta \leq 2, \\
\beta, & \beta \geq 2.
\end{cases}
\end{equation}
The conditional arithmetic estimates obtained from \eqref{an:conditional.bound} in Lemma~\ref{an:conditional} use a spatial event $\smash{\mathcal{A}_{T,\mathrm{sp}}^{\mathrm{A}}}$ supplied by Lemma~\ref{an:spatial.transform} for the frequency cutoffs involved. The constants and the exceptional-probability bound \eqref{an:spatial.exceptional} in Lemma~\ref{an:spatial.transform} are uniform in these cutoffs, and $\smash{\mathcal{A}_{T,\mathrm{sp}}^{\mathrm{A}}}$ may be chosen simultaneously for any collection of at most $t^{B_{0}}$ frequency cutoffs for fixed $B_{0} > 0$. The notation $\mathcal{A}_{T,\mathrm{sp}}$ in a conditional expectation denotes the arithmetic event $\smash{\mathcal{A}_{T,\mathrm{sp}}^{\mathrm{A}}}$ when $\star = \mathrm{A}$ and the whole probability space when $\star = \mathrm{H}$, so that $\ind_{\mathcal{A}_{T,\mathrm{sp}}} = 1$ in the Haar case. This convention uses the untruncated Haar version of \eqref{an:conditional.bound} in Lemma~\ref{an:conditional}.

For the unit-interval estimate \eqref{an:unit.fractional} in Lemma~\ref{an:annular}, the common arithmetic event can be specified explicitly. Let $F_{d}^{\mathrm{A},\mathrm{lo}}(h)$ denote the part of $F_{h}^{\mathrm{A}}$ containing the prime powers with $\log(p^{j}) \leq 1/d$. For a fixed $0 < \e_{0} < 1/4$, define
\begin{equation}\label{an:annular.unit.spatial.event}
\mathcal{A}_{T,\mathrm{sp}}^{\mathrm{A}}
\leqdef \left\{\sup_{h \in I_{\theta}}\max\left\{|F_{h}^{\mathrm{A}}|,\ell_{T}|(F_{h}^{\mathrm{A}})'|\right\} \leq t^{1 + \e_{0}}\right\} \bigcap \bigcap_{0 \leq i \leq \lfloor s_{T}/\log{2}\rfloor}\left\{\sup_{h \in I_{\theta}}|F_{\ell_{T}2^{i}}^{\mathrm{A},\mathrm{lo}}(h)| \leq t^{1 + \e_{0}}\right\}.
\end{equation}
The intersection ranges over the indices $i \in \N_0$ corresponding to the radii $d = \ell_{T}2^{i} \leq 1$ in \eqref{an:annular.detail:unit.decomposition}. These give $\OO(t)$ low-frequency fields. Their coefficient bounds, together with those for the full field and its scaled derivative, are verified in \eqref{an:annular.detail:coefficient.admissibility} in the proof of Lemma~\ref{an:annular}. Thus Lemma~\ref{an:spatial.transform} gives \eqref{an:spatial.exceptional} for the event \eqref{an:annular.unit.spatial.event}.

\begin{lemma}[Annular fractional moments]\label{an:annular}
Let $\theta > 0$. Fix $\star \in \{\mathrm{A},\mathrm{H}\}$, $\beta > 0$ and $0 < a < q$. Let $h_{0} \in I_{\theta}$, let $\gamma$ range over a fixed compact subset of $\R$, and condition on $F_{h_{0}}^{\star} \in [\gamma s_{T},\gamma s_{T} + 1]$. For $d \leqdef \ell_{T}e^{r}$ with $0 \leq r \leq s_{T}$, define
\[
M_{\mathrm{ann}}^{\star,\mathrm{rel}} \leqdef \ell_{T}^{-1}\int_{\{d\leq|h - h_{0}|\leq 2d\}\cap I_{\theta}}\exp(\beta(F_{h}^{\star} - F_{h_{0}}^{\star}))\rd{h}.
\]
Uniformly in these parameters,
\begin{equation}\label{an:annular.bound}
\EE_{\star}[\{M_{\mathrm{ann}}^{\star,\mathrm{rel}}\}^{a}\ind_{\mathcal{A}_{T,\mathrm{sp}}} \mid F_{h_{0}}^{\star} \in [\gamma s_{T},\gamma s_{T} + 1]] \leq C\exp\left\{\left[a\{v(\beta) - \beta\gamma\} + \frac{a^{2}\beta^{2}}{4}\right]r\right\}.
\end{equation}
For the integral of the same normalized exponential over $[h_{0} - 1,h_{0} + 1] \cap I_{\theta}$,
\begin{equation}\label{an:unit.fractional}
\begin{aligned}
&\EE_{\star}\left[\left\{\ell_{T}^{-1}\int_{[h_{0} - 1,h_{0} + 1]\cap I_{\theta}}\exp\{\beta(F_{h}^{\star} - F_{h_{0}}^{\star})\}\rd{h}\right\}^{a}\ind_{\mathcal{A}_{T,\mathrm{sp}}} \mid F_{h_{0}}^{\star} \in [\gamma s_{T},\gamma s_{T} + 1]\right] \\
&\qquad \leq C(1 + s_{T})\exp\left\{\left[a\{v(\beta) - \beta\gamma\} + \frac{a^{2}\beta^{2}}{4}\right]_{+}s_{T}\right\}.
\end{aligned}
\end{equation}
If $a\{v(\beta) - \beta\gamma\} + a^{2}\beta^{2}/4$ is bounded above by a negative constant, the conditional moment on the left-hand side of \eqref{an:unit.fractional} is bounded uniformly in $T$.
\end{lemma}

\begin{proof}[Proof of Lemma~\ref{an:annular}]
Call the conditioning point $h_{0}$ the root. Split the field into a low-frequency part $F_{d}^{\star,\mathrm{lo}}$ and a high-frequency part $F_{d}^{\star,\mathrm{hi}}$ at the cutoff $\log(p^{j}) = 1/d$, assigning a prime power to the low part when $\log(p^{j}) \leq 1/d$. Thus $F_{h}^{\star} = F_{d}^{\star,\mathrm{lo}}(h) + F_{d}^{\star,\mathrm{hi}}(h)$. Define
\[
\widetilde{F}_{d}^{\star,\mathrm{hi}}(h) \leqdef F_{h}^{\star} - F_{d}^{\star,\mathrm{lo}}(h_{0}) = F_{d}^{\star,\mathrm{hi}}(h) + F_{d}^{\star,\mathrm{lo}}(h) - F_{d}^{\star,\mathrm{lo}}(h_{0}).
\]
In particular, the root-relative increment has the exact decomposition
\begin{equation}\label{an:annular.detail:increment}
F_{h}^{\star} - F_{h_{0}}^{\star} = \widetilde{F}_{d}^{\star,\mathrm{hi}}(h) - F_{d}^{\star,\mathrm{hi}}(h_{0}).
\end{equation}
The low and high parts are orthogonal under the Haar law: different powers of a prime have zero covariance, even if those powers fall on opposite sides of the cutoff. Applying the Haar covariance formula \eqref{an:detail:cov.haar.sum} to the two frequency ranges gives the following estimates, whose prime sums are evaluated below. Whenever $d \leq |h - h_{0}| \leq 2d$,
\begin{equation}\label{an:annular.covariances}
\begin{aligned}
&\Var_{\mathrm{H}}(F_{d}^{\mathrm{H},\mathrm{hi}}(h_{0})) = \frac{r}{2} + \OO(1), \qquad
\Var_{\mathrm{H}}(\widetilde{F}_{d}^{\mathrm{H},\mathrm{hi}}(h)) = \frac{r}{2} + \OO(1), \\
&\Cov_{\mathrm{H}}(H(h_{0}),F_{d}^{\mathrm{H},\mathrm{hi}}(h_{0})) = \frac{r}{2} + \OO(1), \\
&|\Cov_{\mathrm{H}}(F_{d}^{\mathrm{H},\mathrm{hi}}(h_{0}),\widetilde{F}_{d}^{\mathrm{H},\mathrm{hi}}(h))| + |\Cov_{\mathrm{H}}(H(h_{0}),\widetilde{F}_{d}^{\mathrm{H},\mathrm{hi}}(h))| \leq C.
\end{aligned}
\end{equation}
The uniformity in $r$ follows from the following explicit estimates. Since $d = \ell_{T}e^{r} = 2\eta_{\sig}e^{r}$ and $0 \leq r \leq s_{T}$, the cutoff obeys $1 \leq 1/d \leq (2\eta_{\sig})^{-1}$. Partial summation using \eqref{analytic:eq:PNT}, followed by the diagonal calculation in \eqref{an:effective.cutoff.estimates}, yields
\begin{equation}\label{an:annular.detail:high.variance}
\begin{aligned}
\Var_{\mathrm{H}}(F_{d}^{\mathrm{H},\mathrm{hi}}(h_{0}))
&= \frac{1}{2}\int_{1/d}^{\log{X}}e^{-2\eta_{\sig}v}\frac{\rd{v}}{v} + \OO(1) \\
&= \frac{1}{2}\int_{1/d}^{(2\eta_{\sig})^{-1}}\frac{\rd{v}}{v} + \OO(1)
= \frac{1}{2}\log\left(\frac{d}{2\eta_{\sig}}\right) + \OO(1)
= \frac{r}{2} + \OO(1).
\end{aligned}
\end{equation}
The covariance with $H(h_{0})$ equals this high-frequency variance by orthogonality.

For the covariance between $F_{d}^{\mathrm{H},\mathrm{hi}}(h_{0})$ and $F_{d}^{\mathrm{H},\mathrm{hi}}(h)$, let $\xi \leqdef (h - h_{0})/d$, so that $1 \leq |\xi| \leq 2$. Its prime contribution is, up to a bounded summation error,
\begin{equation}\label{an:annular.detail:high.covariance.integral}
\frac{1}{2}\int_{1/d}^{\log{X}}e^{-2\eta_{\sig}v}\cos((h - h_{0})v)\frac{\rd{v}}{v} = \frac{1}{2}\int_{1}^{d\log{X}}e^{-(2\eta_{\sig}/d)w}\cos(\xi w)\frac{\rd{w}}{w}.
\end{equation}
The amplitude $b(w) \leqdef e^{-(2\eta_{\sig}/d)w}/w$ is positive and decreasing. With $d\log{X} > 1$, integration by parts makes the bound on \eqref{an:annular.detail:high.covariance.integral} explicit:
\begin{equation}\label{an:annular.detail:high.covariance.integral.next}
\left|\int_{1}^{d\log{X}}b(w)\cos(\xi w)\rd{w}\right| \leq \frac{b(1) + b(d\log{X}) + \int_{1}^{d\log{X}}|b'(w)|\rd{w}}{|\xi|} = \frac{2b(1)}{|\xi|} \leq 2.
\end{equation}
For the low-frequency increment $F_{d}^{\mathrm{H},\mathrm{lo}}(h) - F_{d}^{\mathrm{H},\mathrm{lo}}(h_{0})$, the inequalities $|h - h_{0}| \leq 2d$ and $|e^{\ii x} - 1| \leq |x|$ give
\begin{equation}\label{an:annular.detail:low.increment}
\begin{aligned}
&\Var_{\mathrm{H}}(F_{d}^{\mathrm{H},\mathrm{lo}}(h) - F_{d}^{\mathrm{H},\mathrm{lo}}(h_{0}))
\leq Cd^{2}\sum_{p \leq e^{1/d}}\frac{(\log p)^{2}}{p} + C \leq Cd^{2}\int_{0}^{1/d}v\rd{v} + C \leq C, \\
&|\Cov_{\mathrm{H}}(F_{d}^{\mathrm{H},\mathrm{lo}}(h_{0}),F_{d}^{\mathrm{H},\mathrm{lo}}(h) - F_{d}^{\mathrm{H},\mathrm{lo}}(h_{0}))|
\leq Cd^{2}\int_{0}^{1/d}v\rd{v} + C \leq C.
\end{aligned}
\end{equation}
The second estimate uses $1 - \cos(x) \leq x^{2}/2$. All terms with prime-power exponent at least two contribute bounded quantities, since their variance coefficients are absolutely summable. Partial summation using \eqref{analytic:eq:PNT} gives bounded errors because $|h - h_{0}| \leq 2d \leq 2$. Combining \eqref{an:annular.detail:high.variance}--\eqref{an:annular.detail:low.increment} and orthogonality proves every assertion in \eqref{an:annular.covariances}.

To estimate the annular integral, the proof will bound the supremum of the field on each interval $B$ contained in the annulus and having length at most $\ell_{T}$. Passing from a field value at one point of $B$ to this supremum requires control of the spatial derivative. Define
\[
D_{d}^{\star}(z) \leqdef \ell_{T}(\widetilde{F}_{d}^{\star,\mathrm{hi}})'(z) = \ell_{T}(F_{z}^{\star})',
\]
where the primes denote differentiation with respect to $z$. The second equality holds because $F_{d}^{\star,\mathrm{lo}}(h_{0})$ does not depend on $z$. For $h,h' \in B$, the fundamental theorem of calculus gives
\[
|\widetilde{F}_{d}^{\star,\mathrm{hi}}(h) - \widetilde{F}_{d}^{\star,\mathrm{hi}}(h')|
\leq \int_{B}|(\widetilde{F}_{d}^{\star,\mathrm{hi}})'(z)|\rd{z}
= \ell_{T}^{-1}\int_{B}|D_{d}^{\star}(z)|\rd{z}.
\]
The required estimates concern the variance of $D_{d}^{\mathrm{H}}(z)$ and its covariances with the field values appearing in the conditional exponential moments below.

Differentiating a prime-power term multiplies its coefficient by $-\ii j\log p$. Taking absolute values in the Haar covariance formula \eqref{an:detail:cov.haar.sum} and applying partial summation gives, for $Z^{\mathrm{H}}$ equal to $H(h_{0})$, $F_{d}^{\mathrm{H},\mathrm{hi}}(h_{0})$, or $\widetilde{F}_{d}^{\mathrm{H},\mathrm{hi}}(h)$,
\begin{equation}\label{an:annular.detail:derivative.integrals}
\begin{aligned}
\Var_{\mathrm{H}}(d(F_{d}^{\mathrm{H},\mathrm{lo}})'(z))
&\leq Cd^{2}\sum_{p \leq e^{1/d}}\frac{(\log p)^{2}}{p} + C
\leq Cd^{2}\int_{0}^{1/d}v\rd{v} + C \leq C, \\
\Var_{\mathrm{H}}(\ell_{T}(F_{d}^{\mathrm{H},\mathrm{hi}})'(z))
&\leq C\ell_{T}^{2}\sum_{p \leq X}\frac{(\log p)^{2}}{p^{2\sig}} + C
\leq C\ell_{T}^{2}\int_{0}^{\log{X}}ve^{-2\eta_{\sig}v}\rd{v} + C \leq C, \\
|\Cov_{\mathrm{H}}(d(F_{d}^{\mathrm{H},\mathrm{lo}})'(z),Z^{\mathrm{H}})|
&\leq Cd\sum_{p \leq e^{1/d}}\frac{\log p}{p} + C
\leq Cd\int_{0}^{1/d}\rd{v} + C \leq C, \\
|\Cov_{\mathrm{H}}(\ell_{T}(F_{d}^{\mathrm{H},\mathrm{hi}})'(z),Z^{\mathrm{H}})|
&\leq C\ell_{T}\sum_{p \leq X}\frac{\log p}{p^{2\sig}} + C
\leq C\ell_{T}\int_{0}^{\log{X}}e^{-2\eta_{\sig}v}\rd{v} + C \leq C.
\end{aligned}
\end{equation}
The final inequalities of the second and fourth lines use $\ell_{T} = 2\eta_{\sig}$, and the additive constants absorb the summation errors and the absolutely summable contributions from higher prime powers. Since $\ell_{T} \leq d$, combining the two frequency ranges gives
\begin{equation}\label{an:annular.detail:derivative}
\Var_{\mathrm{H}}(D_{d}^{\mathrm{H}}(z)) \leq C, \qquad |\Cov_{\mathrm{H}}(D_{d}^{\mathrm{H}}(z),Z^{\mathrm{H}})| \leq C,
\end{equation}
uniformly for $z,h \in B$ and all such intervals $B$. In the remainder of the proof, these bounds will be combined with \eqref{an:conditional.bound} of Lemma~\ref{an:conditional} to control the additional derivative term when passing from a fixed point to the supremum over $B$.

For each $\star \in \{\mathrm{A},\mathrm{H}\}$, the field $\{F_{h}^{\star}\}_{h \in I_{\theta}}$, its frequency truncations $\{F_{d}^{\star,\mathrm{lo}}(h)\}_{h \in I_{\theta}}$ and $\smash{\{F_{d}^{\star,\mathrm{hi}}(h)\}_{h \in I_{\theta}}}$, and its first scaled derivative $\{\ell_{T}(F_{h}^{\star})'\}_{h \in I_{\theta}}$ satisfy the coefficient assumptions of Lemma~\ref{an:spatial.transform}. Indeed, the coefficient of the $j$th power of $p$ in $\ell_{T}^{m}(F^{\star})^{(m)}$, for $m \in \{0,1,2\}$, has absolute value
\begin{equation}\label{an:annular.detail:coefficient.admissibility}
\frac{(\ell_{T}j\log{p})^{m}}{jp^{j\sig}}
= \frac{2^{m}}{j}(\eta_{\sig}j\log{p})^{m}e^{-\eta_{\sig}j\log{p}}p^{-j/2}
\leq C_{m}p^{-j/2}.
\end{equation}
A frequency restriction cannot increase the bounds in \eqref{an:annular.detail:coefficient.admissibility}. The case $m = 2$ verifies the derivative assumption for the family member $\ell_{T}(F^{\star})'$. In the arithmetic case, choose the common event $\smash{\mathcal{A}_{T,\mathrm{sp}}^{\mathrm{A}}}$ from Lemma~\ref{an:spatial.transform} for the three families $F^{\star}$, $F_{d}^{\star,\mathrm{lo}}$, and $\ell_{T}(F^{\star})'$. Since $\smash{F_{d}^{\star,\mathrm{hi}} = F^{\star} - F_{d}^{\star,\mathrm{lo}}}$, each coordinate in the exponential estimates \eqref{an:annular.detail:weighted.transforms} is a fixed real linear combination of their evaluations. The arithmetic event $\smash{\mathcal{A}_{T,\mathrm{sp}}^{\mathrm{A}}}$ satisfies \eqref{an:spatial.exceptional}, and the transform estimate \eqref{an:spatial.transform.bound} in Lemma~\ref{an:spatial.transform} applies uniformly over the positions in the annulus and its cells $B$. As before, $\mathcal{A}_{T,\mathrm{sp}}$ denotes the arithmetic event $\smash{\mathcal{A}_{T,\mathrm{sp}}^{\mathrm{A}}}$ when $\star = \mathrm{A}$ and the whole probability space when $\star = \mathrm{H}$, so that $\ind_{\mathcal{A}_{T,\mathrm{sp}}} = 1$ in the Haar case, where the exponential estimates in \eqref{an:annular.detail:weighted.transforms} hold without truncation.

Let $E_{0} \leqdef \{F_{h_{0}}^{\star} \in \gamma s_{T} + [0,1]\}$, and define
\[
\Gamma_{a}(r,\gamma) \leqdef \exp\left\{ - a\beta\gamma r + \frac{a^{2}\beta^{2}r}{4}\right\}, \qquad M_{d}^{\star,\mathrm{inc}} \leqdef \ell_{T}^{-1}\int_{\{d\leq|h - h_{0}|\leq 2d\}\cap I_{\theta}}e^{\beta\widetilde{F}_{d}^{\star,\mathrm{hi}}(h)}\rd{h}.
\]
For single-root conditioning, the mean term $\bb{c}_{T}^{\mathsf{T}}\Sigma_{T}^{-1}\bb{v}$ in \eqref{an:conditional.bound} of Lemma~\ref{an:conditional} is the covariance with $H(h_{0})$ multiplied by $\gamma s_{T}/V_{T}$. Since $V_{T} = s_{T}/2 + \OO(1)$ by \eqref{an:covariance.comparison} in Lemma~\ref{an:covariance} and $r \leq s_{T}$, \eqref{an:annular.covariances} gives
\begin{equation}\label{an:annular.detail:conditional.means}
\frac{\gamma s_{T}}{V_{T}}\Cov_{\mathrm{H}}(H(h_{0}),F_{d}^{\mathrm{H},\mathrm{hi}}(h_{0})) = \gamma r + \OO(1), \qquad
\frac{\gamma s_{T}}{V_{T}}\Cov_{\mathrm{H}}(H(h_{0}),\widetilde{F}_{d}^{\mathrm{H},\mathrm{hi}}(h)) = \OO(1).
\end{equation}
Apply \eqref{an:conditional.bound} in Lemma~\ref{an:conditional} with $k = 1$, conditioning coordinate $Y_{T,1}^{\star} = F_{h_{0}}^{\star}$, and $\lambda_{2} = 1$, first taking the second coordinate to be, respectively,
\[
Y_{T,2}^{\star} = - a\beta F_{d}^{\star,\mathrm{hi}}(h_{0}), \qquad
Y_{T,2}^{\star} = - a\beta F_{d}^{\star,\mathrm{hi}}(h_{0}) + q\beta\widetilde{F}_{d}^{\star,\mathrm{hi}}(h).
\]
Their Haar variances are $a^{2}\beta^{2}r/2 + \OO(1)$ and $(a^{2} + q^{2})\beta^{2}r/2 + \OO(1)$, respectively, because the cross covariance is bounded in \eqref{an:annular.covariances}. By \eqref{an:annular.detail:derivative}, adding $\pm q\beta D_{d}^{\star}(z)$ to the second coordinate changes its Haar variance by $\OO(1)$ and its conditional mean term by $\OO(1)$, since $|\gamma|s_{T}/V_{T} \leq C$. Applying \eqref{an:conditional.bound} also to these two augmented coordinates, combining the variance bounds with \eqref{an:annular.detail:conditional.means}, and dropping the nonnegative variance correction $\bb{c}_{T}^{\mathsf{T}}\Sigma_{T}^{-1}\bb{c}_{T}$ yields the following estimates, with the third valid for both signs and uniformly for $h,z \in B$:
\begin{equation}\label{an:annular.detail:weighted.transforms}
\begin{aligned}
\EE_{\star}[e^{-a\beta F_{d}^{\star,\mathrm{hi}}(h_{0})}\ind_{\mathcal{A}_{T,\mathrm{sp}}} \mid E_{0}] &\leq C\Gamma_{a}(r,\gamma), \\
\EE_{\star}[e^{-a\beta F_{d}^{\star,\mathrm{hi}}(h_{0}) + q\beta\widetilde{F}_{d}^{\star,\mathrm{hi}}(h)}\ind_{\mathcal{A}_{T,\mathrm{sp}}} \mid E_{0}] &\leq C\Gamma_{a}(r,\gamma)e^{q^{2}\beta^{2}r/4}, \\
\EE_{\star}[e^{-a\beta F_{d}^{\star,\mathrm{hi}}(h_{0}) + q\beta\widetilde{F}_{d}^{\star,\mathrm{hi}}(h) \pm q\beta D_{d}^{\star}(z)}\ind_{\mathcal{A}_{T,\mathrm{sp}}} \mid E_{0}] &\leq C\Gamma_{a}(r,\gamma)e^{q^{2}\beta^{2}r/4}.
\end{aligned}
\end{equation}

The estimates in \eqref{an:annular.detail:weighted.transforms} concern field values at fixed positions. To control the supremum over a cell, let $B$ be an interval contained in the annulus, with positive length $L_{B} \leq \ell_{T}$ and left endpoint $h_{B}$. The fundamental theorem of calculus and Jensen's inequality give
\begin{equation}\label{an:annular.detail:cell.supremum}
\begin{aligned}
\exp\left(q\beta\sup_{h \in B}\widetilde{F}_{d}^{\star,\mathrm{hi}}(h)\right)
&\leq e^{q\beta\widetilde{F}_{d}^{\star,\mathrm{hi}}(h_{B})}\exp\left(q\beta\int_{B}|(\widetilde{F}_{d}^{\star,\mathrm{hi}})'(z)|\rd{z}\right) \\
&\leq \frac{1}{L_{B}}\int_{B}e^{q\beta\widetilde{F}_{d}^{\star,\mathrm{hi}}(h_{B}) + q\beta|D_{d}^{\star}(z)|}\rd{z}.
\end{aligned}
\end{equation}
The last inequality in \eqref{an:annular.detail:cell.supremum} uses $L_{B}/\ell_{T} \leq 1$. Multiplying \eqref{an:annular.detail:cell.supremum} by $\smash{e^{-a\beta F_{d}^{\star,\mathrm{hi}}(h_{0})}\ind_{\mathcal{A}_{T,\mathrm{sp}}}}$, taking conditional expectations given $E_{0}$, and using $e^{q\beta|x|} \leq e^{q\beta x} + e^{-q\beta x}$ together with the third estimate in \eqref{an:annular.detail:weighted.transforms} at $h = h_{B}$ give
\begin{equation}\label{an:annular.detail:cell}
\EE_{\star}\left[e^{-a\beta F_{d}^{\star,\mathrm{hi}}(h_{0})}\exp\left(q\beta\sup_{h \in B}\widetilde{F}_{d}^{\star,\mathrm{hi}}(h)\right)\ind_{\mathcal{A}_{T,\mathrm{sp}}} \mid E_{0}\right] \leq C\Gamma_{a}(r,\gamma)e^{q^{2}\beta^{2}r/4}.
\end{equation}

Partition the two components of the annulus, $[h_{0} - 2d,h_{0} - d] \cap I_{\theta}$ and $[h_{0} + d,h_{0} + 2d] \cap I_{\theta}$, into at most $C(d/\ell_{T} + 1) \leq Ce^{r}$ cells of length at most $\ell_{T}$. Since each cell $B$ has length $L_{B} \leq \ell_{T}$,
\begin{equation}\label{an:annular.detail:cell.mass}
\begin{aligned}
M_{d}^{\star,\mathrm{inc}}
&= \sum_{B}\ell_{T}^{-1}\int_{B}e^{\beta\widetilde{F}_{d}^{\star,\mathrm{hi}}(h)}\rd{h} \\
&\leq \sum_{B}\frac{L_{B}}{\ell_{T}}\exp\left(\beta\sup_{h \in B}\widetilde{F}_{d}^{\star,\mathrm{hi}}(h)\right)
\leq \sum_{B}\exp\left(\beta\sup_{h \in B}\widetilde{F}_{d}^{\star,\mathrm{hi}}(h)\right).
\end{aligned}
\end{equation}
Since $q \leq 1$, applying subadditivity of $x \mapsto x^{q}$ to \eqref{an:annular.detail:cell.mass}, multiplying by $e^{-a\beta F_{d}^{\star,\mathrm{hi}}(h_{0})}\ind_{\mathcal{A}_{T,\mathrm{sp}}}$, and using \eqref{an:annular.detail:cell} for each of the at most $Ce^{r}$ cells give
\begin{equation}\label{an:annular.detail:cell.sum}
\EE_{\star}[e^{-a\beta F_{d}^{\star,\mathrm{hi}}(h_{0})}(M_{d}^{\star,\mathrm{inc}})^{q}\ind_{\mathcal{A}_{T,\mathrm{sp}}} \mid E_{0}] \leq C\Gamma_{a}(r,\gamma)e^{(1 + q^{2}\beta^{2}/4)r}.
\end{equation}
By \eqref{an:annular.detail:increment}, $M_{\mathrm{ann}}^{\star,\mathrm{rel}} = e^{-\beta F_{d}^{\star,\mathrm{hi}}(h_{0})}M_{d}^{\star,\mathrm{inc}}$. Let $W \leqdef e^{-a\beta F_{d}^{\star,\mathrm{hi}}(h_{0})}\ind_{\mathcal{A}_{T,\mathrm{sp}}}$. H\"{o}lder's inequality with exponents $q/a$ and $q/(q - a)$, both greater than one, gives
\begin{equation}\label{an:annular.detail:weighted.holder}
\begin{aligned}
\EE_{\star}[\{M_{\mathrm{ann}}^{\star,\mathrm{rel}}\}^{a}\ind_{\mathcal{A}_{T,\mathrm{sp}}} \mid E_{0}]
&= \EE_{\star}[(W(M_{d}^{\star,\mathrm{inc}})^{q})^{a/q}W^{1 - a/q} \mid E_{0}] \\
&\leq \left(\EE_{\star}[W(M_{d}^{\star,\mathrm{inc}})^{q} \mid E_{0}]\right)^{a/q}\left(\EE_{\star}[W \mid E_{0}]\right)^{1 - a/q} \\
&\leq C\Gamma_{a}(r,\gamma)\exp\left\{\frac{a}{q}\left(1 + \frac{q^{2}\beta^{2}}{4}\right)r\right\} = C\Gamma_{a}(r,\gamma)e^{av(\beta)r}.
\end{aligned}
\end{equation}
The last line of \eqref{an:annular.detail:weighted.holder} uses \eqref{an:annular.detail:cell.sum} and the first estimate in \eqref{an:annular.detail:weighted.transforms}; their powers of $\Gamma_{a}(r,\gamma)$ add to one. Substituting the definition of $\Gamma_{a}(r,\gamma)$ proves \eqref{an:annular.bound}.

To prove \eqref{an:unit.fractional}, first bound the part of its integral over $[h_{0} - \ell_{T},h_{0} + \ell_{T}] \cap I_{\theta}$; the remaining region will be covered by dyadic annuli. For $h$ in this central interval, apply the Haar covariance formula \eqref{an:detail:cov.haar.sum} and $1 - \cos(x) \leq x^{2}/2$.
\begin{equation}\label{an:annular.detail:inner.covariances}
\begin{aligned}
&\Var_{\mathrm{H}}(H(h) - H(h_{0})) \leq C|h - h_{0}|^{2}\int_{0}^{\log{X}}ve^{-2\eta_{\sig}v}\rd{v} + C \leq C, \\
&|\Cov_{\mathrm{H}}(H(h_{0}),H(h) - H(h_{0}))| = \frac{1}{2}\Var_{\mathrm{H}}(H(h) - H(h_{0})) \leq C.
\end{aligned}
\end{equation}
The equality uses $\Var_{\mathrm{H}}(H(h)) = \Var_{\mathrm{H}}(H(h_{0})) = V_{T}$. Thus the conditional mean and variance terms in \eqref{an:conditional.bound} of Lemma~\ref{an:conditional} are bounded for the coordinate $a\beta(F_{h}^{\star} - F_{h_{0}}^{\star})$: the mean $\bb{c}_{T}^{\mathsf{T}}\Sigma_{T}^{-1}\bb{v} = \smash{\frac{a\beta\gamma s_{T}}{V_{T}}\Cov_{\mathrm{H}}(H(h_{0}),H(h) - H(h_{0}))}$ is bounded by \eqref{an:annular.detail:inner.covariances} and $|\gamma|s_{T}/V_{T} \leq C$, while the nonnegative variance correction $\bb{c}_{T}^{\mathsf{T}}\Sigma_{T}^{-1}\bb{c}_{T}$ can again be dropped. The derivative estimates \eqref{an:annular.detail:derivative.integrals} give the same conclusion after adding either sign of $a\beta\ell_{T}(F^{\star})'(z)$. Covering the innermost interval by at most two cells, applying \eqref{an:annular.detail:cell.supremum} with $q$ replaced by $a$ and $\widetilde{F}_{d}^{\star,\mathrm{hi}}(h)$ replaced by $F_{h}^{\star} - F_{h_{0}}^{\star}$, and using subadditivity of $x \mapsto x^{a}$ for $a < 1$ therefore give
\begin{equation}\label{an:annular.detail:inner.moment}
\EE_{\star}\left[\left\{\ell_{T}^{-1}\int_{[h_{0} - \ell_{T},h_{0} + \ell_{T}]\cap I_{\theta}}\exp\{\beta(F_{h}^{\star} - F_{h_{0}}^{\star})\}\rd{h}\right\}^{a}\ind_{\mathcal{A}_{T,\mathrm{sp}}} \mid E_{0}\right] \leq C.
\end{equation}

Write $M_{\mathrm{ann}}^{\star,\mathrm{rel}}(d)$ for the annular mass at radius $d$, making this dependence explicit. The interval $[h_{0} - 1,h_{0} + 1] \cap I_{\theta}$ is covered by the central interval $[h_{0} - \ell_{T},h_{0} + \ell_{T}] \cap I_{\theta}$ and the annuli with radii $d = \ell_{T}2^{j} \leq 1$, where $j \in \N_0$. Consequently,
\begin{equation}\label{an:annular.detail:unit.decomposition}
\begin{aligned}
&\ell_{T}^{-1}\int_{[h_{0} - 1,h_{0} + 1]\cap I_{\theta}}\exp\{\beta(F_{h}^{\star} - F_{h_{0}}^{\star})\}\rd{h} \\
&\qquad\leq \ell_{T}^{-1}\int_{[h_{0} - \ell_{T},h_{0} + \ell_{T}]\cap I_{\theta}}\exp\{\beta(F_{h}^{\star} - F_{h_{0}}^{\star})\}\rd{h}
+ \sum_{0 \leq j \leq s_{T}/\log{2}}M_{\mathrm{ann}}^{\star,\mathrm{rel}}(\ell_{T}2^{j}).
\end{aligned}
\end{equation}
Choose the arithmetic event $\smash{\mathcal{A}_{T,\mathrm{sp}}^{\mathrm{A}}}$ simultaneously for these $\OO(t)$ dyadic cutoffs. The estimates \eqref{an:spatial.exceptional} and \eqref{an:spatial.transform.bound} in Lemma~\ref{an:spatial.transform} remain uniform for this enlarged family. In \eqref{an:annular.detail:unit.sum}, the notation $\mathcal{A}_{T,\mathrm{sp}}$ denotes the enlarged arithmetic event $\smash{\mathcal{A}_{T,\mathrm{sp}}^{\mathrm{A}}}$ when $\star = \mathrm{A}$ and still denotes the whole probability space when $\star = \mathrm{H}$, with $\ind_{\mathcal{A}_{T,\mathrm{sp}}} = 1$ in the Haar case. Let $b_{a} \leqdef a\{v(\beta) - \beta\gamma\} + a^{2}\beta^{2}/4$. Applying subadditivity of $x \mapsto x^{a}$ for $a < 1$ to \eqref{an:annular.detail:unit.decomposition}, using \eqref{an:annular.detail:inner.moment} for the central integral and \eqref{an:annular.bound} with $r = j\log{2}$ for each annular mass, yields
\begin{equation}\label{an:annular.detail:unit.sum}
\begin{aligned}
\EE_{\star}\left[\left\{\ell_{T}^{-1}\int_{[h_{0} - 1,h_{0} + 1]\cap I_{\theta}}\exp\{\beta(F_{h}^{\star} - F_{h_{0}}^{\star})\}\rd{h}\right\}^{a}\ind_{\mathcal{A}_{T,\mathrm{sp}}} \mid E_{0}\right]
&\leq C\left(1 + \sum_{0 \leq j \leq s_{T}/\log{2}}e^{b_{a}j\log{2}}\right) \\
&\leq C(1 + s_{T})e^{[b_{a}]_{+}s_{T}}.
\end{aligned}
\end{equation}
This is \eqref{an:unit.fractional}. If $b_{a} \leq - c < 0$, the sum in \eqref{an:annular.detail:unit.sum} is bounded by $\sum_{j = 0}^{\infty}2^{-cj} = (1 - 2^{-c})^{-1}$, proving the final assertion.
\end{proof}

\begin{lemma}[Localization under joint height conditioning]\label{an:joint.annular}
Let $\theta > 0$. Fix $\star \in \{\mathrm{A},\mathrm{H}\}$, $\beta > \beta_{c}(\theta)$, $0 < \e < 1/4$, a compact interval $J \subseteq \R$, and a nonnegative bounded Borel function $\chi$ supported on $J$. For $j \in \Z$, let $B_{j}^{\e} \leqdef [j + \e,j + 1 - \e] \cap I_{\theta}$, omitting empty intervals. Let the representative positions $h_{1},\ldots,h_{k}$ belong to distinct intervals $B_{j}^{\e}$, and write $\bb{h} \leqdef (h_{1},\ldots,h_{k})$. Fix a bin width $b > 0$, independent of $T$, and choose lower endpoints $x_{1},\ldots,x_{k}$ such that $[x_{i},x_{i} + b] \subseteq J$ for every $1 \leq i \leq k$. Define the simultaneous conditioning event $E_{\bb{h}} \leqdef \bigcap_{i = 1}^{k}\{F_{h_{i}}^{\star} - u_{T} \in [x_{i},x_{i} + b]\}$. For $h_{i} \in B_{j_{i}}^{\e}$ and $R \geq 1$, define
\[
M_{\star,i}^{\mathrm{out}}(R) \leqdef \ell_{T}^{-1}\int_{B_{j_{i}}^{\e} \setminus [h_{i} - R\ell_{T},h_{i} + R\ell_{T}]} e^{\beta(F_{h}^{\star} - u_{T})}\chi(F_{h}^{\star} - u_{T})\rd{h}.
\]
There are $a \in (0,2/\beta)$ and $c > 0$, independent of $k$, such that, for each fixed $k$ and $R$, as $T \to \infty$,
\begin{equation}\label{an:joint.annular.bound}
\EE_{\star}[M_{\star,i}^{\mathrm{out}}(R)^{a} \mid E_{\bb{h}}] \leq C_{k,\e,J}R^{-c} + \oo(1),
\end{equation}
uniformly in $\bb{h}$ and $(x_{1},\ldots,x_{k})$ satisfying the stated constraints. The constants may also depend on $\beta$, $\|\chi\|_{\infty}$ and the fixed bin width $b$.
\end{lemma}

\begin{proof}[Proof of Lemma~\ref{an:joint.annular}]
The additional conditioning on $\smash{F_{h_{i'}}^{\star} - u_{T} \in [x_{i'},x_{i'} + b]}$ for $i' \neq i$, beyond the conditioning on the bin at $h_{i}$, changes the mean terms in \eqref{an:conditional.bound} of Lemma~\ref{an:conditional}. The first step is to estimate these changes. Define
\begin{equation}\label{an:annular.detail:joint.mass}
M_{\mathrm{ann}}^{\star,\mathrm{rel}}(d) \leqdef \ell_{T}^{-1}\int_{\{d\leq|h - h_{i}|\leq 2d\}\cap B_{j_{i}}^{\e}}\exp\{\beta(F_{h}^{\star} - F_{h_{i}}^{\star})\}\rd{h}.
\end{equation}
The aim is to bound the conditional fractional moment of $M_{\mathrm{ann}}^{\star,\mathrm{rel}}$ by $Ce^{-cr}$, where $d = \ell_{T}e^{r}$ and the decay rate $c > 0$ is independent of the fixed number $k$ of representatives.

Choose $m_{-} \in (1,m(\theta))$, independently of $k$. Since $h_{i}$ and $h_{i'}$ belong to distinct intervals $B_{j_{i}}^{\e}$ and $B_{j_{i'}}^{\e}$, $|h_{i} - h_{i'}| \geq 2\e$ whenever $i \neq i'$. Thus \eqref{an:overlap.geometry} in Lemma~\ref{an:covariance} gives the covariance matrix of the representative heights as
\begin{equation}\label{an:annular.detail:joint.covariance}
\Sigma_{T} \leqdef (\Cov_{\mathrm{H}}(H(h_{i}),H(h_{i'})))_{1 \leq i,i' \leq k}
= V_{T}\mathrm{Id}_{k} + \mathscr{E}_{T}, \qquad \|\mathscr{E}_{T}\|_{\mathrm{op}} \leq C_{k,\e}\log{t}.
\end{equation}
By \eqref{an:covariance.comparison} in Lemma~\ref{an:covariance}, \eqref{eq:parameter.scales}, and \eqref{eq:effective.cutoff.bounds}, $V_{T} = s_{T}/2 + \OO(1)$ and $s_{T} = t + \OO(\log{t})$. Thus every eigenvalue of $\Sigma_{T}$ lies between $V_{T} - \|\mathscr{E}_{T}\|_{\mathrm{op}}$ and $V_{T} + \|\mathscr{E}_{T}\|_{\mathrm{op}}$, and is comparable to $t$. Moreover, \eqref{an:annular.detail:joint.covariance} implies $\|V_{T}^{-1}\mathscr{E}_{T}\|_{\mathrm{op}} = \OO_{k,\e}((\log{t})/t)$. Since $\Sigma_{T} = V_{T}(\mathrm{Id}_{k} + V_{T}^{-1}\mathscr{E}_{T})$, applying the identity $(\mathrm{Id}_{k} + A)^{-1} - \mathrm{Id}_{k} = - (\mathrm{Id}_{k} + A)^{-1}A$ with $A = V_{T}^{-1}\mathscr{E}_{T}$ and taking operator norms gives
\begin{equation}\label{an:annular.detail:joint.inverse}
\|\Sigma_{T}^{-1} - V_{T}^{-1}\mathrm{Id}_{k}\|_{\mathrm{op}}
\leq V_{T}^{-2}\|(\mathrm{Id}_{k} + V_{T}^{-1}\mathscr{E}_{T})^{-1}\|_{\mathrm{op}}\|\mathscr{E}_{T}\|_{\mathrm{op}}
\leq C_{k,\e}\frac{\log{t}}{t^{2}}.
\end{equation}
Let $\bb{v} \leqdef (u_{T} + x_{1},\ldots,u_{T} + x_{k})$ be the vector of lower endpoints of the bins for the uncentered heights, so that $E_{\bb{h}} = \{(F_{h_{i'}}^{\star})_{1 \leq i' \leq k} \in \bb{v} + [0,b]^{k}\}$. Since $\bb{v} = u_{T}\bb{1}_{k} + \OO_{k,J}(1)$ and $u_{T} = \OO(t)$ by \eqref{an:height.scales} in Lemma~\ref{an:height}, \eqref{an:annular.detail:joint.inverse} yields the tilt vector
\begin{equation}\label{an:annular.detail:joint.tilt}
\bb{\lambda}_{T} \leqdef \Sigma_{T}^{-1}\bb{v} = \frac{u_{T}}{V_{T}}\bb{1}_{k} + \OO_{k,\e,J}\left(\frac{\log{t}}{t}\right), \qquad \frac{u_{T}}{V_{T}} \to 2m(\theta).
\end{equation}
In particular, all coordinates of $\bb{\lambda}_{T}$ are bounded and $\lambda_{T,i} \geq 2m_{-}$ for sufficiently large $T$, uniformly in the representative positions and the allowed bins.

For an annulus of radius $d = \ell_{T}e^{r}$ centered at $h_{i}$, with $0 \leq r \leq s_{T}$, take $h_{i}$ as the root and use the decomposition $F_{h}^{\star} = F_{d}^{\star,\mathrm{lo}}(h) + F_{d}^{\star,\mathrm{hi}}(h)$ from the proof of Lemma~\ref{an:annular}. The low-frequency part contains the prime powers with $\log(p^{j}) \leq 1/d$, and the high-frequency part contains those with $\log(p^{j}) > 1/d$. For $h$ in the domain of integration in \eqref{an:annular.detail:joint.mass}, write
\[
\begin{aligned}
X^{\star} &\leqdef F_{d}^{\star,\mathrm{hi}}(h_{i}), \\
Y^{\star}(h) &\leqdef \widetilde{F}_{d}^{\star,\mathrm{hi}}(h) = F_{h}^{\star} - F_{d}^{\star,\mathrm{lo}}(h_{i})
= F_{d}^{\star,\mathrm{hi}}(h) + F_{d}^{\star,\mathrm{lo}}(h) - F_{d}^{\star,\mathrm{lo}}(h_{i}).
\end{aligned}
\]
Every allowed $h$ and $h_{i}$ is separated from the other representatives by at least $2\e$. The local covariance bounds \eqref{an:annular.covariances} apply to $X^{\mathrm{H}}$ and $Y^{\mathrm{H}}(h)$. The additional covariances with the other representatives require separate estimates at small and large annular radii.

First suppose $d \leq t^{-3}$. The prime terms in $X^{\mathrm{H}}$ and $F_{d}^{\mathrm{H},\mathrm{hi}}(h)$ have $\log{p} > 1/d \geq t^{3}$. For either displacement $\xi = h_{i} - h_{i'}$ or $\xi = h - h_{i'}$, where $i' \neq i$, the separation of $B_{j_{i}}^{\e}$ and $B_{j_{i'}}^{\e}$ gives $|\xi| \geq 2\e$. The integration-by-parts argument used for \eqref{an:high.remote}, now starting at $1/d$, gives
\[
\left|\int_{1/d}^{\log{X}}e^{-2\eta_{\sig}v}\cos(\xi v)\frac{\rd{v}}{v}\right| \leq \frac{2de^{-2\eta_{\sig}/d}}{|\xi|} \leq \frac{2d}{|\xi|} \leq C_{\e}d.
\]
The quantitative prime summation argument in \eqref{an:detail:cov.prime.integral}, with lower endpoint $1/d \geq t^{3}$, bounds the summation error by $C(1 + |I_{\theta}|)t^{C}e^{-ct^{3/2}} = \oo(1)$ uniformly for $|\xi| \leq 2|I_{\theta}|$. The contribution from powers with exponent at least two and $\log(p^{j}) > 1/d$ tends to zero by absolute summability. Hence, for $i' \neq i$ and $h \in B_{j_{i}}^{\e}$ with $d \leq |h - h_{i}| \leq 2d$,
\begin{equation}\label{an:annular.detail:distant.covariances}
\Cov_{\mathrm{H}}(X^{\mathrm{H}},H(h_{i'})) = \oo(1), \quad
\Cov_{\mathrm{H}}(F_{d}^{\mathrm{H},\mathrm{hi}}(h),H(h_{i'})) = \oo(1), \quad
|\Cov_{\mathrm{H}}(Y^{\mathrm{H}}(h),H(h_{i'}))| \leq C.
\end{equation}
For the third bound in \eqref{an:annular.detail:distant.covariances}, the high-frequency part $F_{d}^{\mathrm{H},\mathrm{hi}}(h)$ is covered by the second bound in \eqref{an:annular.detail:distant.covariances}. For the low-frequency increment $\smash{F_{d}^{\mathrm{H},\mathrm{lo}}(h) - F_{d}^{\mathrm{H},\mathrm{lo}}(h_{i})}$, taking absolute values before summation and using $|e^{\ii(h - h_{i})\log{p}} - 1| \leq 2d\log{p}$ gives
\[
|\Cov_{\mathrm{H}}(F_{d}^{\mathrm{H},\mathrm{lo}}(h) - F_{d}^{\mathrm{H},\mathrm{lo}}(h_{i}),H(h_{i'}))|
\leq Cd\sum_{\log{p}\leq 1/d}\frac{\log{p}}{p} + C \leq Cd\left(\frac{1}{d} + 1\right) + C \leq C.
\]
This proves the third bound in \eqref{an:annular.detail:distant.covariances}.

To control the supremum of $Y^{\star}$ on a cell $B \subseteq \{d \leq |h - h_{i}| \leq 2d\} \cap B_{j_{i}}^{\e}$ of length at most $\ell_{T}$, define the scaled derivative
\[
D_{d}^{\star}(z) \leqdef \ell_{T}(Y^{\star})'(z) = \ell_{T}(F_{z}^{\star})', \qquad z \in B,
\]
where the primes denote differentiation with respect to $z$. The second equality holds because $F_{d}^{\star,\mathrm{lo}}(h_{i})$ does not depend on $z$. Equation~\eqref{an:annular.detail:derivative} in the proof of Lemma~\ref{an:annular} bounds the Haar variance of $D_{d}^{\mathrm{H}}(z)$ and its covariances with $X^{\mathrm{H}}$ and $Y^{\mathrm{H}}(h)$, uniformly for $h,z \in B$. The absolute prime-sum estimates \eqref{an:annular.detail:derivative.integrals} in that proof also give
\[
|\Cov_{\mathrm{H}}(D_{d}^{\mathrm{H}}(z),H(h_{i'}))| \leq C, \qquad 1 \leq i' \leq k,
\]
independently of the displacement $z - h_{i'}$. These derivative bounds hold for all $0 \leq r \leq s_{T}$.

Apply \eqref{an:conditional.bound} in Lemma~\ref{an:conditional} with $\smash{Y_{T,\nu}^{\star} = F_{h_{\nu}}^{\star}}$ for $1 \leq \nu \leq k$, and write $\bb{F}^{\star} \leqdef (F_{h_{1}}^{\star},\ldots,F_{h_{k}}^{\star})$. At $\lambda_{k + 1} = 1$, the conditional mean term for $Y_{T,k + 1}^{\star}$ equals $\Cov_{\mathrm{H}}(\bb{F}^{\mathrm{H}},Y_{T,k + 1}^{\mathrm{H}})^{\mathsf{T}}\bb{\lambda}_{T}$. Together, \eqref{an:annular.covariances}, \eqref{an:annular.detail:joint.tilt}, and \eqref{an:annular.detail:distant.covariances} therefore give
\begin{equation}\label{an:annular.detail:joint.means}
\begin{aligned}
\Cov_{\mathrm{H}}(\bb{F}^{\mathrm{H}},X^{\mathrm{H}})^{\mathsf{T}}\bb{\lambda}_{T} &= \lambda_{T,i}\frac{r}{2} + \OO_{k,\e,J}(1) \geq m_{-}r - C_{k,\e,J}, \\
\left|\Cov_{\mathrm{H}}(\bb{F}^{\mathrm{H}},Y^{\mathrm{H}}(h))^{\mathsf{T}}\bb{\lambda}_{T}\right| &\leq C_{k,\e,J}.
\end{aligned}
\end{equation}
Let $q = 2/\beta < 1$, which holds because $\beta > \beta_{c}(\theta) > 2$. For any fixed $0 < a < q$, the local variances in \eqref{an:annular.covariances} give
\[
\Var_{\mathrm{H}}( - a\beta X^{\mathrm{H}} + q\beta Y^{\mathrm{H}}(h)) = \frac{(a^{2} + q^{2})\beta^{2}r}{2} + \OO(1).
\]
Apply \eqref{an:conditional.bound} in Lemma~\ref{an:conditional} first with $Y_{T,k + 1}^{\star} = - a\beta X^{\star}$ and then with $Y_{T,k + 1}^{\star} = - a\beta X^{\star} + q\beta Y^{\star}(h)$, taking $\lambda_{k + 1} = 1$ in both cases. In the arithmetic case, choose the common event $\smash{\mathcal{A}_{T,\mathrm{sp}}^{\mathrm{A}}}$ from Lemma~\ref{an:spatial.transform} for the three families $\{F_{h}^{\star}\}_{h \in I_{\theta}}$, $\{F_{d}^{\star,\mathrm{lo}}(h)\}_{h \in I_{\theta}}$, and $\{\ell_{T}(F_{h}^{\star})'\}_{h \in I_{\theta}}$. As before, $\mathcal{A}_{T,\mathrm{sp}}$ denotes this event when $\star = \mathrm{A}$ and the whole probability space when $\star = \mathrm{H}$, so that $\ind_{\mathcal{A}_{T,\mathrm{sp}}} = 1$ in the Haar case, where the exponential estimates in \eqref{an:annular.detail:joint.transforms} hold without truncation. Their coefficient assumptions follow from \eqref{an:annular.detail:coefficient.admissibility} in the proof of Lemma~\ref{an:annular}, and $X^{\star}$, $Y^{\star}(h)$, and $D_{d}^{\star}(z)$ are fixed real linear combinations of their evaluations. The fixed-width-bin version of \eqref{an:conditional.bound} in Lemma~\ref{an:conditional} therefore applies to $E_{\bb{h}}$ with this common event. For $h,z \in B$, adding $\pm q\beta D_{d}^{\star}(z)$ to the coordinate $ - a\beta X^{\star} + q\beta Y^{\star}(h)$ changes its Haar variance by $\OO(1)$, by \eqref{an:annular.detail:derivative} in the proof of Lemma~\ref{an:annular}. Its conditional mean term changes by $\pm q\beta\Cov_{\mathrm{H}}(\bb{F}^{\mathrm{H}},D_{d}^{\mathrm{H}}(z))^{\mathsf{T}}\bb{\lambda}_{T} = \OO_{k,\e,J}(1)$, by the derivative covariance bounds above and the bounded coordinates of $\bb{\lambda}_{T}$ in \eqref{an:annular.detail:joint.tilt}. Applying \eqref{an:conditional.bound} also to these two augmented coordinates, combining the variance bounds with \eqref{an:annular.detail:joint.means}, and dropping the nonnegative variance correction $\bb{c}_{T}^{\mathsf{T}}\Sigma_{T}^{-1}\bb{c}_{T}$ yields the following estimates, with the third valid for both signs and uniformly for $h,z \in B$:
\begin{equation}\label{an:annular.detail:joint.transforms}
\begin{aligned}
\EE_{\star}[e^{-a\beta X^{\star}}\ind_{\mathcal{A}_{T,\mathrm{sp}}} \mid E_{\bb{h}}] &\leq C\exp\left\{ - a\beta m_{-}r + \frac{a^{2}\beta^{2}r}{4}\right\}, \\
\EE_{\star}[e^{-a\beta X^{\star} + q\beta Y^{\star}(h)}\ind_{\mathcal{A}_{T,\mathrm{sp}}} \mid E_{\bb{h}}] &\leq C\exp\left\{ - a\beta m_{-}r + \frac{(a^{2} + q^{2})\beta^{2}r}{4}\right\}, \\
\EE_{\star}[e^{-a\beta X^{\star} + q\beta Y^{\star}(h) \pm q\beta D_{d}^{\star}(z)}\ind_{\mathcal{A}_{T,\mathrm{sp}}} \mid E_{\bb{h}}] &\leq C\exp\left\{ - a\beta m_{-}r + \frac{(a^{2} + q^{2})\beta^{2}r}{4}\right\}.
\end{aligned}
\end{equation}

For a cell $B$ of positive length at most $\ell_{T}$, let $h_{B}$ denote its left endpoint. Apply \eqref{an:annular.detail:cell.supremum} from the proof of Lemma~\ref{an:annular} with $Y^{\star}$ in place of $\widetilde{F}_{d}^{\star,\mathrm{hi}}$, multiply by $e^{-a\beta X^{\star}}\ind_{\mathcal{A}_{T,\mathrm{sp}}}$, and take conditional expectations given $E_{\bb{h}}$. Using $e^{q\beta|x|} \leq e^{q\beta x} + e^{-q\beta x}$ together with the third estimate in \eqref{an:annular.detail:joint.transforms} at $h = h_{B}$ gives
\begin{equation}\label{an:annular.detail:joint.cell}
\EE_{\star}\left[e^{-a\beta X^{\star}}\exp\left(q\beta\sup_{h \in B}Y^{\star}(h)\right)\ind_{\mathcal{A}_{T,\mathrm{sp}}} \mid E_{\bb{h}}\right] \leq C\exp\left\{ - a\beta m_{-}r + \frac{(a^{2} + q^{2})\beta^{2}r}{4}\right\}.
\end{equation}
Write
\[
M_{d}^{\star,\mathrm{inc}} \leqdef \ell_{T}^{-1}\int_{\{d \leq |h - h_{i}| \leq 2d\}\cap B_{j_{i}}^{\e}}e^{\beta Y^{\star}(h)}\rd{h}, \qquad W \leqdef e^{-a\beta X^{\star}}\ind_{\mathcal{A}_{T,\mathrm{sp}}}.
\]
Partition $[h_{i} - 2d,h_{i} - d] \cap B_{j_{i}}^{\e}$ and $[h_{i} + d,h_{i} + 2d] \cap B_{j_{i}}^{\e}$ into at most $Ce^{r}$ cells of length at most $\ell_{T}$. Subadditivity of $x \mapsto x^{q}$ for $q \leq 1$ and \eqref{an:annular.detail:joint.cell} give
\begin{equation}\label{an:annular.detail:joint.cell.sum}
\begin{aligned}
\EE_{\star}[W(M_{d}^{\star,\mathrm{inc}})^{q} \mid E_{\bb{h}}]
&\leq \sum_{B}\EE_{\star}[We^{q\beta\sup_{h \in B}Y^{\star}(h)} \mid E_{\bb{h}}] \\
&\leq C\exp\left\{\left[1 - a\beta m_{-} + \frac{(a^{2} + q^{2})\beta^{2}}{4}\right]r\right\}.
\end{aligned}
\end{equation}
Since $F_{h}^{\star} - F_{h_{i}}^{\star} = Y^{\star}(h) - X^{\star}$ by \eqref{an:annular.detail:increment} in the proof of Lemma~\ref{an:annular}, the mass in \eqref{an:annular.detail:joint.mass} satisfies $M_{\mathrm{ann}}^{\star,\mathrm{rel}} = e^{-\beta X^{\star}}M_{d}^{\star,\mathrm{inc}}$. Apply the weighted H\"{o}lder calculation \eqref{an:annular.detail:weighted.holder}, using \eqref{an:annular.detail:joint.cell.sum} and the first bound in \eqref{an:annular.detail:joint.transforms}. It yields
\begin{equation}\label{an:joint.annular.exponent}
\begin{aligned}
\EE_{\star}[\{M_{\mathrm{ann}}^{\star,\mathrm{rel}}\}^{a}\ind_{\mathcal{A}_{T,\mathrm{sp}}} \mid E_{\bb{h}}]
&\leq \big(\EE_{\star}[W(M_{d}^{\star,\mathrm{inc}})^{q} \mid E_{\bb{h}}]\big)^{a/q}\big(\EE_{\star}[W \mid E_{\bb{h}}]\big)^{1 - a/q} \\
&\leq C\exp\left\{\left[\frac{a}{q}\left(1 + \frac{q^{2}\beta^{2}}{4}\right) - a\beta m_{-} + \frac{a^{2}\beta^{2}}{4}\right]r\right\} \\
&= C\exp\left\{\left[a\beta(1 - m_{-}) + \frac{a^{2}\beta^{2}}{4}\right]r\right\}.
\end{aligned}
\end{equation}
The last equality in \eqref{an:joint.annular.exponent} uses $q = 2/\beta$, so $q^{-1} + q\beta^{2}/4 = \beta$. Thus the coefficient $(a/q)(1 + q^{2}\beta^{2}/4)$ in the preceding exponent equals $a\beta$.

Choose
\begin{equation}\label{eq:choice.a.Holder}
0 < a < \min\left\{\frac{2}{\beta},\frac{2(m_{-} - 1)}{\beta}\right\}, \qquad c_{0} \leqdef a\beta(m_{-} - 1) - \frac{a^{2}\beta^{2}}{4} > 0.
\end{equation}
Both constants can be fixed independently of $k$. In fact, the restriction on $a$ gives $c_{0} > a\beta(m_{-} - 1)/2$. Equation~\eqref{an:joint.annular.exponent} is consequently bounded by $Ce^{-c_{0}r}$ for all annuli with $d \leq t^{-3}$.

For the remaining annuli, $d > t^{-3}$ and the identity $d = \ell_{T}e^{r} = e^{r - s_{T}}$ give $r = s_{T} + \log{d} > s_{T} - 3\log{t}$. Taking absolute values in the Haar covariance formula \eqref{an:detail:cov.haar.sum}, restricted to the prime powers with $\log(p^{j}) \leq 1/d$, and applying Mertens' estimate \eqref{analytic:eq:Mertens} gives
\begin{equation}\label{an:annular.detail:large.low.covariance}
\begin{aligned}
\sup_{z,z' \in I_{\theta}}|\Cov_{\mathrm{H}}(F_{d}^{\mathrm{H},\mathrm{lo}}(z),H(z'))|
&\leq \Var_{\mathrm{H}}(F_{d}^{\mathrm{H},\mathrm{lo}}(h_{i})) \\
&\leq \frac{1}{2}\sum_{\log{p}\leq 1/d}\frac{1}{p} + C
\leq C\log(2 + 1/d) \leq C\log{t}.
\end{aligned}
\end{equation}
Since $h,h_{i} \in B_{j_{i}}^{\e}$ are separated from each $h_{i'}$ with $i' \neq i$ by at least $2\e$, the full-field covariances $\Cov_{\mathrm{H}}(H(h_{i}),H(h_{i'}))$ and $\Cov_{\mathrm{H}}(H(h),H(h_{i'}))$ are $\OO_{\e}(\log{t})$ by \eqref{an:overlap.geometry} and \eqref{an:covariance.comparison} in Lemma~\ref{an:covariance}. Since $X^{\mathrm{H}} = H(h_{i}) - F_{d}^{\mathrm{H},\mathrm{lo}}(h_{i})$ and $Y^{\mathrm{H}}(h) = H(h) - F_{d}^{\mathrm{H},\mathrm{lo}}(h_{i})$, subtracting the low-frequency covariances bounded in \eqref{an:annular.detail:large.low.covariance} gives, for $i' \neq i$,
\[
|\Cov_{\mathrm{H}}(X^{\mathrm{H}},H(h_{i'}))| + |\Cov_{\mathrm{H}}(Y^{\mathrm{H}}(h),H(h_{i'}))| \leq C_{\e}\log{t}.
\]
The covariances with the root value $H(h_{i})$ and the local variances remain those in \eqref{an:annular.covariances}, with $h_{0} = h_{i}$. Multiplying these covariance estimates by the bounded coordinates of $\bb{\lambda}_{T}$ in \eqref{an:annular.detail:joint.tilt} yields the following replacements for \eqref{an:annular.detail:joint.means}:
\begin{equation}\label{an:annular.detail:large.means}
\begin{aligned}
\Cov_{\mathrm{H}}(\bb{F}^{\mathrm{H}},X^{\mathrm{H}})^{\mathsf{T}}\bb{\lambda}_{T} &\geq m_{-}r - C_{k,\e,J}\log{t}, \\
\left|\Cov_{\mathrm{H}}(\bb{F}^{\mathrm{H}},Y^{\mathrm{H}}(h))^{\mathsf{T}}\bb{\lambda}_{T}\right| &\leq C_{k,\e,J}\log{t}.
\end{aligned}
\end{equation}
Using \eqref{an:annular.detail:large.means} in place of \eqref{an:annular.detail:joint.means} adds at most a factor $\exp(C_{k,\e,J}\log{t})$ to the conditional transform bounds \eqref{an:annular.detail:joint.transforms} and hence to the cell bound \eqref{an:annular.detail:joint.cell}. The cell sum \eqref{an:annular.detail:joint.cell.sum} and the H\"{o}lder calculation \eqref{an:joint.annular.exponent} give
\begin{equation}\label{an:annular.detail:large.radius}
\EE_{\star}[\{M_{\mathrm{ann}}^{\star,\mathrm{rel}}\}^{a}\ind_{\mathcal{A}_{T,\mathrm{sp}}} \mid E_{\bb{h}}] \leq C\exp\{ - c_{0}r + C_{k,\e,J}\log{t}\} \leq Ce^{-c_{0}r/2}
\end{equation}
for sufficiently large $T$, since $r \geq s_{T} - 3\log{t}$ and $s_{T}/t \to 1$. Hence \eqref{an:joint.annular.exponent} and \eqref{an:annular.detail:large.radius} show that every relevant annulus has conditional fractional moment at most $Ce^{-cr}$ with $c \leqdef c_{0}/2$, independent of $k$; the constant and the threshold in $T$ may depend on $k$.

Cover $\smash{B_{j_{i}}^{\e} \setminus [h_{i} - R\ell_{T},h_{i} + R\ell_{T}]}$ by annuli of radii $R\ell_{T}2^{v}$, retaining only those with radius at most one. On $E_{\bb{h}}$, the height $F_{h_{i}}^{\star} - u_{T}$ lies in $J$, so
\begin{equation}\label{an:annular.detail:capped.annulus}
\begin{aligned}
M_{\star,i}^{\mathrm{out}}(R)
&\leq \|\chi\|_{\infty}e^{\beta\sup J}\sum_{\substack{v = 0\\R\ell_{T}2^{v}\leq 1}}^{\infty}\ell_{T}^{-1}\int_{\{R\ell_{T}2^{v}\leq|h - h_{i}|\leq R\ell_{T}2^{v + 1}\}\cap B_{j_{i}}^{\e}}\exp\{\beta(F_{h}^{\star} - F_{h_{i}}^{\star})\}\rd{h} \\
&= \|\chi\|_{\infty}e^{\beta\sup J}\sum_{\substack{v = 0\\R\ell_{T}2^{v}\leq 1}}^{\infty}M_{\mathrm{ann}}^{\star,\mathrm{rel}}(R\ell_{T}2^{v}).
\end{aligned}
\end{equation}
Choose $\smash{\mathcal{A}_{T,\mathrm{sp}}^{\mathrm{A}}}$ simultaneously for $\{F_{h}^{\star}\}_{h \in I_{\theta}}$, $\{\ell_{T}(F_{h}^{\star})'\}_{h \in I_{\theta}}$, and $\{F_{d}^{\star,\mathrm{lo}}(h)\}_{h \in I_{\theta}}$ at their $\OO(t)$ frequency cutoffs $1/d$. This common event satisfies the uniform bounds \eqref{an:spatial.exceptional} and \eqref{an:spatial.transform.bound} in Lemma~\ref{an:spatial.transform}. For the annulus of radius $d = R\ell_{T}2^{v}$, the scale parameter is $r = \log(d/\ell_{T}) = \log{R} + v\log{2}$. Since $a < q < 1$, applying subadditivity of $x \mapsto x^{a}$ to \eqref{an:annular.detail:capped.annulus} and then using \eqref{an:joint.annular.exponent} for annuli with $d \leq t^{-3}$ and \eqref{an:annular.detail:large.radius} for annuli with $d > t^{-3}$ give, after absorbing $(\|\chi\|_{\infty}e^{\beta\sup J})^{a}$ into $C$,
\begin{equation}\label{an:annular.detail:outer.sum}
\begin{aligned}
\EE_{\star}[M_{\star,i}^{\mathrm{out}}(R)^{a}\ind_{\mathcal{A}_{T,\mathrm{sp}}} \mid E_{\bb{h}}]
&\leq C\sum_{\substack{v = 0\\R\ell_{T}2^{v}\leq 1}}^{\infty}\EE_{\star}[M_{\mathrm{ann}}^{\star,\mathrm{rel}}(R\ell_{T}2^{v})^{a}\ind_{\mathcal{A}_{T,\mathrm{sp}}} \mid E_{\bb{h}}] \\
&\leq C\sum_{v = 0}^{\infty}e^{-c(\log{R} + v\log{2})}
= \frac{C}{1 - 2^{-c}}R^{-c}.
\end{aligned}
\end{equation}
If $R\ell_{T} \geq 1$, the outer mass is zero and the inequality holds immediately.

It remains to remove the arithmetic event. Since $\chi$ is supported on $J$ and $|B_{j_{i}}^{\e}| \leq 1$, the deterministic bound
\[
M_{\star,i}^{\mathrm{out}}(R)^{a} \leq \left(\ell_{T}^{-1}\|\chi\|_{\infty}e^{\beta\sup J}\right)^{a} \leq e^{C_{a,J}t}
\]
holds for all sufficiently large $T$, because $\ell_{T}^{-1} = e^{s_{T}}$ and $s_{T} = t + \OO(\log{t})$ by \eqref{eq:effective.cutoff.bounds} and \eqref{eq:parameter.scales}. The fixed-bin lower bound \eqref{an:detail:conditional.bin.lower} in the proof of Lemma~\ref{an:conditional}, with the bin-width extension justified by \eqref{an:conditional.detail:variable.bin.integrals}, gives $\PP_{\star}(E_{\bb{h}}) \geq e^{-C_{k,\e,J}t}$. This is uniform because \eqref{an:annular.detail:joint.covariance} gives eigenvalues comparable to $t$ and $\|\bb{v}\|_{2} = \OO_{k,J}(t)$. Equation~\eqref{an:spatial.exceptional} in Lemma~\ref{an:spatial.transform} therefore implies
\begin{equation}\label{an:annular.detail:event.removal}
\EE[M_{\mathrm{A},i}^{\mathrm{out}}(R)^{a}\ind_{\mathcal{A}_{T,\mathrm{sp}}^{c}} \mid E_{\bb{h}}] \leq e^{C_{a,J}t}\frac{\PP(\mathcal{A}_{T,\mathrm{sp}}^{c})}{\PP(E_{\bb{h}})} \leq \exp\{ - ct^{1 + 2\e_{0}} + \OO_{k,\e,J}(t)\} = \oo(1).
\end{equation}
There is no exceptional-event term analogous to \eqref{an:annular.detail:event.removal} under the Haar law: the conditional estimate \eqref{an:conditional.bound} in Lemma~\ref{an:conditional} holds without truncation, so $\mathcal{A}_{T,\mathrm{sp}}$ denotes the whole probability space when $\star = \mathrm{H}$ and its complement is empty. Combining \eqref{an:annular.detail:event.removal} with \eqref{an:annular.detail:outer.sum} proves the assertion \eqref{an:joint.annular.bound}.
\end{proof}

\subsubsection{Comparison after smoothing a fixed vector}\label{an:smoothing.section}

The following estimate concerns a fixed number of raw field values. The positions may coincide, and no nondegeneracy assumption on their covariance matrix is needed.

\begin{lemma}[Comparison after Gaussian smoothing]\label{an:smoothing}
Let $\theta > 0$. Fix $n_{\mathrm{dim}} \in \N$. Let $h_{1},\ldots,h_{n_{\mathrm{dim}}}$ be arbitrary positions in $I_{\theta}$, and define
\[
\bb{Y}_{T}^{\star} \leqdef (F_{h_{r}}^{\star})_{r = 1}^{n_{\mathrm{dim}}}, \qquad \star \in \{\mathrm{A},\mathrm{H}\},
\]
and let $\bb{g}$ be an independent standard Gaussian vector in $\R^{n_{\mathrm{dim}}}$. There is $c_{n_{\mathrm{dim}}} > 0$ such that, uniformly over the selected positions,
\begin{equation}\label{an:smoothing.tv}
\left\|\mathcal{L}(\bb{Y}_{T}^{\mathrm{A}} + t^{-3/4}\bb{g}) - \mathcal{L}(\bb{Y}_{T}^{\mathrm{H}} + t^{-3/4}\bb{g})\right\|_{\mathrm{TV}} \leq \exp\{ - c_{n_{\mathrm{dim}}}t(\log{t})^{2}\}
\end{equation}
for all sufficiently large $T$.
\end{lemma}

\begin{proof}[Proof of Lemma~\ref{an:smoothing}]
Write $\varphi_{T}^{\star}(\bb{\omega}) \leqdef \EE_{\star}[\exp(\ii\bb{\omega}\cdot\bb{Y}_{T}^{\star})]$ for the characteristic function of the raw field vector. The proof compares these functions on a growing Fourier ball, uses the Gaussian factor $e^{-t^{-3/2}\|\bb{\omega}\|_{2}^{2}/2}$ to control frequencies outside that ball, and then integrates the resulting density bound over a box containing almost all the mass.

For the comparison on $\|\bb{\omega}\|_{2} \leq t^{5/4}\log{t}$, choose an even integer $D_{T}$ satisfying
\begin{equation}\label{an:smoothing.detail:degree}
D_{T} = Kt^{7/2}(\log{t})^{2} + \OO(1), \qquad D_{T} = \oo(t^{\phi}),
\end{equation}
where $K$ is a sufficiently large fixed constant. The second assertion in \eqref{an:smoothing.detail:degree} follows from $\phi > 7/2$ in \eqref{eq:fixed.parameters}, and permits the use of the arithmetic version of the moment bound \eqref{an:haar.moments} in Lemma~\ref{an:polynomial}. For real $x$, the integral remainder in Taylor's formula gives
\[
\left|e^{\ii x} - \sum_{q = 0}^{D_{T} - 1}\frac{(\ii x)^{q}}{q!}\right| \leq \frac{|x|^{D_{T}}}{D_{T}!}.
\]
Define the comparison polynomial by
\[
P_{T,\bb{\omega}}(\bb{y}) \leqdef \sum_{q = 0}^{D_{T} - 1}\frac{(\ii\bb{\omega}\cdot\bb{y})^{q}}{q!}.
\]
Apply the second bound in \eqref{an:haar.moments} of Lemma~\ref{an:polynomial}, or its arithmetic counterpart, to the Taylor remainder with $k = D_{T}/2$. The inequality $D_{T}! \geq (D_{T}/e)^{D_{T}}$ and the Fourier-ball bound $t\|\bb{\omega}\|_{2}^{2} \leq t^{7/2}(\log{t})^{2}$ give, uniformly in the selected positions,
\begin{equation}\label{an:smoothing.detail:taylor.error}
\EE_{\star}\left[\left|e^{\ii\bb{\omega}\cdot\bb{Y}_{T}^{\star}} - P_{T,\bb{\omega}}(\bb{Y}_{T}^{\star})\right|\right]
\leq \frac{(C_{n_{\mathrm{dim}}}D_{T}t\|\bb{\omega}\|_{2}^{2})^{D_{T}/2}}{(D_{T}/e)^{D_{T}}}
\leq \left(\frac{C_{n_{\mathrm{dim}}}t^{7/2}(\log{t})^{2}}{D_{T}}\right)^{D_{T}/2}
\leq e^{-cD_{T}}.
\end{equation}
The last inequality holds after choosing $K$ large enough in \eqref{an:smoothing.detail:degree}.

The polynomial satisfies $|P_{T,\bb{\omega}}(\bb{y})| \leq e^{|\bb{\omega}\cdot\bb{y}|}$. Since $e^{2|x|} \leq e^{2x} + e^{-2x}$, the first bound in \eqref{an:haar.moments} of Lemma~\ref{an:polynomial}, applied with $\bb{\lambda} = 2\bb{\omega}$ and $\bb{\lambda} = - 2\bb{\omega}$, yields
\begin{equation}\label{an:smoothing.detail:polynomial.norm}
\left(\EE_{\mathrm{H}}[|P_{T,\bb{\omega}}(\bb{Y}_{T}^{\mathrm{H}})|^{2}]\right)^{1/2}
\leq \left(\EE_{\mathrm{H}}[e^{2|\bb{\omega}\cdot\bb{Y}_{T}^{\mathrm{H}}|}]\right)^{1/2}
\leq \sqrt{2}\exp(C_{n_{\mathrm{dim}}}t^{7/2}(\log{t})^{2}).
\end{equation}
The real and imaginary parts have degree at most $D_{T} - 1$ and each has Haar $L^{2}$ norm bounded by \eqref{an:smoothing.detail:polynomial.norm}. Applying \eqref{an:polynomial.bound} in Lemma~\ref{an:polynomial} to these two polynomials therefore gives
\begin{equation}\label{an:smoothing.detail:polynomial.error}
\left|\EE[P_{T,\bb{\omega}}(\bb{Y}_{T}^{\mathrm{A}})] - \EE_{\mathrm{H}}[P_{T,\bb{\omega}}(\bb{Y}_{T}^{\mathrm{H}})]\right| \leq \frac{CX^{2D_{T}}}{T}\exp(C_{n_{\mathrm{dim}}}t^{7/2}(\log{t})^{2}) = T^{-1 + \oo(1)},
\end{equation}
by \eqref{eq:fixed.parameters} and \eqref{an:smoothing.detail:degree}. Since $D_{T} = \oo(\log{T})$, this error is smaller than $e^{-cD_{T}}$ for all sufficiently large $T$. Combining \eqref{an:smoothing.detail:taylor.error} under the two laws with \eqref{an:smoothing.detail:polynomial.error} proves
\begin{equation}\label{an:characteristic.ball}
\begin{aligned}
|\varphi_{T}^{\mathrm{A}}(\bb{\omega}) - \varphi_{T}^{\mathrm{H}}(\bb{\omega})|
&\leq \sum_{\star \in \{\mathrm{A},\mathrm{H}\}}\EE_{\star}\left[\left|e^{\ii\bb{\omega}\cdot\bb{Y}_{T}^{\star}} - P_{T,\bb{\omega}}(\bb{Y}_{T}^{\star})\right|\right]
+ \left|\EE[P_{T,\bb{\omega}}(\bb{Y}_{T}^{\mathrm{A}})] - \EE_{\mathrm{H}}[P_{T,\bb{\omega}}(\bb{Y}_{T}^{\mathrm{H}})]\right| \\
&\leq e^{-cD_{T}}, \qquad \|\bb{\omega}\|_{2} \leq t^{5/4}\log{t}.
\end{aligned}
\end{equation}

The smoothed vectors $\bb{Y}_{T}^{\star} + t^{-3/4}\bb{g}$ have bounded continuous densities $p_{T}^{\star}$. Independence of $\bb{g}$ identifies their characteristic functions with $\smash{e^{-t^{-3/2}\|\bb{\omega}\|_{2}^{2}/2}\varphi_{T}^{\star}(\bb{\omega})}$. These functions are integrable because $|\varphi_{T}^{\star}| \leq 1$, so Fourier inversion applies even when $(\Cov_{\star}(F_{h_{r}}^{\star},F_{h_{s}}^{\star}))_{1 \leq r,s \leq n_{\mathrm{dim}}}$ is singular.

On the Fourier ball $\|\bb{\omega}\|_{2} \leq t^{5/4}\log{t}$, use \eqref{an:characteristic.ball}; outside it, use $|\varphi_{T}^{\mathrm{A}} - \varphi_{T}^{\mathrm{H}}| \leq 2$. The change of variables $\bb{z} \leqdef t^{-3/4}\bb{\omega}$ sends the boundary of the ball to radius $t^{1/2}\log{t}$. The Gaussian tail integral is therefore bounded by
\begin{equation}\label{an:smoothing.detail:fourier.tail}
\begin{aligned}
\int_{\|\bb{\omega}\|_{2} > t^{5/4}\log{t}}e^{-t^{-3/2}\|\bb{\omega}\|_{2}^{2}/2}\rd{\bb{\omega}}
&= t^{3n_{\mathrm{dim}}/4}\int_{\|\bb{z}\|_{2} > t^{1/2}\log{t}}e^{-\|\bb{z}\|_{2}^{2}/2}\rd{\bb{z}} \\
&\leq t^{3n_{\mathrm{dim}}/4}e^{-t(\log{t})^{2}/4}\int_{\R^{n_{\mathrm{dim}}}}e^{-\|\bb{z}\|_{2}^{2}/4}\rd{\bb{z}} \\
&\leq C_{n_{\mathrm{dim}}}t^{3n_{\mathrm{dim}}/4}e^{-t(\log{t})^{2}/4}.
\end{aligned}
\end{equation}
Combining \eqref{an:characteristic.ball} and \eqref{an:smoothing.detail:fourier.tail} in the Fourier inversion formula gives
\begin{equation}\label{an:smoothing.detail:density.bound}
\begin{aligned}
\|p_{T}^{\mathrm{A}} - p_{T}^{\mathrm{H}}\|_{\infty}
&\leq \frac{1}{(2\pi)^{n_{\mathrm{dim}}}}\int_{\R^{n_{\mathrm{dim}}}}e^{-t^{-3/2}\|\bb{\omega}\|_{2}^{2}/2}|\varphi_{T}^{\mathrm{A}}(\bb{\omega}) - \varphi_{T}^{\mathrm{H}}(\bb{\omega})|\rd{\bb{\omega}} \\
&\leq C_{n_{\mathrm{dim}}}t^{5n_{\mathrm{dim}}/4}(\log{t})^{n_{\mathrm{dim}}}e^{-cD_{T}} + C_{n_{\mathrm{dim}}}t^{3n_{\mathrm{dim}}/4}e^{-ct(\log{t})^{2}} \\[1mm]
&\leq \exp\{ - ct(\log{t})^{2}\}.
\end{aligned}
\end{equation}
The final inequality absorbs the fixed powers of $t$ and $\log{t}$ by reducing $c$.

To convert the supremum bound \eqref{an:smoothing.detail:density.bound} into an $L^{1}$ bound, it remains to control both laws outside a box of polynomial volume. Choose a sufficiently large fixed $A_{0}$ and apply the moment bound in \eqref{an:haar.moments} of Lemma~\ref{an:polynomial}, or its arithmetic counterpart, with $k = \lfloor t^{2}\rfloor = \oo(t^{\phi})$ to each coordinate. Markov's inequality and a union bound yield
\[
\PP_{\star}\left(\max_{1 \leq q \leq n_{\mathrm{dim}}}|F_{h_{q}}^{\star}| > A_{0}t^{3/2}/2\right)
\leq n_{\mathrm{dim}}\frac{(C_{n_{\mathrm{dim}}}\lfloor t^{2}\rfloor t)^{\lfloor t^{2}\rfloor}}{(A_{0}t^{3/2}/2)^{2\lfloor t^{2}\rfloor}}
\leq n_{\mathrm{dim}}\left(\frac{4C_{n_{\mathrm{dim}}}}{A_{0}^{2}}\right)^{\lfloor t^{2}\rfloor} \leq e^{-ct^{2}}.
\]
The independent Gaussian perturbation $\bb{g}$ also satisfies
\[
\PP\left(\max_{1 \leq q \leq n_{\mathrm{dim}}}|t^{-3/4}g_{q}| > A_{0}t^{3/2}/2\right) \leq C_{n_{\mathrm{dim}}}\exp( - cA_{0}^{2}t^{9/2}) \leq e^{-ct^{2}}.
\]
Thus the box $Q_{T} \leqdef [ - A_{0}t^{3/2},A_{0}t^{3/2}]^{n_{\mathrm{dim}}}$ satisfies
\begin{equation}\label{an:smoothing.detail:box.tail}
\PP_{\star}(\bb{Y}_{T}^{\star} + t^{-3/4}\bb{g} \notin Q_{T}) \leq 2e^{-ct^{2}}.
\end{equation}
Integrate the density difference over $Q_{T}$ using \eqref{an:smoothing.detail:density.bound}, and over its complement using \eqref{an:smoothing.detail:box.tail}. This gives
\[
\begin{aligned}
\int_{\R^{n_{\mathrm{dim}}}}|p_{T}^{\mathrm{A}}(\bb{x}) - p_{T}^{\mathrm{H}}(\bb{x})|\rd{\bb{x}}
&\leq (2A_{0}t^{3/2})^{n_{\mathrm{dim}}}\|p_{T}^{\mathrm{A}} - p_{T}^{\mathrm{H}}\|_{\infty} + \sum_{\star \in \{\mathrm{A},\mathrm{H}\}}\PP_{\star}(\bb{Y}_{T}^{\star} + t^{-3/4}\bb{g} \notin Q_{T}) \\
&\leq \exp\{ - c_{n_{\mathrm{dim}}}t(\log{t})^{2}\}.
\end{aligned}
\]
The volume factor $(2A_{0}t^{3/2})^{n_{\mathrm{dim}}}$ is absorbed by reducing $c_{n_{\mathrm{dim}}}$, while $e^{-ct^{2}}$ is smaller than the claimed error bound for all sufficiently large $T$. Total variation is bounded by the resulting $L^{1}$ distance, proving \eqref{an:smoothing.tv}. The coefficient bounds used in \eqref{an:polynomial.bound} and \eqref{an:haar.moments} of Lemma~\ref{an:polynomial} do not depend on the positions, and no covariance inverse appears. Thus all estimates remain uniform when some selected positions coincide.
\end{proof}

\subsection{Mass bounds and Gibbs mass spreading}\label{sec:temperature}

Let $\mathcal{B}_{T}$ be the partition of $I_{\theta}$ into consecutive intervals of length one, with possibly shorter intervals at the endpoints. For $B \in \mathcal{B}_{T}$, $\star \in \{\mathrm{A},\mathrm{H}\}$, and a Borel set $A_{\mathrm{ht}} \subseteq \R$, define
\begin{equation}\label{temp:interval.masses}
Z_{B,\beta}^{\star}(A_{\mathrm{ht}}) \leqdef \ell_{T}^{-1}\int_{B}e^{\beta(F_{h}^{\star} - u_{T})}\ind_{\{F_{h}^{\star} - u_{T}\in A_{\mathrm{ht}}\}}\rd{h}, \qquad M_{\beta,T}^{\star} \leqdef \sum_{B \in \mathcal{B}_{T}}Z_{B,\beta}^{\star}(\R).
\end{equation}
Two intervals are called distant when their distance is at least two. Each interval has at most a fixed number of intervals that are not distant from it.

\subsubsection{A moment bound at critical temperature}

\begin{lemma}\label{temp:critical.moment}
Let $\theta > 0$. There are constants $a > 0$ and $C < \infty$ such that, for every $\star \in \{\mathrm{A},\mathrm{H}\}$ and every fixed $L \geq 0$,
\begin{equation}\label{temp:critical.moment.tail}
\limsup_{T \to \infty}\EE_{\star}\left[\sum_{B \in \mathcal{B}_{T}}Z_{B,\beta_{c}(\theta)}^{\star}(( - \infty, - L])^{1 + a}\right] \leq Ce^{-a\beta_{c}(\theta) L}.
\end{equation}
\end{lemma}

\begin{proof}[Proof of Lemma~\ref{temp:critical.moment}]
The proof conditions on the height of a root $h_{0} \in B$ in the identity \eqref{temp:detail:root.factorization} for $Z_{B,\beta_{c}(\theta)}^{\star}(( - \infty, - L])^{1 + a}$ and separates four contributions in \eqref{temp:detail:critical.decomposition}. Roots whose field heights are near $u_{T}$ give the dependence on $L$ in \eqref{temp:critical.moment.tail}; contributions from the remaining roots have a smaller exponential rate. Choose $m_{0} \in (1,m(\theta))$ and then choose $a > 0$ so small that
\begin{equation}\label{temp:detail:critical.parameters}
a < 1 \wedge m(\theta)^{-1}, \qquad 1 + \frac{am(\theta)}{2} < m_{0}, \qquad am(\theta) < m(\theta) - 1.
\end{equation}
These conditions are compatible because $m(\theta) > 1$. Abbreviate $Z_{B} \leqdef Z_{B,\beta_{c}(\theta)}^{\star}(( - \infty, - L])$. Multiplying the defining integral for $Z_{B}$ in \eqref{temp:interval.masses} by $Z_{B}^{a}$ gives
\begin{equation}\label{temp:detail:root.factorization}
Z_{B}^{1 + a} = \ell_{T}^{-1}\int_{B}\exp\{\beta_{c}(\theta)(F_{h_{0}}^{\star} - u_{T})\}\ind_{\{F_{h_{0}}^{\star}\leq u_{T} - L\}}Z_{B}^{a}\rd{h_{0}}.
\end{equation}

For the arithmetic law, take $\mathcal{A}_{T,\mathrm{sp}}$ to be the event $\smash{\mathcal{A}_{T,\mathrm{sp}}^{\mathrm{A}}}$ defined in \eqref{an:annular.unit.spatial.event}, which is used in \eqref{an:unit.fractional} of Lemma~\ref{an:annular}. The Haar conditional bounds in \eqref{an:unit.fractional} of Lemma~\ref{an:annular} hold without an exceptional event, so $\mathcal{A}_{T,\mathrm{sp}}$ denotes the whole probability space under the Haar law.

For every $j \in \N_0$ satisfying $u_{T} - j - 1 \geq m_{0}s_{T}$, denote the contribution of the bin $[u_{T} - j - 1,u_{T} - j]$ to the identity \eqref{temp:detail:root.factorization}, restricted to the event $\mathcal{A}_{T,\mathrm{sp}}$, by
\begin{equation}\label{temp:detail:critical.bin.contribution}
\mathcal{C}_{T,j}^{\star} \leqdef \EE_{\star}\Bigg[\sum_{B \in \mathcal{B}_{T}}\ell_{T}^{-1}\int_{B}\exp\{\beta_{c}(\theta)(F_{h_{0}}^{\star} - u_{T})\}Z_{B}^{a} \ind_{\{F_{h_{0}}^{\star}\in [u_{T} - j - 1,u_{T} - j],\,F_{h_{0}}^{\star}\leq u_{T} - L\}}\ind_{\mathcal{A}_{T,\mathrm{sp}}}\rd{h_{0}}\Bigg].
\end{equation}
Denote the contributions of the lower nonnegative roots, the negative roots, and the exceptional event by $\mathcal{C}_{T,\mathrm{low}}^{\star}$, $\mathcal{C}_{T,\mathrm{neg}}^{\star}$, and $\mathcal{C}_{T,\mathrm{exc}}^{\star}$, respectively, where
\begin{equation}\label{temp:detail:critical.other.contributions}
\begin{aligned}
\mathcal{C}_{T,\mathrm{low}}^{\star}
&\leqdef \EE_{\star}\Bigg[\sum_{B \in \mathcal{B}_{T}}\ell_{T}^{-1}\int_{B}\exp\{\beta_{c}(\theta)(F_{h_{0}}^{\star} - u_{T})\}Z_{B}^{a} \ind_{\{0 \leq F_{h_{0}}^{\star} < m_{0}s_{T} + 1,\,F_{h_{0}}^{\star}\leq u_{T} - L\}}\ind_{\mathcal{A}_{T,\mathrm{sp}}}\rd{h_{0}}\Bigg], \\
\mathcal{C}_{T,\mathrm{neg}}^{\star}
&\leqdef \EE_{\star}\Bigg[\sum_{B \in \mathcal{B}_{T}}\ell_{T}^{-1}\int_{B}\exp\{\beta_{c}(\theta)(F_{h_{0}}^{\star} - u_{T})\}Z_{B}^{a}\ind_{\{F_{h_{0}}^{\star} < 0\}}\rd{h_{0}}\Bigg], \\
\mathcal{C}_{T,\mathrm{exc}}^{\star}
&\leqdef \EE_{\star}\left[\sum_{B \in \mathcal{B}_{T}}Z_{B}^{1 + a}\ind_{\mathcal{A}_{T,\mathrm{sp}}^{c}}\right].
\end{aligned}
\end{equation}

By \eqref{temp:detail:root.factorization}, separating $\mathcal{A}_{T,\mathrm{sp}}$ from its complement and splitting the root integral according to the field value gives, for all sufficiently large $T$,
\begin{equation}\label{temp:detail:critical.decomposition}
\EE_{\star}\left[\sum_{B \in \mathcal{B}_{T}}Z_{B}^{1 + a}\right] \leq \sum_{\substack{j \geq \lfloor L\rfloor \\ u_{T} - j - 1 \geq m_{0}s_{T}}}\mathcal{C}_{T,j}^{\star} + \mathcal{C}_{T,\mathrm{low}}^{\star} + \mathcal{C}_{T,\mathrm{neg}}^{\star} + \mathcal{C}_{T,\mathrm{exc}}^{\star}.
\end{equation}
The smallest lower endpoint among the high bins lies in $[m_{0}s_{T},m_{0}s_{T} + 1)$ for all sufficiently large $T$. Thus $\mathcal{C}_{T,\mathrm{low}}^{\star}$ in \eqref{temp:detail:critical.decomposition} includes all roots with nonnegative field values outside the high bins. It will be bounded by summing the left-hand side of \eqref{temp:detail:critical.lower.bin} over the unit bins covering $[0,m_{0}s_{T} + 1)$. Any overlap with the high bins cannot reduce this upper bound. In the definition of $\mathcal{C}_{T,\mathrm{neg}}^{\star}$ in \eqref{temp:detail:critical.other.contributions}, the height cap and the event $\mathcal{A}_{T,\mathrm{sp}}$ have been removed for an upper bound. Under the Haar law, $\mathcal{C}_{T,\mathrm{exc}}^{\mathrm{H}} = 0$.

The four terms on the right-hand side of \eqref{temp:detail:critical.decomposition} will be controlled, respectively, by \eqref{temp:detail:critical.high.roots}, \eqref{temp:detail:critical.lower.sum}, \eqref{temp:detail:critical.negative.roots}, and \eqref{temp:detail:cap.exceptional}.

For each fixed root $h_{0} \in B$, write $F_{h}^{\star} - u_{T} = (F_{h_{0}}^{\star} - u_{T}) + (F_{h}^{\star} - F_{h_{0}}^{\star})$ in the defining integral for $Z_{B}$ in \eqref{temp:interval.masses}. The factor $\smash{\exp\{\beta_{c}(\theta)(F_{h_{0}}^{\star} - u_{T})\}}$ does not depend on the integration variable $h$ and can be taken outside the integral. Raising this identity to the power $a > 0$ and then bounding the height indicator by one yields
\begin{equation}\label{temp:detail:root.relative.bound}
Z_{B}^{a} \leq \exp\{a\beta_{c}(\theta)(F_{h_{0}}^{\star} - u_{T})\}\left(\ell_{T}^{-1}\int_{B}\exp\{\beta_{c}(\theta)(F_{h}^{\star} - F_{h_{0}}^{\star})\}\rd{h}\right)^{a}.
\end{equation}
Every $h \in B$ satisfies $|h - h_{0}| \leq 1$. Thus the integral $\ell_{T}^{-1}\int_{B}\exp\{\beta_{c}(\theta)(F_{h}^{\star} - F_{h_{0}}^{\star})\}\rd{h}$ in \eqref{temp:detail:root.relative.bound} is bounded by the integral over $[h_{0} - 1,h_{0} + 1] \cap I_{\theta}$ in \eqref{an:unit.fractional} of Lemma~\ref{an:annular}. This permits conditioning on each unit root-height bin before integrating in $h_{0}$.

To bound $\mathcal{C}_{T,\mathrm{exc}}^{\mathrm{A}}$ in \eqref{temp:detail:critical.decomposition}, note that $L \geq 0$ and the height restriction give $Z_{B} \leq \ell_{T}^{-1}|B| \leq e^{s_{T}}$. Moreover, $|\mathcal{B}_{T}| \leq e^{\theta t + \OO(1)}$. Hence, by \eqref{an:spatial.exceptional} of Lemma~\ref{an:spatial.transform},
\begin{equation}\label{temp:detail:cap.exceptional}
\begin{aligned}
\mathcal{C}_{T,\mathrm{exc}}^{\mathrm{A}}
&\leq |\mathcal{B}_{T}|e^{(1 + a)s_{T}}\PP(\mathcal{A}_{T,\mathrm{sp}}^{c}) \\
&\leq \exp\{(1 + a)s_{T} + \theta t + \OO(1)\}\exp\{ - ct^{1 + 2\e_{0}} + Ct\} \\
&\leq \exp\{ - ct^{1 + 2\e_{0}} + \OO_{a}(t)\}.
\end{aligned}
\end{equation}

First consider the root-height bins $[u_{T} - j - 1,u_{T} - j]$ indexed by $j \in \N_0$ satisfying $u_{T} - j - 1 \geq m_{0}s_{T}$; such bins exist for all sufficiently large $T$ because $u_{T}/s_{T} \to m(\theta)$ and $m_{0} < m(\theta)$. Set $\gamma \leqdef (u_{T} - j - 1)/s_{T}$. Since $\beta_{c}(\theta) = 2m(\theta) > 2$, the coefficient of $r$ in the exponent on the right-hand side of \eqref{an:annular.bound} of Lemma~\ref{an:annular} satisfies
\begin{equation}\label{temp:detail:critical.annular.exponent}
a\beta_{c}(\theta)(1 - \gamma) + \frac{a^{2}\beta_{c}(\theta)^{2}}{4} \leq 2am(\theta)(1 - m_{0}) + a^{2}m(\theta)^{2} \leq -c < 0,
\end{equation}
for some fixed $c > 0$, by \eqref{temp:detail:critical.parameters}. Therefore the geometric sum in \eqref{an:annular.detail:unit.sum} gives
\begin{equation}\label{temp:detail:critical.conditional.unit}
\begin{aligned}
&\EE_{\star}\left[\left(\ell_{T}^{-1}\int_{B}\exp\{\beta_{c}(\theta)(F_{h}^{\star} - F_{h_{0}}^{\star})\}\rd{h}\right)^{a}\ind_{\mathcal{A}_{T,\mathrm{sp}}} \mid F_{h_{0}}^{\star} \in [u_{T} - j - 1,u_{T} - j]\right] \\
&\qquad\leq C\left(1 + \sum_{i = 0}^{\infty}2^{-ci}\right)
\leq C'.
\end{aligned}
\end{equation}
The constants in \eqref{temp:detail:critical.conditional.unit} are independent of $j$, $T$, and the root position $h_{0} \in B$.

For each such bin, the factors $\exp\{\beta_{c}(\theta)(F_{h_{0}}^{\star} - u_{T})\}$ in \eqref{temp:detail:root.factorization} and $\exp\{a\beta_{c}(\theta)(F_{h_{0}}^{\star} - u_{T})\}$ in \eqref{temp:detail:root.relative.bound} have product at most $e^{-(1 + a)\beta_{c}(\theta)j}$. After dropping the height restriction $F_{h_{0}}^{\star}\leq u_{T} - L$, conditioning on the root-height bin allows the remaining relative integral to be bounded by \eqref{temp:detail:critical.conditional.unit}. Bounding the root-height probability uniformly in $h_{0}$ and summing the integrals over $B$ in \eqref{temp:detail:critical.bin.contribution} then contributes $\sum_{B \in \mathcal{B}_{T}}\ell_{T}^{-1}|B| = \ell_{T}^{-1}|I_{\theta}| = N_{T}$ by \eqref{eq:N.T}. Using also the depth estimate \eqref{an:depth.bound} in Lemma~\ref{an:height} therefore gives
\begin{equation}\label{temp:detail:critical.depth.sum}
\begin{aligned}
\mathcal{C}_{T,j}^{\star}
&\leq Ce^{-(1 + a)\beta_{c}(\theta)j}\sum_{B \in \mathcal{B}_{T}}\ell_{T}^{-1}\int_{B}\PP_{\star}(F_{h_{0}}^{\star} \in [u_{T} - j - 1,u_{T} - j])\rd{h_{0}} \\
&\leq Ce^{-(1 + a)\beta_{c}(\theta)j}N_{T}\sup_{h_{0} \in I_{\theta}}\PP_{\star}(F_{h_{0}}^{\star} \in [u_{T} - j - 1,u_{T} - j]) \\
&\leq C\exp\left\{ - a\beta_{c}(\theta)j + \left(\frac{u_{T}}{V_{T}} - \beta_{c}(\theta)\right)j - \frac{j^{2}}{2V_{T}}\right\}.
\end{aligned}
\end{equation}
Equation \eqref{an:height.scales} in Lemma~\ref{an:height} gives $u_{T}/V_{T} \to \beta_{c}(\theta)$, so $u_{T}/V_{T} - \beta_{c}(\theta) \leq a\beta_{c}(\theta)/2$ for all sufficiently large $T$. Since $ - j^{2}/(2V_{T}) \leq 0$, the right-hand side of \eqref{temp:detail:critical.depth.sum} is then at most $Ce^{-a\beta_{c}(\theta)j/2}$, uniformly in the allowed indices $j$. For each fixed $j$, its limiting upper value is at most $Ce^{-a\beta_{c}(\theta)j}$.

Only roots with $F_{h_{0}}^{\star} \leq u_{T} - L$ contribute, so the bin $[u_{T} - j - 1,u_{T} - j]$ can contribute only if $u_{T} - j - 1 < u_{T} - L$, equivalently $j \geq \lfloor L\rfloor$, apart from endpoints of zero probability. The additional restriction $u_{T} - j - 1 \geq m_{0}s_{T}$ ensures $\gamma \geq m_{0}$ and hence the uniform conditional bound \eqref{temp:detail:critical.conditional.unit}, by \eqref{temp:detail:critical.annular.exponent}. Thus the bins contributing to this part of the proof satisfy both conditions $m_{0}s_{T} \leq u_{T} - j - 1 < u_{T} - L$.

To apply domination on the fixed index set $j \geq \lfloor L\rfloor$, extend $\mathcal{C}_{T,j}^{\star}$ by zero when $u_{T} - j - 1 < m_{0}s_{T}$. This extension records only the contributions from bins with lower endpoint at least $m_{0}s_{T}$; the contributions of the remaining root heights, including the possible bin crossing $m_{0}s_{T}$, are estimated in \eqref{temp:detail:critical.lower.sum} and \eqref{temp:detail:critical.negative.roots}. Domination by the summable sequence $\smash{Ce^{-a\beta_{c}(\theta)j/2}}$ therefore gives
\begin{equation}\label{temp:detail:critical.high.roots}
\limsup_{T \to \infty}\sum_{\substack{j \geq \lfloor L\rfloor \\ u_{T} - j - 1 \geq m_{0}s_{T}}}\mathcal{C}_{T,j}^{\star} \leq C\sum_{j = \lfloor L\rfloor}^{\infty}e^{-a\beta_{c}(\theta)j} \leq C'e^{-a\beta_{c}(\theta)L}.
\end{equation}
The constant $C'$ is independent of $L$.

The smallest lower endpoint among the bins treated in \eqref{temp:detail:critical.high.roots} lies in $[m_{0}s_{T},m_{0}s_{T} + 1)$ for all sufficiently large $T$. Hence any root $h_{0}$ with $0 \leq F_{h_{0}}^{\star} \leq u_{T} - L$ whose field value lies outside the bins treated in \eqref{temp:detail:critical.high.roots} must satisfy $F_{h_{0}}^{\star} < m_{0}s_{T} + 1$. The interval $[0,m_{0}s_{T} + 1)$ can be covered by $\OO(s_{T})$ unit bins $[\gamma s_{T},\gamma s_{T} + 1]$ with $0 \leq \gamma \leq m_{0} + s_{T}^{-1}$. The additional $s_{T}^{-1}$ in $\gamma$ allows a height enlargement of at most one. Bounding the contribution of each full covering bin also bounds the contribution of the remaining roots within that bin, even if part of the covering bin overlaps a previously treated bin.

To estimate these contributions, first establish a uniform exponential bound on scalar bins. For $\gamma$ in any fixed compact subset of $[0,\infty)$, apply exponential Markov on $\mathcal{A}_{T,\mathrm{sp}}$ with tilt $2\gamma$. The quadratic cumulant estimate \eqref{an:joint.cumulants} with the bounded remainder \eqref{an:conditional.detail:remainder}, the truncated transform comparison \eqref{an:spatial.transform.bound} in Lemma~\ref{an:spatial.transform}, and $V_{T} = s_{T}/2 + \OO(1)$ from \eqref{an:detail:cov.variance} give
\begin{equation}\label{temp:detail:scalar.markov}
\begin{aligned}
\PP_{\star}(\{F_{h}^{\star} \geq \gamma s_{T}\} \cap \mathcal{A}_{T,\mathrm{sp}})
&\leq e^{-2\gamma^{2}s_{T}}\EE_{\star}[e^{2\gamma F_{h}^{\star}}\ind_{\mathcal{A}_{T,\mathrm{sp}}}] \\
&\leq e^{-2\gamma^{2}s_{T}}\left(\exp\{2\gamma^{2}V_{T} + \OO(1)\} + \exp\{ - ct^{1 + 2\e_{0}} + Ct\}\right) \\
&\leq \exp\{ - \gamma^{2}s_{T} + \OO(1)\}.
\end{aligned}
\end{equation}
For the Haar law, the transform comparison error in \eqref{temp:detail:scalar.markov} can be omitted. For the arithmetic law, adding $\PP(\mathcal{A}_{T,\mathrm{sp}}^{c})$ removes the event $\mathcal{A}_{T,\mathrm{sp}}$ from the probability on the left-hand side of \eqref{temp:detail:scalar.markov}: the bound \eqref{an:spatial.exceptional} in Lemma~\ref{an:spatial.transform} is smaller than $\exp( - Cs_{T})$ for every fixed $C > 0$. Since $\log{N_{T}} = m(\theta)^{2}s_{T} + \oo(s_{T})$ by \eqref{eq:effective.cutoff.bounds} and \eqref{eq:parameter.scales}, \eqref{temp:detail:scalar.markov} consequently implies
\begin{equation}\label{temp:detail:bin.exponent}
N_{T}\sup_{h \in I_{\theta}}\PP_{\star}(F_{h}^{\star} \in [\gamma s_{T},\gamma s_{T} + 1]) \leq \exp\{(m(\theta)^{2} - \gamma^{2})s_{T} + \oo(s_{T})\}.
\end{equation}
The remainder is uniform in $\gamma$ on every such compact set.

For a fixed covering bin $[\gamma s_{T},\gamma s_{T} + 1]$, consider the expected sum of the root integrals in \eqref{temp:detail:root.factorization}, restricted to this bin and to the event $\mathcal{A}_{T,\mathrm{sp}}$. Its contribution satisfies
\begin{align}\label{temp:detail:critical.lower.bin}
&\EE_{\star}\Bigg[\sum_{B \in \mathcal{B}_{T}}\ell_{T}^{-1}\int_{B}\exp\{\beta_{c}(\theta)(F_{h_{0}}^{\star} - u_{T})\}Z_{B}^{a} \ind_{\{F_{h_{0}}^{\star}\in [\gamma s_{T},\gamma s_{T} + 1],\,F_{h_{0}}^{\star}\leq u_{T} - L\}}\ind_{\mathcal{A}_{T,\mathrm{sp}}}\rd{h_{0}}\Bigg] \notag \\
&\qquad\leq e^{(1 + a)\beta_{c}(\theta)(\gamma s_{T} + 1 - u_{T})}\sum_{B \in \mathcal{B}_{T}}\ell_{T}^{-1}\int_{B}\PP_{\star}(F_{h_{0}}^{\star} \in [\gamma s_{T},\gamma s_{T} + 1]) \notag \\
&\qquad\qquad \times \EE_{\star}\left[\left(\ell_{T}^{-1}\int_{B}e^{\beta_{c}(\theta)(F_{h}^{\star} - F_{h_{0}}^{\star})}\rd{h}\right)^{a}\ind_{\mathcal{A}_{T,\mathrm{sp}}} \mid F_{h_{0}}^{\star} \in [\gamma s_{T},\gamma s_{T} + 1]\right]\rd{h_{0}} \notag \\
&\qquad\leq C(1 + s_{T})e^{(1 + a)\beta_{c}(\theta)(\gamma s_{T} + 1 - u_{T})}\exp\left\{\left[a\beta_{c}(\theta)(1 - \gamma) + \frac{a^{2}\beta_{c}(\theta)^{2}}{4}\right]_{+}s_{T}\right\} \notag \\
&\qquad\qquad \times \sum_{B \in \mathcal{B}_{T}}\ell_{T}^{-1}\int_{B}\PP_{\star}(F_{h_{0}}^{\star} \in [\gamma s_{T},\gamma s_{T} + 1])\rd{h_{0}} \notag \\
&\qquad\leq C(1 + s_{T})e^{(1 + a)\beta_{c}(\theta)(\gamma s_{T} + 1 - u_{T})}N_{T}\sup_{h_{0} \in I_{\theta}}\PP_{\star}(F_{h_{0}}^{\star} \in [\gamma s_{T},\gamma s_{T} + 1]) \notag \\
&\qquad\qquad \times \exp\left\{\left[a\beta_{c}(\theta)(1 - \gamma) + \frac{a^{2}\beta_{c}(\theta)^{2}}{4}\right]_{+}s_{T}\right\}.
\end{align}
The first inequality in \eqref{temp:detail:critical.lower.bin} follows from \eqref{temp:detail:root.relative.bound}, the bound $F_{h_{0}}^{\star} \leq \gamma s_{T} + 1$ on the bin, and conditioning on the root-height bin after dropping the restriction $F_{h_{0}}^{\star} \leq u_{T} - L$. The second inequality uses \eqref{an:unit.fractional} in Lemma~\ref{an:annular}, since $B \subseteq [h_{0} - 1,h_{0} + 1] \cap I_{\theta}$. The last inequality bounds the root-height probability uniformly in $h_{0}$ and uses $\sum_{B \in \mathcal{B}_{T}}\ell_{T}^{-1}|B| = N_{T}$ from \eqref{eq:N.T}.

Substituting \eqref{temp:detail:bin.exponent} into \eqref{temp:detail:critical.lower.bin} and using $u_{T}/s_{T} \to m(\theta)$ from \eqref{an:height.scales} in Lemma~\ref{an:height} gives an upper bound of $C(1 + s_{T})$ multiplied by an exponential whose exponent is
\begin{equation}\label{temp:detail:critical.lower.exponent}
\left\{m(\theta)^{2} - \gamma^{2} + (1 + a)\beta_{c}(\theta)(\gamma - m(\theta)) + \left[a\beta_{c}(\theta)(1 - \gamma) + \frac{a^{2}\beta_{c}(\theta)^{2}}{4}\right]_{+} + \oo(1)\right\}s_{T}.
\end{equation}
The $\oo(1)$ term is uniform over the covering bins. First suppose that $0 \leq \gamma \leq m_{0}$. If the positive part in \eqref{temp:detail:critical.lower.exponent} is zero, the expression in braces, excluding the $\oo(1)$ term, equals
\begin{equation}\label{temp:detail:critical.zero.branch}
 - (m(\theta) - \gamma)^{2} - 2am(\theta)(m(\theta) - \gamma) \leq - (m(\theta) - m_{0})^{2},
\end{equation}
because $\beta_{c}(\theta) = 2m(\theta)$, $a > 0$, and $m(\theta) - \gamma \geq m(\theta) - m_{0} > 0$. If the positive part in \eqref{temp:detail:critical.lower.exponent} is nonzero, the expression in braces, excluding the $\oo(1)$ term, equals
\begin{equation}\label{temp:detail:critical.positive.branch}
 - (m(\theta) - \gamma)^{2} - 2am(\theta)(m(\theta) - 1) + a^{2}m(\theta)^{2} \leq - am(\theta)(m(\theta) - 1),
\end{equation}
where the last inequality uses \eqref{temp:detail:critical.parameters}. Equations \eqref{temp:detail:critical.zero.branch} and \eqref{temp:detail:critical.positive.branch} therefore bound the expression in braces in \eqref{temp:detail:critical.lower.exponent} by $ - c_{0} + \oo(1)$, where $c_{0} \leqdef \min\{(m(\theta) - m_{0})^{2},am(\theta)(m(\theta) - 1)\} > 0$. Allowing $m_{0} < \gamma \leq m_{0} + s_{T}^{-1}$ changes the terms other than $\oo(1)$ in these braces by $\OO(s_{T}^{-1})$, uniformly over this enlarged range, so the bound $ - c_{0} + \oo(1)$ remains valid. Since there are $\OO(s_{T})$ bins and the prefactor for each bin is $C(1 + s_{T})$, summing the left-hand side of \eqref{temp:detail:critical.lower.bin} over these bins bounds $\mathcal{C}_{T,\mathrm{low}}^{\star}$ in \eqref{temp:detail:critical.other.contributions} and gives
\begin{equation}\label{temp:detail:critical.lower.sum}
\mathcal{C}_{T,\mathrm{low}}^{\star} \leq Cs_{T}(1 + s_{T})\exp\{ - c_{0}s_{T} + \oo(s_{T})\} \leq C\exp( - c_{0}s_{T}/2).
\end{equation}
The error in the exponent is uniform over the bins by \eqref{temp:detail:bin.exponent} and the uniformity of the preceding estimates.

Finally, to bound $\mathcal{C}_{T,\mathrm{neg}}^{\star}$, which corresponds to roots with $F_{h_{0}}^{\star} < 0$, retain the deterministic cap $Z_{B}^{a} \leq e^{as_{T}}$ instead of using a conditional estimate. Indeed, since $L \geq 0$, the exponential factor $\smash{e^{\beta_{c}(\theta)(F_{h}^{\star} - u_{T})}}$ is at most one on $\{F_{h}^{\star} \leq u_{T} - L\}$. Thus \eqref{temp:interval.masses} gives $Z_{B} \leq \ell_{T}^{-1}|B| \leq \ell_{T}^{-1} = e^{s_{T}}$, using $|B| \leq 1$ and \eqref{eq:N.T}, and raising this bound to the power $a > 0$ proves the claimed cap. The definition of $\mathcal{C}_{T,\mathrm{neg}}^{\star}$ in \eqref{temp:detail:critical.other.contributions} then gives
\begin{equation}\label{temp:detail:critical.negative.roots}
\mathcal{C}_{T,\mathrm{neg}}^{\star}
\leq e^{as_{T}}\ell_{T}^{-1}\int_{I_{\theta}}e^{-\beta_{c}(\theta)u_{T}}\rd{h_{0}} = N_{T}e^{-\beta_{c}(\theta)u_{T} + as_{T}}
= \exp\{( - m(\theta)^{2} + a)s_{T} + \oo(s_{T})\}.
\end{equation}
Here $a < 1 < m(\theta)^{2}$, $u_{T}/s_{T} \to m(\theta)$ by \eqref{an:height.scales} in Lemma~\ref{an:height}, and $\log{N_{T}} = m(\theta)^{2}s_{T} + \oo(s_{T})$ by \eqref{eq:effective.cutoff.bounds} and \eqref{eq:parameter.scales}. Adding the high-root sum \eqref{temp:detail:critical.high.roots}, the lower-root bounds \eqref{temp:detail:critical.lower.sum} and \eqref{temp:detail:critical.negative.roots}, and the exceptional contribution \eqref{temp:detail:cap.exceptional} proves \eqref{temp:critical.moment.tail}.
\end{proof}

\subsubsection{Mass spreading and lower mass bounds}

\begin{proposition}[Mass spreading at critical temperature]\label{temp:critical}
Let $\theta > 0$. For $\star \in \{\mathrm{A},\mathrm{H}\}$, as $T \to \infty$,
\begin{equation}\label{temp:critical.mass}
M_{\beta_{c}(\theta),T}^{\star} \to \infty \quad\text{in probability}, \qquad \max_{B \in \mathcal{B}_{T}}Z_{B,\beta_{c}(\theta)}^{\star}(\R) = \OO_{\PP_{\star}}(1).
\end{equation}
Consequently,
\begin{equation}\label{temp:critical.spreading}
\max_{B \in \mathcal{B}_{T}}G_{\beta_{c}(\theta),T}^{\star}(B) = \frac{\max_{B}Z_{B,\beta_{c}(\theta)}^{\star}(\R)}{M_{\beta_{c}(\theta),T}^{\star}} \to 0 \quad\text{in probability}.
\end{equation}
\end{proposition}

\begin{proof}[Proof of Proposition~\ref{temp:critical}]
First obtain a lower bound for the total mass $\smash{M_{\beta_{c}(\theta),T}^{\star}}$ from the height window $F_{h}^{\star} - u_{T} \in [ - L,0]$. Truncating the mass $Z_{B,\beta_{c}(\theta)}^{\star}([ - L,0])$ of each interval $B \in \mathcal{B}_{T}$ controls the contribution of nearby intervals to the second moment, while \eqref{an:remote.bound} in Lemma~\ref{an:remote} controls the distant pairs. Fix $L > 0$ and define
\begin{equation}\label{temp:detail:window.definition}
\begin{aligned}
Z_{B} &\leqdef Z_{B,\beta_{c}(\theta)}^{\star}([ - L,0]), \\
M_{T,L}^{\mathrm{win}}
&\leqdef \ell_{T}^{-1}\int_{I_{\theta}}e^{\beta_{c}(\theta)(F_{h}^{\star} - u_{T})}\ind_{\{F_{h}^{\star} - u_{T} \in [ - L,0]\}}\rd{h}
= \sum_{B \in \mathcal{B}_{T}}Z_{B}.
\end{aligned}
\end{equation}
Let $f(x) \leqdef e^{\beta_{c}(\theta)x}\ind_{[ - L,0]}(x)$. This function is bounded and compactly supported, and its only discontinuities, at $ - L$ and $0$, have zero mass under the limiting measure in \eqref{an:height.measure} of Lemma~\ref{an:height}. Integrating $f$ against the measures in that weak convergence therefore gives, uniformly in $h \in I_{\theta}$,
\[
N_{T}\EE_{\star}[f(F_{h}^{\star} - u_{T})] \to \int_{\R}f(x)\beta_{c}(\theta)e^{-\beta_{c}(\theta)x}\rd{x} = \beta_{c}(\theta)L.
\]
Averaging this convergence over $h \in I_{\theta}$ and using $N_{T} = \ell_{T}^{-1}|I_{\theta}|$ from \eqref{eq:N.T} gives the first line below. The second follows from \eqref{temp:critical.moment.tail} in Lemma~\ref{temp:critical.moment} with $L = 0$:
\begin{equation}\label{temp:detail:window.moments}
\begin{aligned}
\EE_{\star}[M_{T,L}^{\mathrm{win}}]
&= \frac{1}{|I_{\theta}|}\int_{I_{\theta}}N_{T}\EE_{\star}[f(F_{h}^{\star} - u_{T})]\rd{h} \to \beta_{c}(\theta)L, \\
\limsup_{T \to \infty}\sum_{B}\EE_{\star}[Z_{B}^{1 + a}]
&\leq \limsup_{T \to \infty}\sum_{B}\EE_{\star}[Z_{B,\beta_{c}(\theta)}^{\star}(( - \infty,0])^{1 + a}] \leq C.
\end{aligned}
\end{equation}
The constant in the second estimate of \eqref{temp:detail:window.moments} is independent of $L$.

For distant intervals $B,B' \in \mathcal{B}_{T}$, the Cartesian product $B \times B'$ is contained in $\mathscr{D}_{2,T}(2) = \{(h,h') \in I_{\theta}^{\times 2}: |h - h'| \geq 2\}$. Thus the integrated estimate \eqref{an:remote.bound} in Lemma~\ref{an:remote}, with $k = 2$, $d_{0} = 2$, and the same function $f$, gives
\begin{equation}\label{temp:detail:window.distant}
\begin{aligned}
\limsup_{T \to \infty}\sum_{(B,B')\,\mathrm{distant}}\EE_{\star}[Z_{B}Z_{B'}]
&\leq \limsup_{T \to \infty}\EE_{\star}\left[\ell_{T}^{-2}\int_{\mathscr{D}_{2,T}(2)}f(F_{h}^{\star} - u_{T})f(F_{h'}^{\star} - u_{T})\rd{h}\rd{h'}\right] \\
&\leq (\beta_{c}(\theta)L)^{2}.
\end{aligned}
\end{equation}

Let $\widetilde{Z}_{B} \leqdef Z_{B}\ind_{\{Z_{B}\leq 1\}}$. On $\{Z_{B} > 1\}$, the inequality $Z_{B} \leq Z_{B}^{1 + a}$ and \eqref{temp:detail:window.moments} imply
\begin{equation}\label{temp:detail:window.lost}
\limsup_{T \to \infty}\EE_{\star}\left[M_{T,L}^{\mathrm{win}} - \sum_{B}\widetilde{Z}_{B}\right]
\leq \limsup_{T \to \infty}\sum_{B \in \mathcal{B}_{T}}\EE_{\star}[Z_{B}^{1 + a}]
\leq C.
\end{equation}
To control the second moment, split the ordered interval pairs into distant pairs and the remaining pairs, including diagonal pairs. Each interval occurs in only a fixed number of remaining pairs. The inequalities $2xy \leq x^{2} + y^{2}$ and $\widetilde{Z}_{B}^{2} \leq \widetilde{Z}_{B}$ give
\begin{equation}\label{temp:detail:window.near}
\begin{aligned}
\sum_{(B,B')\,\mathrm{not\ distant}}\EE_{\star}[\widetilde{Z}_{B}\widetilde{Z}_{B'}]
&\leq \frac12\sum_{(B,B')\,\mathrm{not\ distant}}\left(\EE_{\star}[\widetilde{Z}_{B}^{2}] + \EE_{\star}[\widetilde{Z}_{B'}^{2}]\right) \\
&\leq C\sum_{B}\EE_{\star}[\widetilde{Z}_{B}^{2}] \leq C\EE_{\star}[M_{T,L}^{\mathrm{win}}].
\end{aligned}
\end{equation}
For the distant pairs, use $\widetilde{Z}_{B}\widetilde{Z}_{B'} \leq Z_{B}Z_{B'}$ and \eqref{temp:detail:window.distant}. Combining \eqref{temp:detail:window.moments}, \eqref{temp:detail:window.lost}, and \eqref{temp:detail:window.near} yields, whenever $\beta_{c}(\theta)L > C$,
\begin{equation}\label{temp:detail:window.variance}
\liminf_{T \to \infty}\EE_{\star}\left[\sum_{B}\widetilde{Z}_{B}\right] \geq \beta_{c}(\theta)L - C, \qquad
\limsup_{T \to \infty}\Var_{\star}\left(\sum_{B}\widetilde{Z}_{B}\right) \leq 3C\beta_{c}(\theta)L.
\end{equation}
For the variance bound in \eqref{temp:detail:window.variance}, the distant and nearby pair bounds \eqref{temp:detail:window.distant} and \eqref{temp:detail:window.near} give
\[
\limsup_{T \to \infty}\Var_{\star}\left(\sum_{B}\widetilde{Z}_{B}\right)
\leq (\beta_{c}(\theta)L)^{2} + C\beta_{c}(\theta)L - (\beta_{c}(\theta)L - C)^{2}
\leq 3C\beta_{c}(\theta)L.
\]
The condition $\beta_{c}(\theta)L > C$ ensures that the lower mean bound can be squared; the last line drops the nonpositive term $ - C^{2}$.

Assume that $\beta_{c}(\theta)L/2 > C$. The first bound in \eqref{temp:detail:window.variance} then ensures that $\EE_{\star}[\sum_{B}\widetilde{Z}_{B}] > \beta_{c}(\theta)L/2$ for all sufficiently large $T$. Since $\smash{M_{T,L}^{\mathrm{win}} \geq \sum_{B}\widetilde{Z}_{B}}$, Chebyshev's inequality and the two bounds in \eqref{temp:detail:window.variance} give
\begin{equation}\label{temp:detail:critical.chebyshev}
\begin{aligned}
&\limsup_{T \to \infty}\PP_{\star}(M_{T,L}^{\mathrm{win}} \leq \beta_{c}(\theta)L/2) \\
&\qquad\leq \limsup_{T \to \infty}\frac{\Var_{\star}(\sum_{B}\widetilde{Z}_{B})}{(\EE_{\star}[\sum_{B}\widetilde{Z}_{B}] - \beta_{c}(\theta)L/2)^{2}}
\leq \frac{3C\beta_{c}(\theta)L}{(\beta_{c}(\theta)L/2 - C)^{2}}.
\end{aligned}
\end{equation}
For each such fixed $L$, first let $T \to \infty$; the upper bound in \eqref{temp:detail:critical.chebyshev} is then $\OO(L^{-1})$ as $L \to \infty$. Hence
\begin{equation}\label{eq:temp:window.mass}
\lim_{L \to \infty}\limsup_{T \to \infty}\PP_{\star}(M_{T,L}^{\mathrm{win}} \leq \beta_{c}(\theta)L/2) = 0.
\end{equation}
For every fixed $R > 0$, choose $L$ with $\beta_{c}(\theta)L/2 > R$. Since $\smash{M_{\beta_{c}(\theta),T}^{\star} \geq M_{T,L}^{\mathrm{win}}}$, \eqref{eq:temp:window.mass} gives $\PP_{\star}(M_{\beta_{c}(\theta),T}^{\star} \leq R) \to 0$ after letting $L \to \infty$. This proves that, as $T \to \infty$,
\begin{equation}\label{temp:detail:critical.mass.divergence}
M_{\beta_{c}(\theta),T}^{\star} \to \infty \quad\text{in probability},
\end{equation}
as claimed in \eqref{temp:critical.mass}.

It remains to prove tightness of the largest unit-interval mass $\smash{\max_{B \in \mathcal{B}_{T}}Z_{B,\beta_{c}(\theta)}^{\star}(\R)}$ in \eqref{temp:critical.mass}. The bound \eqref{temp:critical.moment.tail} in Lemma~\ref{temp:critical.moment} with $L = 0$, followed by a union bound and Markov's inequality, gives, for every $R > 0$,
\begin{equation}\label{temp:detail:critical.lower.maximum}
\begin{aligned}
\limsup_{T \to \infty}\PP_{\star}\left(\max_{B}Z_{B,\beta_{c}(\theta)}^{\star}(( - \infty,0]) > R\right)
&\leq R^{-1 - a}\limsup_{T \to \infty}\sum_{B \in \mathcal{B}_{T}}\EE_{\star}\left[Z_{B,\beta_{c}(\theta)}^{\star}(( - \infty,0])^{1 + a}\right] \\
&\leq CR^{-1 - a}.
\end{aligned}
\end{equation}
For each fixed $K > 0$, \eqref{an:height.measure} in Lemma~\ref{an:height} also gives
\begin{equation}\label{temp:detail:critical.upper.window}
\EE_{\star}\left[\sum_{B}Z_{B,\beta_{c}(\theta)}^{\star}([0,K])\right] \to \beta_{c}(\theta)K.
\end{equation}
On the height-cap event $\mathcal{H}_{T,K}^{\star} \leqdef \{\sup_{h \in I_{\theta}}(F_{h}^{\star} - u_{T}) \leq K\}$, the full interval masses satisfy
\begin{equation}\label{temp:detail:critical.full.decomposition}
\max_{B \in \mathcal{B}_{T}}Z_{B,\beta_{c}(\theta)}^{\star}(\R)
\leq \max_{B \in \mathcal{B}_{T}}Z_{B,\beta_{c}(\theta)}^{\star}(( - \infty,0]) + \sum_{B \in \mathcal{B}_{T}}Z_{B,\beta_{c}(\theta)}^{\star}([0,K]).
\end{equation}
Apply \eqref{temp:detail:critical.lower.maximum} to the first term in \eqref{temp:detail:critical.full.decomposition}, Markov's inequality with \eqref{temp:detail:critical.upper.window} to the second term, and the second bound in \eqref{an:height.tails} of Lemma~\ref{an:height} to the complementary event $(\mathcal{H}_{T,K}^{\star})^{c}$ to obtain
\begin{equation}\label{temp:detail:critical.full.maximum}
\limsup_{T \to \infty}\PP_{\star}\left(\max_{B}Z_{B,\beta_{c}(\theta)}^{\star}(\R) > 2R\right) \leq CR^{-1 - a} + \frac{\beta_{c}(\theta)K}{R} + Ce^{-\beta_{c}(\theta)K}.
\end{equation}
In \eqref{temp:detail:critical.full.maximum}, first let $R \to \infty$ with $K$ fixed, and then let $K \to \infty$. This proves the tightness assertion in \eqref{temp:critical.mass}. Finally, for every $\e > 0$ and $R > 0$, this ratio satisfies
\[
\PP_{\star}\left(\max_{B}G_{\beta_{c}(\theta),T}^{\star}(B) > \e\right)
\leq \PP_{\star}\left(\max_{B}Z_{B,\beta_{c}(\theta)}^{\star}(\R) > R\right) + \PP_{\star}\left(M_{\beta_{c}(\theta),T}^{\star} \leq R/\e\right).
\]
Letting $T \to \infty$ and then $R \to \infty$, using \eqref{temp:critical.mass}, proves the assertion in \eqref{temp:critical.spreading}.
\end{proof}

\begin{corollary}[Lower mass bound at low temperature]\label{temp:lower.mass}
Let $\theta > 0$. For every $\beta > \beta_{c}(\theta)$ and $\star \in \{\mathrm{A},\mathrm{H}\}$,
\[
\lim_{\delta \downarrow 0}\limsup_{T \to \infty}\PP_{\star}(M_{\beta,T}^{\star} < \delta) = 0.
\]
\end{corollary}

\begin{proof}[Proof of Corollary~\ref{temp:lower.mass}]
Fix $L > 0$. Since $\beta - \beta_{c}(\theta) > 0$, every $x \in [ - L,0]$ satisfies
\begin{equation}\label{temp:detail:lower.height.factor}
e^{\beta x} = e^{(\beta - \beta_{c}(\theta))x}e^{\beta_{c}(\theta)x} \geq e^{-(\beta - \beta_{c}(\theta))L}e^{\beta_{c}(\theta)x}.
\end{equation}
In the integral defining $M_{\beta,T}^{\star}$, retain only the positions $h \in I_{\theta}$ satisfying $F_{h}^{\star} - u_{T} \in [ - L,0]$, and apply \eqref{temp:detail:lower.height.factor} with $x = F_{h}^{\star} - u_{T}$ to obtain
\begin{equation}\label{temp:detail:lower.window}
M_{\beta,T}^{\star}
= \ell_{T}^{-1}\int_{I_{\theta}}e^{\beta(F_{h}^{\star} - u_{T})}\rd{h}
\geq e^{-(\beta - \beta_{c}(\theta))L}M_{T,L}^{\mathrm{win}},
\end{equation}
where $M_{T,L}^{\mathrm{win}}$ is defined in \eqref{temp:detail:window.definition}. Given $\e > 0$, choose a fixed $L$ large enough that \eqref{eq:temp:window.mass} gives
\begin{equation}\label{temp:detail:lower.window.probability}
\limsup_{T \to \infty}\PP_{\star}(M_{T,L}^{\mathrm{win}} \leq \beta_{c}(\theta)L/2) \leq \e.
\end{equation}
For every $0 < \delta < (\beta_{c}(\theta)L/2)e^{-(\beta - \beta_{c}(\theta))L}$, \eqref{temp:detail:lower.window} implies the inclusion
\[
\{M_{\beta,T}^{\star} < \delta\} \subseteq \{M_{T,L}^{\mathrm{win}} < \beta_{c}(\theta)L/2\}.
\]
The inclusion and \eqref{temp:detail:lower.window.probability} give
\[
\lim_{\delta \downarrow 0}\limsup_{T \to \infty}\PP_{\star}(M_{\beta,T}^{\star} < \delta)
\leq \limsup_{T \to \infty}\PP_{\star}(M_{T,L}^{\mathrm{win}} \leq \beta_{c}(\theta)L/2)
\leq \e.
\]
Here $L$ is fixed before $T \to \infty$ and $\delta \downarrow 0$. Letting $\e \downarrow 0$ proves the assertion.
\end{proof}

\subsubsection{Mass spreading and vanishing off-diagonal overlaps}

\begin{proposition}[Mass spreading at high temperature]\label{temp:high}
Let $\theta > 0$. For every $0 < \beta < \beta_{c}(\theta)$ and $\star \in \{\mathrm{A},\mathrm{H}\}$, as $T \to \infty$,
\[
\max_{B \in \mathcal{B}_{T}}G_{\beta,T}^{\star}(B) = \frac{\max_{B}Z_{B,\beta}^{\star}(\R)}{M_{\beta,T}^{\star}} \to 0 \quad\text{in probability}.
\]
\end{proposition}

\begin{proof}[Proof of Proposition~\ref{temp:high}]
Fix $K > 0$ and abbreviate $Z_{B} \leqdef Z_{B,\beta}^{\star}(( - \infty,K])$. The proof estimates $\EE_{\star}[\sum_{B \in \mathcal{B}_{T}}Z_{B}]$, the mean total capped mass, and compares the $(1 + a)$ moments of the interval masses $Z_{B}$ with the corresponding power of its asymptotic scale. The probability decomposition \eqref{temp:detail:high.probability.decomposition} reduces the assertion to concentration of the total capped mass and a bound on the largest capped interval mass. The moment decomposition \eqref{temp:detail:high.moment.decomposition} separates the nonnegative root-height bins, the negative roots, and the arithmetic exceptional event that enter the required moment estimate.

Define
\begin{equation}\label{temp:detail:high.scale.definition}
A_{T}^{\mathrm{res}}(\beta) \leqdef N_{T}e^{-\beta u_{T}}\EE_{\mathrm{H}}[e^{\beta H(0)}].
\end{equation}
On the event $\{\sup_{h \in I_{\theta}}(F_{h}^{\star} - u_{T}) \leq K\}$, the capped masses $Z_{B}$ equal the full masses $Z_{B,\beta}^{\star}(\R)$. If also $\sum_{B}Z_{B} \geq A_{T}^{\mathrm{res}}(\beta)/2$ and $\max_{B}Z_{B} \leq \e A_{T}^{\mathrm{res}}(\beta)/2$, then $\max_{B}G_{\beta,T}^{\star}(B) \leq \e$. Hence, for every $\e > 0$,
\begin{equation}\label{temp:detail:high.probability.decomposition}
\begin{aligned}
\PP_{\star}\left(\max_{B \in \mathcal{B}_{T}}G_{\beta,T}^{\star}(B) > \e\right)
&\leq \PP_{\star}\left(\sup_{h \in I_{\theta}}(F_{h}^{\star} - u_{T}) > K\right)
+ \PP_{\star}\left(\sum_{B \in \mathcal{B}_{T}}Z_{B} < A_{T}^{\mathrm{res}}(\beta)/2\right) \\
&\qquad + \PP_{\star}\left(\max_{B \in \mathcal{B}_{T}}Z_{B} > \e A_{T}^{\mathrm{res}}(\beta)/2\right).
\end{aligned}
\end{equation}
The first probability is controlled by the second estimate in \eqref{an:height.tails} of Lemma~\ref{an:height}. The other two probabilities will tend to zero for each fixed $K$, by \eqref{temp:detail:high.total} and \eqref{temp:detail:high.maximum}, respectively.

Initially choose $0 < a < 1 \wedge (2/\beta)$, with $a$ to be decreased when comparing the exponential rates. For a fixed $0 < \e_{0} < 1/4$, use the common arithmetic event $\smash{\mathcal{A}_{T,\mathrm{sp}}^{\mathrm{A}}}$ defined in \eqref{an:annular.unit.spatial.event}, which includes the dyadic cutoffs used in \eqref{an:unit.fractional} of Lemma~\ref{an:annular}. Write $\mathcal{A}_{T,\mathrm{sp}}$ for this event under the arithmetic law and for the whole probability space under the Haar law. Multiplying the definition of $Z_{B}$ in \eqref{temp:interval.masses} by $Z_{B}^{a}$ gives
\begin{equation}\label{temp:detail:high.root.factorization}
Z_{B}^{1 + a} = \ell_{T}^{-1}\int_{B}e^{\beta(F_{h_{0}}^{\star} - u_{T})}\ind_{\{F_{h_{0}}^{\star}\leq u_{T} + K\}}Z_{B}^{a}\rd{h_{0}}.
\end{equation}
For each integer $0 \leq j \leq \lfloor u_{T} + K\rfloor$, define the contribution of the nonnegative root-height bin $[j,j + 1)$ by
\begin{equation}\label{temp:detail:high.bin.contribution}
\begin{aligned}
\widetilde{\mathcal{C}}_{T,j}^{\star}
&\leqdef \EE_{\star}\Bigg[\sum_{B \in \mathcal{B}_{T}}\ell_{T}^{-1}\int_{B}e^{\beta(F_{h_{0}}^{\star} - u_{T})}Z_{B}^{a} \ind_{\{F_{h_{0}}^{\star}\in [j,j + 1),\,F_{h_{0}}^{\star}\leq u_{T} + K\}}\ind_{\mathcal{A}_{T,\mathrm{sp}}}\rd{h_{0}}\Bigg].
\end{aligned}
\end{equation}
Define the negative-root and exceptional contributions by
\begin{equation}\label{temp:detail:high.other.contributions}
\begin{aligned}
\widetilde{\mathcal{C}}_{T,\mathrm{neg}}^{\star}
&\leqdef \EE_{\star}\Bigg[\sum_{B \in \mathcal{B}_{T}}\ell_{T}^{-1}\int_{B}e^{\beta(F_{h_{0}}^{\star} - u_{T})}Z_{B}^{a} \ind_{\{F_{h_{0}}^{\star} < 0,\,F_{h_{0}}^{\star}\leq u_{T} + K\}}\ind_{\mathcal{A}_{T,\mathrm{sp}}}\rd{h_{0}}\Bigg], \\
\widetilde{\mathcal{C}}_{T,\mathrm{exc}}^{\star}
&\leqdef \EE_{\star}\left[\sum_{B \in \mathcal{B}_{T}}Z_{B}^{1 + a}\ind_{\mathcal{A}_{T,\mathrm{sp}}^{c}}\right].
\end{aligned}
\end{equation}
For all sufficiently large $T$, $u_{T} + K > 0$, and each height $0 \leq F_{h_{0}}^{\star} \leq u_{T} + K$ lies in exactly one half-open bin in \eqref{temp:detail:high.bin.contribution}, with the height cutoff included in the integrand. Thus \eqref{temp:detail:high.root.factorization} gives the exact decomposition
\begin{equation}\label{temp:detail:high.moment.decomposition}
\EE_{\star}\left[\sum_{B \in \mathcal{B}_{T}}Z_{B}^{1 + a}\right]
= \sum_{j = 0}^{\lfloor u_{T} + K\rfloor}\widetilde{\mathcal{C}}_{T,j}^{\star} + \widetilde{\mathcal{C}}_{T,\mathrm{neg}}^{\star} + \widetilde{\mathcal{C}}_{T,\mathrm{exc}}^{\star}.
\end{equation}
The contributions of the nonnegative root-height bins $[j,j + 1)$ and of roots with $F_{h_{0}}^{\star} < 0$ will be bounded together in \eqref{temp:detail:high.retained.moment}, using \eqref{temp:detail:high.bin.bound} and \eqref{temp:detail:high.negative.roots}, respectively. The exceptional contribution will be bounded in \eqref{temp:detail:high.exceptional.moment}; under the Haar law, $\widetilde{\mathcal{C}}_{T,\mathrm{exc}}^{\mathrm{H}} = 0$.

First identify the mean scale in \eqref{temp:detail:high.scale.definition}. The cumulant estimate \eqref{an:joint.cumulants}, its remainder bound \eqref{an:conditional.detail:remainder}, and the variance estimate \eqref{an:detail:cov.variance} give $\log\EE_{\mathrm{H}}[e^{\beta H(0)}] = \beta^{2}V_{T}/2 + \OO(1) = \beta^{2}s_{T}/4 + \OO(1)$. Together with \eqref{eq:effective.cutoff.bounds} and $u_{T}/s_{T} \to m(\theta)$ from \eqref{an:height.scales} in Lemma~\ref{an:height}, this yields
\[
\begin{aligned}
\log A_{T}^{\mathrm{res}}(\beta)
&= \log{N_{T}} - \beta u_{T} + \log\EE_{\mathrm{H}}[e^{\beta H(0)}] \\
&= \left(m(\theta)^{2} - \beta m(\theta) + \frac{\beta^{2}}{4}\right)s_{T} + \oo(s_{T})
= (m(\theta) - \beta/2)^{2}s_{T} + \oo(s_{T}).
\end{aligned}
\]
By \eqref{eq:N.T}, $A_{T}^{\mathrm{res}}(\beta) = \ell_{T}^{-1}e^{-\beta u_{T}}|I_{\theta}|\EE_{\mathrm{H}}[e^{\beta H(0)}]$. Summing the defining integrals for $Z_{B}$ and applying \eqref{an:remote.tilt.mean} in Lemma~\ref{an:capped.subcritical} therefore give
\begin{equation}\label{eq:temp:annealed.scale}
\begin{aligned}
\log A_{T}^{\mathrm{res}}(\beta)
&= (m(\theta) - \beta/2)^{2}s_{T} + \oo(s_{T}), \\
\EE_{\star}\left[\sum_{B \in \mathcal{B}_{T}}Z_{B}\right]
&= \ell_{T}^{-1}e^{-\beta u_{T}}\EE_{\star}\left[\int_{I_{\theta}}e^{\beta F_{h}^{\star}}\ind_{\{F_{h}^{\star}\leq u_{T} + K\}}\rd{h}\right] \\
&= \ell_{T}^{-1}e^{-\beta u_{T}}|I_{\theta}|\EE_{\mathrm{H}}[e^{\beta H(0)}]\{1 + \oo(1)\}
= A_{T}^{\mathrm{res}}(\beta)\{1 + \oo(1)\}.
\end{aligned}
\end{equation}
The exponential rate in \eqref{eq:temp:annealed.scale} is positive because $\beta < 2m(\theta)$. Under the Haar tilt with density $e^{\beta H(h)}/\EE_{\mathrm{H}}[e^{\beta H(h)}]$, differentiating \eqref{an:joint.cumulants} and using the derivative bound in \eqref{an:conditional.detail:remainder} gives mean $\beta V_{T} + \OO(1) = \beta s_{T}/2 + \OO(1)$. Thus the cap $u_{T} + K$ lies above the tilted mean of $H(h)$ by $(m(\theta) - \beta/2)s_{T} + \oo(s_{T})$. The exponential tail estimate \eqref{an:remote.detail:tilted.tail} in the proof of Lemma~\ref{an:capped.subcritical} shows that the Haar exponential moment lost above this cap is $\oo(\EE_{\mathrm{H}}[e^{\beta H(h)}])$. Integrating over $I_{\theta}$ gives the Haar case of \eqref{an:remote.tilt.mean}; the arithmetic case follows from the truncated-transform comparison in the same proof.

The next step is to use \eqref{temp:detail:high.moment.decomposition} to prove that, for sufficiently small $a > 0$ and some $c > 0$,
\begin{equation}\label{eq:temp:high.moment}
\EE_{\star}\left[\sum_{B \in \mathcal{B}_{T}}Z_{B}^{1 + a}\right] \leq A_{T}^{\mathrm{res}}(\beta)^{1 + a}e^{-cs_{T}}
\end{equation}
for all sufficiently large $T$. Recall from \eqref{an:annular.parameters} that
\[
v(\beta) = \begin{cases}1 + \beta^{2}/4, & \beta \leq 2, \\ \beta, & \beta \geq 2.\end{cases}
\]
For every root $h_{0} \in B$, separating the root height from the increment in the defining integral for $Z_{B}$ gives
\begin{equation}\label{temp:detail:high.relative.bound}
\begin{aligned}
Z_{B}^{a}
&= e^{a\beta(F_{h_{0}}^{\star} - u_{T})}\left(\ell_{T}^{-1}\int_{B}e^{\beta(F_{h}^{\star} - F_{h_{0}}^{\star})}\ind_{\{F_{h}^{\star}\leq u_{T} + K\}}\rd{h}\right)^{a} \\
&\leq e^{a\beta(F_{h_{0}}^{\star} - u_{T})}\left(\ell_{T}^{-1}\int_{B}e^{\beta(F_{h}^{\star} - F_{h_{0}}^{\star})}\rd{h}\right)^{a}.
\end{aligned}
\end{equation}
Fix an integer $0 \leq j \leq \lfloor u_{T} + K\rfloor$ and set $\gamma \leqdef j/s_{T}$. Enlarge the half-open bin $[j,j + 1)$ to the closed bin $[\gamma s_{T},\gamma s_{T} + 1]$ and remove the root-height cap for an upper bound. On this bin, the product of the root factors in \eqref{temp:detail:high.bin.contribution} and \eqref{temp:detail:high.relative.bound} is at most $e^{(1 + a)\beta(\gamma s_{T} + 1 - u_{T})}$. Conditioning on the closed bin therefore gives
\begin{align}\label{temp:detail:high.bin.bound}
\widetilde{\mathcal{C}}_{T,j}^{\star}
&\leq e^{(1 + a)\beta(\gamma s_{T} + 1 - u_{T})}\sum_{B \in \mathcal{B}_{T}}\ell_{T}^{-1}\int_{B}\PP_{\star}(F_{h_{0}}^{\star} \in [\gamma s_{T},\gamma s_{T} + 1]) \notag \\
&\qquad \times \EE_{\star}\left[\left(\ell_{T}^{-1}\int_{B}e^{\beta(F_{h}^{\star} - F_{h_{0}}^{\star})}\rd{h}\right)^{a}\ind_{\mathcal{A}_{T,\mathrm{sp}}} \mid F_{h_{0}}^{\star} \in [\gamma s_{T},\gamma s_{T} + 1]\right]\rd{h_{0}} \notag \\
&\leq C(1 + s_{T})e^{(1 + a)\beta(\gamma s_{T} + 1 - u_{T})}\exp\left\{\left[av(\beta) - a\beta\gamma + \frac{a^{2}\beta^{2}}{4}\right]_{+}s_{T}\right\} \notag \\
&\qquad \times \sum_{B \in \mathcal{B}_{T}}\ell_{T}^{-1}\int_{B}\PP_{\star}(F_{h_{0}}^{\star} \in [\gamma s_{T},\gamma s_{T} + 1])\rd{h_{0}} \notag \\
&\leq C(1 + s_{T})e^{(1 + a)\beta(\gamma s_{T} + 1 - u_{T})}N_{T}\sup_{h_{0} \in I_{\theta}}\PP_{\star}(F_{h_{0}}^{\star} \in [\gamma s_{T},\gamma s_{T} + 1]) \notag \\
&\qquad \times \exp\left\{\left[av(\beta) - a\beta\gamma + \frac{a^{2}\beta^{2}}{4}\right]_{+}s_{T}\right\} \notag \\
&\leq C(1 + s_{T})\exp\{E(\gamma)s_{T} + \oo(s_{T})\},
\end{align}
where
\begin{equation}\label{temp:detail:high.exponent}
E(\gamma) \leqdef m(\theta)^{2} - \gamma^{2} + (1 + a)\beta(\gamma - m(\theta)) + \left[av(\beta) - a\beta\gamma + \frac{a^{2}\beta^{2}}{4}\right]_{+}.
\end{equation}
The second inequality in \eqref{temp:detail:high.bin.bound} uses \eqref{an:unit.fractional} in Lemma~\ref{an:annular}, since $B \subseteq [h_{0} - 1,h_{0} + 1] \cap I_{\theta}$ and $a < 1 \wedge (2/\beta)$. The third inequality in \eqref{temp:detail:high.bin.bound} bounds the bin probability uniformly in $h_{0}$ and uses $\sum_{B}\ell_{T}^{-1}|B| = N_{T}$ from \eqref{eq:N.T}. The last inequality in \eqref{temp:detail:high.bin.bound} uses \eqref{temp:detail:bin.exponent} from the proof of Lemma~\ref{temp:critical.moment} and $u_{T}/s_{T} \to m(\theta)$ from \eqref{an:height.scales} in Lemma~\ref{an:height}. Since $0 \leq \gamma \leq (u_{T} + K)/s_{T} \to m(\theta)$, all slopes lie in one fixed compact interval, and the $\oo(s_{T})$ term in the exponent on the last line of \eqref{temp:detail:high.bin.bound} is uniform over the bins.

When the positive part in \eqref{temp:detail:high.exponent} is zero, completing the square gives the exact identity
\begin{equation}\label{temp:detail:high.zero.branch}
E(\gamma) - (1 + a)(m(\theta) - \beta/2)^{2} = - \left(\gamma - \frac{(1 + a)\beta}{2}\right)^{2} + a\left\{ - m(\theta)^{2} + \frac{(1 + a)\beta^{2}}{4}\right\}.
\end{equation}
Since $\beta^{2}/4 < \beta_{c}(\theta)^{2}/4 = m(\theta)^{2}$, the last term in \eqref{temp:detail:high.zero.branch} is negative for sufficiently small $a$. On the other branch of \eqref{temp:detail:high.exponent}, the corresponding identity is
\begin{equation}\label{temp:detail:high.positive.branch}
\begin{aligned}
E(\gamma) - (1 + a)(m(\theta) - \beta/2)^{2}
&= - (\gamma - \beta/2)^{2} + a\left\{v(\beta) - m(\theta)^{2} - \frac{\beta^{2}}{4}\right\} + \frac{a^{2}\beta^{2}}{4} \\
&\leq \begin{cases} - a\theta + a^{2}\beta^{2}/4, & \beta \leq 2, \\ - a\{\theta + (\beta - 2)^{2}/4\} + a^{2}\beta^{2}/4, & \beta \geq 2.\end{cases}
\end{aligned}
\end{equation}
Decrease $a > 0$ so that $(1 + a)\beta^{2}/4 < m(\theta)^{2}$ and $a\beta^{2}/4 < \theta/2$. Since the square term is nonpositive, the right-hand side of \eqref{temp:detail:high.zero.branch} is at most $ - a\{m(\theta)^{2} - (1 + a)\beta^{2}/4\}$. Both upper bounds in \eqref{temp:detail:high.positive.branch} are at most $ - a\theta + a^{2}\beta^{2}/4 < - a\theta/2$. Consequently, uniformly over all permitted $\gamma$,
\begin{equation}\label{temp:detail:high.exponent.margin}
E(\gamma) - (1 + a)(m(\theta) - \beta/2)^{2}
\leq - \min\left\{a\left(m(\theta)^{2} - \frac{(1 + a)\beta^{2}}{4}\right),\frac{a\theta}{2}\right\} < 0.
\end{equation}

To bound the negative-root and exceptional contributions in \eqref{temp:detail:high.other.contributions}, first note the deterministic cap
\begin{equation}\label{temp:detail:high.deterministic.cap}
Z_{B} = \ell_{T}^{-1}\int_{B}e^{\beta(F_{h}^{\star} - u_{T})}\ind_{\{F_{h}^{\star}\leq u_{T} + K\}}\rd{h}
\leq e^{\beta K}\ell_{T}^{-1}|B| \leq e^{s_{T} + \beta K},
\end{equation}
which uses $|B| \leq 1$ and $\ell_{T}^{-1} = e^{s_{T}}$ from \eqref{eq:N.T}. For a negative root, $F_{h_{0}}^{\star} < 0$, so $e^{\beta F_{h_{0}}^{\star}} \leq 1$. Applying \eqref{temp:detail:high.deterministic.cap} to $Z_{B}^{a}$ and then dropping the remaining indicators gives
\begin{equation}\label{temp:detail:high.negative.roots}
\begin{aligned}
\widetilde{\mathcal{C}}_{T,\mathrm{neg}}^{\star}
&\leq e^{as_{T} + a\beta K}\EE_{\star}\left[\sum_{B \in \mathcal{B}_{T}}\ell_{T}^{-1}\int_{B}e^{\beta(F_{h_{0}}^{\star} - u_{T})}\ind_{\{F_{h_{0}}^{\star} < 0\}}\rd{h_{0}}\right] \\
&\leq e^{as_{T} + a\beta K}\ell_{T}^{-1}\int_{I_{\theta}}e^{-\beta u_{T}}\rd{h_{0}}
= N_{T}e^{-\beta u_{T} + as_{T} + a\beta K} \\
&= \exp\{(m(\theta)^{2} - \beta m(\theta) + a)s_{T} + \oo(s_{T})\}.
\end{aligned}
\end{equation}
The last line uses \eqref{eq:effective.cutoff.bounds}, \eqref{eq:parameter.scales}, and \eqref{an:height.scales} in Lemma~\ref{an:height}, with $K$ fixed. Dividing \eqref{temp:detail:high.negative.roots} by the scale in \eqref{eq:temp:annealed.scale} gives
\begin{equation}\label{temp:detail:high.negative.rate}
\frac{\widetilde{\mathcal{C}}_{T,\mathrm{neg}}^{\star}}{A_{T}^{\mathrm{res}}(\beta)^{1 + a}}
\leq \exp\left\{\left[ - \frac{\beta^{2}}{4} + a\{1 - (m(\theta) - \beta/2)^{2}\} + \oo(1)\right]s_{T}\right\}.
\end{equation}
The coefficient of $s_{T}$ in \eqref{temp:detail:high.negative.rate}, excluding $\oo(1)$, is negative after reducing $a$, since $\beta > 0$ is fixed. Fix this $a$ and choose
\begin{equation}\label{temp:detail:high.rate.margin}
0 < c_{0} < \min\left\{a\left(m(\theta)^{2} - \frac{(1 + a)\beta^{2}}{4}\right),\frac{a\theta}{2},\frac{\beta^{2}}{4} - a\{1 - (m(\theta) - \beta/2)^{2}\}\right\}.
\end{equation}
All three entries in the minimum in \eqref{temp:detail:high.rate.margin} are positive by the restrictions on $a$. Combining \eqref{temp:detail:high.exponent.margin} with the choice of $c_{0}$ in \eqref{temp:detail:high.rate.margin} gives $E(\gamma) - (1 + a)(m(\theta) - \beta/2)^{2} \leq - c_{0}$ uniformly over the bins. There are $\lfloor u_{T} + K\rfloor + 1 = \OO(s_{T})$ bins, and each bound in \eqref{temp:detail:high.bin.bound} has prefactor $C(1 + s_{T})$. Summing these bounds and adding \eqref{temp:detail:high.negative.rate} yields
\begin{equation}\label{temp:detail:high.retained.moment}
\begin{aligned}
\frac{\EE_{\star}[\sum_{B}Z_{B}^{1 + a}\ind_{\mathcal{A}_{T,\mathrm{sp}}}]}{A_{T}^{\mathrm{res}}(\beta)^{1 + a}}
&= \frac{\sum_{j = 0}^{\lfloor u_{T} + K\rfloor}\widetilde{\mathcal{C}}_{T,j}^{\star} + \widetilde{\mathcal{C}}_{T,\mathrm{neg}}^{\star}}{A_{T}^{\mathrm{res}}(\beta)^{1 + a}} \\
&\leq \{Cs_{T}(1 + s_{T}) + C\}\exp\{ - c_{0}s_{T} + \oo(s_{T})\} \\
&\leq Ce^{-c_{0}s_{T}/2}.
\end{aligned}
\end{equation}
The error in the exponent is uniform over the bins by \eqref{temp:detail:high.bin.bound} and \eqref{eq:temp:annealed.scale}.

For the arithmetic exceptional contribution, apply \eqref{temp:detail:high.deterministic.cap} and the interval count $|\mathcal{B}_{T}| \leq e^{\theta t + \OO(1)}$. Together with \eqref{an:spatial.exceptional} in Lemma~\ref{an:spatial.transform}, these give
\begin{equation}\label{temp:detail:high.exceptional.moment}
\begin{aligned}
\widetilde{\mathcal{C}}_{T,\mathrm{exc}}^{\mathrm{A}}
= \EE\left[\sum_{B \in \mathcal{B}_{T}}Z_{B}^{1 + a}\ind_{\mathcal{A}_{T,\mathrm{sp}}^{c}}\right]
&\leq |\mathcal{B}_{T}|e^{(1 + a)(s_{T} + \beta K)}\PP(\mathcal{A}_{T,\mathrm{sp}}^{c}) \\[-3mm]
&\leq \exp\{(1 + a)(s_{T} + \beta K) + \theta t + \OO(1)\}\exp\{ - ct^{1 + 2\e_{0}} + Ct\} \\[1mm]
&\leq \exp\{ - ct^{1 + 2\e_{0}} + \OO_{a,\beta,K}(t)\}.
\end{aligned}
\end{equation}
By \eqref{eq:temp:annealed.scale}, dividing the bound in \eqref{temp:detail:high.exceptional.moment} by $A_{T}^{\mathrm{res}}(\beta)^{1 + a}$ still gives $\exp\{ - ct^{1 + 2\e_{0}} + \OO_{a,\beta,K}(t)\}$, which is at most $e^{-c_{0}s_{T}}$ for all sufficiently large $T$, with the same $c_{0}$ chosen in \eqref{temp:detail:high.rate.margin}. This weaker bound suffices to combine the exceptional contribution with \eqref{temp:detail:high.retained.moment}. The same bound holds for $\smash{\widetilde{\mathcal{C}}_{T,\mathrm{exc}}^{\mathrm{H}} = 0}$. Returning to \eqref{temp:detail:high.moment.decomposition} and using \eqref{temp:detail:high.retained.moment} therefore gives
\[
\frac{\EE_{\star}[\sum_{B}Z_{B}^{1 + a}]}{A_{T}^{\mathrm{res}}(\beta)^{1 + a}}
= \frac{\sum_{j = 0}^{\lfloor u_{T} + K\rfloor}\widetilde{\mathcal{C}}_{T,j}^{\star}
+ \widetilde{\mathcal{C}}_{T,\mathrm{neg}}^{\star}}{A_{T}^{\mathrm{res}}(\beta)^{1 + a}}
+ \frac{\widetilde{\mathcal{C}}_{T,\mathrm{exc}}^{\star}}{A_{T}^{\mathrm{res}}(\beta)^{1 + a}}
\leq Ce^{-c_{0}s_{T}/2} + e^{-c_{0}s_{T}} \leq e^{-cs_{T}}
\]
for some $c > 0$ and all sufficiently large $T$, proving \eqref{eq:temp:high.moment}.

To prove concentration of $\sum_{B}Z_{B}$, let $\e_{T}^{\mathrm{mass}} \leqdef \exp\{ - cs_{T}/(2(1 + a))\}$ and define
\begin{equation}\label{temp:detail:high.truncation}
\widetilde{Z}_{B} \leqdef Z_{B}\ind_{\{Z_{B}\leq\e_{T}^{\mathrm{mass}}A_{T}^{\mathrm{res}}(\beta)\}}.
\end{equation}
Then $\e_{T}^{\mathrm{mass}} \to 0$ and $(\e_{T}^{\mathrm{mass}})^{-1 - a}e^{-cs_{T}} = e^{-cs_{T}/2} \to 0$. On $\{Z_{B} > \e_{T}^{\mathrm{mass}}A_{T}^{\mathrm{res}}(\beta)\}$, the inequality $Z_{B} \leq (\e_{T}^{\mathrm{mass}}A_{T}^{\mathrm{res}}(\beta))^{-a}Z_{B}^{1 + a}$ and \eqref{eq:temp:high.moment} give
\begin{equation}\label{temp:detail:high.lost}
\begin{aligned}
\EE_{\star}\left[\sum_{B}(Z_{B} - \widetilde{Z}_{B})\right]
&= \EE_{\star}\left[\sum_{B}Z_{B}\ind_{\{Z_{B}>\e_{T}^{\mathrm{mass}}A_{T}^{\mathrm{res}}(\beta)\}}\right] \\
&\leq (\e_{T}^{\mathrm{mass}}A_{T}^{\mathrm{res}}(\beta))^{-a}\EE_{\star}\left[\sum_{B}Z_{B}^{1 + a}\right] \\
&\leq A_{T}^{\mathrm{res}}(\beta)(\e_{T}^{\mathrm{mass}})^{-a}e^{-cs_{T}}
= \oo(A_{T}^{\mathrm{res}}(\beta)).
\end{aligned}
\end{equation}
In particular, \eqref{eq:temp:annealed.scale} and \eqref{temp:detail:high.lost} imply that
\begin{equation}\label{temp:detail:high.truncated.mean}
\EE_{\star}\left[\sum_{B \in \mathcal{B}_{T}}\widetilde{Z}_{B}\right] = A_{T}^{\mathrm{res}}(\beta)\{1 + \oo(1)\}.
\end{equation}

For its second moment, split the ordered pairs of intervals into distant pairs and the remaining pairs, including the diagonal pairs:
\begin{equation}\label{temp:detail:high.second.moment.decomposition}
\EE_{\star}\left[\left(\sum_{B}\widetilde{Z}_{B}\right)^{2}\right]
= \sum_{(B,B')\,\mathrm{not\ distant}}\EE_{\star}[\widetilde{Z}_{B}\widetilde{Z}_{B'}]
+ \sum_{(B,B')\,\mathrm{distant}}\EE_{\star}[\widetilde{Z}_{B}\widetilde{Z}_{B'}].
\end{equation}
As in \eqref{temp:detail:window.near} in the proof of Proposition~\ref{temp:critical}, each interval has only a bounded number of intervals that are not distant from it (two intervals are called distant when their distance is at least two). The inequality $2xy \leq x^{2} + y^{2}$ and the bound $\widetilde{Z}_{B}^{2} \leq \e_{T}^{\mathrm{mass}}A_{T}^{\mathrm{res}}(\beta)\widetilde{Z}_{B}$ from \eqref{temp:detail:high.truncation} therefore give
\begin{equation}\label{temp:detail:high.near.pairs}
\begin{aligned}
\hspace{-3mm} \sum_{(B,B')\,\mathrm{not\ distant}} \hspace{-3mm} \EE_{\star}[\widetilde{Z}_{B}\widetilde{Z}_{B'}]
&\leq \frac{1}{2}\sum_{(B,B')\,\mathrm{not\ distant}}\left(\EE_{\star}[\widetilde{Z}_{B}^{2}] + \EE_{\star}[\widetilde{Z}_{B'}^{2}]\right) \\
&\leq C\sum_{B}\EE_{\star}[\widetilde{Z}_{B}^{2}]
\leq C\e_{T}^{\mathrm{mass}}A_{T}^{\mathrm{res}}(\beta)\EE_{\star}\left[\sum_{B}\widetilde{Z}_{B}\right]
= \oo(A_{T}^{\mathrm{res}}(\beta)^{2}),
\end{aligned}
\end{equation}
where the last equality uses \eqref{temp:detail:high.truncated.mean} and $\e_{T}^{\mathrm{mass}} \to 0$.

For distant pairs, $B \times B' \subseteq \mathscr{D}_{2,T}(2)$, with $\mathscr{D}_{2,T}(2)$ defined in \eqref{eq:D.k.T}, and these rectangles are disjoint up to boundaries of Lebesgue measure zero. Since $\smash{\widetilde{Z}_{B}\widetilde{Z}_{B'} \leq Z_{B}Z_{B'}}$, the defining integrals for $Z_{B}$ in \eqref{temp:interval.masses} give
\begin{equation}\label{temp:detail:high.distant.integral}
\begin{aligned}
\sum_{(B,B')\,\mathrm{distant}} \EE_{\star}[\widetilde{Z}_{B}\widetilde{Z}_{B'}]
&\leq \ell_{T}^{-2}e^{-2\beta u_{T}} \hspace{-3mm} \sum_{(B,B')\,\mathrm{distant}}\EE_{\star}\left[\int_{B \times B'}e^{\beta(F_{h}^{\star} + F_{h'}^{\star})}\ind_{\{F_{h}^{\star}\leq u_{T} + K,\,F_{h'}^{\star}\leq u_{T} + K\}}\rd{h}\rd{h'}\right] \\
&\leq \ell_{T}^{-2}e^{-2\beta u_{T}}\EE_{\star}\left[\int_{\mathscr{D}_{2,T}(2)}e^{\beta(F_{h}^{\star} + F_{h'}^{\star})}\ind_{\{F_{h}^{\star}\leq u_{T} + K,\,F_{h'}^{\star}\leq u_{T} + K\}}\rd{h}\rd{h'}\right].
\end{aligned}
\end{equation}
The normalization in \eqref{temp:detail:high.scale.definition} satisfies
\begin{equation}\label{temp:detail:high.scale.square}
A_{T}^{\mathrm{res}}(\beta)^{2} = \ell_{T}^{-2}e^{-2\beta u_{T}}|I_{\theta}|^{2}\EE_{\mathrm{H}}[e^{\beta H(0)}]^{2}.
\end{equation}
Dividing \eqref{temp:detail:high.distant.integral} by $A_{T}^{\mathrm{res}}(\beta)^{2}$ and using \eqref{temp:detail:high.scale.square} and applying \eqref{an:remote.tilt.pairs} in Lemma~\ref{an:capped.subcritical} with $d_{0} = 2$ yields
\begin{equation}\label{temp:detail:high.distant.pairs}
\limsup_{T \to \infty}\frac{\sum_{(B,B')\,\mathrm{distant}}\EE_{\star}[\widetilde{Z}_{B}\widetilde{Z}_{B'}]}{A_{T}^{\mathrm{res}}(\beta)^{2}} \leq 1.
\end{equation}
Equations \eqref{temp:detail:high.second.moment.decomposition}, \eqref{temp:detail:high.near.pairs}, and \eqref{temp:detail:high.distant.pairs} imply
\[
\limsup_{T \to \infty}\frac{\EE_{\star}[(\sum_{B}\widetilde{Z}_{B})^{2}]}{A_{T}^{\mathrm{res}}(\beta)^{2}} \leq 1.
\]
The normalized mean $\EE_{\star}[\sum_{B}\widetilde{Z}_{B}]/A_{T}^{\mathrm{res}}(\beta)$ tends to one by \eqref{temp:detail:high.truncated.mean}. Consequently,
\[
0 \leq \limsup_{T \to \infty}\frac{\Var_{\star}(\sum_{B}\widetilde{Z}_{B})}{A_{T}^{\mathrm{res}}(\beta)^{2}}
= \limsup_{T \to \infty}\frac{\EE_{\star}[(\sum_{B}\widetilde{Z}_{B})^{2}]}{A_{T}^{\mathrm{res}}(\beta)^{2}} - 1
\leq 1 - 1 = 0.
\]
This proves
\begin{equation}\label{temp:detail:high.variance}
\EE_{\star}\left[\sum_{B}\widetilde{Z}_{B}\right] = A_{T}^{\mathrm{res}}(\beta)\{1 + \oo(1)\}, \qquad \Var_{\star}\left(\sum_{B}\widetilde{Z}_{B}\right) = \oo(A_{T}^{\mathrm{res}}(\beta)^{2}).
\end{equation}
To recover the original sum $\sum_{B}Z_{B}$ from its truncation $\sum_{B}\widetilde{Z}_{B}$, the triangle inequality gives
\begin{equation}\label{temp:detail:high.total.decomposition}
\left|\frac{\sum_{B}Z_{B}}{A_{T}^{\mathrm{res}}(\beta)} - 1\right|
\leq \frac{|\sum_{B}\widetilde{Z}_{B} - \EE_{\star}[\sum_{B}\widetilde{Z}_{B}]|}{A_{T}^{\mathrm{res}}(\beta)}
+ \left|\frac{\EE_{\star}[\sum_{B}\widetilde{Z}_{B}]}{A_{T}^{\mathrm{res}}(\beta)} - 1\right|
+ \frac{\sum_{B}(Z_{B} - \widetilde{Z}_{B})}{A_{T}^{\mathrm{res}}(\beta)}.
\end{equation}
For each fixed $\delta > 0$, the middle term in \eqref{temp:detail:high.total.decomposition} is at most $\delta/3$ for all sufficiently large $T$, by \eqref{temp:detail:high.variance}. Chebyshev's inequality for the first term and Markov's inequality for the last then give, as $T \to \infty$,
\begin{equation}\label{temp:detail:high.concentration.bound}
\PP_{\star}\left(\left|\frac{\sum_{B}Z_{B}}{A_{T}^{\mathrm{res}}(\beta)} - 1\right| > \delta\right)
\leq \frac{9\Var_{\star}(\sum_{B}\widetilde{Z}_{B})}{\delta^{2}A_{T}^{\mathrm{res}}(\beta)^{2}} + \frac{3\EE_{\star}[\sum_{B}(Z_{B} - \widetilde{Z}_{B})]}{\delta A_{T}^{\mathrm{res}}(\beta)} \to 0,
\end{equation}
using \eqref{temp:detail:high.variance} and \eqref{temp:detail:high.lost}. The probability bound \eqref{temp:detail:high.concentration.bound} proves
\begin{equation}\label{temp:detail:high.total}
\frac{\sum_{B}Z_{B}}{A_{T}^{\mathrm{res}}(\beta)} \to 1 \quad\text{in probability}.
\end{equation}

For the largest capped interval mass $\max_{B \in \mathcal{B}_{T}}Z_{B}$, a union bound and Markov's inequality, followed by \eqref{eq:temp:high.moment}, give, for every fixed $\e > 0$,
\begin{equation}\label{temp:detail:high.maximum}
\PP_{\star}\left(\max_{B}Z_{B} > \e A_{T}^{\mathrm{res}}(\beta)\right)
\leq \sum_{B}\PP_{\star}(Z_{B} > \e A_{T}^{\mathrm{res}}(\beta))
\leq \frac{\EE_{\star}[\sum_{B}Z_{B}^{1 + a}]}{(\e A_{T}^{\mathrm{res}}(\beta))^{1 + a}}
\leq \e^{-1 - a}e^{-cs_{T}} \to 0.
\end{equation}
Finally, return to the probability decomposition \eqref{temp:detail:high.probability.decomposition}. Equation \eqref{temp:detail:high.total} makes its second term tend to zero, while \eqref{temp:detail:high.maximum} with $\e/2$ in place of $\e$ bounds its third term by $(2/\e)^{1 + a}e^{-cs_{T}}$. Hence the second estimate in \eqref{an:height.tails} of Lemma~\ref{an:height} gives
\[
\limsup_{T \to \infty}\PP_{\star}\left(\max_{B \in \mathcal{B}_{T}}G_{\beta,T}^{\star}(B) > \e\right)
\leq \limsup_{T \to \infty}\PP_{\star}\left(\sup_{h \in I_{\theta}}(F_{h}^{\star} - u_{T}) > K\right)
\leq Ce^{-\beta_{c}(\theta)K}.
\]
Here $K$ is fixed throughout the limit $T \to \infty$, and the last constant is independent of $K$. Letting $K \to \infty$ proves the assertion.
\end{proof}

\begin{corollary}[Vanishing off-diagonal overlaps]\label{temp:overlap}
Let $\theta > 0$. For every $0 < \beta \leq \beta_{c}(\theta)$, as $T \to \infty$,
\[
\EE (G_{\beta,T}^{\mathrm{A}})^{\times 2}[|R_{T}^{\mathrm{A}}(h_{1},h_{2})|] \to 0.
\]
Consequently, all the conclusions of Theorem~\ref{expand:main} hold in this range, including the uniformity of the $\oo(1)$ term in \eqref{expand:GG} over $\|\Psi_{T}\|_{\infty} \leq 1$.
\end{corollary}

\begin{proof}[Proof of Corollary~\ref{temp:overlap}]
Fix $0 < \beta \leq \beta_{c}(\theta)$. Two positions at distance less than one belong to the same or adjacent intervals of $\mathcal{B}_{T}$. Summing over $B \in \mathcal{B}_{T}$ containing $h_{1}$ and, for each $B$, over $B$ and its at most two adjacent intervals containing $h_{2}$ gives
\begin{equation}\label{temp:detail:overlap.mass}
(G_{\beta,T}^{\mathrm{A}})^{\times 2}[|h_{1} - h_{2}| < 1]
\leq \sum_{B \in \mathcal{B}_{T}}G_{\beta,T}^{\mathrm{A}}(B) \hspace{-3mm}\sum_{\substack{B' = B\ \mathrm{or}\\B'\,\mathrm{adjacent\ to}\,B}} \hspace{-3mm} G_{\beta,T}^{\mathrm{A}}(B')
\leq 3\max_{B \in \mathcal{B}_{T}}G_{\beta,T}^{\mathrm{A}}(B).
\end{equation}
Proposition~\ref{temp:high} for $\beta < \beta_{c}(\theta)$, and \eqref{temp:critical.spreading} in Proposition~\ref{temp:critical} for $\beta = \beta_{c}(\theta)$, show that this maximum tends to zero in probability as $T \to \infty$. Since it lies in $[0,1]$, the convergence also holds in $L^{1}$. Taking expectation in \eqref{temp:detail:overlap.mass} therefore yields the following separation estimate under the annealed two-replica law $\EE (G_{\beta,T}^{\mathrm{A}})^{\times 2}$:
\begin{equation}\label{temp:detail:overlap.separation}
\EE (G_{\beta,T}^{\mathrm{A}})^{\times 2}[|h_{1} - h_{2}| < 1]
\leq 3\EE\left[\max_{B \in \mathcal{B}_{T}}G_{\beta,T}^{\mathrm{A}}(B)\right] \to 0, \qquad T \to \infty.
\end{equation}

Since $|R_{T}^{\mathrm{A}}(h,h')| \leq 1$ by Cauchy--Schwarz in \eqref{analytic:eq:overlaps}, and $r_{T}(h,h') = 0$ when $|h - h'| \geq 1$, \eqref{an:overlap.geometry} in Lemma~\ref{an:covariance} yields, uniformly for $h,h' \in I_{\theta}$,
\begin{equation}\label{temp:detail:overlap.localization}
|R_{T}^{\mathrm{A}}(h,h')| \leq \ind_{\{|h - h'| < 1\}} + C\frac{\log{t}}{t}.
\end{equation}
Integrating \eqref{temp:detail:overlap.localization} against the two-replica Gibbs measure, taking expectation, and applying \eqref{temp:detail:overlap.separation} gives
\begin{equation}\label{temp:detail:overlap.absolute.mean}
\EE (G_{\beta,T}^{\mathrm{A}})^{\times 2}[|R_{T}^{\mathrm{A}}(h_{1},h_{2})|]
\leq \EE (G_{\beta,T}^{\mathrm{A}})^{\times 2}[|h_{1} - h_{2}| < 1] + C\frac{\log{t}}{t}
\to 0, \qquad T \to \infty.
\end{equation}
This proves the stated convergence and, by Markov's inequality, the two-overlap limit $\delta_{0}$.

For every fixed number $n \geq 2$ of replicas, each off-diagonal entry has the same two-replica marginal. A union bound and Markov's inequality therefore give, for every $\delta > 0$,
\[
\EE (G_{\beta,T}^{\mathrm{A}})^{\times n}\left[\max_{1 \leq i < j \leq n}|R_{T}^{\mathrm{A}}(h_{i},h_{j})| > \delta\right]
\leq \frac{n(n - 1)}{2\delta}\EE (G_{\beta,T}^{\mathrm{A}})^{\times 2}[|R_{T}^{\mathrm{A}}(h_{1},h_{2})|]
\to 0,
\]
where the convergence follows from \eqref{temp:detail:overlap.absolute.mean}. The diagonal entries equal one, so the finite overlap array converges in probability under its annealed law to $\mathrm{Id}_{n}$. For $n = 1$, the array is identically $\mathrm{Id}_{1}$. Every continuous test function on $[ - 1,1]^{n \times n}$ is bounded and uniformly continuous, so its expectation under the annealed overlap law converges to its value at $\mathrm{Id}_{n}$ as $T \to \infty$, proving the joint overlap convergence asserted in Theorem~\ref{expand:main}.

Finally, fix $n \in \N$ and $1 \leq k \leq n$, and suppose $\|\Psi_{T}\|_{\infty} \leq 1$. Integrating out the unused replicas bounds the absolute value of each expectation containing an off-diagonal overlap in \eqref{expand:GG} of Theorem~\ref{expand:main} by $\EE (G_{\beta,T}^{\mathrm{A}})^{\times 2}[|R_{T}^{\mathrm{A}}(h_{1},h_{2})|]$. Also, $|\EE (G_{\beta,T}^{\mathrm{A}})^{\times n}[\Psi_{T}]| \leq 1$. The triangle inequality therefore gives the full residual bound
\begin{equation}\label{temp:detail:gg.residual}
\begin{aligned}
&\Bigg|\EE (G_{\beta,T}^{\mathrm{A}})^{\times(n + 1)}[R_{T}^{\mathrm{A}}(h_{k},h_{n + 1})\Psi_{T}(\bb{h})] - \frac{1}{n}\EE (G_{\beta,T}^{\mathrm{A}})^{\times 2}[R_{T}^{\mathrm{A}}(h_{1},h_{2})]\EE (G_{\beta,T}^{\mathrm{A}})^{\times n}[\Psi_{T}(\bb{h})] \\
&\hspace{80mm}- \frac{1}{n}\sum_{\substack{1 \leq j \leq n \\ j \neq k}}\EE (G_{\beta,T}^{\mathrm{A}})^{\times n}[R_{T}^{\mathrm{A}}(h_{k},h_{j})\Psi_{T}(\bb{h})]\Bigg| \\
&\leq \left(1 + \frac{1}{n} + \frac{n - 1}{n}\right)\EE (G_{\beta,T}^{\mathrm{A}})^{\times 2}[|R_{T}^{\mathrm{A}}(h_{1},h_{2})|]
= 2\EE (G_{\beta,T}^{\mathrm{A}})^{\times 2}[|R_{T}^{\mathrm{A}}(h_{1},h_{2})|] \to 0.
\end{aligned}
\end{equation}
The convergence follows from \eqref{temp:detail:overlap.absolute.mean}, and the bound in \eqref{temp:detail:gg.residual} is independent of the measurable positional test $\Psi_{T}$ with $\|\Psi_{T}\|_{\infty} \leq 1$. It also holds when $n = 1$, in which case the sum over the other sampled replicas is empty.
\end{proof}

\subsection{Normalization of random measures}\label{expand:normalization}

Dividing by the random total mass requires care because that mass may be close to zero. The next lemma handles this difficulty through a polynomial approximation, using bounds on $\EE[M^{r}]$ without requiring $\EE[M^{-r}]$ to be finite.

\begin{lemma}[Polynomial inverse-mass approximation]\label{norm:inverse}
Let $M$ be a nonnegative random variable such that
\begin{equation}\label{norm:detail:inverse.moments}
\EE[M^{r}] \leq C_{\mathrm{mom}}^{r}r!, \qquad r \in \N,
\end{equation}
where $C_{\mathrm{mom}} < \infty$ is independent of $r$. Then the following assertions hold.
\begin{enumerate}[label=(\roman*),itemsep=0mm]
\item The moments of $M$ determine its distribution.
\item For every $n \in \N$ and every $\e > 0$, there is a real polynomial $P_{\mathrm{inv}}$ such that
\begin{equation}\label{norm:detail:inverse.approximation}
\EE[\{1 - M^{n}P_{\mathrm{inv}}(M)\}^{2}] \leq \PP(M = 0) + \e^{2}.
\end{equation}
\end{enumerate}
\end{lemma}

\begin{proof}[Proof of Lemma~\ref{norm:inverse}]
If $M = 0$ almost surely, assertions~(i) and~(ii) hold, with $P_{\mathrm{inv}} = 0$ in~(ii). Otherwise, $C_{\mathrm{mom}} > 0$. Fix $0 < c < C_{\mathrm{mom}}^{-1}$. Monotone convergence and \eqref{norm:detail:inverse.moments} give
\begin{equation}\label{norm:detail:exponential.moment}
\EE[e^{cM}] = \sum_{r = 0}^{\infty}\frac{c^{r}\EE[M^{r}]}{r!} \leq \sum_{r = 0}^{\infty}(c \, C_{\mathrm{mom}})^{r} = \frac{1}{1 - c \, C_{\mathrm{mom}}} < \infty.
\end{equation}
Since $M$ is nonnegative, \eqref{norm:detail:exponential.moment} implies that its moment generating function is finite in a neighborhood of zero. Hence its moments determine its distribution, proving~(i).

To prove~(ii), fix $n \in \N$ and $\e > 0$, and define the measure $\Xi_{n}(\rd{x}) \leqdef x^{2n}\PP(M \in \rd{x})$ on $[0,\infty)$. The hypotheses of \citet[Theorem~1.1]{deJeu2003} in dimension one, with $R = 1$ and $\rho(s) = cs/2$ for $s > 1$, are verified as follows.
\begin{enumerate}[label=(\arabic*),itemsep=0mm]
\item The extension of $\Xi_{n}$ by zero to $(-\infty,0)$ is a nonnegative Borel measure on $\mathbb R$.

\item The function $\rho : (1,\infty) \to [0,\infty)$ is nondecreasing and of class $C^{1}$, and
\[
\int_{1}^{\infty}\frac{\rho(s)}{s^{2}}\rd{s} = \frac{c}{2}\int_{1}^{\infty}\frac{1}{s}\rd{s} = \infty.
\]

\item Since $\Xi_{n}$ is supported on $[0,\infty)$, \eqref{norm:detail:exponential.moment} yields the required integrability:
\[
\begin{aligned}
\int_{|x| > 1}\exp\left(\int_{1}^{|x|}\frac{\rho(s)}{s}\rd{s}\right)\Xi_{n}(\rd{x})
&= e^{-c/2}\int_{(1,\infty)}e^{cx/2}\Xi_{n}(\rd{x}) \\
&\leq e^{-c/2}\EE[M^{2n}e^{cM/2}] \\
&\leq e^{-c/2}\left(\sup_{x \geq 0}x^{2n}e^{-cx/2}\right)\EE[e^{cM}] < \infty.
\end{aligned}
\]
The supremum is finite because exponential decay dominates polynomial growth.
\end{enumerate}
The theorem with $p = 2$ shows that polynomials with complex coefficients are dense in $L^{2}_{\mathbb C}(\Xi_{n})$, the space of square-integrable complex-valued functions with respect to $\Xi_{n}$.

Set $x^{-n} = 0$ at $x = 0$ for its $L^{2}_{\mathbb C}(\Xi_{n})$ representative. Since $\Xi_{n}(\{0\}) = 0$,
\[
\int_{(0,\infty)}(x^{-n})^{2}\Xi_{n}(\rd{x}) = \int_{(0,\infty)}x^{-2n}x^{2n}\PP(M \in \rd{x}) = \PP(M > 0) \leq 1.
\]
Polynomial density supplies a polynomial $Q$ with complex coefficients whose $L^{2}_{\mathbb C}(\Xi_{n})$ distance from $x\mapsto x^{-n}$ is at most $\e$. Let $P_{\mathrm{inv}}$ be the polynomial obtained by taking the real parts of the coefficients of $Q$. Thus $P_{\mathrm{inv}}(x) = \Re Q(x)$ for real $x$, and $|x^{-n} - P_{\mathrm{inv}}(x)| \leq |x^{-n} - Q(x)|$ for $x > 0$. Separating the atom at zero now yields
\[
\begin{aligned}
\EE[\{1 - M^{n}P_{\mathrm{inv}}(M)\}^{2}]
&= \PP(M = 0) + \int_{(0,\infty)}|1 - x^{n}P_{\mathrm{inv}}(x)|^{2}\PP(M \in \rd{x}) \\
&\leq \PP(M = 0) + \int_{(0,\infty)}|x^{-n} - Q(x)|^{2}\Xi_{n}(\rd{x}) \\
&\leq \PP(M = 0) + \e^{2},
\end{aligned}
\]
which proves \eqref{norm:detail:inverse.approximation} and hence~(ii).
\end{proof}

\begin{proposition}[Normalization principle]\label{norm:transfer}
For every $T$, let $\mu_{T}^{\mathrm{A}}$ and $\mu_{T}^{\mathrm{H}}$ be nonnegative finite random measures on a common measurable space $\mathscr{S}_{T}$, with almost surely positive total masses. For each $b \in \N$, let $\nu_{T,b}^{\mathrm{A}}$ and $\nu_{T,b}^{\mathrm{H}}$ be nonnegative finite random measures on $\mathscr{S}_{T}$. For $\star \in \{\mathrm{A},\mathrm{H}\}$, define
\[
G_{T}^{\star} \leqdef \frac{\mu_{T}^{\star}}{\mu_{T}^{\star}(\mathscr{S}_{T})}, \qquad M_{T,b}^{\star} \leqdef \nu_{T,b}^{\star}(\mathscr{S}_{T}), \qquad G_{T,b}^{\star} \leqdef \frac{\nu_{T,b}^{\star}}{M_{T,b}^{\star}} \quad\text{on }\{M_{T,b}^{\star} > 0\}.
\]
On the event $\{M_{T,b}^{\star} = 0\}$, let $G_{T,b}^{\star}$ be an arbitrary fixed probability measure on $\mathscr{S}_{T}$. Suppose that the following conditions hold:
\begin{enumerate}[label=(\roman*),itemsep=0mm]
\item For every fixed $b,r \in \N$, as $T \to \infty$,
\begin{equation}\label{norm:detail:moment.measure.assumption}
\left\|\EE_{\mathrm{A}}[(\nu_{T,b}^{\mathrm{A}})^{\times r}] - \EE_{\mathrm{H}}[(\nu_{T,b}^{\mathrm{H}})^{\times r}]\right\|_{\mathrm{TV}} \to 0.
\end{equation}
\item For every fixed $b$, there is $C_{\mathrm{mom},b} < \infty$, independent of $r$, such that
\begin{equation}\label{norm:detail:mass.moment.assumption}
\limsup_{T\to\infty}\EE_{\mathrm{H}}[(M_{T,b}^{\mathrm{H}})^{r}] \leq C_{\mathrm{mom},b}^{r}r!, \qquad r \in \N.
\end{equation}
\item Every subsequential distributional limit $M_{b}$ of $M_{T,b}^{\mathrm{H}}$ as $T \to \infty$ satisfies
\begin{equation}\label{norm:detail:zero.atom}
\PP(M_{b} = 0) \leq \e_{b}^{\mathrm{atom}},
\end{equation}
where $\e_{b}^{\mathrm{atom}}$ is deterministic and tends to zero as $b \to \infty$.
\item For $\star \in \{\mathrm{A},\mathrm{H}\}$,
\begin{equation}\label{norm:detail:capture.assumption}
\lim_{b\to\infty}\limsup_{T\to\infty}\EE_{\star}[\|G_{T,b}^{\star} - G_{T}^{\star}\|_{\mathrm{TV}}] = 0.
\end{equation}
\end{enumerate}
Then, for every fixed $n \in \N$, as $T \to \infty$,
\begin{equation}\label{norm:detail:normalization.conclusion}
\left\|\EE_{\mathrm{A}}[(G_{T}^{\mathrm{A}})^{\times n}] - \EE_{\mathrm{H}}[(G_{T}^{\mathrm{H}})^{\times n}]\right\|_{\mathrm{TV}} \to 0.
\end{equation}
\end{proposition}

\begin{proof}[Proof of Proposition~\ref{norm:transfer}]
Fix $b,n \in \N$ and an arbitrary sequence of values of $T$ tending to infinity. The first step extracts a subsequence along which the total mass $M_{T,b}^{\mathrm{H}}$ of $\nu_{T,b}^{\mathrm{H}}$ converges in distribution to a nonnegative random variable $M_{b}$. Lemma~\ref{norm:inverse} then supplies a real polynomial $P_{\mathrm{inv}}$ such that $M_{b}^{n}P_{\mathrm{inv}}(M_{b})$ approximates $\ind_{\{M_{b}>0\}}$ in $L^{2}$ to any prescribed accuracy. By Markov's inequality and the first-moment bound in \eqref{norm:detail:mass.moment.assumption},
\[
\limsup_{T\to\infty}\PP_{\mathrm{H}}(M_{T,b}^{\mathrm{H}} > L) \leq L^{-1}\limsup_{T\to\infty}\EE_{\mathrm{H}}[M_{T,b}^{\mathrm{H}}] \leq C_{\mathrm{mom},b}L^{-1}.
\]
Thus, for fixed $b$, the Haar masses $M_{T,b}^{\mathrm{H}}$ form a tight family along the chosen sequence. After passing to a subsequence, they converge in distribution to a nonnegative random variable $M_{b}$. For each fixed $r \in \N$, \eqref{norm:detail:mass.moment.assumption} at order $r + 1$ gives, for all sufficiently large $T$,
\begin{equation}\label{norm:detail:mass.uniform.integrability}
\EE_{\mathrm{H}}[(M_{T,b}^{\mathrm{H}})^{r}\ind_{\{M_{T,b}^{\mathrm{H}}>L\}}] \leq L^{-1}\EE_{\mathrm{H}}[(M_{T,b}^{\mathrm{H}})^{r + 1}] \leq C_{b,r}L^{-1}.
\end{equation}
By \eqref{norm:detail:mass.uniform.integrability}, for each fixed $r$, the $r$th powers are uniformly integrable after omitting finitely many initial terms of this subsequence. Convergence in distribution, followed by truncation at $L$ and then $L \to \infty$, yields
\begin{equation}\label{norm:detail:limiting.moments}
\EE_{\mathrm{H}}[(M_{T,b}^{\mathrm{H}})^{r}] \to \EE[M_{b}^{r}] \leq C_{\mathrm{mom},b}^{r}r!, \qquad r \in \N.
\end{equation}
The same subsequence works for every fixed $r$, since the distributional convergence has already been obtained and \eqref{norm:detail:mass.uniform.integrability} applies separately at each order $r \in \N$, without requiring any further subsequence extraction. Since $(M_{T,b}^{\star})^{r} = (\nu_{T,b}^{\star})^{\times r}(\mathscr{S}_{T}^{r})$, the definition of the total variation norm and \eqref{norm:detail:moment.measure.assumption} give
\begin{equation}\label{norm:detail:arithmetic.mass.moments}
\left|\EE_{\mathrm{A}}[(M_{T,b}^{\mathrm{A}})^{r}] - \EE_{\mathrm{H}}[(M_{T,b}^{\mathrm{H}})^{r}]\right| \leq \left\|\EE_{\mathrm{A}}[(\nu_{T,b}^{\mathrm{A}})^{\times r}] - \EE_{\mathrm{H}}[(\nu_{T,b}^{\mathrm{H}})^{\times r}]\right\|_{\mathrm{TV}} \to 0.
\end{equation}
Consequently, \eqref{norm:detail:arithmetic.mass.moments} and \eqref{norm:detail:limiting.moments} give the same limits for all fixed arithmetic and Haar moments.

Let $\e > 0$, and apply \eqref{norm:detail:inverse.approximation} in Lemma~\ref{norm:inverse} to $M_{b}$, using \eqref{norm:detail:limiting.moments}. This supplies a polynomial $P_{\mathrm{inv}}$, fixed along the chosen subsequence. The function $\{1 - x^{n}P_{\mathrm{inv}}(x)\}^{2}$ is a polynomial of fixed degree. Its expectation therefore converges by \eqref{norm:detail:limiting.moments} and \eqref{norm:detail:arithmetic.mass.moments}. Using \eqref{norm:detail:inverse.approximation} in Lemma~\ref{norm:inverse} and the zero-atom bound \eqref{norm:detail:zero.atom} gives
\begin{equation}\label{norm:detail:polynomial.error}
\lim_{T\to\infty}\EE_{\star}[\{1 - (M_{T,b}^{\star})^{n}P_{\mathrm{inv}}(M_{T,b}^{\star})\}^{2}] = \EE[\{1 - M_{b}^{n}P_{\mathrm{inv}}(M_{b})\}^{2}] \leq \e_{b}^{\mathrm{atom}} + \e^{2}.
\end{equation}
Here and below the limits in the subsequence argument are taken along the subsequence just selected.

The polynomial error in \eqref{norm:detail:polynomial.error} controls normalization uniformly over the measurable test functions $\Psi_{T}:\mathscr{S}_{T}^{n} \to \R$ with $\|\Psi_{T}\|_{\infty} \leq 1$. Indeed, for a nonnegative finite measure $\nu$ of mass $M > 0$, let $G \leqdef \nu/M$. For every such $\Psi_{T}$, the identity $\nu^{\times n}[\Psi_{T}] = M^{n}G^{\times n}[\Psi_{T}]$ gives
\begin{equation}\label{norm:detail:normalization.pointwise}
|G^{\times n}[\Psi_{T}] - \nu^{\times n}[\Psi_{T}]P_{\mathrm{inv}}(M)| = |G^{\times n}[\Psi_{T}]|\,|1 - M^{n}P_{\mathrm{inv}}(M)| \leq |1 - M^{n}P_{\mathrm{inv}}(M)|.
\end{equation}
When $M = 0$, the nonnegative measure $\nu$ is zero, whereas the prescribed $G$ is a probability measure. The inequality \eqref{norm:detail:normalization.pointwise} then holds because its right-hand side is one. Taking expectations in \eqref{norm:detail:normalization.pointwise} and applying the Cauchy--Schwarz inequality therefore gives, uniformly over these test functions $\Psi_{T}$,
\begin{equation}\label{norm:detail:normalization.error}
\left|\EE_{\star}[(G_{T,b}^{\star})^{\times n}[\Psi_{T}]] - \EE_{\star}[(\nu_{T,b}^{\star})^{\times n}[\Psi_{T}]P_{\mathrm{inv}}(M_{T,b}^{\star})]\right| \leq \left(\EE_{\star}[\{1 - (M_{T,b}^{\star})^{n}P_{\mathrm{inv}}(M_{T,b}^{\star})\}^{2}]\right)^{1/2}.
\end{equation}
The limiting upper value in \eqref{norm:detail:normalization.error} is at most $\sqrt{\e_{b}^{\mathrm{atom}} + \e^{2}}$ by \eqref{norm:detail:polynomial.error}.

Write $P_{\mathrm{inv}}(x) = \sum_{j = 0}^{D_{\mathrm{inv}}}a_{j}^{\mathrm{inv}}x^{j}$, and extend the test function by
\[
\Psi_{T,j}(h_{1},\ldots,h_{n + j}) \leqdef \Psi_{T}(h_{1},\ldots,h_{n}).
\]
Integration over each additional position contributes one factor of the mass $M_{T,b}^{\star}$. In particular,
\[
(\nu_{T,b}^{\star})^{\times(n + j)}[\Psi_{T,j}] = (\nu_{T,b}^{\star})^{\times n}[\Psi_{T}](M_{T,b}^{\star})^{j}.
\]
Expanding the polynomial $P_{\mathrm{inv}}$ therefore gives
\[
\EE_{\star}[(\nu_{T,b}^{\star})^{\times n}[\Psi_{T}]P_{\mathrm{inv}}(M_{T,b}^{\star})] = \sum_{j = 0}^{D_{\mathrm{inv}}}a_{j}^{\mathrm{inv}}\EE_{\star}[(\nu_{T,b}^{\star})^{\times(n + j)}[\Psi_{T,j}]].
\]
The expectation of the absolute value of each summand is at most $|a_{j}^{\mathrm{inv}}|\EE_{\star}[(M_{T,b}^{\star})^{n + j}]$, which is finite for all sufficiently large $T$ by \eqref{norm:detail:limiting.moments} and \eqref{norm:detail:arithmetic.mass.moments}. Moreover, $\|\Psi_{T,j}\|_{\infty} \leq 1$. Hence, uniformly over these test functions $\Psi_{T}$,
\begin{equation}\label{norm:detail:polynomial.comparison}
\begin{aligned}
&\left|\EE_{\mathrm{A}}[(\nu_{T,b}^{\mathrm{A}})^{\times n}[\Psi_{T}]P_{\mathrm{inv}}(M_{T,b}^{\mathrm{A}})] - \EE_{\mathrm{H}}[(\nu_{T,b}^{\mathrm{H}})^{\times n}[\Psi_{T}]P_{\mathrm{inv}}(M_{T,b}^{\mathrm{H}})]\right| \\
&\qquad \leq \sum_{j = 0}^{D_{\mathrm{inv}}}|a_{j}^{\mathrm{inv}}|\left\|\EE_{\mathrm{A}}[(\nu_{T,b}^{\mathrm{A}})^{\times(n + j)}] - \EE_{\mathrm{H}}[(\nu_{T,b}^{\mathrm{H}})^{\times(n + j)}]\right\|_{\mathrm{TV}} \to 0,
\end{aligned}
\end{equation}
by \eqref{norm:detail:moment.measure.assumption}. Combining \eqref{norm:detail:normalization.error} for the two laws with \eqref{norm:detail:polynomial.comparison}, and taking the supremum over $\|\Psi_{T}\|_{\infty} \leq 1$, gives, along the chosen subsequence,
\[
\limsup_{T\to\infty}\left\|\EE_{\mathrm{A}}[(G_{T,b}^{\mathrm{A}})^{\times n}] - \EE_{\mathrm{H}}[(G_{T,b}^{\mathrm{H}})^{\times n}]\right\|_{\mathrm{TV}} \leq 2\sqrt{\e_{b}^{\mathrm{atom}} + \e^{2}}.
\]
Letting $\e \downarrow 0$ removes the approximation error. Every original sequence tending to infinity admits such a subsequence, so the same bound holds for the full limiting upper value:
\begin{equation}\label{norm:detail:retained.comparison}
\limsup_{T\to\infty}\left\|\EE_{\mathrm{A}}[(G_{T,b}^{\mathrm{A}})^{\times n}] - \EE_{\mathrm{H}}[(G_{T,b}^{\mathrm{H}})^{\times n}]\right\|_{\mathrm{TV}} \leq 2\sqrt{\e_{b}^{\mathrm{atom}}}.
\end{equation}
Indeed, a sequence violating \eqref{norm:detail:retained.comparison} by a fixed positive amount would have a subsequence to which the preceding argument applies, giving a contradiction.

For probability measures $G,G'$, the product difference can be written as
\[
G^{\times n} - (G')^{\times n} = \sum_{i = 1}^{n}G^{\times(i - 1)} \otimes (G - G') \otimes (G')^{\times(n - i)}.
\]
Each factor outside $G - G'$ has total mass one, so integration over the remaining $n - 1$ coordinates preserves the bound one on a test function. The triangle inequality therefore gives
\begin{equation}\label{norm:detail:product.variation}
\|G^{\times n} - (G')^{\times n}\|_{\mathrm{TV}}
\leq \sum_{i = 1}^{n}\|G^{\times(i - 1)} \otimes (G - G') \otimes (G')^{\times(n - i)}\|_{\mathrm{TV}}
\leq n\|G - G'\|_{\mathrm{TV}}.
\end{equation}
The triangle inequality, \eqref{norm:detail:retained.comparison}, and \eqref{norm:detail:product.variation} imply
\[
\limsup_{T\to\infty}\left\|\EE_{\mathrm{A}}[(G_{T}^{\mathrm{A}})^{\times n}] - \EE_{\mathrm{H}}[(G_{T}^{\mathrm{H}})^{\times n}]\right\|_{\mathrm{TV}}
\leq 2\sqrt{\e_{b}^{\mathrm{atom}}} + n\sum_{\star \in \{\mathrm{A},\mathrm{H}\}}\limsup_{T\to\infty}\EE_{\star}[\|G_{T,b}^{\star} - G_{T}^{\star}\|_{\mathrm{TV}}].
\]
Finally let $b \to \infty$, using $\e_{b}^{\mathrm{atom}} \to 0$ and \eqref{norm:detail:capture.assumption} to obtain \eqref{norm:detail:normalization.conclusion}.
\end{proof}

\subsection{Comparison of annealed sampling laws at low temperature}\label{sec:capped}

Throughout this subsection, assume that $\theta > 0$ and $\beta > \beta_c(\theta)$. This subsection shows that, at low temperature, the Euler and Haar Gibbs measures give asymptotically the same probabilities for events involving any fixed number of sampled positions, after averaging over the field randomness. These events may depend on $T$.

\begin{theorem}[Comparison of annealed sampling laws]\label{cap:sampling}
Let $\theta > 0$ and $\beta > \beta_{c}(\theta)$. For every fixed $n \in \N$, as $T \to \infty$,
\[
\left\|\EE (G_{\beta,T}^{\mathrm{A}})^{\times n} - \EE_{\mathrm{H}}(G_{\beta,T}^{\mathrm{H}})^{\times n}\right\|_{\mathrm{TV}} \to 0.
\]
The total variation norm is taken on $I_{\theta}^{\times n}$.
\end{theorem}

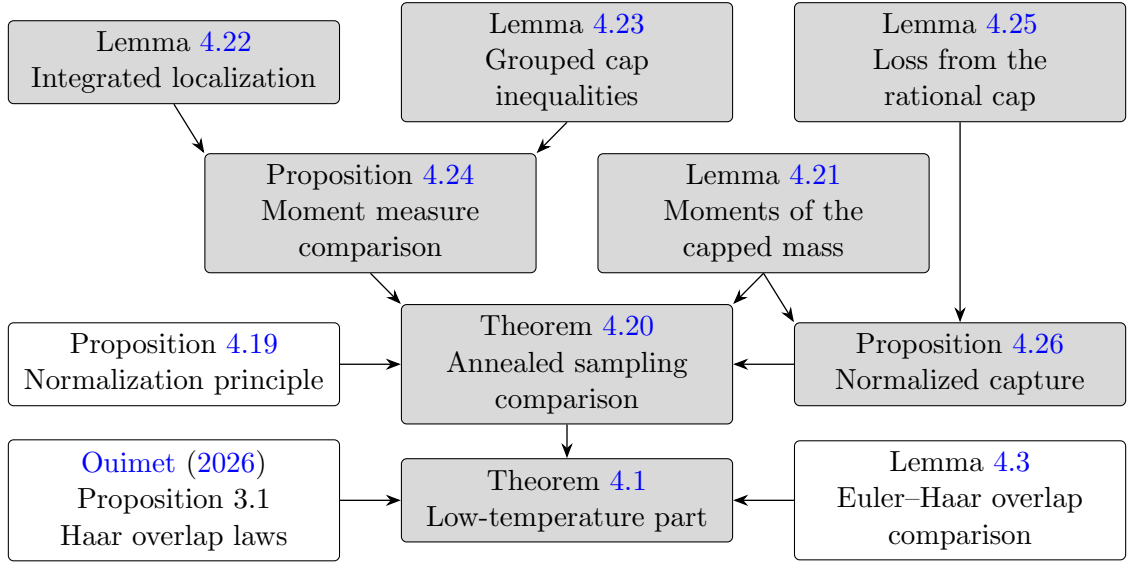
\begin{figure}[H]
\centering
\begin{tikzpicture}[x = 1cm, y = 1cm, every node/.style = {font = \small, align = center, execute at begin node = {\hyphenpenalty = 10000\exhyphenpenalty = 10000\relax}}, result/.style = {draw, rounded corners = 2pt, fill = black!15, text width = 4.1cm, minimum height = 1.05cm, inner sep = 4pt}, input/.style = {result, fill = white}, dependency/.style = {-{Stealth[length = 2mm]}, line width = 0.5pt}]
\node[result] (localization) at (0,0) {Lemma~\ref{cap:localization.bound} \\Integrated localization};
\node[result] (inequalities) at (5.2,0) {Lemma~\ref{cap:inequalities} \\Grouped cap inequalities};
\node[result] (loss) at (10.4,0) {Lemma~\ref{cap:loss} \\Loss from the rational cap};

\node[result] (comparison) at (2.6, - 1.998) {Proposition~\ref{cap:moment.measures} \\Moment measure comparison};
\node[result] (moments) at (7.8, - 1.998) {Lemma~\ref{cap:moments} \\Moments of the capped mass};

\node[input] (normalization) at (0, - 3.996) {Proposition~\ref{norm:transfer} \\Normalization principle};
\node[result] (sampling) at (5.2, - 3.996) {Theorem~\ref{cap:sampling} \\Annealed sampling comparison};
\node[result] (capture) at (10.4, - 3.996) {Proposition~\ref{cap:normalized.capture} \\Normalized capture};

\node[input] (haar) at (0, - 5.796) {\mbox{\citet{OuimetHaar2026}} \\Proposition~3.1 \\Haar overlap laws};
\node[result] (main) at (5.2, - 5.796) {Theorem~\ref{expand:main} \\Low-temperature part};
\node[input] (overlaps) at (10.4, - 5.796) {Lemma~\ref{an:covariance} \\Euler--Haar overlap comparison};

\draw[dependency] (localization.south) to (comparison.north west);
\draw[dependency] (inequalities.south) to (comparison.north east);
\draw[dependency] (loss.south) to (capture.north);
\draw[dependency] (moments.south) to (capture.north west);
\draw[dependency] (comparison.south) to (sampling.north west);
\draw[dependency] (moments.south) to (sampling.north east);
\draw[dependency] (normalization.east) to (sampling.west);
\draw[dependency] (capture.west) to (sampling.east);
\draw[dependency] (sampling.south) to (main.north);
\draw[dependency] (haar.east) to (main.west);
\draw[dependency] (overlaps.west) to (main.east);
\end{tikzpicture}
\caption{Dependencies among the results of Subsection~\ref{sec:capped} and their use in the low-temperature proof of Theorem~\ref{expand:main}. Arrows point from inputs to conclusions. White boxes denote previously established results. Other inputs from Subsections~\ref{expand:arithmetic} and~\ref{sec:temperature} are omitted.}\label{fig:capped.dependencies}
\end{figure}

The proof of Theorem~\ref{cap:sampling} is given at the end of this subsection, after Lemmas~\ref{cap:moments}, \ref{cap:localization.bound}, \ref{cap:inequalities}, and~\ref{cap:loss} and Propositions~\ref{cap:moment.measures} and~\ref{cap:normalized.capture}. These results establish the moment bounds, moment measure comparison, and normalized capture needed to apply Proposition~\ref{norm:transfer}. Theorem~\ref{cap:sampling} is the key input for the low-temperature part of the proof of Theorem~\ref{expand:main}, given in Subsection~\ref{sec:main.proof}, where it allows the Haar overlap laws and approximate GG identities to be transferred to the Euler Gibbs measure. Figure~\ref{fig:capped.dependencies} summarizes these proof dependencies.

\subsubsection{Capped measures and their moments}\label{cap:retained}

Fix a smooth function $\chi : \R \to [0,1]$ with compact support, a number $0 < \e < 1/4$, and a cap parameter $C_{0} > 0$. Define
\begin{equation}\label{cap:retained.definition}
f(x) \leqdef e^{\beta x}\chi(x), \qquad \mu_{\star}^{\chi}(\rd{h}) \leqdef \ell_{T}^{-1}f(F_{h}^{\star} - u_{T})\rd{h}.
\end{equation}
For $j \in \Z$, let $B_{j}^{\e} \leqdef [j + \e,j + 1 - \e] \cap I_{\theta}$, omitting empty intervals, and let $\mathcal{B}_{T}^{\e}$ denote the resulting collection. Distinct members of this collection are separated by at least $2\e$. Define
\begin{equation}\label{cap:definition}
\varpi_{C_{0}}(x) \leqdef \frac{C_{0}}{C_{0} + x}, \quad x\geq 0, \qquad \nu_{\star}(\rd{h}) \leqdef \sum_{B \in \mathcal{B}_{T}^{\e}}\varpi_{C_{0}}(\mu_{\star}^{\chi}(B))\mu_{\star}^{\chi}|_{B}(\rd{h}).
\end{equation}
The same measure $\nu_{\star}$ is used to form $\nu_{\star}^{\times r}$ for every $r \in \N$. Its mass on each interval $B \in \mathcal{B}_{T}^{\e}$ is at most $C_{0}$.

\begin{lemma}[Moments of the capped mass]\label{cap:moments}
Let $\theta > 0$ and $\beta > \beta_{c}(\theta)$. There is a constant $C_{\mathrm{mom}} = C_{\mathrm{mom}}(\beta,\chi,C_{0}) < \infty$, independent of $r$, such that, for every $r \in \N$ and every $\star \in \{\mathrm{A},\mathrm{H}\}$,
\[
\limsup_{T \to \infty}\EE_{\star}[\nu_{\star}(I_{\theta})^{r}] \leq C_{\mathrm{mom}}^{r}r!.
\]
\end{lemma}

\begin{proof}[Proof of Lemma~\ref{cap:moments}]
Since $\nu_{\star}(B) \leq C_{0}$ and $\nu_{\star}(B) \leq \mu_{\star}^{\chi}(B)$, each power $\nu_{\star}(B)^{q}$ is bounded by $C_{0}^{q - 1}\mu_{\star}^{\chi}(B)$. For $B \in \mathcal{B}_{T}^{\e}$, \eqref{cap:definition} gives
\begin{equation}\label{cap:local.detail:interval.mass}
0 \leq \nu_{\star}(B) = \frac{C_{0}\mu_{\star}^{\chi}(B)}{C_{0} + \mu_{\star}^{\chi}(B)} \leq \min\{C_{0},\mu_{\star}^{\chi}(B)\}.
\end{equation}
By \eqref{cap:local.detail:interval.mass}, for every $q \in \N$,
\begin{equation}\label{cap:local.detail:single.interval.moment}
\nu_{\star}(B)^{q} = \nu_{\star}(B)^{q - 1}\nu_{\star}(B) \leq C_{0}^{q - 1}\mu_{\star}^{\chi}(B).
\end{equation}
The collection $\mathcal{B}_{T}^{\e}$ is finite for fixed $T$, and
\[
\nu_{\star}(I_{\theta})^{r} = \sum_{(B_{1},\ldots,B_{r}) \in (\mathcal{B}_{T}^{\e})^{r}}\prod_{i = 1}^{r}\nu_{\star}(B_{i}).
\]
Group these ordered $r$-tuples according to the partition of $\{1,\ldots,r\}$ determined by equality of their interval labels. Fix a partition with $k$ nonempty groups, order its groups by their least elements, and let their sizes be $n_{1},\ldots,n_{k}$, where $n_{1} + \cdots + n_{k} = r$. Each assignment of $k$ distinct intervals to these ordered groups occurs exactly once, so the partition contributes
\begin{equation}\label{cap:local.detail:partition.moment}
\sum_{\substack{B_{1},\ldots,B_{k} \in \mathcal{B}_{T}^{\e} \\ B_{1},\ldots,B_{k}\text{ distinct}}}\prod_{i = 1}^{k}\nu_{\star}(B_{i})^{n_{i}} \leq C_{0}^{r - k}\sum_{\substack{B_{1},\ldots,B_{k} \in \mathcal{B}_{T}^{\e} \\ B_{1},\ldots,B_{k}\text{ distinct}}}\prod_{i = 1}^{k}\mu_{\star}^{\chi}(B_{i}),
\end{equation}
by \eqref{cap:local.detail:single.interval.moment}. The definition \eqref{cap:retained.definition}, positivity, and the separation of distinct intervals in $\mathcal{B}_{T}^{\e}$ give
\[
\begin{aligned}
\sum_{\substack{B_{1},\ldots,B_{k} \in \mathcal{B}_{T}^{\e} \\ B_{1},\ldots,B_{k}\text{ distinct}}}\prod_{i = 1}^{k}\mu_{\star}^{\chi}(B_{i})
&= \ell_{T}^{-k}\sum_{\substack{B_{1},\ldots,B_{k} \in \mathcal{B}_{T}^{\e} \\ B_{1},\ldots,B_{k}\text{ distinct}}}\int_{B_{1}\times\cdots\times B_{k}}\prod_{i = 1}^{k}f(F_{h_{i}}^{\star} - u_{T})\rd{\bb{h}} \\
&\leq \ell_{T}^{-k}\int_{\min_{i\neq i'}|h_{i} - h_{i'}|\geq 2\e}\prod_{i = 1}^{k}f(F_{h_{i}}^{\star} - u_{T})\rd{\bb{h}}.
\end{aligned}
\]
For $k = 1$, the separation condition imposes no restriction. The height estimate at separated positions now controls the sum over distinct intervals. Tonelli's theorem and \eqref{an:remote.bound} in Lemma~\ref{an:remote}, with $d_{0} = 2\e$ and $f(x) = e^{\beta x}\chi(x)$, therefore give
\begin{equation}\label{cap:local.detail:distinct.interval.masses}
\limsup_{T\to\infty}\EE_{\star}\left[\sum_{\substack{B_{1},\ldots,B_{k} \in \mathcal{B}_{T}^{\e} \\ B_{1},\ldots,B_{k}\text{ distinct}}}\prod_{i = 1}^{k}\mu_{\star}^{\chi}(B_{i})\right] \leq \Theta(f)^{k}, \qquad \Theta(f) = \int_{\R}f(x)\beta_{c}(\theta)e^{-\beta_{c}(\theta)x}\rd{x}.
\end{equation}
There are $\left\{\begin{smallmatrix}r\\k\end{smallmatrix}\right\}$ partitions into $k$ groups, where the braces denote Stirling numbers of the second kind. Summing \eqref{cap:local.detail:partition.moment} over these finitely many partitions and using \eqref{cap:local.detail:distinct.interval.masses} gives
\begin{equation}\label{cap:stirling}
\limsup_{T \to \infty}\EE_{\star}[\nu_{\star}(I_{\theta})^{r}] \leq \sum_{k = 1}^{r}\left\{\begin{matrix}r\\k\end{matrix}\right\}C_{0}^{r - k}\Theta(f)^{k}.
\end{equation}

A Poisson moment identity bounds the right-hand side of \eqref{cap:stirling} uniformly in $r$. Indeed, let $N_{\mathrm{Pois}}$ have the Poisson distribution with mean $\Theta(f)/C_{0}$. Its $k$th falling factorial moment is $(\Theta(f)/C_{0})^{k}$. Expanding the ordinary power in falling factorials therefore gives
\[
C_{0}^{r}\EE[N_{\mathrm{Pois}}^{r}]
= C_{0}^{r}\sum_{k = 1}^{r}\left\{\begin{matrix}r\\k\end{matrix}\right\}\EE[N_{\mathrm{Pois}}(N_{\mathrm{Pois}} - 1)\cdots(N_{\mathrm{Pois}} - k + 1)]
= \sum_{k = 1}^{r}\left\{\begin{matrix}r\\k\end{matrix}\right\}C_{0}^{r - k}\Theta(f)^{k}.
\]
For $x \geq 0$, the term $x^{r}/r!$ in the Taylor series for $e^{x}$ gives $x^{r} \leq r!e^{x}$. The Poisson exponential moment consequently yields
\begin{equation}\label{cap:local.detail:poisson.moment}
C_{0}^{r}\EE[N_{\mathrm{Pois}}^{r}] \leq C_{0}^{r}r!\EE[e^{N_{\mathrm{Pois}}}]
= C_{0}^{r}r!\exp\left(\frac{\Theta(f)}{C_{0}}(e - 1)\right)
\leq C_{\mathrm{mom}}^{r}r!,
\end{equation}
with $C_{\mathrm{mom}} \leqdef \max\{1,C_{0}\}\exp((\Theta(f)/C_{0})(e - 1))$. Combining \eqref{cap:local.detail:poisson.moment} with \eqref{cap:stirling} proves the lemma.
\end{proof}

\subsubsection{Localization within individual intervals}\label{cap:localization}

For a tuple $\bb{h} \in I_{\theta}^{\times k}$ whose coordinates lie in distinct intervals of $\mathcal{B}_{T}^{\e}$, let $B(h_{i})$ denote the interval containing $h_{i}$. Let $\mathscr{S}_{k,T}^{\e}$ be the set of all such tuples. For $R \geq 1$, define
\begin{equation}\label{cap:localization.outer.mass}
M_{\star,i}^{\mathrm{out}}(R) \leqdef \mu_{\star}^{\chi}\big(B(h_{i}) \setminus [h_{i} - R\ell_{T},h_{i} + R\ell_{T}]\big).
\end{equation}
This is the $\mu_{\star}^{\chi}$-mass of the portion of $B(h_{i})$ lying at distance greater than $R\ell_{T}$ from $h_{i}$. The next lemma quantifies how the interval mass localizes near its representative position.

\begin{lemma}[Integrated localization estimate]\label{cap:localization.bound}
Let $\theta > 0$ and $\beta > \beta_{c}(\theta)$. For every fixed $k \in \N$, there exist $0 < a \leq 1$, $c > 0$, and $C < \infty$ such that, for $\star \in \{\mathrm{A},\mathrm{H}\}$, $1 \leq i \leq k$ and $R \geq 1$,
\[
\limsup_{T \to \infty}\EE_{\star}\left[\int_{\mathscr{S}_{k,T}^{\e}}M_{\star,i}^{\mathrm{out}}(R)^{a}\prod_{q = 1}^{k}\mu_{\star}^{\chi}(\rd{h}_{q})\right] \leq CR^{-c}.
\]
The constants may depend on $k$, $\e$, $\beta$, and $\chi$.
\end{lemma}

\begin{proof}[Proof of Lemma~\ref{cap:localization.bound}]
First, the conditional estimate \eqref{an:joint.annular.bound} in Lemma~\ref{an:joint.annular} is multiplied by $\PP_{\star}(E_{\bb{h},\bb{v}})$, the probability of the joint rare-height event defined below in \eqref{eq:rare.height.events}, and then integrated over the representative positions. Choose a compact interval $J$ containing the support of $f(x) = e^{\beta x}\chi(x)$, cover it by finitely many closed bins $J_{1},\ldots,J_{L}$ of one fixed positive width, and let $J'$ be a compact interval containing all these bins. For $\bb{v} \leqdef (v_{1},\ldots,v_{k}) \in \{1,\ldots,L\}^{k}$, define the joint bin event
\begin{equation}\label{eq:rare.height.events}
E_{\bb{h},\bb{v}} \leqdef \bigcap_{q = 1}^{k}\{F_{h_{q}}^{\star} - u_{T} \in J_{v_{q}}\}.
\end{equation}
Only tuples for which all centered representative heights $F_{h_{q}}^{\star} - u_{T}$ belong to $J$ contribute to the integral in Lemma~\ref{cap:localization.bound}. On these tuples,
\begin{equation}\label{cap:local.detail:bin.domination}
\prod_{q = 1}^{k}f(F_{h_{q}}^{\star} - u_{T}) \leq \|f\|_{\infty}^{k}\sum_{\bb{v} \in \{1,\ldots,L\}^{k}}\ind_{E_{\bb{h},\bb{v}}}.
\end{equation}
Using a cover allows a height at a bin endpoint to occur in more than one event; the resulting finite multiplicity is harmless in \eqref{cap:local.detail:bin.domination}.

The hypotheses of Lemma~\ref{an:joint.annular} hold for $E_{\bb{h},\bb{v}}$: the representatives lie in distinct intervals of $\mathcal{B}_{T}^{\e}$, their separation is at least $2\e$, and all bins have the same fixed width and lie in $J'$. The decay in \eqref{an:joint.annular.bound} therefore applies uniformly in the positions and bin choices. To make the fractional exponent explicit, choose $1 < m_{-} < m(\theta)$ and then
\[
0 < a < \min\{1,2/\beta,2(m_{-} - 1)/\beta\}.
\]
The restriction $a < 2/\beta$ ensures that the H\"{o}lder exponents used to derive \eqref{an:joint.annular.exponent} exceed one; see \eqref{eq:choice.a.Holder}. For this choice, the exponent in \eqref{an:joint.annular.exponent} satisfies
\begin{equation}\label{cap:local.detail:negative.exponent}
a\beta(1 - m_{-}) + \frac{a^{2}\beta^{2}}{4} \leq - \frac{a\beta(m_{-} - 1)}{2} < 0.
\end{equation}
For $\ell_{T} \leq d \leq 1$, the definitions \eqref{cap:retained.definition} and \eqref{cap:localization.outer.mass} identify the annular contribution through
\begin{equation}\label{cap:local.detail:annular.mass}
\begin{aligned}
&\mu_{\star}^{\chi}(\{y \in B(h_{i}) : d \leq |y - h_{i}| \leq 2d\}) \\
&\qquad = \ell_{T}^{-1}\int_{\substack{y \in B(h_{i}) \\ d \leq |y - h_{i}| \leq 2d}}e^{\beta(F_{y}^{\star} - u_{T})}\chi(F_{y}^{\star} - u_{T})\rd{y} \\
&\qquad \leq \|\chi\|_{\infty}\exp\{\beta(F_{h_{i}}^{\star} - u_{T})\}\ell_{T}^{-1}\int_{\substack{y \in B(h_{i}) \\ d \leq |y - h_{i}| \leq 2d}}\exp\{\beta(F_{y}^{\star} - F_{h_{i}}^{\star})\}\rd{y}.
\end{aligned}
\end{equation}
On $E_{\bb{h},\bb{v}}$, the prefactor in \eqref{cap:local.detail:annular.mass} is bounded by $\|\chi\|_{\infty}e^{\beta\sup J'}$. Apply \eqref{an:joint.annular.exponent} for $\ell_{T} \leq d \leq t^{-3}$ and \eqref{an:annular.detail:large.radius} for $t^{-3} < d \leq 1$, both in the proof of Lemma~\ref{an:joint.annular}, with $r = \log(d/\ell_{T})$. Together with \eqref{cap:local.detail:negative.exponent}, these estimates give, after possibly reducing $c > 0$,
\begin{equation}\label{cap:local.detail:annular.moment}
\EE_{\star}[\mu_{\star}^{\chi}(\{y \in B(h_{i}) : d \leq |y - h_{i}| \leq 2d\})^{a}\ind_{\mathcal{A}_{T,\mathrm{sp}}} \mid E_{\bb{h},\bb{v}}] \leq C(d/\ell_{T})^{-c}.
\end{equation}
Here $\mathcal{A}_{T,\mathrm{sp}}$ is the common spatial event under the arithmetic law and the whole probability space under the Haar law.

The part of $B(h_{i})$ lying at distance greater than $R\ell_{T}$ from $h_{i}$ is covered by annuli centered at $h_{i}$ with inner radius $2^{v}R\ell_{T}$ and outer radius $2^{v + 1}R\ell_{T}$, for $v \in \N_0$ satisfying $2^{v}R\ell_{T} \leq 1$. Thus \eqref{cap:localization.outer.mass} gives
\begin{equation}\label{cap:local.detail:outer.decomposition}
M_{\star,i}^{\mathrm{out}}(R) \leq \sum_{\substack{v = 0\\2^{v}R\ell_{T}\leq 1}}^{\infty}\mu_{\star}^{\chi}(\{y \in B(h_{i}) : 2^{v}R\ell_{T} \leq |y - h_{i}| \leq 2^{v + 1}R\ell_{T}\}).
\end{equation}
Choose the arithmetic event $\smash{\mathcal{A}_{T,\mathrm{sp}}^{\mathrm{A}}}$ simultaneously for the fields $\{F_{h}^{\star}\}_{h \in I_{\theta}}$, $\{\ell_{T}(F_{h}^{\star})'\}_{h \in I_{\theta}}$, and $\{F_{d}^{\star,\mathrm{lo}}(h)\}_{h \in I_{\theta}}$ at the $\OO(t)$ radii $d = 2^{v}R\ell_{T} \leq 1$, as in the proof of Lemma~\ref{an:joint.annular}. This event satisfies \eqref{an:spatial.exceptional} and \eqref{an:spatial.transform.bound} in Lemma~\ref{an:spatial.transform}. In \eqref{cap:local.detail:outer.moment}, $\mathcal{A}_{T,\mathrm{sp}}$ denotes this arithmetic event when $\star = \mathrm{A}$ and the whole probability space when $\star = \mathrm{H}$. Applying subadditivity of $x \mapsto x^{a}$ to \eqref{cap:local.detail:outer.decomposition} and then using \eqref{cap:local.detail:annular.moment} give
\begin{equation}\label{cap:local.detail:outer.moment}
\EE_{\star}[M_{\star,i}^{\mathrm{out}}(R)^{a}\ind_{\mathcal{A}_{T,\mathrm{sp}}} \mid E_{\bb{h},\bb{v}}] \leq C\sum_{v = 0}^{\infty}(2^{v}R)^{-c} = \frac{C}{1 - 2^{-c}}R^{-c}.
\end{equation}
If $R\ell_{T} \geq 1$, the outer mass is zero and \eqref{cap:local.detail:outer.moment} holds immediately.

For completeness, the arithmetic contribution outside $\mathcal{A}_{T,\mathrm{sp}}$ is negligible after conditioning. The support bound on $f$ and $|B(h_{i})| \leq 1$ give $M_{\mathrm{A},i}^{\mathrm{out}}(R)^{a} \leq (\|f\|_{\infty}\ell_{T}^{-1})^{a}$. Combining the exceptional probability \eqref{an:spatial.exceptional} in Lemma~\ref{an:spatial.transform} with the bin lower bound \eqref{an:detail:conditional.bin.lower} and the fixed-width extension \eqref{an:conditional.detail:variable.bin.integrals} in the proof of Lemma~\ref{an:conditional} gives
\begin{equation}\label{cap:local.detail:event.removal}
\begin{aligned}
\EE[M_{\mathrm{A},i}^{\mathrm{out}}(R)^{a}\ind_{\mathcal{A}_{T,\mathrm{sp}}^{c}} \mid E_{\bb{h},\bb{v}}]
&\leq (\|f\|_{\infty}\ell_{T}^{-1})^{a}\frac{\PP(\mathcal{A}_{T,\mathrm{sp}}^{c})}{\PP(E_{\bb{h},\bb{v}})} \\
&\leq \exp\{ - ct^{1 + 2\e_{0}} + \OO_{k,J'}(t)\} = \oo(1),
\end{aligned}
\end{equation}
uniformly in the positions and bins. The scale bound $\log(\ell_{T}^{-1}) = \OO(t)$ follows from \eqref{eq:effective.cutoff.bounds} and \eqref{eq:parameter.scales}. Combining \eqref{cap:local.detail:outer.moment} and \eqref{cap:local.detail:event.removal} removes the arithmetic event; under the Haar law, \eqref{cap:local.detail:outer.moment} already holds without an indicator. Multiplying the resulting conditional estimate by the joint-bin probability gives
\begin{equation}\label{cap:local.detail:joint.bin.moment}
\begin{aligned}
\EE_{\star}[M_{\star,i}^{\mathrm{out}}(R)^{a}\ind_{E_{\bb{h},\bb{v}}}]
&= \PP_{\star}(E_{\bb{h},\bb{v}})\EE_{\star}[M_{\star,i}^{\mathrm{out}}(R)^{a} \mid E_{\bb{h},\bb{v}}] \\
&\leq \{CR^{-c} + \oo(1)\}\PP_{\star}(E_{\bb{h},\bb{v}}),
\end{aligned}
\end{equation}
uniformly in the representative positions $\bb{h} \in \mathscr{S}_{k,T}^{\e}$ and the finitely many bin choices, for every fixed $R$.

It remains to integrate the inequality \eqref{cap:local.detail:joint.bin.moment} over the representative positions. Tonelli's theorem, \eqref{cap:local.detail:bin.domination}, and the density in \eqref{cap:retained.definition} give
\begin{equation}\label{cap:local.detail:integrated.bins}
\EE_{\star}\left[\int_{\mathscr{S}_{k,T}^{\e}}M_{\star,i}^{\mathrm{out}}(R)^{a}\prod_{q = 1}^{k}\mu_{\star}^{\chi}(\rd{h}_{q})\right]
\leq \|f\|_{\infty}^{k}\{CR^{-c} + \oo(1)\}\ell_{T}^{-k}\sum_{\bb{v} \in \{1,\ldots,L\}^{k}}\int_{\mathscr{S}_{k,T}^{\e}}\PP_{\star}(E_{\bb{h},\bb{v}})\rd{\bb{h}}.
\end{equation}
Each bin lies in $J'$, and there are $L^{k}$ choices of $\bb{v}$. Every tuple in $\mathscr{S}_{k,T}^{\e}$ has separation at least $2\e$, so \eqref{an:remote.bound} in Lemma~\ref{an:remote}, applied to $\ind_{J'}$, yields
\begin{equation}\label{cap:local.detail:bin.probability}
\begin{aligned}
&\limsup_{T\to\infty}\ell_{T}^{-k}\sum_{\bb{v} \in \{1,\ldots,L\}^{k}}\int_{\mathscr{S}_{k,T}^{\e}}\PP_{\star}(E_{\bb{h},\bb{v}})\rd{\bb{h}} \\
&\qquad \leq L^{k}\limsup_{T\to\infty}\ell_{T}^{-k}\int_{\min_{q\neq q'}|h_{q} - h_{q'}|\geq 2\e}\EE_{\star}\left[\prod_{q = 1}^{k}\ind_{J'}(F_{h_{q}}^{\star} - u_{T})\right]\rd{\bb{h}} \\
&\qquad \leq L^{k}\Theta(\ind_{J'})^{k} < \infty.
\end{aligned}
\end{equation}
The indicator of the compact interval $J'$ is admissible in \eqref{an:remote.bound} of Lemma~\ref{an:remote} because its discontinuities are confined to its endpoints. Applying \eqref{cap:local.detail:bin.probability} in \eqref{cap:local.detail:integrated.bins} proves the lemma, with a constant that absorbs $L^{k}$, $\|f\|_{\infty}^{k}$ and $\Theta(\ind_{J'})^{k}$.
\end{proof}

The following purely algebraic lemma is used in the proof of Proposition~\ref{cap:moment.measures} to control the errors from restricting sampled positions to local windows and replacing interval masses by the masses within these windows.

\begin{lemma}[Inequalities for grouped cap factors]\label{cap:inequalities}
Let $C_{0} > 0$ and $0 < a \leq 1$ be real numbers. Recall from \eqref{cap:definition} that $\varpi_{C_{0}}(x) = C_{0}/(C_{0} + x)$ for real $x \geq 0$. For real numbers $M \geq M^{\mathrm{out}} \geq 0$ and an integer $n \geq 2$,
\begin{equation}\label{cap:discard.group}
\varpi_{C_{0}}(M)^{n}M^{n - 2}M^{\mathrm{out}} \leq C_{0}^{n - 1 - a}(M^{\mathrm{out}})^{a}.
\end{equation}
For real numbers $x,y \geq 0$ and $n \in \N$,
\begin{equation}\label{cap:replace.group}
x^{n - 1}|\varpi_{C_{0}}(x)^{n} - \varpi_{C_{0}}(x + y)^{n}| \leq nC_{0}^{n - 1}\min\{1,y/C_{0}\} \leq nC_{0}^{n - 1 - a}y^{a}.
\end{equation}
\end{lemma}

\begin{proof}[Proof of Lemma~\ref{cap:inequalities}]
For every $M \geq 0$, the rational cap $\varpi_{C_{0}}(M)$ satisfies
\begin{equation}\label{cap:local.detail:basic.cap.bounds}
0 < \varpi_{C_{0}}(M) \leq 1, \qquad M\varpi_{C_{0}}(M) = \frac{C_{0}M}{C_{0} + M} \leq C_{0}.
\end{equation}
Since $M^{\mathrm{out}} \leq M$, the inequalities in \eqref{cap:local.detail:basic.cap.bounds} give
\begin{equation}\label{cap:local.detail:cap.square}
\begin{aligned}
\varpi_{C_{0}}(M)^{2}M^{\mathrm{out}} &\leq M^{\mathrm{out}} \\
\varpi_{C_{0}}(M)^{2}M^{\mathrm{out}} \leq \varpi_{C_{0}}(M)^{2}M \leq C_{0}\varpi_{C_{0}}(M) &\leq C_{0}.
\end{aligned}
\end{equation}
For every $z \geq 0$ and $0 < a \leq 1$, splitting into $z \leq C_{0}$ and $z \geq C_{0}$ shows that $\min\{z,C_{0}\} \leq C_{0}^{1 - a}z^{a}$. Therefore, using \eqref{cap:local.detail:basic.cap.bounds} and \eqref{cap:local.detail:cap.square}, for $n \geq 2$,
\[
\begin{aligned}
\varpi_{C_{0}}(M)^{n}M^{n - 2}M^{\mathrm{out}} &= \{M\varpi_{C_{0}}(M)\}^{n - 2}\varpi_{C_{0}}(M)^{2}M^{\mathrm{out}} \\
&\leq C_{0}^{n - 2}\min\{M^{\mathrm{out}},C_{0}\} \\
&\leq C_{0}^{n - 1 - a}(M^{\mathrm{out}})^{a},
\end{aligned}
\]
which proves \eqref{cap:discard.group}.

For \eqref{cap:replace.group}, the function $\varpi_{C_{0}}$ is decreasing, so $\varpi_{C_{0}}(x)^{n} - \varpi_{C_{0}}(x + y)^{n} \geq 0$. For $z \geq 0$,
\[
0 \leq 1 - (1 + z)^{-n} = n\int_{0}^{z}(1 + u)^{-n - 1}\rd{u} \leq nz.
\]
Combining this bound with $1 - (1 + z)^{-n} \leq 1$ gives
\begin{equation}\label{cap:local.detail:power.difference}
0 \leq 1 - (1 + z)^{-n} \leq \min\{1,nz\} \leq n\min\{1,z\}.
\end{equation}
Applying \eqref{cap:local.detail:power.difference} with $z = y/(C_{0} + x)$ gives
\begin{equation}\label{cap:local.detail:cap.power.difference}
\begin{aligned}
0 \leq \varpi_{C_{0}}(x)^{n} - \varpi_{C_{0}}(x + y)^{n} &= \varpi_{C_{0}}(x)^{n}\left\{1 - \left(1 + \frac{y}{C_{0} + x}\right)^{-n}\right\} \\
&\leq n\varpi_{C_{0}}(x)^{n}\min\left\{1,\frac{y}{C_{0} + x}\right\}.
\end{aligned}
\end{equation}
Multiplying \eqref{cap:local.detail:cap.power.difference} by $x^{n - 1}$, and then using \eqref{cap:local.detail:basic.cap.bounds} and $C_{0} + x \geq C_{0}$, gives
\[
\begin{aligned}
x^{n - 1}|\varpi_{C_{0}}(x)^{n} - \varpi_{C_{0}}(x + y)^{n}| &\leq n\{x\varpi_{C_{0}}(x)\}^{n - 1}\varpi_{C_{0}}(x)\min\left\{1,\frac{y}{C_{0} + x}\right\} \\
&\leq nC_{0}^{n - 1}\min\{1,y/C_{0}\} \\
&\leq nC_{0}^{n - 1 - a}y^{a}.
\end{aligned}
\]
The last line uses $\min\{1,z\} \leq z^{a}$ for $z \geq 0$ and $0 < a \leq 1$, with $z = y/C_{0}$. This proves both inequalities in \eqref{cap:replace.group}. At $n = 1$, the factors with exponent zero are one, so the calculation also covers $x = 0$.
\end{proof}

\subsubsection{Comparison of moment measures}\label{cap:comparison}

\begin{proposition}[Moment measure comparison]\label{cap:moment.measures}
Let $\theta > 0$ and $\beta > \beta_{c}(\theta)$. Fix a smooth function $\chi : \R \to [0,1]$ with compact support and real numbers $0 < \e < 1/4$ and $C_{0} > 0$. For the measures $\nu_{\mathrm{A}}$ and $\nu_{\mathrm{H}}$ defined by \eqref{cap:definition}, and every fixed $r \in \N$, as $T \to \infty$,
\[
\left\|\EE\nu_{\mathrm{A}}^{\times r} - \EE_{\mathrm{H}}\nu_{\mathrm{H}}^{\times r}\right\|_{\mathrm{TV}} \to 0.
\]
\end{proposition}

\begin{proof}[Proof of Proposition~\ref{cap:moment.measures}]
Fix $r\in \N$. Partition the position space according to which coordinates belong to the same interval of $\mathcal{B}_{T}^{\e}$. Fix one partition of $\{1,\ldots,r\}$ into $k$ groups of sizes $n_{1},\ldots,n_{k}$, with $\sum_{i = 1}^{k}n_{i} = r$, and choose $a \in (0,1]$ and $c > 0$ from Lemma~\ref{cap:localization.bound} for this value of $k$. Select the smallest label in each group as its representative, and denote the representative positions by $y_{1},\ldots,y_{k}$. The representatives belong to $\mathscr{S}_{k,T}^{\e}$, and every other position in group $i$ belongs to $B(y_{i})$. By the definition of the capped measure $\nu_{\star}$ in \eqref{cap:definition}, the contribution of this partition to $\nu_{\star}^{\times r}$ has density
\begin{equation}\label{cap:comparison.detail:partition.density}
\frac{\rd{\nu_{\star}^{\times r}}}{\rd{h}_{1}\cdots\rd{h}_{r}}(h_{1},\ldots,h_{r})
= \prod_{i = 1}^{k}\left(\varpi_{C_{0}}(\mu_{\star}^{\chi}(B(y_{i})))^{n_{i}}\prod_{q \text{ in group }i}\ell_{T}^{-1}f(F_{h_{q}}^{\star} - u_{T})\right).
\end{equation}
All errors below are integrated absolute errors in \eqref{cap:comparison.detail:partition.density}. For an arbitrary deterministic measurable function $\Psi_{T}$ satisfying $\|\Psi_{T}\|_{\infty} \leq 1$, multiplication of an error density by $\Psi_{T}$ cannot increase its integrated absolute value. Thus the resulting bounds are uniform in the test function required by the total variation norm.

\emph{Step 1. Localization of the sampled positions.} Localize the nonrepresentative positions, leaving each cap factor $\varpi_{C_{0}}(\mu_{\star}^{\chi}(B(y_{i})))^{n_{i}}$ unchanged. For $R \geq 1$ and each group of size $n_{i} \geq 2$, restrict them to $B(y_{i}) \cap [y_{i} - R\ell_{T},y_{i} + R\ell_{T}]$. The error in this step is the mass of the configurations discarded by these restrictions. If the restriction fails, at least one of the $n_{i} - 1$ nonrepresentative positions lies outside this local window, whose complement in $B(y_{i})$ has mass $M_{\star,i}^{\mathrm{out}}(R) = \mu_{\star}^{\chi}(B(y_{i}) \setminus [y_{i} - R\ell_{T},y_{i} + R\ell_{T}])$. After integrating the nonrepresentative positions, the loss from group $i$ is the first expression below. A union bound and \eqref{cap:discard.group} in Lemma~\ref{cap:inequalities}, applied with $M = \mu_{\star}^{\chi}(B(y_{i}))$, $M^{\mathrm{out}} = M_{\star,i}^{\mathrm{out}}(R)$, and $n = n_{i}$, give
\begin{equation}\label{cap:comparison.detail:group.loss}
\begin{aligned}
&\mu_{\star}^{\chi}(\rd{y}_{i})\varpi_{C_{0}}(\mu_{\star}^{\chi}(B(y_{i})))^{n_{i}}
\left[\mu_{\star}^{\chi}(B(y_{i}))^{n_{i} - 1} - \left\{\mu_{\star}^{\chi}(B(y_{i})) - M_{\star,i}^{\mathrm{out}}(R)\right\}^{n_{i} - 1}\right] \\
&\qquad \leq (n_{i} - 1)\mu_{\star}^{\chi}(\rd{y}_{i})\varpi_{C_{0}}(\mu_{\star}^{\chi}(B(y_{i})))^{n_{i}}\mu_{\star}^{\chi}(B(y_{i}))^{n_{i} - 2}M_{\star,i}^{\mathrm{out}}(R) \\
&\qquad \leq (n_{i} - 1)C_{0}^{n_{i} - 1 - a}\mu_{\star}^{\chi}(\rd{y}_{i})M_{\star,i}^{\mathrm{out}}(R)^{a}.
\end{aligned}
\end{equation}
For every group $q \in \{1,\ldots,k\} \setminus \{i\}$, integrating its nonrepresentative positions over the whole interval $B(y_{q})$ gives the bound
\begin{equation}\label{cap:comparison.detail:other.groups}
\varpi_{C_{0}}(\mu_{\star}^{\chi}(B(y_{q})))^{n_{q}}\mu_{\star}^{\chi}(B(y_{q}))^{n_{q} - 1}\mu_{\star}^{\chi}(\rd{y}_{q}) \leq C_{0}^{n_{q} - 1}\mu_{\star}^{\chi}(\rd{y}_{q}),
\end{equation}
by \eqref{cap:local.detail:basic.cap.bounds}. Multiplying \eqref{cap:comparison.detail:group.loss} and \eqref{cap:comparison.detail:other.groups}, summing over the groups with nonrepresentative positions, and applying Lemma~\ref{cap:localization.bound} give the following upper bound for the limiting expected loss from localization:
\begin{equation}\label{cap:comparison.detail:localization.error}
\limsup_{T \to \infty}\sum_{\substack{1 \leq i \leq k\\n_{i}\geq 2}}(n_{i} - 1)C_{0}^{r - k - a}\EE_{\star}\left[\int_{\mathscr{S}_{k,T}^{\e}}M_{\star,i}^{\mathrm{out}}(R)^{a}\prod_{q = 1}^{k}\mu_{\star}^{\chi}(\rd{y}_{q})\right]
\leq C_{r}R^{-c}.
\end{equation}
There are finitely many partitions and groups for fixed $r$, so $c > 0$ may be chosen as the minimum of their localization exponents. Groups of size one cause no loss in this step.

\emph{Step 2. Localization of the cap masses.} Localize the cap factors $\varpi_{C_{0}}(\mu_{\star}^{\chi}(B(y_{i})))^{n_{i}}$ themselves. All nonrepresentative positions are now restricted to their local windows, and the error in this step comes from changing the weights of the configurations that remain. Define
\[
M_{\star,i}^{\mathrm{loc}}(R) \leqdef \mu_{\star}^{\chi}\big(B(y_{i}) \cap [y_{i} - R\ell_{T},y_{i} + R\ell_{T}]\big).
\]
Then $\mu_{\star}^{\chi}(B(y_{i})) = M_{\star,i}^{\mathrm{loc}}(R) + M_{\star,i}^{\mathrm{out}}(R)$. In each group, replace the factor $\varpi_{C_{0}}(\mu_{\star}^{\chi}(B(y_{i})))^{n_{i}}$ by its localized counterpart $\varpi_{C_{0}}(M_{\star,i}^{\mathrm{loc}}(R))^{n_{i}}$. Integrating the $n_{i} - 1$ localized nonrepresentative positions and applying \eqref{cap:replace.group} in Lemma~\ref{cap:inequalities} with $x = M_{\star,i}^{\mathrm{loc}}(R)$, $y = M_{\star,i}^{\mathrm{out}}(R)$, and $n = n_{i}$ bounds the error for this group by
\begin{equation}\label{cap:comparison.detail:cap.factor.loss}
\begin{aligned}
&\mu_{\star}^{\chi}(\rd{y}_{i})M_{\star,i}^{\mathrm{loc}}(R)^{n_{i} - 1}\left|\varpi_{C_{0}}(M_{\star,i}^{\mathrm{loc}}(R))^{n_{i}} - \varpi_{C_{0}}(\mu_{\star}^{\chi}(B(y_{i})))^{n_{i}}\right| \\
&\qquad \leq n_{i}C_{0}^{n_{i} - 1 - a}\mu_{\star}^{\chi}(\rd{y}_{i})M_{\star,i}^{\mathrm{out}}(R)^{a}.
\end{aligned}
\end{equation}
The bound \eqref{cap:comparison.detail:cap.factor.loss} also holds when $n_{i} = 1$, with no nonrepresentative integration. For every group $q \in \{1,\ldots,k\} \setminus \{i\}$, monotonicity of $\varpi_{C_{0}}$ and \eqref{cap:local.detail:basic.cap.bounds} give
\begin{equation}\label{cap:comparison.detail:localized.group.bound}
\varpi_{C_{0}}(\mu_{\star}^{\chi}(B(y_{q})))^{n_{q}}M_{\star,q}^{\mathrm{loc}}(R)^{n_{q} - 1}
\leq \varpi_{C_{0}}(M_{\star,q}^{\mathrm{loc}}(R))^{n_{q}}M_{\star,q}^{\mathrm{loc}}(R)^{n_{q} - 1}
\leq C_{0}^{n_{q} - 1}.
\end{equation}
Thus \eqref{cap:comparison.detail:localized.group.bound} bounds the factor from group $q$ by $C_{0}^{n_{q} - 1}$ whether its cap uses the mass of the whole interval $B(y_{q})$ or only the mass in the local window. In the telescoping identity
\begin{equation}\label{cap:comparison.detail:product.difference}
\prod_{i = 1}^{k}a_{i} - \prod_{i = 1}^{k}b_{i} = \sum_{i = 1}^{k}(a_{i} - b_{i})\prod_{q < i}a_{q}\prod_{q > i}b_{q},
\end{equation}
take $a_{i} = \varpi_{C_{0}}(M_{\star,i}^{\mathrm{loc}}(R))^{n_{i}}$ and $b_{i} = \varpi_{C_{0}}(\mu_{\star}^{\chi}(B(y_{i})))^{n_{i}}$. In the $i$th term of \eqref{cap:comparison.detail:product.difference}, \eqref{cap:comparison.detail:cap.factor.loss} bounds the change in the cap factor of group $i$ when the full interval mass is replaced by the local mass, while \eqref{cap:comparison.detail:localized.group.bound} bounds the factors from all groups $q \neq i$. Integrating against $\smash{\prod_{q = 1}^{k}\mu_{\star}^{\chi}(\rd{y}_{q})}$ and applying Lemma~\ref{cap:localization.bound} therefore give the following upper bound for the limiting expected error from localizing the cap factors:
\begin{equation}\label{cap:comparison.detail:cap.localization.error}
\limsup_{T \to \infty}\sum_{i = 1}^{k}n_{i}C_{0}^{r - k - a}\EE_{\star}\left[\int_{\mathscr{S}_{k,T}^{\e}}M_{\star,i}^{\mathrm{out}}(R)^{a}\prod_{q = 1}^{k}\mu_{\star}^{\chi}(\rd{y}_{q})\right] \leq C_{r}R^{-c}.
\end{equation}
Since $r$ is fixed, there are finitely many partitions and groups. Taking the minimum of their localization exponents therefore gives a common $c > 0$ in \eqref{cap:comparison.detail:localization.error} and \eqref{cap:comparison.detail:cap.localization.error}.

\emph{Step 3. Quadrature approximation of the local masses.} By \eqref{cap:comparison.detail:localization.error} and \eqref{cap:comparison.detail:cap.localization.error}, it suffices to compare the Euler and Haar expressions after both localization steps, first for fixed $R$ and then letting $R \to \infty$. Approximate each local mass by a finite sum so that the resulting density depends on finitely many field values. Express each local window in the rescaled coordinate $v = (h - y_{i})/\ell_{T}$ and write
\[
I_{i,T}^{\mathrm{loc}}(y_{i}) \leqdef \{v \in [ - R,R] : y_{i} + \ell_{T}v \in B(y_{i})\}.
\]
The substitution $h = y_{i} + \ell_{T}v$ cancels the factor $\ell_{T}^{-1}$ in the definition of $\mu_{\star}^{\chi}$ and gives
\begin{equation}\label{cap:local.integral}
M_{\star,i}^{\mathrm{loc}}(R) = \int_{I_{i,T}^{\mathrm{loc}}(y_{i})}f(F_{y_{i} + \ell_{T}v}^{\star} - u_{T})\rd{v},
\end{equation}
where $f(x) = e^{\beta x}\chi(x)$ as defined in \eqref{cap:retained.definition}. For $\ell_{\mathrm{quad}} > 0$, partition $[ - R,R]$ into intervals of length at most $\ell_{\mathrm{quad}}$, intersect them with $I_{i,T}^{\mathrm{loc}}(y_{i})$, and select a node from each intersection of positive length. The nodes can, for example, be the midpoints of the intersections, which makes them deterministic measurable functions of $y_{i}$. Let $w_{i,v}$ be the length of the intersection associated with node $v$. Replace the integral in \eqref{cap:local.integral} by the following quadrature approximation, that is, a weighted sum of values of its integrand:
\begin{equation}\label{cap:comparison.detail:quadrature}
\widetilde{M}_{\star,i}^{\mathrm{loc}}(R,\ell_{\mathrm{quad}}) \leqdef \sum_{v}w_{i,v}f(F_{y_{i} + \ell_{T}v}^{\star} - u_{T}), \qquad w_{i,v} \geq 0, \qquad \sum_{v}w_{i,v} \leq 2R.
\end{equation}
The number of nodes is bounded solely in terms of $R$ and $\ell_{\mathrm{quad}}$, and does not depend on $T$ or the representative positions.

For an absolutely continuous function $f_{\mathrm{quad}}$ on a partition interval $I_{\mathrm{quad}}$, let $v_{I_{\mathrm{quad}}} \in I_{\mathrm{quad}}$ denote the chosen node. The fundamental theorem of calculus gives
\begin{equation}\label{cap:comparison.detail:one.interval.quadrature}
\begin{aligned}
\left|\int_{I_{\mathrm{quad}}}f_{\mathrm{quad}}(v)\rd{v} - |I_{\mathrm{quad}}|f_{\mathrm{quad}}(v_{I_{\mathrm{quad}}})\right| &\leq \int_{I_{\mathrm{quad}}}|f_{\mathrm{quad}}(v) - f_{\mathrm{quad}}(v_{I_{\mathrm{quad}}})|\rd{v} \\
&\leq |I_{\mathrm{quad}}|\int_{I_{\mathrm{quad}}}|f_{\mathrm{quad}}'(v)|\rd{v}.
\end{aligned}
\end{equation}
Apply \eqref{cap:comparison.detail:one.interval.quadrature} to $f_{\mathrm{quad}}(v) \leqdef f(F_{y_{i} + \ell_{T}v}^{\star} - u_{T})$. The chain rule gives $f_{\mathrm{quad}}'(v) = f'(F_{y_{i} + \ell_{T}v}^{\star} - u_{T})\ell_{T}(F^{\star})'(y_{i} + \ell_{T}v)$. Summing over the positive-length intersections $I_{\mathrm{quad}}$ used in \eqref{cap:comparison.detail:quadrature} yields
\begin{equation}\label{cap:comparison.detail:quadrature.error}
\begin{aligned}
|M_{\star,i}^{\mathrm{loc}}(R) - \widetilde{M}_{\star,i}^{\mathrm{loc}}(R,\ell_{\mathrm{quad}})|
&\leq \sum_{I_{\mathrm{quad}}}\left|\int_{I_{\mathrm{quad}}}f_{\mathrm{quad}}(v)\rd{v} - |I_{\mathrm{quad}}|f_{\mathrm{quad}}(v_{I_{\mathrm{quad}}})\right| \\
&\leq \ell_{\mathrm{quad}}\int_{I_{i,T}^{\mathrm{loc}}(y_{i})}|f'(F_{y_{i} + \ell_{T}v}^{\star} - u_{T})|\,|\ell_{T}(F^{\star})'(y_{i} + \ell_{T}v)|\rd{v} \\
&\leq \ell_{\mathrm{quad}}\|f'\|_{\infty}\int_{I_{i,T}^{\mathrm{loc}}(y_{i})}|\ell_{T}(F^{\star})'(y_{i} + \ell_{T}v)|\rd{v}.
\end{aligned}
\end{equation}
Also, direct differentiation gives
\begin{equation}\label{cap:comparison.detail:cap.derivative}
\left|\frac{\rd}{\rd{x}}\varpi_{C_{0}}(x)^{n}\right| = \frac{nC_{0}^{n}}{(C_{0} + x)^{n + 1}} \leq \frac{n}{C_{0}}, \qquad x \geq 0.
\end{equation}
The mean value theorem and \eqref{cap:comparison.detail:cap.derivative}, followed by \eqref{cap:comparison.detail:quadrature.error}, give
\begin{equation}\label{cap:comparison.detail:cap.quadrature.error}
\begin{aligned}
\left|\varpi_{C_{0}}(M_{\star,i}^{\mathrm{loc}}(R))^{n_{i}} - \varpi_{C_{0}}(\widetilde{M}_{\star,i}^{\mathrm{loc}}(R,\ell_{\mathrm{quad}}))^{n_{i}}\right|
&\leq \frac{n_{i}}{C_{0}}|M_{\star,i}^{\mathrm{loc}}(R) - \widetilde{M}_{\star,i}^{\mathrm{loc}}(R,\ell_{\mathrm{quad}})| \\
&\leq \frac{n_{i}\ell_{\mathrm{quad}}}{C_{0}}\|f'\|_{\infty}\int_{I_{i,T}^{\mathrm{loc}}(y_{i})}|\ell_{T}(F^{\star})'(y_{i} + \ell_{T}v)|\rd{v}.
\end{aligned}
\end{equation}
For the fixed partition, consider the density obtained from \eqref{cap:comparison.detail:partition.density} after both localization steps. Replace every cap factor $\varpi_{C_{0}}(M_{\star,i}^{\mathrm{loc}}(R))^{n_{i}}$ by $\smash{\varpi_{C_{0}}(\widetilde{M}_{\star,i}^{\mathrm{loc}}(R,\ell_{\mathrm{quad}}))^{n_{i}}}$, keeping the factors $\ell_{T}^{-1}f(F_{h_{q}}^{\star} - u_{T})$ unchanged. Let $\Delta_{T,\mathrm{quad}}^{\star}$ denote the integral of the absolute difference between the densities before and after this replacement. The integral is over all $r$ sampled positions: the representatives $(y_{1},\ldots,y_{k})$ belong to $\mathscr{S}_{k,T}^{\e}$, and each nonrepresentative position in group $i$ belongs to $B(y_{i}) \cap [y_{i} - R\ell_{T},y_{i} + R\ell_{T}]$. The cap factors do not vary with the nonrepresentative positions, so integrating those positions gives
\begin{equation}\label{eq:Delta.T.quad}
\begin{aligned}
\Delta_{T,\mathrm{quad}}^{\star}
&= \int_{\mathscr{S}_{k,T}^{\e}}
\left|\prod_{i = 1}^{k}\varpi_{C_{0}}(M_{\star,i}^{\mathrm{loc}}(R))^{n_{i}}
- \prod_{i = 1}^{k}\varpi_{C_{0}}(\widetilde{M}_{\star,i}^{\mathrm{loc}}(R,\ell_{\mathrm{quad}}))^{n_{i}}\right| \\
&\qquad\qquad \times \prod_{i = 1}^{k}\left\{M_{\star,i}^{\mathrm{loc}}(R)^{n_{i} - 1}\mu_{\star}^{\chi}(\rd{y}_{i})\right\}.
\end{aligned}
\end{equation}
Each factor $M_{\star,i}^{\mathrm{loc}}(R)^{n_{i} - 1}$ comes from integrating the $n_{i} - 1$ nonrepresentative positions of group $i$ over their local window, leaving only the $k$ representative positions to be integrated in the displayed expression.

To bound $\EE_{\star}[\Delta_{T,\mathrm{quad}}^{\star}]$, integration over the representative positions $\bb{y} = (y_{1},\ldots,y_{k}) \in \mathscr{S}_{k,T}^{\e}$ must retain the small probability that all centered heights $F_{y_{q}}^{\star} - u_{T}$ lie in the support of $f$. Fix a bounded interval $J$ with integer endpoints containing the support of $f$, and define
\[
E_{\bb{y},J}^{\star} \leqdef \{F_{y_{q}}^{\star} - u_{T} \in J \text{ for every }1 \leq q \leq k\}.
\]
Since the representatives lie in distinct intervals of $\mathcal{B}_{T}^{\e}$, their pairwise distances are at least $2\e$. Apply \eqref{an:remote.bound} in Lemma~\ref{an:remote} with $d_{0} = 2\e$, $k$ equal to the number of groups, and the test function $\ind_{J}$ in place of $f$. Since $\mathscr{S}_{k,T}^{\e} \subseteq \mathscr{D}_{k,T}(2\e)$, this gives
\begin{equation}\label{cap:comparison.detail:representative.probability}
\limsup_{T \to \infty}\ell_{T}^{-k}\int_{\mathscr{S}_{k,T}^{\e}}\PP_{\star}(E_{\bb{y},J}^{\star})\rd{\bb{y}} \leq \left(\int_{J}\beta_{c}(\theta)e^{-\beta_{c}(\theta)x}\rd{x}\right)^{k} < \infty.
\end{equation}
The bound \eqref{cap:comparison.detail:representative.probability} also holds for any fixed bounded enlargement of $J$.

To include the scaled derivative $\ell_{T}(F^{\star})'(y_{i} + \ell_{T}v)$ in this probability bound, use the common spatial event of Lemma~\ref{an:spatial.transform} for the two functions $F^{\star}(h)$ and $\ell_{T}(F^{\star})'(h)$. The required coefficient bounds follow from \eqref{an:annular.detail:coefficient.admissibility} in the proof of Lemma~\ref{an:annular}. Thus \eqref{an:spatial.exceptional} and \eqref{an:spatial.transform.bound} in Lemma~\ref{an:spatial.transform} hold with one event $\mathcal{A}_{T,\mathrm{sp}}^{\star}$ for all representative positions and all $v \in I_{i,T}^{\mathrm{loc}}(y_{i})$.

Partition $J$ into closed unit bins and condition on the bin of every representative height. For the application of Lemma~\ref{an:conditional}, take $\smash{Z_{T,q}^{\star} = F_{y_{q}}^{\star}}$, $1 \leq q \leq k$, and $\smash{Z_{T,k + 1}^{\star} \!= \ell_{T}(F^{\star})'(y_{i} + \ell_{T}v)}$. The derivative variance $\sig_{T}^{2} = \Var_{\mathrm{H}}(Z_{T,k + 1}^{\mathrm{H}})$ and its covariances $c_{T,q} = \Cov_{\mathrm{H}}(Z_{T,k + 1}^{\mathrm{H}},H(y_{q}))$ satisfy
\begin{equation}\label{cap:comparison.detail:derivative.covariances}
\sig_{T}^{2} \leq C + C\ell_{T}^{2}\int_{\log 2}^{\log X}w e^{-\ell_{T}w}\rd{w} \leq C, \qquad
|c_{T,q}| \leq C + C\ell_{T}\int_{\log 2}^{\log X}e^{-\ell_{T}w}\rd{w} \leq C.
\end{equation}
Here the bounded terms include the higher prime powers. Partial summation using \eqref{analytic:eq:PNT}, followed by \eqref{an:effective.cutoff.estimates} and $\ell_{T}\int_{0}^{\infty}e^{-\ell_{T}w}\rd{w} = 1$, bounds the prime contributions. Since $|y_{q} - y_{q'}| \geq 2\e$ for $q \neq q'$, \eqref{an:overlap.geometry} and \eqref{an:covariance.comparison} in Lemma~\ref{an:covariance}, together with \eqref{eq:parameter.scales} and \eqref{eq:effective.cutoff.bounds}, give
\begin{equation}\label{cap:comparison.detail:representative.covariance}
\Sigma_{T} = \big(\Cov_{\mathrm{H}}(H(y_{q}),H(y_{q'}))\big)_{1 \leq q,q' \leq k} = V_{T}\mathrm{Id}_{k} + \OO_{k,\e}(\log{t}) = \frac{t}{2}\mathrm{Id}_{k} + \OO_{k,\e}(\log{t}).
\end{equation}
In particular, the eigenvalues of $\Sigma_{T}$ lie between $t/4$ and $t$ for all sufficiently large $T$. The lower endpoints of the conditioning bins form a vector $\bb{\alpha} = u_{T}\bb{1}_{k} + \OO_{J}(1)$. Since $u_{T} = \OO(t)$ by \eqref{an:height.scales} in Lemma~\ref{an:height}, \eqref{cap:comparison.detail:derivative.covariances} and \eqref{cap:comparison.detail:representative.covariance} imply
\begin{equation}\label{cap:comparison.detail:conditional.parameters}
\|\Sigma_{T}^{-1}\bb{\alpha}\|_{2} \leq C_{k,J,\e}, \qquad |\bb{c}_{T}^{\mathsf{T}}\Sigma_{T}^{-1}\bb{\alpha}| \leq C_{k,J,\e}, \qquad 0 \leq \sig_{T}^{2} - \bb{c}_{T}^{\mathsf{T}}\Sigma_{T}^{-1}\bb{c}_{T} \leq C.
\end{equation}
The nonnegativity in \eqref{cap:comparison.detail:conditional.parameters} follows because the expression is the Haar variance of $Z_{T,k + 1}^{\mathrm{H}} - \bb{c}_{T}^{\mathsf{T}}\Sigma_{T}^{-1}(H(y_{1}),\ldots,H(y_{k}))^{\mathsf{T}}$. Define the joint-bin event
\[
E_{\bb{y},\bb{\alpha}}^{\star} \leqdef \{(F_{y_{1}}^{\star},\ldots,F_{y_{k}}^{\star}) \in \bb{\alpha} + [0,1]^{k}\}.
\]
Apply \eqref{an:conditional.bound} in Lemma~\ref{an:conditional} at $\lambda_{k + 1} = 1$ and $\lambda_{k + 1} = - 1$, using \eqref{cap:comparison.detail:conditional.parameters}. Since $|z| \leq e^{z} + e^{-z}$, this gives
\begin{equation}\label{cap:comparison.detail:derivative.bin}
\EE_{\star}[\ind_{\mathcal{A}_{T,\mathrm{sp}}^{\star}}|Z_{T,k + 1}^{\star}| \mid E_{\bb{y},\bb{\alpha}}^{\star}]
\leq \sum_{\lambda \in \{ - 1,1\}}\EE_{\star}[\ind_{\mathcal{A}_{T,\mathrm{sp}}^{\star}}e^{\lambda Z_{T,k + 1}^{\star}} \mid E_{\bb{y},\bb{\alpha}}^{\star}]
\leq C_{k,J,\e}.
\end{equation}
In the sums below, $\bb{\alpha}$ ranges over all vectors of lower endpoints of the closed unit bins in $u_{T} + J$. The events $E_{\bb{y},\bb{\alpha}}^{\star}$ cover $E_{\bb{y},J}^{\star}$ and are contained in it. Each representative height belongs to at most two bins, so each outcome belongs to at most $2^{k}$ joint-bin events. Multiplying \eqref{cap:comparison.detail:derivative.bin} by $\PP_{\star}(E_{\bb{y},\bb{\alpha}}^{\star})$, integrating over $v \in I_{i,T}^{\mathrm{loc}}(y_{i})$, and summing over the bin choices therefore give
\begin{equation}\label{cap:comparison.detail:conditional.derivative}
\begin{aligned}
&\EE_{\star}\left[\ind_{\mathcal{A}_{T,\mathrm{sp}}^{\star}}\ind_{E_{\bb{y},J}^{\star}}\int_{I_{i,T}^{\mathrm{loc}}(y_{i})}|\ell_{T}(F^{\star})'(y_{i} + \ell_{T}v)|\rd{v}\right] \\
&\qquad \leq \sum_{\bb{\alpha}}\int_{I_{i,T}^{\mathrm{loc}}(y_{i})}\EE_{\star}\left[\ind_{\mathcal{A}_{T,\mathrm{sp}}^{\star}}\ind_{E_{\bb{y},\bb{\alpha}}^{\star}}|\ell_{T}(F^{\star})'(y_{i} + \ell_{T}v)|\right]\rd{v} \\
&\qquad = \sum_{\bb{\alpha}}\PP_{\star}(E_{\bb{y},\bb{\alpha}}^{\star})\int_{I_{i,T}^{\mathrm{loc}}(y_{i})}\EE_{\star}\left[\ind_{\mathcal{A}_{T,\mathrm{sp}}^{\star}}|\ell_{T}(F^{\star})'(y_{i} + \ell_{T}v)| \mid E_{\bb{y},\bb{\alpha}}^{\star}\right]\rd{v} \\
&\qquad \leq C_{k,J,\e}|I_{i,T}^{\mathrm{loc}}(y_{i})|\sum_{\bb{\alpha}}\PP_{\star}(E_{\bb{y},\bb{\alpha}}^{\star}) \\
&\qquad \leq 2^{k + 1}RC_{k,J,\e}\PP_{\star}(E_{\bb{y},J}^{\star}) \\
&\qquad \leq C_{r,R}\PP_{\star}(E_{\bb{y},J}^{\star}).
\end{aligned}
\end{equation}
The first inequality uses the covering by joint-bin events and Tonelli's theorem. The next inequality uses \eqref{cap:comparison.detail:derivative.bin}, uniformly over the bins and $v$. To bound the sum of the joint-bin probabilities, fix a realization of the field. Each height $F_{y_{q}}^{\star}$ belongs to at most two closed unit bins in $u_{T} + J$, since two adjacent bins share an endpoint but otherwise do not overlap. A joint-bin event specifies one bin for each of the $k$ heights, so there are at most $2^{k}$ vectors $\bb{\alpha}$ for which $E_{\bb{y},\bb{\alpha}}^{\star}$ occurs. If any height lies outside $u_{T} + J$, none of these events occurs. Thus $\sum_{\bb{\alpha}}\ind_{E_{\bb{y},\bb{\alpha}}^{\star}} \leq 2^{k}\ind_{E_{\bb{y},J}^{\star}}$ pointwise, and taking expectations gives $\smash{\sum_{\bb{\alpha}}\PP_{\star}(E_{\bb{y},\bb{\alpha}}^{\star}) \leq 2^{k}\PP_{\star}(E_{\bb{y},J}^{\star})}$. Combining this probability bound with $|I_{i,T}^{\mathrm{loc}}(y_{i})| \leq 2R$ gives the factor $2^{k + 1}R$ in \eqref{cap:comparison.detail:conditional.derivative}. The constants are then absorbed into $C_{r,R}$ because $k \leq r$ and $J$ and $\e$ are fixed. Under the Haar law, the untruncated version of \eqref{an:conditional.bound} in Lemma~\ref{an:conditional} also proves \eqref{cap:comparison.detail:conditional.derivative} without $\smash{\ind_{\mathcal{A}_{T,\mathrm{sp}}^{\mathrm{H}}}}$.

Since every cap factor $\varpi_{C_{0}}(\cdot)$ lies in $[0,1]$, \eqref{cap:comparison.detail:product.difference} and \eqref{cap:comparison.detail:cap.quadrature.error} give
\[
\begin{aligned}
&\left|\prod_{i = 1}^{k}\varpi_{C_{0}}(M_{\star,i}^{\mathrm{loc}}(R))^{n_{i}} - \prod_{i = 1}^{k}\varpi_{C_{0}}(\widetilde{M}_{\star,i}^{\mathrm{loc}}(R,\ell_{\mathrm{quad}}))^{n_{i}}\right| \\
&\qquad \leq \sum_{i = 1}^{k}\left|\varpi_{C_{0}}(M_{\star,i}^{\mathrm{loc}}(R))^{n_{i}} - \varpi_{C_{0}}(\widetilde{M}_{\star,i}^{\mathrm{loc}}(R,\ell_{\mathrm{quad}}))^{n_{i}}\right| \\
&\qquad \leq \frac{\ell_{\mathrm{quad}}\|f'\|_{\infty}}{C_{0}}\sum_{i = 1}^{k}n_{i}\int_{I_{i,T}^{\mathrm{loc}}(y_{i})}|\ell_{T}(F^{\star})'(y_{i} + \ell_{T}v)|\rd{v}.
\end{aligned}
\]
Each sampled factor $f$ is bounded by $\|f\|_{\infty}$, and their product vanishes unless every representative height belongs to $u_{T} + J$. Hence $\prod_{q = 1}^{r}f(F_{h_{q}}^{\star} - u_{T}) \leq \|f\|_{\infty}^{r}\ind_{E_{\bb{y},J}^{\star}}$. For fixed representatives, the preceding bound on the cap difference does not depend on the nonrepresentative positions. Integrating those positions therefore contributes at most $(2R\ell_{T})^{r - k}$, since each of the $r - k$ positions ranges over an interval of length at most $2R\ell_{T}$. The definition \eqref{eq:Delta.T.quad}, the preceding bounds, Tonelli's theorem, and \eqref{cap:comparison.detail:conditional.derivative} give
\begin{equation}\label{cap:comparison.detail:integrated.error.bound}
\begin{aligned}
\EE_{\star}[\ind_{\mathcal{A}_{T,\mathrm{sp}}^{\star}}\Delta_{T,\mathrm{quad}}^{\star}]
&\leq \frac{\ell_{\mathrm{quad}}\|f\|_{\infty}^{r}\|f'\|_{\infty}}{C_{0}}\ell_{T}^{-r}(2R\ell_{T})^{r - k} \\
&\qquad \times \sum_{i = 1}^{k}n_{i}\int_{\mathscr{S}_{k,T}^{\e}}\EE_{\star}\left[\ind_{\mathcal{A}_{T,\mathrm{sp}}^{\star}}\ind_{E_{\bb{y},J}^{\star}}\int_{I_{i,T}^{\mathrm{loc}}(y_{i})}|\ell_{T}(F^{\star})'(y_{i} + \ell_{T}v)|\rd{v}\right]\rd{\bb{y}} \\
&\leq \frac{C_{r,R}\ell_{\mathrm{quad}}\|f\|_{\infty}^{r}\|f'\|_{\infty}}{C_{0}}\ell_{T}^{-r}(2R\ell_{T})^{r - k}\left(\sum_{i = 1}^{k}n_{i}\right)\int_{\mathscr{S}_{k,T}^{\e}}\PP_{\star}(E_{\bb{y},J}^{\star})\rd{\bb{y}} \\
&\leq C_{r,R}\ell_{\mathrm{quad}}\ell_{T}^{-r}(2R\ell_{T})^{r - k}\int_{\mathscr{S}_{k,T}^{\e}}\PP_{\star}(E_{\bb{y},J}^{\star})\rd{\bb{y}} \\
&\leq C_{r,R}\ell_{\mathrm{quad}}\ell_{T}^{-k}\int_{\mathscr{S}_{k,T}^{\e}}\PP_{\star}(E_{\bb{y},J}^{\star})\rd{\bb{y}}.
\end{aligned}
\end{equation}
Here $\sum_{i = 1}^{k}n_{i} = r$, and the factors involving the fixed function $f$ and cap parameter $C_{0}$ are absorbed into $C_{r,R}$. The last inequality also absorbs $(2R)^{r - k}$ into this constant, since $R$ is fixed and $k \leq r$. Under the Haar law, the same calculation bounds $\EE_{\mathrm{H}}[\Delta_{T,\mathrm{quad}}^{\mathrm{H}}]$ with both occurrences of $\smash{\ind_{\mathcal{A}_{T,\mathrm{sp}}^{\mathrm{H}}}}$ omitted, using the version of \eqref{cap:comparison.detail:conditional.derivative} without that event. By \eqref{cap:comparison.detail:representative.probability}, the normalized integral $\smash{\ell_{T}^{-k}\int_{\mathscr{S}_{k,T}^{\e}}\PP_{\star}(E_{\bb{y},J}^{\star})\rd{\bb{y}}}$ has bounded limsup. Taking the limsup in \eqref{cap:comparison.detail:integrated.error.bound} and absorbing this additional bound into $C_{r,R}$ therefore give
\begin{equation}\label{cap:comparison.detail:integrated.error.bound.next}
\limsup_{T \to \infty}\EE_{\mathrm{A}}[\ind_{\mathcal{A}_{T,\mathrm{sp}}^{\mathrm{A}}}\Delta_{T,\mathrm{quad}}^{\mathrm{A}}] \leq C_{r,R}\ell_{\mathrm{quad}}, \qquad
\limsup_{T \to \infty}\EE_{\mathrm{H}}[\Delta_{T,\mathrm{quad}}^{\mathrm{H}}] \leq C_{r,R}\ell_{\mathrm{quad}}.
\end{equation}

In order to remove the restriction to the arithmetic event, note that the cap factors $\varpi_{C_{0}}(M_{\star,i}^{\mathrm{loc}}(R))^{n_{i}}$ and $\smash{\varpi_{C_{0}}(\widetilde{M}_{\star,i}^{\mathrm{loc}}(R,\ell_{\mathrm{quad}}))^{n_{i}}}$ lie in $[0,1]$ by \eqref{cap:local.detail:basic.cap.bounds}. The error $\Delta_{T,\mathrm{quad}}^{\star}$ defined in \eqref{eq:Delta.T.quad} satisfies the deterministic bound $\Delta_{T,\mathrm{quad}}^{\star} \leq C_{r}\ell_{T}^{-r}|I_{\theta}|^{r} = \exp(\OO_{r}(t))$, using \eqref{eq:parameter.scales} and \eqref{eq:effective.cutoff.bounds}. The exceptional probability \eqref{an:spatial.exceptional} in Lemma~\ref{an:spatial.transform} gives
\begin{equation}\label{cap:comparison.detail:exceptional.error}
\EE_{\mathrm{A}}[\ind_{(\mathcal{A}_{T,\mathrm{sp}}^{\mathrm{A}})^{c}}\Delta_{T,\mathrm{quad}}^{\mathrm{A}}] \leq \exp\{ - ct^{1 + 2\e_{0}} + \OO_{r}(t)\} = \oo(1).
\end{equation}
Together with \eqref{cap:comparison.detail:integrated.error.bound.next}, \eqref{cap:comparison.detail:exceptional.error} proves
\begin{equation}\label{cap:comparison.detail:integrated.quadrature.error}
\limsup_{T \to \infty}\EE_{\star}[\Delta_{T,\mathrm{quad}}^{\star}] \leq C_{r,R}\ell_{\mathrm{quad}}.
\end{equation}

\emph{Step 4. Comparison after Gaussian smoothing.} With the quadrature error controlled by \eqref{cap:comparison.detail:integrated.quadrature.error}, it remains to compare the Euler and Haar expectations of the density in \eqref{cap:comparison.detail:partition.density} after localization and quadrature approximation, for fixed $r$, $R$, and $\ell_{\mathrm{quad}}$. Its weight is a smooth bounded function of finitely many field values, to which Lemma~\ref{an:smoothing} applies after adding a Gaussian perturbation. Let $\bb{Y}_{T}^{\star}$ be the vector consisting of the field values at the $r$ sampled positions and at the quadrature nodes:
\[
\bb{Y}_{T}^{\star} \leqdef \left(F_{h_{1}}^{\star},\ldots,F_{h_{r}}^{\star},\big(F_{y_{i} + \ell_{T}v}^{\star}\big)_{i,v}\right).
\]
Here $(i,v)$ ranges over the groups $1 \leq i \leq k$ and their quadrature nodes in \eqref{cap:comparison.detail:quadrature}, ordered first by group and then by increasing $v$. Coordinates corresponding to coincident positions are kept separately. If necessary, append unused copies of $F_{h_{1}}^{\star}$ so that the dimension is independent of the representative positions. This dimension is bounded in terms of $r$, $R$, and $\ell_{\mathrm{quad}}$ only. Denote the sampled coordinates of an argument $\bb{x}$ by $x_{q}$, $1 \leq q \leq r$, and its quadrature coordinates by $x_{i,v}$. The unused repeated coordinates do not enter the weight. Write
\begin{equation}\label{cap:comparison.detail:weight.definition}
W_{T,\bb{h}}(\bb{x}) \leqdef \prod_{i = 1}^{k}\varpi_{C_{0}}\left(\sum_{v}w_{i,v}f(x_{i,v} - u_{T})\right)^{n_{i}}\prod_{q = 1}^{r}f(x_{q} - u_{T}).
\end{equation}
Multiplying the value of \eqref{cap:comparison.detail:weight.definition} at $\bb{Y}_{T}^{\star}$ by $\ell_{T}^{-r}$ gives the localized density obtained from \eqref{cap:comparison.detail:partition.density}, with the cap masses replaced by the quadrature sums in \eqref{cap:comparison.detail:quadrature}. Each derivative in a sampled coordinate $x_{q}$, $1 \leq q \leq r$, differentiates one factor $f$, and each derivative in a quadrature coordinate $x_{i,v}$ contributes $w_{i,v}f'(x_{i,v} - u_{T})$ times a cap derivative. Therefore \eqref{cap:comparison.detail:quadrature}, \eqref{cap:comparison.detail:cap.derivative}, and $0 \leq \varpi_{C_{0}} \leq 1$ from \eqref{cap:local.detail:basic.cap.bounds} imply
\begin{equation}\label{cap:comparison.detail:smooth.weight}
0 \leq W_{T,\bb{h}} \leq \|f\|_{\infty}^{r}, \qquad \|\nabla W_{T,\bb{h}}\|_{\infty} \leq C_{r,R,\ell_{\mathrm{quad}}}.
\end{equation}
The bounds in \eqref{cap:comparison.detail:smooth.weight} are uniform in $T$ and in every allowed tuple of sampled positions.

Let $\mathscr{D}_{T,R} \subseteq I_{\theta}^{r}$ denote the integration domain for the fixed partition after localization: the representatives satisfy $\bb{y} \in \mathscr{S}_{k,T}^{\e}$, and each nonrepresentative position in group $i$ belongs to $B(y_{i}) \cap [y_{i} - R\ell_{T},y_{i} + R\ell_{T}]$. Let $\bb{g}$ be a standard Gaussian vector of the same dimension as $\bb{Y}_{T}^{\star}$, independent of the field. The smoothing estimate \eqref{an:smoothing.tv} in Lemma~\ref{an:smoothing} gives, uniformly over $\bb{h} = (h_{1},\ldots,h_{r}) \in \mathscr{D}_{T,R}$,
\begin{equation}\label{cap:smoothed.tv}
\left\|\mathcal{L}(\bb{Y}_{T}^{\mathrm{A}} + t^{-3/4}\bb{g}) - \mathcal{L}(\bb{Y}_{T}^{\mathrm{H}} + t^{-3/4}\bb{g})\right\|_{\mathrm{TV}} \leq \exp( - ct(\log{t})^{2}).
\end{equation}
For each fixed tuple $\bb{h}$, this total variation bound controls the difference between the expectations of the nonnegative bounded function $W_{T,\bb{h}}$ under the two laws. Together with \eqref{cap:comparison.detail:smooth.weight}, it gives
\begin{equation}\label{cap:comparison.detail:smoothed.weight.comparison}
\begin{aligned}
&\left|\EE_{\mathrm{A}}\left[W_{T,\bb{h}}(\bb{Y}_{T}^{\mathrm{A}} + t^{-3/4}\bb{g})\right] - \EE_{\mathrm{H}}\left[W_{T,\bb{h}}(\bb{Y}_{T}^{\mathrm{H}} + t^{-3/4}\bb{g})\right]\right| \\
&\qquad \leq \|W_{T,\bb{h}}\|_{\infty}\left\|\mathcal{L}(\bb{Y}_{T}^{\mathrm{A}} + t^{-3/4}\bb{g}) - \mathcal{L}(\bb{Y}_{T}^{\mathrm{H}} + t^{-3/4}\bb{g})\right\|_{\mathrm{TV}} \\
&\qquad \leq \|f\|_{\infty}^{r}\exp( - ct(\log{t})^{2}).
\end{aligned}
\end{equation}
Integrating \eqref{cap:comparison.detail:smoothed.weight.comparison} over the sampled positions and restoring the density factor $\ell_{T}^{-r}$ yield
\begin{equation}\label{cap:comparison.detail:smoothed.error}
\begin{aligned}
&\ell_{T}^{-r}\int_{\mathscr{D}_{T,R}}\left|\EE_{\mathrm{A}}\left[W_{T,\bb{h}}(\bb{Y}_{T}^{\mathrm{A}} + t^{-3/4}\bb{g})\right] - \EE_{\mathrm{H}}\left[W_{T,\bb{h}}(\bb{Y}_{T}^{\mathrm{H}} + t^{-3/4}\bb{g})\right]\right|\rd{\bb{h}} \\
&\qquad \leq \ell_{T}^{-r}|I_{\theta}|^{r}\|f\|_{\infty}^{r}\exp( - ct(\log{t})^{2})
\leq \exp\{ - ct(\log{t})^{2} + \OO_{r}(t)\} = \oo(1).
\end{aligned}
\end{equation}
Since \eqref{cap:comparison.detail:smoothed.error} bounds the integral of the absolute difference, it also bounds the absolute difference of the spatial integrals after multiplication of both integrands by any measurable function $\Psi_{T}(\bb{h})$ with $\|\Psi_{T}\|_{\infty} \leq 1$.

\emph{Step 5. Removal of the Gaussian perturbation.} Remove the Gaussian perturbation from the two expectations in \eqref{cap:comparison.detail:smoothed.error}, treating separately perturbations with sup norm at most one and greater than one. Let $J^{+}$ be a fixed bounded interval containing all points at distance at most one from the supports of $f$ and $f'$. On $\{\|t^{-3/4}\bb{g}\|_{\infty} \leq 1\}$, a nonzero derivative of $W_{T,\bb{h}}$ anywhere along the segment between $\bb{Y}_{T}^{\star}$ and $\bb{Y}_{T}^{\star} + t^{-3/4}\bb{g}$ requires every unperturbed representative height to lie in $u_{T} + J^{+}$. Indeed, differentiating a representative factor replaces $f$ by $f'$, while differentiating any other factor leaves all representative factors $f$ present. Integrating the derivative along this segment and using \eqref{cap:comparison.detail:smooth.weight} therefore gives, on $\{\|t^{-3/4}\bb{g}\|_{\infty} \leq 1\}$,
\begin{equation}\label{cap:comparison.detail:smoothing.support}
\begin{aligned}
|W_{T,\bb{h}}(\bb{Y}_{T}^{\star} + t^{-3/4}\bb{g}) - W_{T,\bb{h}}(\bb{Y}_{T}^{\star})|
&\leq t^{-3/4}\int_{0}^{1}|\nabla W_{T,\bb{h}}(\bb{Y}_{T}^{\star} + zt^{-3/4}\bb{g})\cdot\bb{g}|\rd{z} \\
&\leq C_{r,R,\ell_{\mathrm{quad}}}t^{-3/4}\|\bb{g}\|_{1}\ind_{E_{\bb{y},J^{+}}^{\star}}.
\end{aligned}
\end{equation}
Independence of $\bb{g}$ from the field $\{F_{h}^{\star}\}_{h \in I_{\theta}}$ gives $\smash{\EE_{\star}[\|\bb{g}\|_{1}\ind_{E_{\bb{y},J^{+}}^{\star}}] = \EE\|\bb{g}\|_{1}\PP_{\star}(E_{\bb{y},J^{+}}^{\star})}$. Integrating \eqref{cap:comparison.detail:smoothing.support} over the sampled positions and using the volume bound $(2R\ell_{T})^{r - k}$ for the $r - k$ local nonrepresentative positions, established in the derivation of \eqref{cap:comparison.detail:integrated.error.bound}, therefore gives
\begin{equation}\label{cap:comparison.detail:integrated.smoothing}
\begin{aligned}
&\ell_{T}^{-r}\int_{\mathscr{D}_{T,R}}\EE_{\star}\left[\ind_{\{\|t^{-3/4}\bb{g}\|_{\infty} \leq 1\}}\left|W_{T,\bb{h}}(\bb{Y}_{T}^{\star} + t^{-3/4}\bb{g}) - W_{T,\bb{h}}(\bb{Y}_{T}^{\star})\right|\right]\rd{\bb{h}} \\[-1mm]
&\qquad \leq C_{r,R,\ell_{\mathrm{quad}}}t^{-3/4}\EE\|\bb{g}\|_{1}\ell_{T}^{-r}(2R\ell_{T})^{r - k}\int_{\mathscr{S}_{k,T}^{\e}}\PP_{\star}(E_{\bb{y},J^{+}}^{\star})\rd{\bb{y}} \\[-1mm]
&\qquad \leq C_{r,R,\ell_{\mathrm{quad}}}t^{-3/4}.
\end{aligned}
\end{equation}
The last inequality holds for all sufficiently large $T$ by \eqref{cap:comparison.detail:representative.probability}, applied to $J^{+}$ (see the comment just below that equation). The factor $\EE\|\bb{g}\|_{1}$ is absorbed into the constant $C_{r,R,\ell_{\mathrm{quad}}}$ because the dimension of $\bb{g}$ is bounded in terms of $r$, $R$, and $\ell_{\mathrm{quad}}$.

A union bound over the fixed number of coordinates of $\bb{g}$ gives $\PP(\|t^{-3/4}\bb{g}\|_{\infty} > 1) \leq C_{r,R,\ell_{\mathrm{quad}}}\exp( - ct^{3/2})$. Since $W_{T,\bb{h}}(\bb{Y}_{T}^{\star} + t^{-3/4}\bb{g})$ and $W_{T,\bb{h}}(\bb{Y}_{T}^{\star})$ both lie in $[0,\|f\|_{\infty}^{r}]$ by \eqref{cap:comparison.detail:smooth.weight} and $\mathscr{D}_{T,R} \subseteq I_{\theta}^{r}$, the expected integrated error on this complementary Gaussian event satisfies
\begin{equation}\label{cap:comparison.detail:gaussian.tail}
\begin{aligned}
&\ell_{T}^{-r}\int_{\mathscr{D}_{T,R}}\EE_{\star}\left[\ind_{\{\|t^{-3/4}\bb{g}\|_{\infty} > 1\}}\left|W_{T,\bb{h}}(\bb{Y}_{T}^{\star} + t^{-3/4}\bb{g}) - W_{T,\bb{h}}(\bb{Y}_{T}^{\star})\right|\right]\rd{\bb{h}} \\
&\qquad \leq \ell_{T}^{-r}|I_{\theta}|^{r}\|f\|_{\infty}^{r}\PP(\|t^{-3/4}\bb{g}\|_{\infty} > 1) \\
&\qquad \leq C_{r,R,\ell_{\mathrm{quad}}}\exp\{ - ct^{3/2} + \OO_{r}(t)\} = \oo(1).
\end{aligned}
\end{equation}
Combining \eqref{cap:comparison.detail:integrated.smoothing} and \eqref{cap:comparison.detail:gaussian.tail} with \eqref{cap:comparison.detail:smoothed.error}, the triangle inequality gives
\[
\begin{aligned}
&\sup_{\|\Psi_{T}\|_{\infty} \leq 1}\left|\ell_{T}^{-r}\int_{\mathscr{D}_{T,R}}\Psi_{T}(\bb{h})\left\{\EE_{\mathrm{A}}\left[W_{T,\bb{h}}(\bb{Y}_{T}^{\mathrm{A}})\right] - \EE_{\mathrm{H}}\left[W_{T,\bb{h}}(\bb{Y}_{T}^{\mathrm{H}})\right]\right\}\rd{\bb{h}}\right| \\
&\qquad \leq \ell_{T}^{-r}\int_{\mathscr{D}_{T,R}}\left|\EE_{\mathrm{A}}\left[W_{T,\bb{h}}(\bb{Y}_{T}^{\mathrm{A}})\right] - \EE_{\mathrm{H}}\left[W_{T,\bb{h}}(\bb{Y}_{T}^{\mathrm{H}})\right]\right|\rd{\bb{h}} = \oo(1),
\end{aligned}
\]
where the supremum is over deterministic measurable functions $\Psi_{T}:I_{\theta}^{r} \to \R$. This proves the comparison of the two quadrature expressions uniformly over these test functions.

\emph{Step 6. Combination of the estimates.} Sum over the finitely many partitions $\pi$ of $\{1,\ldots,r\}$. For each partition, let $p_{\star,\pi}^{(0)}$ denote the original density in \eqref{cap:comparison.detail:partition.density}. Let $\smash{p_{\star,\pi}^{(1)}}$ be the density after restricting the nonrepresentative positions to their local windows, $\smash{p_{\star,\pi}^{(2)}}$ the density after also localizing the cap masses, and $\smash{p_{\star,\pi}^{(3)}}$ the density after replacing those local masses by quadrature sums. Finally, let $\smash{p_{\star,\pi}^{(4)}}$ be obtained from $\smash{p_{\star,\pi}^{(3)}}$ by replacing $\bb{Y}_{T}^{\star}$ inside $W_{T,\bb{h}}$ by $\bb{Y}_{T}^{\star} + t^{-3/4}\bb{g}_{\pi}$, where $\bb{g}_{\pi}$ is the independent standard Gaussian vector of the dimension required for this partition. All densities are extended by zero outside their respective integration domains, and $\|p\|_{1} \leqdef \int_{I_{\theta}^{r}}|p(\bb{h})|\rd{\bb{h}}$. The triangle inequality gives the full decomposition below. Steps~1--3 control the first three errors, Step~5 handles the two Gaussian events, and Step~4 gives the final comparison:
\[
\begin{aligned}
&\limsup_{T \to \infty}\left\|\EE\nu_{\mathrm{A}}^{\times r} - \EE_{\mathrm{H}}\nu_{\mathrm{H}}^{\times r}\right\|_{\mathrm{TV}} \\
&\quad \leq \sum_{\pi}\sum_{\star \in \{\mathrm{A},\mathrm{H}\}}
\Bigg[
\underbrace{\limsup_{T \to \infty}\EE_{\star}\left[\left\|p_{\star,\pi}^{(0)} - p_{\star,\pi}^{(1)}\right\|_{1}\right]}_{\leq C_{r}R^{-c}\text{ by }\eqref{cap:comparison.detail:localization.error}}
+ \underbrace{\limsup_{T \to \infty}\EE_{\star}\left[\left\|p_{\star,\pi}^{(1)} - p_{\star,\pi}^{(2)}\right\|_{1}\right]}_{\leq C_{r}R^{-c}\text{ by }\eqref{cap:comparison.detail:cap.localization.error}} \\
&\qquad
+ \underbrace{\limsup_{T \to \infty}\EE_{\star}\left[\left\|p_{\star,\pi}^{(2)} - p_{\star,\pi}^{(3)}\right\|_{1}\right]}_{\leq C_{r,R}\ell_{\mathrm{quad}}\text{ by }\eqref{cap:comparison.detail:integrated.quadrature.error}}
+ \underbrace{\limsup_{T \to \infty}\EE_{\star}\left[\ind_{\{\|t^{-3/4}\bb{g}_{\pi}\|_{\infty} \leq 1\}}\left\|p_{\star,\pi}^{(3)} - p_{\star,\pi}^{(4)}\right\|_{1}\right]}_{= 0\text{ by }\eqref{cap:comparison.detail:integrated.smoothing}} \\
&\qquad
+ \underbrace{\limsup_{T \to \infty}\EE_{\star}\left[\ind_{\{\|t^{-3/4}\bb{g}_{\pi}\|_{\infty} > 1\}}\left\|p_{\star,\pi}^{(3)} - p_{\star,\pi}^{(4)}\right\|_{1}\right]}_{= 0\text{ by }\eqref{cap:comparison.detail:gaussian.tail}}
\Bigg]
+ \sum_{\pi}\underbrace{\limsup_{T \to \infty}\left\|\EE_{\mathrm{A}}p_{\mathrm{A},\pi}^{(4)} - \EE_{\mathrm{H}}p_{\mathrm{H},\pi}^{(4)}\right\|_{1}}_{= 0\text{ by }\eqref{cap:comparison.detail:smoothed.error}} \\
&\quad \leq C_{r}R^{-c} + C_{r,R}\ell_{\mathrm{quad}}.
\end{aligned}
\]
The last inequality absorbs the finite number of partitions and the two choices of $\star$ into the constants. For fixed $R$, let $\ell_{\mathrm{quad}} \downarrow 0$, and then let $R \to \infty$. This proves the proposition. The order of these limits keeps the Gaussian vector dimension fixed during the limit $T \to \infty$. The quadrature approximates only the local masses inside the cap factors and leaves every sampled position, and hence every measurable positional test function, unchanged.
\end{proof}

\subsubsection{Capture of the normalized mass}\label{cap:capture}

Write
\[
\mu_{\star}(\rd{h}) \leqdef \ell_{T}^{-1}e^{\beta(F_{h}^{\star} - u_{T})}\rd{h}, \qquad M_{\beta,T}^{\star} = \mu_{\star}(I_{\theta}),
\]
so that the normalization of $\mu_{\mathrm{A}}$ is $G_{\beta,T}^{\mathrm{A}}$ and the normalization of $\mu_{\mathrm{H}}$ is $G_{\beta,T}^{\mathrm{H}}$.

By \eqref{cap:definition}, on each interval $B \in \mathcal{B}_{T}^{\e}$, the rational cap $x \mapsto x\varpi_{C_{0}}(x)$ replaces the mass $\mu_{\star}^{\chi}(B)$ by $\nu_{\star}(B) = C_{0}\mu_{\star}^{\chi}(B)/(C_{0} + \mu_{\star}^{\chi}(B))$. The next lemma shows that choosing $C_{0}$ large enough makes the total expected mass lost arbitrarily small for all sufficiently large $T$.

\begin{lemma}[Loss from the rational cap]\label{cap:loss}
Let $\theta > 0$ and $\beta > \beta_{c}(\theta)$. Fix a smooth function $\chi : \R \to [0,1]$ with compact support. Then there exist $0 < a \leq 1$ and $C_{\chi} < \infty$ such that, for each $\star \in \{\mathrm{A},\mathrm{H}\}$,
\begin{equation}\label{cap:uncapped.moment}
\limsup_{T \to \infty}\sum_{B \in \mathcal{B}_{T}^{\e}}\EE_{\star}[\mu_{\star}^{\chi}(B)^{1 + a}] \leq C_{\chi}.
\end{equation}
Consequently,
\begin{equation}\label{cap:loss.bound}
\limsup_{T \to \infty}\EE_{\star}\left[\sum_{B \in \mathcal{B}_{T}^{\e}}(\mu_{\star}^{\chi}(B) - \nu_{\star}(B))\right] \leq C_{\chi}C_{0}^{-a}.
\end{equation}
\end{lemma}

\begin{proof}[Proof of Lemma~\ref{cap:loss}]
The identity $\mu_{\star}^{\chi}(B)^{1 + a} = \int_{B}\mu_{\star}^{\chi}(B)^{a}\mu_{\star}^{\chi}(\rd{h})$ expresses the required moment as an integral over a position $h \in B$. Conditioning on the height at $h$ will then control the remaining interval mass. Choose $a \in (0,1]$ for the case $k = 1$ of \eqref{an:joint.annular.bound} in Lemma~\ref{an:joint.annular}. The definition of $\mu_{\star}^{\chi}(\rd{h})$ in \eqref{cap:retained.definition} and Tonelli's theorem give
\begin{equation}\label{cap:capture.detail:rooted.moment}
\sum_{B \in \mathcal{B}_{T}^{\e}}\EE_{\star}[\mu_{\star}^{\chi}(B)^{1 + a}]
= \ell_{T}^{-1}\sum_{B \in \mathcal{B}_{T}^{\e}}\int_{B}\EE_{\star}\left[f(F_{h}^{\star} - u_{T})\mu_{\star}^{\chi}(B)^{a}\right]\rd{h}.
\end{equation}
Let $J$ contain the support of $f(x) = e^{\beta x}\chi(x)$. Only the heights $F_{h}^{\star} - u_{T} \in J$ contribute to \eqref{cap:capture.detail:rooted.moment}. For $B \in \mathcal{B}_{T}^{\e}$ and $h \in B$, use $h_{1} = h$ and the interval $B$ in the definition of $M_{\star,1}^{\mathrm{out}}(1)$ from Lemma~\ref{an:joint.annular}, so that
\[
\begin{aligned}
M_{\star,1}^{\mathrm{out}}(1) &= \mu_{\star}^{\chi}(B \setminus [h - \ell_{T},h + \ell_{T}]), \\
\mu_{\star}^{\chi}(B \cap [h - \ell_{T},h + \ell_{T}]) &\leq \ell_{T}^{-1}(2\ell_{T})\|f\|_{\infty} = 2\|f\|_{\infty}.
\end{aligned}
\]
Subadditivity of $x \mapsto x^{a}$ for $a \leq 1$ therefore yields
\begin{equation}\label{cap:capture.detail:local.mass}
\mu_{\star}^{\chi}(B)^{a}
= \left\{\mu_{\star}^{\chi}(B \cap [h - \ell_{T},h + \ell_{T}]) + M_{\star,1}^{\mathrm{out}}(1)\right\}^{a}
\leq (2\|f\|_{\infty})^{a} + M_{\star,1}^{\mathrm{out}}(1)^{a}.
\end{equation}

Cover $J$ by finitely many bins $J_{1},\ldots,J_{L}$ of one fixed positive width, and choose a compact interval $J'$ containing them all. For each bin, \eqref{an:joint.annular.bound} in Lemma~\ref{an:joint.annular}, with one representative and $R = 1$, gives
\begin{equation}\label{cap:capture.detail:conditional.outside.mass}
\EE_{\star}[M_{\star,1}^{\mathrm{out}}(1)^{a} \mid F_{h}^{\star} - u_{T} \in J_{v}] \leq C_{\chi} + \oo(1).
\end{equation}
The error is uniform in $h$ and in the finitely many bins, and is absorbed in $C_{\chi}$ for all sufficiently large $T$. Multiply \eqref{cap:capture.detail:local.mass} by $f(F_{h}^{\star} - u_{T})$, bound this factor by $\|f\|_{\infty}\sum_{v = 1}^{L}\ind_{\{F_{h}^{\star} - u_{T} \in J_{v}\}}$, and apply \eqref{cap:capture.detail:conditional.outside.mass}. This gives
\begin{equation}\label{cap:capture.detail:conditional.mass}
\EE_{\star}\left[f(F_{h}^{\star} - u_{T})\mu_{\star}^{\chi}(B)^{a}\right]
\leq C_{\chi}\sum_{v = 1}^{L}\PP_{\star}(F_{h}^{\star} - u_{T} \in J_{v})
\leq C_{\chi}\PP_{\star}(F_{h}^{\star} - u_{T} \in J').
\end{equation}
The last constant absorbs the finite number of bins, including any endpoint multiplicity. The conditional estimate \eqref{an:joint.annular.bound} in Lemma~\ref{an:joint.annular} holds without truncation under both laws. In its proof, the arithmetic contribution from $\mathcal{A}_{T,\mathrm{sp}}^{c}$ is bounded by $\oo(1)$ in \eqref{an:annular.detail:event.removal}, so no separate exceptional-event estimate is needed in \eqref{cap:capture.detail:conditional.mass}. For $k = 1$, the proof of Lemma~\ref{an:joint.annular} uses only $|h - h'| \leq 1$ within $B$ and no separation between different intervals. Thus $a$ and $C_{\chi}$ here can be chosen independently of $\e$.

The weak convergence \eqref{an:height.measure} in Lemma~\ref{an:height}, applied to $\ind_{J'}$ whose discontinuities have zero limiting mass, bounds the probability in \eqref{cap:capture.detail:conditional.mass} by $C_{\chi}N_{T}^{-1}$ for all sufficiently large $T$, uniformly in $h \in I_{\theta}$. Inserting this bound into \eqref{cap:capture.detail:rooted.moment}, and using $N_{T} = |I_{\theta}|/\ell_{T}$ from \eqref{eq:N.T}, gives
\[
\sum_{B \in \mathcal{B}_{T}^{\e}}\EE_{\star}[\mu_{\star}^{\chi}(B)^{1 + a}]
\leq C_{\chi}\ell_{T}^{-1}\sum_{B \in \mathcal{B}_{T}^{\e}}|B|N_{T}^{-1}
\leq C_{\chi}\ell_{T}^{-1}|I_{\theta}|N_{T}^{-1} = C_{\chi}.
\]
Taking the limit superior as $T \to \infty$ proves \eqref{cap:uncapped.moment}.

For the cap loss, \eqref{cap:definition} gives $\nu_{\star}(B) = C_{0}\mu_{\star}^{\chi}(B)/(C_{0} + \mu_{\star}^{\chi}(B))$. If $x \geq 0$ and $0 < a \leq 1$, then $x/(1 + x) \leq \min\{x,1\} \leq x^{a}$. With $x = M/C_{0}$, this yields
\begin{equation}\label{cap:capture.detail:pointwise.cap.loss}
M - \frac{C_{0}M}{C_{0} + M} = M\frac{M/C_{0}}{1 + M/C_{0}} \leq M\min\{M/C_{0},1\} \leq C_{0}^{-a}M^{1 + a}.
\end{equation}
Apply \eqref{cap:capture.detail:pointwise.cap.loss} with $M = \mu_{\star}^{\chi}(B)$ and sum over $B \in \mathcal{B}_{T}^{\e}$. Equation~\eqref{cap:uncapped.moment} then gives
\[
\limsup_{T \to \infty}\EE_{\star}\left[\sum_{B \in \mathcal{B}_{T}^{\e}}(\mu_{\star}^{\chi}(B) - \nu_{\star}(B))\right] \leq C_{0}^{-a}\limsup_{T \to \infty}\sum_{B \in \mathcal{B}_{T}^{\e}}\EE_{\star}[\mu_{\star}^{\chi}(B)^{1 + a}] \leq C_{\chi}C_{0}^{-a},
\]
which is \eqref{cap:loss.bound}.
\end{proof}

\begin{proposition}[Normalized capture]\label{cap:normalized.capture}
Let $\theta > 0$ and $\beta > \beta_{c}(\theta)$. Recall $M_{\beta,T}^{\star} = \mu_{\star}(I_{\theta})$. There are deterministic choices of $\chi_{b}$, $\e_{b}$ and $C_{0,b}$, indexed by $b \in \N$, for which the measures $\nu_{\star,b}$ defined by \eqref{cap:definition} satisfy, for each $\star \in \{\mathrm{A},\mathrm{H}\}$,
\begin{equation}\label{cap:normalized.capture:tv}
\lim_{b \to \infty}\limsup_{T \to \infty}\EE_{\star}\left[\left\|\frac{\mu_{\star}}{M_{\beta,T}^{\star}} - \frac{\nu_{\star,b}}{\nu_{\star,b}(I_{\theta})}\right\|_{\mathrm{TV}}\right] = 0.
\end{equation}
An arbitrary deterministic probability measure is used when $\nu_{\star,b}(I_{\theta}) = 0$. There are numbers $\e_{b}^{\mathrm{atom}} \downarrow 0$ as $b \to \infty$ such that every subsequential distributional limit of $\nu_{\mathrm{H},b}(I_{\theta})$ as $T \to \infty$ assigns mass at most $\e_{b}^{\mathrm{atom}}$ to zero, i.e., $\lim_{\delta \downarrow 0}\limsup_{T \to \infty}\PP_{\mathrm{H}}(\nu_{\mathrm{H},b}(I_{\theta}) \leq \delta) \leq \e_{b}^{\mathrm{atom}} ~\forall b \in \N$.
\end{proposition}

\begin{proof}[Proof of Proposition~\ref{cap:normalized.capture}]
The mass loss from $\mu_{\star}$ to $\nu_{\star}$ comes from the height cutoff $\chi$, the restriction to the intervals in $\mathcal{B}_{T}^{\e}$, and the rational cap $x \mapsto x\varpi_{C_{0}}(x) = C_{0}x/(C_{0} + x)$. First bound these three losses before normalization, and then use the lower-tail bound in Corollary~\ref{temp:lower.mass} for $M_{\beta,T}^{\star}$ to control their normalized effect.

Choose $L,K > 0$, and let $\chi$ be smooth, take values in $[0,1]$, equal one on $[ - L,K]$, and vanish outside $[ - L - 1,K + 1]$. Let $D_{T}^{\e} \leqdef I_{\theta} \setminus \bigcup_{B \in \mathcal{B}_{T}^{\e}}B$ be the omitted spatial set. Its Lebesgue measure satisfies $|D_{T}^{\e}| \leq 2\e|I_{\theta}| + C$, where the constant accounts for the at most two partial unit intervals at the endpoints. The uniform one-point estimate \eqref{an:height.measure} in Lemma~\ref{an:height}, tested against $f(x) = e^{\beta x}\chi(x)$, implies
\[
\sup_{h \in I_{\theta}}\EE_{\star}[f(F_{h}^{\star} - u_{T})] \leq C_{\chi}N_{T}^{-1}
\]
for all sufficiently large $T$. Since $\ell_{T}N_{T} = |I_{\theta}|$ by \eqref{eq:N.T} and $|I_{\theta}| \to \infty$, integration over $D_{T}^{\e}$ gives
\begin{equation}\label{cap:capture.detail:boundary.loss}
\begin{aligned}
\EE_{\star}[\mu_{\star}^{\chi}(D_{T}^{\e})]
&= \ell_{T}^{-1}\int_{D_{T}^{\e}}\EE_{\star}[f(F_{h}^{\star} - u_{T})]\rd{h} \\
&\leq C_{\chi}\ell_{T}^{-1}|D_{T}^{\e}|N_{T}^{-1} = C_{\chi}\frac{|D_{T}^{\e}|}{|I_{\theta}|} \leq C_{\chi}\e + \oo(1).
\end{aligned}
\end{equation}
The cap definition \eqref{cap:definition} and $0 \leq \chi \leq 1$ give $0 \leq \nu_{\star} \leq \mu_{\star}^{\chi} \leq \mu_{\star}$. Given that $M_{\beta,T}^{\star} = \mu_{\star}(I_{\theta})$ and that $\nu_{\star}(D_{T}^{\e}) = 0$ by \eqref{cap:definition}, the unnormalized mass loss splits exactly as
\begin{equation}\label{cap:capture.detail:loss.decomposition}
M_{\beta,T}^{\star} - \nu_{\star}(I_{\theta})
= (\mu_{\star} - \mu_{\star}^{\chi})(I_{\theta}) + \mu_{\star}^{\chi}(D_{T}^{\e}) + \sum_{B \in \mathcal{B}_{T}^{\e}}(\mu_{\star}^{\chi}(B) - \nu_{\star}(B)).
\end{equation}
On the event $\sup_{h \in I_{\theta}}(F_{h}^{\star} - u_{T}) \leq K$, the first term in \eqref{cap:capture.detail:loss.decomposition} satisfies
\[
(\mu_{\star} - \mu_{\star}^{\chi})(I_{\theta})
\leq \ell_{T}^{-1}\int_{I_{\theta}}e^{\beta(F_{h}^{\star} - u_{T})}\ind_{\{F_{h}^{\star} - u_{T}< - L\}}\rd{h},
\]
because $\chi = 1$ on $[ - L,K]$. The first estimate in \eqref{an:height.tails} of Lemma~\ref{an:height} controls this low-height mass. Together with the boundary estimate \eqref{cap:capture.detail:boundary.loss}, and the cap-loss estimate \eqref{cap:loss.bound} in Lemma~\ref{cap:loss}, it gives
\begin{equation}\label{cap:unnormalized.loss}
\limsup_{T \to \infty}\EE_{\star}\left[(M_{\beta,T}^{\star} - \nu_{\star}(I_{\theta}))\ind_{( - \infty,K]}\left(\sup_{h \in I_{\theta}}(F_{h}^{\star} - u_{T})\right)\right] \leq Ce^{-c_{\beta}L} + C_{\chi}\e + C_{\chi}C_{0}^{-a},
\end{equation}
where $c_{\beta}$ is any fixed number in $(0,\beta - \beta_{c}(\theta))$. The supremum estimate in \eqref{an:height.tails} of Lemma~\ref{an:height} also gives
\begin{equation}\label{cap:capture.detail:upper.height}
\limsup_{T \to \infty}\PP_{\star}\left(\sup_{h \in I_{\theta}}(F_{h}^{\star} - u_{T}) > K\right) \leq Ce^{-\beta_{c}(\theta)K}.
\end{equation}

Let $0 \leq \nu \leq \mu$ be finite measures on a common measurable space $\mathscr{S}$, with masses $M_{\nu} \leqdef \nu(\mathscr{S})$ and $M_{\mu} \leqdef \mu(\mathscr{S}) > 0$. If $M_{\nu} > 0$, write
\[
\frac{\mu}{M_{\mu}} - \frac{\nu}{M_{\nu}} = \frac{\mu - \nu}{M_{\mu}} + \left(\frac{1}{M_{\mu}} - \frac{1}{M_{\nu}}\right)\nu.
\]
Since $\mu - \nu$ is a nonnegative measure, the triangle inequality gives
\begin{equation}\label{cap:normalization.inequality}
\begin{aligned}
\left\|\frac{\mu}{M_{\mu}} - \frac{\nu}{M_{\nu}}\right\|_{\mathrm{TV}}
&\leq \frac{(\mu - \nu)(\mathscr{S})}{M_{\mu}} + \left(\frac{1}{M_{\nu}} - \frac{1}{M_{\mu}}\right)\nu(\mathscr{S}) \\
&= \frac{M_{\mu} - M_{\nu}}{M_{\mu}} + \frac{M_{\mu} - M_{\nu}}{M_{\mu}} = \frac{2(M_{\mu} - M_{\nu})}{M_{\mu}}.
\end{aligned}
\end{equation}
If $M_{\nu} = 0$, interpret $\nu/M_{\nu}$ as the arbitrary deterministic probability measure specified in the statement of the proposition. Then \eqref{cap:normalization.inequality} still holds, since its right-hand side equals two and the total variation distance between probability measures is at most two.

For $\delta_{\mathrm{mass}} > 0$, apply \eqref{cap:normalization.inequality} on the event $\{M_{\beta,T}^{\star} \geq \delta_{\mathrm{mass}}\} \cap \{\sup_{h \in I_{\theta}}(F_{h}^{\star} - u_{T}) \leq K\}$. On the complement use the deterministic bound two. Equations \eqref{cap:unnormalized.loss} and \eqref{cap:capture.detail:upper.height} then imply
\begin{equation}\label{cap:capture.detail:normalized.error}
\begin{aligned}
\limsup_{T \to \infty}\EE_{\star}\left[\left\|\frac{\mu_{\star}}{M_{\beta,T}^{\star}} - \frac{\nu_{\star}}{\nu_{\star}(I_{\theta})}\right\|_{\mathrm{TV}}\right]
&\leq 2\limsup_{T \to \infty}\PP_{\star}(M_{\beta,T}^{\star} < \delta_{\mathrm{mass}}) + 2Ce^{-\beta_{c}(\theta)K} \\
&\qquad + \frac{2}{\delta_{\mathrm{mass}}}\left(Ce^{-c_{\beta}L} + C_{\chi}\e + C_{\chi}C_{0}^{-a}\right).
\end{aligned}
\end{equation}
The lower-mass estimate in Corollary~\ref{temp:lower.mass} makes the first limiting upper probability arbitrarily small as $\delta_{\mathrm{mass}} \downarrow 0$.

The parameters can now be chosen explicitly in the necessary order. For each $b \in \N$, first choose $\delta_{\mathrm{mass},b} > 0$ such that the limiting upper probability in \eqref{cap:capture.detail:normalized.error} is at most $2^{-b}/8$ under both laws. Choose $K_{b}$ so large that $Ce^{-\beta_{c}(\theta)K_{b}} \leq 2^{-b}/8$. Next choose $L_{b}$ so large that $Ce^{-c_{\beta}L_{b}} \leq \delta_{\mathrm{mass},b}2^{-b}/24$, and fix a smooth cutoff $\chi_{b}$ for this pair $(L_{b},K_{b})$. With $C_{\chi_{b}}$ and the exponent $a$ in \eqref{cap:loss.bound} of Lemma~\ref{cap:loss} now fixed, choose $0 < \e_{b} < 1/4$ and $C_{0,b} > 0$ so that
\[
C_{\chi_{b}}\e_{b} \leq \frac{\delta_{\mathrm{mass},b}2^{-b}}{24}, \qquad C_{\chi_{b}}C_{0,b}^{-a} \leq \frac{\delta_{\mathrm{mass},b}2^{-b}}{24}.
\]
The constants are chosen to work under both laws. These are deterministic choices, made before letting $T \to \infty$. Substitution into \eqref{cap:capture.detail:normalized.error} gives a bound at most $3\cdot 2^{-b}/4$, which proves the asserted convergence in expected total variation in \eqref{cap:normalized.capture:tv}.

It remains to prove the asserted bound on the mass at zero of subsequential distributional limits of $\nu_{\mathrm{H},b}(I_{\theta})$. The same parameter choices give this bound. Indeed,
\[
\{\nu_{\star,b}(I_{\theta}) \leq \delta_{\mathrm{mass},b}/2\} \subseteq \{M_{\beta,T}^{\star} < \delta_{\mathrm{mass},b}\} \cup \{M_{\beta,T}^{\star} - \nu_{\star,b}(I_{\theta}) \geq \delta_{\mathrm{mass},b}/2\}.
\]
For the second event, split off $\{\sup_{h \in I_{\theta}}(F_{h}^{\star} - u_{T}) > K_{b}\}$ and use Markov's inequality on its complement. By \eqref{cap:unnormalized.loss} and \eqref{cap:capture.detail:upper.height},
\begin{equation}\label{cap:capture.detail:small.retained.mass}
\begin{aligned}
\limsup_{T \to \infty}\PP_{\star}(\nu_{\star,b}(I_{\theta}) \leq \delta_{\mathrm{mass},b}/2)
&\leq \limsup_{T \to \infty}\PP_{\star}(M_{\beta,T}^{\star} < \delta_{\mathrm{mass},b}) + Ce^{-\beta_{c}(\theta)K_{b}} \\
&\qquad + \frac{2}{\delta_{\mathrm{mass},b}}\left(Ce^{-c_{\beta}L_{b}} + C_{\chi_{b}}\e_{b} + C_{\chi_{b}}C_{0,b}^{-a}\right)
\leq 2^{-b}/2.
\end{aligned}
\end{equation}
For fixed $b$, the first-moment ($r=1$) case of Lemma~\ref{cap:moments} and Markov's inequality give tightness of $\nu_{\mathrm{H},b}(I_{\theta})$. Let $M_{b}$ be a distributional limit along any subsequence $T_{k} \to \infty$. Since $( - \infty,\delta_{\mathrm{mass},b}/2)$ is open, the open-set form of the Portmanteau theorem gives
\begin{equation}\label{cap:capture.detail:zero.atom}
\PP(M_{b} = 0)
\leq \PP(M_{b} < \delta_{\mathrm{mass},b}/2)
\leq \liminf_{k \to \infty}\PP_{\mathrm{H}}(\nu_{\mathrm{H},b}(I_{\theta}) < \delta_{\mathrm{mass},b}/2) \leq 2^{-b}/2,
\end{equation}
where the last inequality uses \eqref{cap:capture.detail:small.retained.mass}. For each fixed $b \in \N$ and every $0 < \delta \leq \delta_{\mathrm{mass},b}/2$, monotonicity and the same estimate give $\limsup_{T \to \infty}\PP_{\mathrm{H}}(\nu_{\mathrm{H},b}(I_{\theta}) \leq \delta) \leq 2^{-b}/2$. Letting $\delta \downarrow 0$ proves the last assertion of the proposition with $\e_{b}^{\mathrm{atom}} \leqdef 2^{-b}/2$.
\end{proof}

The moment bounds in Lemma~\ref{cap:moments} and the comparison and approximation results in Propositions~\ref{cap:moment.measures} and~\ref{cap:normalized.capture} verify the hypotheses of Proposition~\ref{norm:transfer}, yielding Theorem~\ref{cap:sampling}, as shown next.

\begin{proof}[Proof of Theorem~\ref{cap:sampling}]
Apply Proposition~\ref{norm:transfer} on $\mathscr{S}_{T} = I_{\theta}$ with the full measures $\mu_{T}^{\star} = \mu_{\star}$ and the capped measures $\nu_{T,b}^{\star} = \nu_{\star,b}$ constructed in Proposition~\ref{cap:normalized.capture}. The full density $\ell_{T}^{-1}e^{\beta(F_{h}^{\star} - u_{T})}$ is continuous and positive, so its integral is finite and positive. Normalization cancels its common deterministic multiplier $\ell_{T}^{-1}e^{-\beta u_{T}}$ and gives
\[
\frac{\mu_{\mathrm{A}}(\rd{h})}{\mu_{\mathrm{A}}(I_{\theta})} = \frac{e^{\beta V_{h}}\rd{h}}{\int_{I_{\theta}}e^{\beta V_{h'}}\rd{h'}} = G_{\beta,T}^{\mathrm{A}}(\rd{h}), \qquad
\frac{\mu_{\mathrm{H}}(\rd{h})}{\mu_{\mathrm{H}}(I_{\theta})} = \frac{e^{\beta H(h)}\rd{h}}{\int_{I_{\theta}}e^{\beta H(h')}\rd{h'}} = G_{\beta,T}^{\mathrm{H}}(\rd{h}).
\]

For each fixed $b$, the parameters $\chi_{b}$, $\e_{b}$ and $C_{0,b}$ are independent of $T$ and of the moment order $r$. Consequently, Proposition~\ref{cap:moment.measures} verifies the moment-measure condition \eqref{norm:detail:moment.measure.assumption} of Proposition~\ref{norm:transfer} for every fixed $r$. The mass-moment bound in Lemma~\ref{cap:moments}, with these same parameters, gives
\[
\limsup_{T \to \infty}\EE_{\mathrm{H}}[\nu_{\mathrm{H},b}(I_{\theta})^{r}] \leq C_{\mathrm{mom},b}^{r}r!, \qquad r \in \N,
\]
where $C_{\mathrm{mom},b}$ is independent of $r$. This is the mass-moment condition \eqref{norm:detail:mass.moment.assumption} in Proposition~\ref{norm:transfer}. The zero-atom estimate \eqref{cap:capture.detail:zero.atom} proves its condition \eqref{norm:detail:zero.atom}, and the normalized capture estimate in Proposition~\ref{cap:normalized.capture} proves its condition \eqref{norm:detail:capture.assumption}. When $\nu_{\star,b}(I_{\theta}) = 0$, use the same deterministic probability measure in both propositions.

All four conditions in Proposition~\ref{norm:transfer} therefore hold for a single sequence $(\nu_{\star,b})_{b \in \N}$. The conclusion \eqref{norm:detail:normalization.conclusion} of Proposition~\ref{norm:transfer} yields
\[
\left\|\EE (G_{\beta,T}^{\mathrm{A}})^{\times n} - \EE_{\mathrm{H}}(G_{\beta,T}^{\mathrm{H}})^{\times n}\right\|_{\mathrm{TV}} \to 0
\]
for every fixed $n \in \N$, as asserted.
\end{proof}

\subsection{Proof of the Euler overlap theorem}\label{sec:main.proof}

\begin{proof}[Proof of Theorem~\ref{expand:main}]
For $0 < \beta \leq \beta_{c}(\theta)$, all the conclusions follow from Corollary~\ref{temp:overlap}.

Now suppose that $\beta > \beta_{c}(\theta)$. The parameters in \eqref{eq:fixed.parameters} satisfy the range in \citet[Eq.~(1.1)]{OuimetHaar2026}, with the same $\Phi$, $\phi$, and $C_{\sig}$. Proposition~3.1 of that paper therefore supplies the Haar overlap laws and approximate GG identities used below; it does not require the last condition in \eqref{eq:fixed.parameters}. Theorem~\ref{cap:sampling} gives, for every fixed $r \in \N$,
\begin{equation}\label{expand:detail:sampling.distance}
\Delta_{r,T} \leqdef \left\|\EE (G_{\beta,T}^{\mathrm{A}})^{\times r} - \EE_{\mathrm{H}}(G_{\beta,T}^{\mathrm{H}})^{\times r}\right\|_{\mathrm{TV}} \to 0.
\end{equation}
The uniform approximations of $R_{T}^{\mathrm{A}}$ and $R_{T}^{\mathrm{H}}$ by $r_{T}$ in \eqref{an:overlap.geometry} of Lemma~\ref{an:covariance} imply
\begin{equation}\label{expand:detail:kernel.distance}
e_{T}^{\mathrm{H}} \leqdef \sup_{h,h' \in I_{\theta}}|R_{T}^{\mathrm{A}}(h,h') - R_{T}^{\mathrm{H}}(h,h')| = \OO\left(\frac{\log{t}}{t}\right) \to 0.
\end{equation}
Cauchy--Schwarz also yields
\begin{equation}\label{expand:detail:overlap.bounded}
|R_{T}^{\mathrm{A}}(h,h')| \leq 1, \qquad |R_{T}^{\mathrm{H}}(h,h')| \leq 1.
\end{equation}
Equations \eqref{expand:detail:sampling.distance}, \eqref{expand:detail:kernel.distance}, and \eqref{expand:detail:overlap.bounded} together give, for $1 \leq i,j \leq r$ and $\|\Psi_{T}\|_{\infty} \leq 1$,
\begin{equation}\label{expand:detail:one.overlap.comparison}
\begin{aligned}
&\left|\EE (G_{\beta,T}^{\mathrm{A}})^{\times r}[R_{T}^{\mathrm{A}}(h_{i},h_{j})\Psi_{T}(\bb{h})] - \EE_{\mathrm{H}}(G_{\beta,T}^{\mathrm{H}})^{\times r}[R_{T}^{\mathrm{H}}(h_{i},h_{j})\Psi_{T}(\bb{h})]\right| \\
&\qquad \leq \|R_{T}^{\mathrm{A}}(h_{i},h_{j})\Psi_{T}(\bb{h})\|_{\infty}\Delta_{r,T} + \|(R_{T}^{\mathrm{A}} - R_{T}^{\mathrm{H}})(h_{i},h_{j})\Psi_{T}(\bb{h})\|_{\infty}
\leq \Delta_{r,T} + e_{T}^{\mathrm{H}}.
\end{aligned}
\end{equation}

Fix $1 \leq k \leq n$, and write the GG residual in \eqref{expand:GG} as
\[
\begin{aligned}
\Delta_{\mathrm{A}}^{\mathrm{GG}}(\Psi_{T})
&\leqdef \EE (G_{\beta,T}^{\mathrm{A}})^{\times(n + 1)}[R_{T}^{\mathrm{A}}(h_{k},h_{n + 1})\Psi_{T}(\bb{h})] - \frac{1}{n}\EE (G_{\beta,T}^{\mathrm{A}})^{\times 2}[R_{T}^{\mathrm{A}}(h_{1},h_{2})]\EE (G_{\beta,T}^{\mathrm{A}})^{\times n}[\Psi_{T}(\bb{h})] \\
&\qquad - \frac{1}{n}\sum_{\substack{1 \leq \ell \leq n \\ \ell \neq k}}\EE (G_{\beta,T}^{\mathrm{A}})^{\times n}[R_{T}^{\mathrm{A}}(h_{k},h_{\ell})\Psi_{T}(\bb{h})],
\end{aligned}
\]
and denote by $\Delta_{\mathrm{H}}^{\mathrm{GG}}(\Psi_{T})$ the Haar counterpart. Using \eqref{expand:detail:sampling.distance}, \eqref{expand:detail:overlap.bounded}, and \eqref{expand:detail:one.overlap.comparison}, it follows that, uniformly over $\|\Psi_{T}\|_{\infty} \leq 1$,
\begin{equation}\label{expand:detail:GG.comparison}
\begin{aligned}
|\Delta_{\mathrm{A}}^{\mathrm{GG}}(\Psi_{T}) - \Delta_{\mathrm{H}}^{\mathrm{GG}}(\Psi_{T})|
&\leq (\Delta_{n + 1,T} + e_{T}^{\mathrm{H}}) + \frac{(\Delta_{2,T} + e_{T}^{\mathrm{H}}) + \Delta_{n,T}}{n} + \frac{n - 1}{n}(\Delta_{n,T} + e_{T}^{\mathrm{H}}).
\end{aligned}
\end{equation}
The right-hand side tends to zero by \eqref{expand:detail:sampling.distance} and \eqref{expand:detail:kernel.distance}, uniformly in $\Psi_{T}$. The Haar GG identities of \citet[Proposition~3.1]{OuimetHaar2026} give $\sup_{\|\Psi_{T}\|_{\infty} \leq 1}|\Delta_{\mathrm{H}}^{\mathrm{GG}}(\Psi_{T})| \to 0$. Combining these convergences with the triangle inequality gives
\[
\begin{aligned}
\sup_{\|\Psi_{T}\|_{\infty} \leq 1}|\Delta_{\mathrm{A}}^{\mathrm{GG}}(\Psi_{T})|
\leq \sup_{\|\Psi_{T}\|_{\infty} \leq 1}|\Delta_{\mathrm{A}}^{\mathrm{GG}}(\Psi_{T}) - \Delta_{\mathrm{H}}^{\mathrm{GG}}(\Psi_{T})|
+ \sup_{\|\Psi_{T}\|_{\infty} \leq 1}|\Delta_{\mathrm{H}}^{\mathrm{GG}}(\Psi_{T})|
\to 0.
\end{aligned}
\]
This proves \eqref{expand:GG} for the Euler field, uniformly over $\|\Psi_{T}\|_{\infty} \leq 1$.

To transfer the limiting joint distributions of the overlaps under the annealed Haar sampling law to those under the annealed Euler sampling law, let $f$ be continuous on $[ - 1,1]^{n \times n}$ and define its modulus of continuity by
\[
\omega_{f}(\delta) \leqdef \sup\{|f(\bb{x}) - f(\bb{y})| : \bb{x},\bb{y} \in [ - 1,1]^{n \times n},\ \max_{i,j}|x_{i,j} - y_{i,j}| \leq \delta\}.
\]
Compactness of $[ - 1,1]^{n \times n}$ gives $\omega_{f}(\delta) \to 0$ as $\delta \downarrow 0$. Write $\bb{Q}_{T}^{\mathrm{A}} \leqdef (R_{T}^{\mathrm{A}}(h_{i},h_{j}))_{1 \leq i,j \leq n}$ and $\bb{Q}_{T}^{\mathrm{H}} \leqdef (R_{T}^{\mathrm{H}}(h_{i},h_{j}))_{1 \leq i,j \leq n}$. Add and subtract the Haar sampling expectation of $f(\bb{Q}_{T}^{\mathrm{A}})$. The first difference is controlled by \eqref{expand:detail:sampling.distance}; the second is pointwise bounded by $\omega_{f}(e_{T}^{\mathrm{H}})$ through \eqref{expand:detail:kernel.distance}. Consequently,
\begin{equation}\label{expand:detail:array.comparison}
\begin{aligned}
&\left|\EE (G_{\beta,T}^{\mathrm{A}})^{\times n}[f(\bb{Q}_{T}^{\mathrm{A}})] - \EE_{\mathrm{H}}(G_{\beta,T}^{\mathrm{H}})^{\times n}[f(\bb{Q}_{T}^{\mathrm{H}})]\right| \\
&\qquad \leq \left|\EE (G_{\beta,T}^{\mathrm{A}})^{\times n}[f(\bb{Q}_{T}^{\mathrm{A}})] - \EE_{\mathrm{H}}(G_{\beta,T}^{\mathrm{H}})^{\times n}[f(\bb{Q}_{T}^{\mathrm{A}})]\right|
+ \EE_{\mathrm{H}}(G_{\beta,T}^{\mathrm{H}})^{\times n}[|f(\bb{Q}_{T}^{\mathrm{A}}) - f(\bb{Q}_{T}^{\mathrm{H}})|] \\
&\qquad \leq \|f\|_{\infty}\Delta_{n,T} + \omega_{f}(e_{T}^{\mathrm{H}}) \to 0.
\end{aligned}
\end{equation}
By \citet[Proposition~3.1]{OuimetHaar2026}, the Haar overlap array has the limiting law in \eqref{eq:main.joint.overlap} of Theorem~\ref{thm:limiting.joint.overlap.distribution}. Since $\theta > 0$, its limiting matrix has entries $\ind_{\{K_{i} = K_{j}\}}$, with the group labels sampled from the Poisson--Dirichlet weights of parameter $\beta_{c}(\theta)/\beta$. Combining that Haar limit with \eqref{expand:detail:array.comparison} gives
\[
\lim_{T \to \infty}\EE (G_{\beta,T}^{\mathrm{A}})^{\times n}[f(\bb{Q}_{T}^{\mathrm{A}})] = \EE_{\bb{\xi}}\big[\EE[f((\ind_{\{K_{i} = K_{j}\}})_{1 \leq i,j \leq n}) \mid \bb{\xi}]\big].
\]
This proves the joint overlap assertion. For $n = 2$, take $f(\bb{x}) \leqdef g(x_{1,2})$ with $g$ continuous on $[ - 1,1]$. Equation~\eqref{expand:detail:array.comparison} and the Haar two-overlap law in \citet[Proposition~3.1]{OuimetHaar2026} give the stated mixture of $\delta_{0}$ and $\delta_{1}$.
\end{proof}

\section*{Statement of AI use}
\addcontentsline{toc}{section}{Statement of AI use}

A rough draft of this paper, containing all the main results stated in Section~\ref{sec:main.results}, was written in 2019--2021 during my postdoctoral appointment at the California Institute of Technology. ChatGPT 6 Astra assisted in fixing and writing some of the missing mathematical arguments. The portions written by AI have been completely reworked and validated by hand. I take full responsibility for the accuracy and rigor of all proofs.

\section*{Funding}
\addcontentsline{toc}{section}{Funding}

F.\ Ouimet is supported by the Fonds de recherche du Québec -- Nature et technologies (FRQNT) through the Relève professorale grant doi:10.69777/377977, and by the Natural Sciences and Engineering Research Council of Canada (NSERC) through Discovery Grant RGPIN-2026-04471 and Discovery Launch Supplement DGECR-2026-00449.

\renewcommand{\bibsection}{\section*{References}\addcontentsline{toc}{section}{References}}

\setlength{\bibsep}{0pt plus 0ex}

\bibliographystyle{authordate1-doi-mr}
\bibliography{bib}

\end{document}